\documentclass{amsart}
\usepackage{amssymb,amscd, mathrsfs}
\newcommand{\nc}{\newcommand}
\newcommand{\longhookrightarrow}{\DOTSB\lhook\joinrel\longrightarrow} 
\renewcommand{\r}{\mathfrak{r}}
\nc{\ban}[1]{\rm{(#1)}}%{\rm{(}#1\rm{)}}
\nc{\longra}{\longrightarrow}
\nc{\Um}{\U^-(\Wg)}
\nc{\CClam}{C}
\nc{\zero}{0}
\nc{\I}{\mathsf{I}}
\nc{\nmap}{\Psi}
\nc{\WW}{\mathsf{W}}
\nc{\vQ}{\mathsf{Q}}
\nc{\vchi}{\mathsf{\chi}}
\nc{\sig}{\sigma}
\nc{\TT}{T}
\nc{\sE}{\bar{E}}
\nc{\Lone}{T^*}
\nc{\Lnew}{T_1}
\nc{\Lold}{T_1^{\old}}
\nc{\nno}{\nonumber}
\nc{\omu}{\bar{v}_{\mu}}
\nc{\gU}{\a_{U''}}
\nc{\olam}{\overline{|\lam\ket}}
\nc{\bCmu}{K(\mu)}%{\overline{C_+(\mu)}}
\nc{\bClam}{K(\lam)}%{\overline{C_+(\lam)}}
\nc{\Clam}{C_+(\lam)}
\nc{\ttp}{\tp}
\nc{\FF}{\mathsf{F}}
\nc{\nW}{\bar{\W}_+}
\nc{\newWU}{\wt{U_k(\bg')\* \bCl}_{\new}}
\nc{\oldWU}{\wt{U_k(\bg')\* \bCl}_{\old}}
\nc{\newU}{U_k(\bg')\* \bCl_{\new}}
\nc{\ra}{\rightarrow}
\nc{\Hp}[1]{H^{#1}}
\nc{\UV}{\U(\V)}
\nc{\wJ}{\widehat{J}}
\nc{\wBB}{\wt{\B}}
\nc{\UN}[1]{\U_N(#1)}
\nc{\QO}{Q}
\nc{\AN}[1]{\mathscr{A}_N(#1)}
\nc{\adMod}{\text{-}\mathsf{adMod}}
\nc{\GrMod}{\text{-}\mathsf{grMod}}%{\text{-}\mathsf{Mod^{ad}}}
\nc{\bCgoo}[1]{C_k^{#1}(\sg)''}
\nc{\bCgo}{\bCg''}
\nc{\bCGo}{\bCg''}
\nc{\wUU}{\wt{\UU}}
\nc{\QN}{\U(\V)/\wt{\II}_N}
\nc{\QNS}[1]{\U(\V_{#1})/\wt{\II}_N^{(#1)}}
\nc{\LV}{L(\V)}
\nc{\wh}{\widehat}
\nc{\adQ}{\ad Q }
\nc{\dV}{Q_{(0)}}
\nc{\wdV}{\wt{\QO}}
\nc{\Basis}{\mathcal{B}}
\nc{\UU}{\mathbf{U}}
\nc{\bU}{\mathbb{U}}
\nc{\II}{\mathscr{I}}
\nc{\wt}{\widetilde}
\nc{\Fin}{\mathsf{Fin}}
\nc{\Vkg}{V_k(\sg)}
\nc{\wpar}{\partial_{-1}}%{\widetilde{\partial}_{-1}}
\nc{\wU}{\widetilde{U}}
\nc{\Wg}{\W_k(\sg)}
\nc{\wc}{\widetilde{c}}
\nc{\twlam}{\overline{t_{-\srho\che}\circ \lam}}
\nc{\bge}{\bg^e}
\nc{\bgf}{\bg^f}
\nc{\sge}{\sg^e}
\nc{\sgf}{\sg^f}
\nc{\WN}{\mathbf{N}}
\nc{\bt}{\bh}
\nc{\bCl}{\Cl}
\nc{\bCg}{C_k(\sg)}
\nc{\sCg}{\sC(\sg)}
\nc{\sN}{\bar{N}}
\nc{\sG}{\bar{G}}
\nc{\smu}{\bar{\mu}}
\nc{\salpha}{\bar{\alpha}}
\nc{\hatJ}{\widehat{J}}
\nc{\sA}{\bar{A}}
\nc{\Cat}{{\mathcal{C}}}
\nc{\Lie}[1]{{\mathfrak{L}(#1)}}%{{\mathscr{L}(#1)}}
\nc{\sw}{\bar{w}}
\nc{\sC}{\bar{C}}
\nc{\dBGG}{\dot{\BGG}}
\nc{\tri}{\triangle}
\nc{\GG}{\mathsf{Y}}
\nc{\Prnongen}{{\textsl{Pr}}_{\kappa,\mathrm{nondeg}}}
\nc{\Center}{{\mathcal{Z}}}
\nc{\Wcat}{\mathscr{O}(\Wg)}%{{\mathsf{Rep}}\Wg}%{\Wg\text{-}\mathsf{grmod}}
\nc{\WCat}{\Wcat}
\nc{\fmap}{\rightsquigarrow }
\nc{\fdomain}{{\mathcal{C}}}
\nc{\D}{{\mathsf{D}}}

\nc{\whJ}{\widehat{J}}
\nc{\bd}{Q}
\nc{\blam}{\lam}
\nc{\dchi}{\chi}
\nc{\sV}{\bar{V}}
\nc{\sCl}{\bar{\Cl}}
\nc{\sH}{\bar{H}}
\nc{\schi}{\bar{\chi}}
\nc{\sL}{\bar{L}}
\nc{\sM}{\bar{M}}
\nc{\sd}{\bar{Q}}
\nc{\slam}{\bar{\lam}}

\newcommand{\teigi}{:=} %{\underset{\mathrm{def}}{=}}

\newcommand{\V}{{\mathbb V}}

\newcommand{\why}{{\bf L}}

\newcommand{\M}{{\bf{M}}}
\newcommand{\B}{{\mathbf{B}}}

\newcommand{\N}{{\mathbb{N}}}
\newcommand{\U}{{\mathfrak{U}}}%{{\mathscr{U}}}
\newcommand{\DW}{\Dg_{\new}}

\newcommand{\Dg}{{\bf D}}

\newcommand{\isomap}{\overset{\sim}{\rightarrow} }
\newcommand{\W}{{\mathscr{W}}}

\newcommand{\sI}{\bar{I}}
\newcommand{\1}{{\mathbf{1}}}
\newcommand{\Ln}[1]{L\bar{\mathfrak{n}}_{#1}}

\newcommand{\sQ}{\bar{Q}}
\newcommand{\eW}{\widetilde{W}}
\renewcommand{\Pr}{\textsl{{Pr}}_{k}}

\newcommand{\op}{\rm{op}}

\newcommand{\gen}{{\rm{gen}}}
\newcommand{\sroots}{\bar{\Delta}}
\newcommand{\roots}{\Delta}
\newcommand{\proots}{\Delta_+}

\newcommand{\rroots}{\Delta^{\rm re}}
\newcommand{\iroots}{\Delta^{\rm im}}
\newcommand{\prroots}{\Delta_+^{\rm re}}

\newcommand{\nrroots}{\Delta_-^{\rm re}}
\newcommand{\F}{{\mathcal{F}}}

\renewcommand{\a}{{\mathfrak{a}}}
\newcommand{\sa}{\bar{\mathfrak{a}}}
\renewcommand{\sl}{{\mathfrak{sl}}}

\newcommand{\sPi}{\bar{\Pi}}

\newcommand{\Cl}{{\mathcal{C}l}}
\newcommand{\Tg}{L^{\g}}
\newcommand{\Tf}{L^f}

\newcommand{\bW}{{W}}

\newcommand{\semiinf}{\frac{\infty}{2}}

\renewcommand{\theenumi}{\roman{enumi}}

\newcommand{\C}{{\mathbb C}}
\newcommand{\Z}{{\mathbb Z}}
\newcommand{\Q}{{\mathbb Q}}
\newcommand{\inv}{^{-1}}
\newcommand{\dual}[1]{{#1}^*}
\newcommand{\lam}{\lambda}
\newcommand{\Lam}{\Lambda}
\renewcommand{\P}{{\mathcal{P}}}
\renewcommand{\*}{{\otimes}}
\newcommand{\+}{\mathop{\oplus}}
\newcommand{\h}{ {\mathfrak h}}
\newcommand{\g}{ {\mathfrak g}}

\renewcommand{\sb}{\bar{\mathfrak{b}}}
\newcommand{\che}{^{\vee}}
\renewcommand{\S}{{\mathbb{S}}} %{{\mathbb{S}}}

\newcommand{\bra}{{\langle}}
\newcommand{\ket}{{\rangle}}

\DeclareMathOperator{\Vac}{Vac}

\DeclareMathOperator{\Zhu}{\mathfrak{Zh}}

\DeclareMathOperator{\Wh}{Wh}

\DeclareMathOperator{\nondeg}{nondeg}
\DeclareMathOperator{\tp}{top}
\DeclareMathOperator{\vac}{vac_k}

\DeclareMathOperator{\tr}{tr}

\DeclareMathOperator{\Gr}{gr}
\DeclareMathOperator{\End}{End}

\DeclareMathOperator{\Hom}{Hom}

\DeclareMathOperator{\id}{id}
\DeclareMathOperator{\ad}{ad}
\DeclareMathOperator{\haru}{span}

\DeclareMathOperator{\ch}{ch}

\DeclareMathOperator{\rank}{rank}
\DeclareMathOperator{\gr}{gr}
\DeclareMathOperator{\height}{ht}

\DeclareMathOperator{\new}{new}
\DeclareMathOperator{\old}{old}

\DeclareMathOperator{\even}{even}
\DeclareMathOperator{\odd}{odd}

\newcommand{\sg}{  \bar{\mathfrak g}}
\newcommand{\sh}{\bar{ \h}}
\newcommand{\sn}{\bar{\mathfrak{n}}}
\newcommand{\snp}{\bar{\mathfrak{n}}_+}
\newcommand{\snn}{\bar{\mathfrak{n}}_-}

\newcommand{\sP}{\bar P}
\newcommand{\sW}{\bar{W}}
\newcommand{\srho}{ \bar{\rho}}

\newcommand{\sproots}{\bar{\Delta}_+}
\newcommand{\snroots}{\bar{\Delta}_-}

\newcommand{\bg}{{\g}}
\newcommand{\bh}{ \h}

\newcommand{\bQ}{ Q}

\DeclareMathOperator{\re}{re}

\DeclareMathOperator{\im}{Im}
\DeclareMathOperator{\Ker}{ker}

\newcommand{\sBGG}{ {\BGG}(\sg)}

\newcommand{\ud}[2]{{\genfrac{}{}{0pt}{}{#1}{#2}}}

\newcommand{\n}{{{\mathfrak{n}}}}

\newcommand{\BGG}{{\mathcal O}}

\newcommand{\Obj}{Obj}

\title{Representation theory of $\W$-algebras}
\author{Tomoyuki Arakawa}
\address{Department of Mathematics, 
Nara Women's University, Nara 630-8506, JAPAN}
\email{arakawa@cc.nara-wu.ac.jp}
\thanks{The author is partially supported 
by the JSPS Grant-in-Aid for Young Scientists (B)
No.\ 17740006}
\begin{document}
%\today
%%%%%%%% Thoremstyles %%%%%%%%%%
\theoremstyle{plain}
\newtheorem{Th}{Theorem}[subsection]
\newtheorem{MainTh}{Main Theorem}
\newtheorem{Pro}[Th]{Proposition}
\newtheorem{ProDef}[Th]{Proposition and Definition}
\newtheorem{DefPro}[Th]{Definition and Proposition}
\newtheorem{Lem}[Th]{Lemma}
\newtheorem{Co}[Th]{Corollary}
\newtheorem{Facts}[Th]{Facts}
\theoremstyle{definition}
\theoremstyle{remark}
\newtheorem{Def}[Th]{Definition}
\newtheorem{Rem}[Th]{Remark}
\newtheorem{Conj}{Conjecture}
\newtheorem{Claim}{Claim}
\newtheorem*{ClaimNN}{Claim}
\newtheorem{Notation}{Notation}
\newtheorem{Ex}[Th]{Example}
\newcommand{\st}{{\mathrm{st}}}
%\numberwithin{equation}{subsection}
%\begin{flushright}
%personal preprint, \today
%\end{flushright}
\maketitle
\begin{abstract}
We study  the representation theory of the $\W$-algebra
$\Wg$ associated with a simple Lie algebra
$\sg$ at level $k$. We show that the ``$-$'' reduction 
functor
is exact and sends an irreducible module  to
zero or an irreducible module at any  level $k\in \C$.
Moreover,  we show that the character  of
each irreducible highest weight representation
of $\W_k(\sg)$ is completely determined by that of 
the corresponding  irreducible highest weight representation
of affine Lie algebra $\g$ of $\sg$.
As a consequence
we complete (for the  ^^ ^^ $-$" reduction)
the proof of the 
conjecture  of E. Frenkel, V. Kac and M. Wakimoto
on
the existence and the construction of 
the modular invariant representations
of $\W$-algebras.
%This article is the detailed  version of
%\cite{A-irr} with extended results.
\end{abstract}
\tableofcontents
\section{Introduction}
%\subsection{}
This paper is a continuation of the  author's earlier work \cite{A04}
in which we  proved (completely for  the ``$-$'' reduction and partially for
the  ``$+$'' reduction)
the vanishing conjecture of
E. Frenkel, V. Kac and M. Wakimoto \cite{FKW}.

\smallskip

Let $\sg$ be a finite-dimensional 
Complex simple Lie algebra,
$\sg=\sn_-\+\sh\+\sn_+$ its triangular decomposition,
$\bg$ the Kac-Moody affinization of $\sg$,
 $k$ the level of $\bg$.
Then
the corresponding  {\em $\W$-algebra}
 $\W_k(\sg)$
(\cite{Za,FZa,  Luk,FL, LF, FF_W}),
which is a M\"obius conformal vertex algebra,
can be constructed %as a vertex algebra
by the method of 
B.  Feigin and E. Frenkel  \cite{FF_W} (see Section
\ref{Section:definition-of-W-algebra}).
Let $h\che$ be the dual 
Coxeter number
of $\sg$.
It is known \cite{FF_W}
that
for $k\ne -h\che$
the vertex algebra $\W_k(\sg)$
is conformal (or a vertex operator algebra)
of central charge $c(k)$,
where
\begin{align*}
 c(k)
=\rank \sg-12\left((k+h\che)|\srho\che|^2
-2 \bra \srho,
\srho\che\ket+ \frac{1}{k+h\che}|\srho|^2
\right),
\end{align*} 
while $\W_{-h\che}(\sg)$
 coincides with 
the center 
of the universal affine vertex algebra
$V_{-h\che}(\sg)$
associated with $\sg$ at level $-h\che$
(see Section \ref{subsection:univ-affine}).
Here
$\srho$ is the half sum of positive roots of $\sg$
and
$\srho\che$
is the half sum of positive coroots of $\sg$.

The simplest $\W$-algebra
$\W_k(\mathfrak{sl}_2)$ 
is the Virasoro vertex algebra,
%at central charge $c(k)$
%(for $k\ne -2$),
whose representation theory
is well-studied
(\cite{KacDet, FeiFuch-Det,
KR,
FeiFuch}).
However,
$\W$-algebras associated with other simple Lie algebras
are not generated by Lie algebras,
and hence one cannot apply the Lie algebra theory
to $\W$-algebras directly.
As a consequence,
 surprisingly little is known about 
their representation theory
(see eg.\ \cite{BMP}).

%The purpose of this paper is to determine
%the character of each 
%irreducible highest weight representation
%of  $\W_k(\sg)$ at each  $k\in \C$.
%Applied
%for an admissible number \cite{KW2} $k\in \Q$,
%this in particular completes %we complete 
%(for the  ``$-$'' reduction)
%the proof the conjecture of
%E. Frenkel, V. Kac and M. Wakimoto (\cite{FKW})
%on
%the existence and the construction of the minimal series representations
%($=$ the modular invariant representations)
%of $\W_k(\sg)$.
%

\smallskip

Let $\Zhu(\W_k(\sg))$ be the {\em Zhu algebra}
\cite{FZ, Zhu}
of $\W_k(\sg)$ (see Section \ref{subsection:MOdules-and-Zhu-algberas}).
Then,
by construction \cite{FF_W},
$\Zhu(\W_k(\sg))$
is isomorphic to 
the center $\Center(\sg)$
of the universal enveloping algebra $U(\sg)$
of $\sg$ (Theorem \ref{Th:2006-10-19-30-31} (ii)).
Therefore,
by Zhu's theorem \cite{Zhu},
irreducible highest weight representations
%(=irreducible admissible representation in the sense of \cite{ABD})
of $\W_k(\sg)$
are parameterized\footnote{%
%An irreducible highest weight representation
%is a irreducible admissible 
%representation in the sense of \cite{ABD}.
To be precise,
in the case that
the level $k$ is critical,
``irreducible highest weight 
representations''
should be replaced by ``graded irreducible highest weight
representations'', see Section 
\ref{subSection:The Virasoro Field at Non-Critical Level}.}
by the central characters 
of $\Center(\sg)\cong S(\sh)^{\sW}$.
Let $\gamma_{\slam}:\Center(\sg)\ra \C$
be the evaluation
at the Verma module $\sM(\slam)$ of $\sg$ 
with highest weight $\bar{\lam}\in \dual{\sh}$.
Denote by
$\why(\gamma_{\slam})$,
with $\slam\in -\srho+\sW\backslash \dual{\sh}$,
be the corresponding
irreducible highest weight representation
of $\W_k(\sg)$.
Here $\sW$ is the Weyl group of
$\sg$.  %and $\srho$ is the half sum of positive roots of $\sg$.
Then
$\why(\gamma_{\slam})$
is the simple quotient of
the Verma module
$\M(\gamma_{\slam})$
with highest weight $\gamma_{\slam}$
(see Section \ref{section:Simple module over}).
As in the case of Lie algebras, 
one has
 a simple character formula for $\M(\gamma_{\slam})$
 (see Proposition \ref{Pro:character-formula-of-Verma}).

One of the aim of 
this paper %problems of  the representation theory
%of  $\Wg$
% is of course
is to find the character formula of each irreducible highest weight
representation $\why(\gamma_{\slam})$,
that is,
to determine the integer $m_{\slam,\smu}$
in the expression
$\ch \why(\gamma_{\slam})=\sum_{\smu}m_{\slam \smu}\ch \M(\gamma_{\smu})$.

%\smallskip

To this end 
we study %Our main tool for studying
%representations
%of $\Wg$
%is
 the method of quantum reduction \cite{FF_W, FKW},
 in particular
the
``$-$'' reduction functor  \cite{FKW},
which is the modified version of the reduction functor
originally introduced in \cite{FF_W}
(the original reduction functor is referred to as
the ``$+$'' reduction functor):
Fix
the 
triangular decomposition 
$\bg=\bg_-\+\bh\+\bg_+$
of
 the affine Lie algebra
$\g$ %=\sg\* \C[t,t\inv]\+\C K\+ \C \Dg$
in the standard manner.
Then,
one has the
corresponding Bernstein-Gelfand-Gelfand category
$\BGG_{k}$
of $\bg$
 at level $k$ (see Section \ref{subsection:Category}).
Let $\dual{\bh}_k\subset \dual{\bh}$
be the set of weights of level $k$,
$M(\lam)\in \BGG_k$
the Verma module of $\bg$ with highest weight $\lam\in \dual{\bh}_k$,
$L(\lam)\in \BGG_k$ the unique simple quotient of $M(\lam)$.
Denote by $H^{\bullet}_{-}(M)$ 
 the cohomology with coefficient in $M\in \BGG_k$
associated with the
quantized Drinfeld-Sokolov ``$-$'' reduction 
\cite{FF_W,FKW} (see Section \ref{subsection:H_-}).
Then
the correspondence
\begin{align*}
M\longmapsto H^0_-(M)
\end{align*} %with $i\in \Z$,
defines a functor from $\BGG_k$ to
the  category of $\W_k(\sg)$-modules\footnote{The action
of $\Wg$ on $H^{\bullet}_-(M)$ is slightly changed,
see Section \ref{subsection:H_+}.}.
%We call this functor
%the ``$-$'' reduction functor in this paper.

\smallskip

Now recall that
if $\lam\in \dual{\bh}$ is admissible
\cite{KW2}
then $L(\lam)$ is called
an {\em Kac-Wakimoto admissible representation}\footnote{This should 
not confused with the notion
of admissible modules in the theory of
vertex algebras \cite{ABD},
which also appears  in this paper.}.
Conjecturally,
Kac-Wakimoto admissible representations exhaust all modular invariant
representations of $\bg$.
In \cite{FKW},
E. Frenkel, V. Kac and M. Wakimoto conjectured
that,
for an admissible weight $\lam\in \dual{\bh}$,
one has
$H^i_-(L(\lam))=\zero $
with $i\ne 0$ and $H^0_-(L(\lam))$
is an irreducible $\W_k(\sg)$-module
if $\lam$ is non-degenerate \cite{KW90}.
(They showed that $H^0_-(L(\lam))$
is zero if $\lam$ is a degenerate admissible weight\footnote{
A
dominant 
integral
weights of $\bg$ is a degenerate admissible weight.
Therefore,
an integrable representation
 of $\bg$ goes to zero by the functor $H_-^0(?)$ .}.)
The remarkable consequence  of this conjecture is that
%from this
%it  follows that
the existence and the construction of modular
invariant representations of $\W$-algebras.
Namely,
one can obtain
(non-trivial) modular invariant representations of $\W$-algebras
as
the image of  non-degenerate
Kac-Wakimoto
admissible representations of $\bg$
by the functor $H_-^0(?)$.
%of $\W_k(\sg)$, see 
%because admissible representations of $\bg$
%are modular invariant
In the case that $\sg=\mathfrak{sl}_2$,
this conjecture is  known to be true (\cite{FKW}),
%fact which follows from earlier works
%\cite{BO, FF-boson}
and the corresponding modular invariant 
representations are the minimal series representations
\cite{BPZ}
of the Virasoro (vertex) algebra.
%as explained in 

\smallskip

Of this conjecture
the part concerning the vanishing of cohomology
is proved in \cite{A04}.
Therefore,
it is remains to prove
%the remaining part to be proved
%is 
the irreducibility.
It turns out that the method of quantum reduction
produces not only 
modular invariant representations
but also 
 {\em all}
 the highest weight irreducible representations
of $\W$-algebras,
determining their characters.

\medskip

Let $\dual{\h}\ni \lam \mapsto \slam
\in \dual{\sh}$ be the restriction.
We refer to $\slam$ as the {\em classical part}
 of $\lam$.
The weight $\slam\in \dual{\sh}$
is called {\em{anti-dominant}}
if
$\sM(\slam)$ is irreducible over $\sg$.
 One knows that
this condition is equivalent
to %the condition that
\begin{align}
 \bra \slam+\srho,\salpha\che\ket \not \in \{1,2,\dots\}
\text{ for each element
$\salpha$ of  positive roots $\sroots_+$ of $\sg$}.
\end{align}
It is clear that any central character of $\Center(\sg)$
is of the form $\gamma_{\slam}$
for some anti-dominant $\slam\in \dual{\sh}$.
Note  that
if $\lam\in \dual{\bh}$
an admissible weight %\cite{KW2} ,
then by definition
it
is non-degenerate %\cite{KW90} 
if and only if its  classical part $\slam$
is anti-dominant\footnote{But 
a non-degenerate admissible weight $\lam$ is by no means 
``anti-dominant'' as a weight of affine Lie algebra $\bg$,
see \cite{KW90,FKW}.
(Indeed it is ``regular dominant'' so that the Weyl-Kac type character
formula holds for $L(\lam)$ (\cite{KW1}).)}.

Let $D(\why(\gamma_{\slam}))$ be the 
restricted dual of
$\why(\gamma_{\slam})$ in the sense of \cite{FHL}.
Then 
\begin{align*}
D(\why(\gamma_{\slam}))\cong \why(\gamma_{-w_0(\slam)})
\end{align*}
as $\Wg$-modules,
where $w_0$ is the longest element 
of $\sW$ (Theorem \ref{Th:Duality} (i)).
\begin{MainTh}\label{MainTh1}
Let $k$ be an arbitrary complex number.
\begin{enumerate}
 \item
{\rm{(Theorem \ref{Th:vanighing--})}}
For any object $M$ of $\BGG_k$,
the cohomology 
 $H^i_-(M)$ is zero unless $i=0$.
\item 
{\rm{(Theorem \ref{Th:irr--})}}
For any $\lam\in \dual{\bh}_k$,
 there is an  isomorphism
%of $\Wg$-modules:
\begin{align*}
 H^0_-(L(\lam))\cong \begin{cases}
		 D(\why(\gamma_{\slam}))&\text{if the classical part
		 $\slam$
is anti-dominant,}\\
\zero &\text{otherwise.}
		\end{cases}
\end{align*}
\end{enumerate}
\end{MainTh}
By Main Theorem 1 (1),
which 
 generalizes  the result obtained in \cite{A04},
it follows that
the  correspondence 
$
M\mapsto H^0_-(M)
$
defines an  exact functor from $\BGG_k$
to the category of $\W_k(\sg)$-modules at any $k\in \C$.
Also,
it holds that
$H_-^0(M(\lam))\cong \M(\gamma_{-w_0(\bar{\lam})})$
 (Theorem \ref{Th:image-of-Verma--}).
Therefore,
if we choose $\lam$ so that its classical part $\slam$ is
anti-dominant,
then
 we have
\begin{align*}
 \ch \why(\gamma_{\slam})=\sum_{\mu}
[L(\lam): M(\mu)]\ch \M(\gamma_{\bar{\mu}})
\end{align*}
(Theorem \ref{Th:ch-formula}).
Here
$[L(\lam): M(\mu)]\in \Z$
is determined by the formula
$
 \ch L(\lam)=\sum_{\mu}
[L(\lam): M(\mu)]\ch M(\mu)
$.
This integer
$[L(\lam):M(\mu)]$
is known\footnote{
The character formula of 
$L(\lam)$ is not known
if 
the level of $\lam$ is critical.
 However,
since
$\W_{-h\che}(\sg)$ is a commutative 
vertex algebra \cite{FF_W},
%Hence 
in this case 
one knows that
 each $\why(\gamma_{\slam})$ is one-dimensional,
see Section 
\ref{subSection:The Virasoro Field at Non-Critical Level}.} 
and written in terms of
Kazhdan-Lusztig polynomials
\cite{KTI, Cas, KTII, KT1, KT}.
Therefore
the above formula gives the  character
of each irreducible highest weight  representation 
$\why(\gamma_{\bar{\lam}})$.

%It is  clear 
%that 
The 
above explained 
conjecture of E. Frenkel, V. Kac and M. Wakimoto \cite{FKW}
follows immediately 
from   Main Theorem 1
%to 
%Kac-Wakimoto admissible %(conjecturally $=$ all modular invariant)
% representations  
%of $\g$ 
(Corollary \ref{Co:FKW-COng--}).

\medskip
Next 
 we 
describe 
 our results
for the
``$+$'' reduction.
In \cite{FKW},
E. Frenkel, V. Kac and M. Wakimoto 
gave the similar conjecture 
also for
the ``$+$'' reduction functor.
The 
``$+$'' reduction
can be defined  more directly using the vertex algebra theory
and 
%in
% particular,
in fact
$\W$-algebras themselves are defined
through the ``$+$'' reduction
(see Section \ref{subsection:def-of-W-algebras}
and Section \ref{subsection:H_+})).
However,
as observed in
  \cite{FKW},
the ``$+$'' reduction
is in general not 
as simple as  the ``$-$" reduction
(in particular if we consider the whole category $\BGG_k$).
%though
%$\W$-algebras themselves are defined through
%``$+$'' reduction.
This is basically due to the fact
that the
 ``$+$'' reduction is not compatible with the 
degree operator $\Dg$ of $\bg$, as opposed to 
the ``$-$'' reduction 
(however see \cite{A05Duke} for the
case that $\sg={\mathfrak{sl}}_2$).
Nevertheless, 
we have the following %(weaker)
 result:
Let $H^i_+(V)$
be the cohomology associated to the
quantized Drinfeld-Sokolov ``$+$'' reduction (\cite{FF_W, FB}, see
Section \ref{subsection:H_+}).
Denote by $\BGG_k^{[\lam]}$ the block of $\BGG_k$
corresponding to a weight $\lam\in \dual{\bh}_k$.
\begin{MainTh}\label{MainTh2}
 Suppose that $\lam\in \dual{\h}$ is non-critical and
satisfies the following condition:
\begin{align*}
 \bra \lam,\alpha\che\ket \not \in \Z\quad
\text{for all $\alpha\in \{-\bar{\beta}+n\delta;
\bar{\beta}\in \sproots, 1\leq n\leq \height \bar{\beta}\}$},
\end{align*}
where $\height \bar\beta$ is the height of $\bar \beta$.
\begin{enumerate}
 \item {$($\cite{A04}$)$}
Let $M$ be any object of  $\BGG_k^{[\lam]}$.
Then
the cohomology 
 $H^i_+(M)$ is zero
for all $i\ne 0$.
\item  {\rm{(Theorem \ref{Th:Main-results-for-+})}}
There is an isomorphism
%of $\Wg$-modules
$H_+^0(L(\lam))\cong \why(\gamma_{\slam-(k+h\che)\srho\che}
)$.
\end{enumerate}
%Here, $t_{\srho\che}$ is the element of extended Weyl group of $\bg$
%corresponding to $\srho\che$.
\end{MainTh}
Let $k$ be a rational number
such that
there exists a
non-degenerated principal admissible
weight \cite{KW2} at level $k$.
Then
the condition of
Main Theorem 2
is satisfied
for $\lam =k\Lam_0$ %(Corollary \ref{Co:2006-11-08-21-11} (ii))
and all the minimal series representations
of $\Wg$ at level $k$
can be obtained as the image 
of the functor $H^0_+(?)$
(see Remark \ref{Rem:2006-11-07-01-07}).
This
 proves the conjecture
of E. Frenkel, V. Kac and M. Wakimoto 
for the ``$+$'' reduction partially 
(Corollary \ref{Co:CONG-FKW-PLUS}).

\medskip

This paper is organized as follows.
In Section  \ref{section:quantum reduction for finite-dimensional Lie
algebras},
we collect some 
fundamental results
 concerning 
the
quantum reduction for finite-dimensional Lie algebras \cite{Kostant-Wh,
KS,
Ma}.
We remark that the finite-dimensional version of the functor $H^0_-(?)$
is identical to Soergel's functor \cite{Soe}
(\cite{Bac}, see Theorem \ref{TH:Soergel} 
and Remark \ref{Rem:-cohomology-as-semi-inf}).
In Section  3,
we  prepare some general theory about filtration 
of
vertex superalgebras
and BRST cohomology.
In particular the notion of strict filtration 
of vertex algebras
is introduced (Section \ref{subsecton-new:StrictFiltration}).
The definition of
the current algebra \cite{MNT}(=the universal enveloping
algebra in the sense of I. Frenkel and Y. Zhu \cite{FZ}) and 
the Zhu algebra
of a vertex algebra is also recalled.
Some technical notion needed for this is summarized 
in Appendix \ref{section:The Current Algebra of Vertex (Super)algebras}.
In section 4,
we recall 
the definition of 
$\W$-algebras %$\W_k(\sg)$
and collect necessary information about its structure,
following \cite{FB}.
In particular
we define an important strict filtration of $\Wg$
arising from the argument of \cite{FB}
(Section \ref{subsection:Poorf-vanihing-FBB}).
We also describe the current algebra
and the Zhu algebra of $\Wg$ (Theorems \ref{Th:vanishing-UBCGO}
\ref{Th:2006-08-29-09-53} and \ref{Th:2006-10-19-30-31}).
Theorem \ref{Th:2006-10-19-30-31} should be compared with 
Kostant's Theorem \cite{Kostant-Wh} (Theorem \ref{Th:Kostant-WH})
through Kostant-Sternberg's Theorem \cite{KS}
(Theorem \ref{Th:zhu1}).
In Section  \ref{section:Simple module over},
we define the Verma module $\M(\gamma_{\bar{\lam}})$ 
of $\W_{\kappa}(\sg)$ with highest weight $\gamma_{\bar{\lam}}$ and
its simple quotient $\why(\gamma_{\bar{\lam}})$.
The duality structure of modules over $\W$-algebras
is also discussed (Theorem \ref{Th:Duality}).
In Section  \ref{section:FunctorH},
we recall the definition of
 the category $\BGG_k$
and two reduction functors $H^{\bullet}_{+}(?)$
 and $H^{\bullet}_{-}(?)$.
We also recall some  of the structure of
the category $\BGG_k$
for the later purpose.
 Section \ref{section:RepresentationTheoryofW-algebrasThRSRroughtheFunctor}
is the main part of this paper:
we study
the representation theory of $\W$-algebras through the  
``$-$'' reduction functor
$H_-^0(?)$.
%In particular we  prove 
Section  \ref{section:Proof-Vanishing}
is devoted to the 
proof of Theorem \ref{Th:vanishing-dual-Verma--},
which is
the most crucial assertion in the proof
of Main Theorem \ref{MainTh1}.
 In Section  \ref{section:results-for-+}
we prove Main Theorem \ref{MainTh2}.

\medskip

In view of \cite{dBT, KRW, KW2003},
the $\W$-algebras $\Wg$ considered in this article
are the $\W$-algebras associated with principal nilpotent
orbits of $\sg$.
There has been  rapidly growing interest  in
$\W$-algebras associated with other nilpotent orbits
(see eg.\ \cite{Premet, GG, BK})\footnote{%
While the present article was being refereed, the paper
 [De Sole, A., Kac, V., {\em 
Finite vs affine W-algebras},
Japanese J.  Math. 1, No.\ 1,
 April, 2006,
137--261]  appeared, 
which studies the relationship between affine 
$\W$-algebras and finite 
$\W$-algebras in full generality.
}.
The method developed in the present paper works 
for general $\W$-algebras  as well (cf.\ \cite{A05Duke}).
The detail will appear in our forthcoming papers.

\medskip

The author thanks Profs.\ V. Kac, M. Wakimoto and I. Frenkel % and E. Frenkel
for conversation and valuable comments. 
Results in this paper were presented in part in
^^ Quantum Integral Systems and Infinite Dimensional Algebras',
 Kyoto, February 2004,
^^ International Conference on Infinite Dimensional Lie Theory, Beijing,
 July 2004',
^^ Perspectives Arising From Vertex Algebra Theory, Osaka, November
2004',
^^ International Conference on Algebra in Memory of Kostia Beidar',
Tainan, March 2005.
He thanks the organizers of these conferences. 
Finally, the author
 thanks the  referee and the editor
whose comments greatly improved
the paper.

\bigskip

\noindent {\em Notation. }
Throughout this paper the ground field
is the complex number $\C$
and
tensor products and dimensions are always meant to be as vector spaces
over $\C$ if not otherwise stated.

\section{Quantum Reduction for Finite-Dimensional Simple Lie Algebras}
\label{section:quantum reduction for finite-dimensional Lie algebras}
\subsection{The Setting}\label{subsection:Setting}
Let $\sg$ be 
a complex
finite dimensional   simple  Lie algebra,
$l=\rank \sg$.
Let
$(\cdot |\cdot)$
 the  normalized invariant 
inner product of $\sg$,
that is,
$(\cdot |\cdot)=\frac{1}{2h\che}\times $Killing form
where $h\che$ is the dual Coxeter number of $\sg$.

For $x\in \sg$
we  
write $\sg^x$ for the centralizer of $x$ in $\sg$.

Let $e$ be a principal nilpotent element of $\sg$,
so that
$\dim \sg^e=l$.
%One knows that $\dim \sg^e=l$
One knows that
$\sg^e$ is commutative.
By the Jacobson-Morozov theorem,
there exists an ${\mathfrak{sl}}_2$-triple
$\{e, h_0, f\}$ associated $e$,
that is,
\begin{align*}
[h_0,e]=2e,\quad  [h_0,f]=-2f, \quad [e,f]=h_0.
\end{align*}

Set
  \begin{align}
\sg_j\teigi \{x\in \sg;  [h_0,x]=2j x\}
\quad \text{ for $j\in \Z$}.
  \end{align}
This gives 
a  triangular decomposition 
$	       \sg=\snn\+\sh\+\snp$,
where
\begin{align}
\sh=\sg_0,\quad
\sn_+=\bigoplus_{j>0}\sg_j,\quad 
\sn_-=\bigoplus_{j<0}\sg_j.
\end{align}

We often identify the Cartan subalgebra $\sh$ 
with its dual  $\dual{\sh}$ using the form $(\cdot|\cdot)$.

Let $\sb_-=\sn_-\+\sh$.
Then there is an exact sequence
\begin{align}
 0\longrightarrow \sg^f\longhookrightarrow \sb_-
\overset{\ad f}\longrightarrow 
\sn_-
\longrightarrow 0.
\label{eq:2006-08-03-13-31}
\end{align}

Let 
$\sroots$ be the set of roots of $\sg$.
We have the decomposition 
$\sroots=\bigsqcup_{j\in \Z}\sroots_{j}$,
where 
$\sroots_{j}=\{\alpha\in \sroots; 
\bra \alpha,h_0\ket=2j\}$.
The set
\begin{align}
\sPi\teigi \sroots_1
\end{align}
is a basis of $\sroots$,
$\sroots_+\teigi \sqcup_{j>0}\sroots_j$
is the set of positive roots
and 
$\sroots_+\teigi \sqcup_{j<0}\sroots_j$
is the set of negative roots.
Let $\sQ$
be the roots lattice,
$\sP$ the weight lattice,
$\sQ\che$ the coroot lattice,
$\sP\che$ the coweight lattice.
We write $\sW$ for the Weyl groups of $\sg$
and
$w_0$ for the longest element of $\sW$.

Let $\srho$ be the half  sum of
positive roots,
$\srho\che$ the half  sum of positive coroots.
For $\alpha\in \sproots$,
the number
$\bra \alpha,\srho\che\ket$
is called the
{\em height} of
$\alpha$ and denote by
$\height \alpha$.
We have
\begin{align}
 h_0=2\srho\che,
\end{align}and so $\sroots_j=\{\alpha\in \sroots_+
;  \height \alpha=j\}$
for $j>0$.

Fix an anti-Lie algebra involution
\begin{align}
\sg\ni X\mapsto X^t \in \sg
\label{eq:2006-10-10-10-14}
\end{align}
such that
$e^t=f$, $f^t=e$
and $h^t=h$ for all $h\in \sh$.
Choose  root vectors $\{J_{\alpha}\}$ such that
$J_{\alpha}^t=J_{-\alpha}$ and $(J_{\alpha},J_{-\alpha})=1$.
Let
$\{J_i;  i\in \sI\}$,
where 
\begin{align}
\sI=\{1,2,\dots,l\},
\end{align}
 be a basis of $\sh$.
Then
the set
$\{J_a;  a\in \sI\sqcup \sroots\}$
forms a basis of $\sg$.
Let $c_{a,b}^d$ be the
structure constant with respect to
this basis:
$[J_a,J_b]=\sum_{d}c_{a,b}^d J_d$.
We have
$c_{\alpha,\beta}^{\gamma}
=-c_{-\alpha,-\beta}^{-\gamma}$
for $\alpha,\beta,\gamma\in \sproots$.

Let $d_1,d_2,\dots ,d_l$
be the exponents of $\sg$.
There exists a basis
 $\{P_i;  i\in \sI\}$ of $\sg^f$
such that 
that is,
\begin{align}
P_i\in \sg_{{-d_i}},\quad \text{or equivalently,}\quad
[\srho\che,P_i]=-d_i P_i.
\end{align}
By \eqref{eq:2006-08-03-13-31}
for $\alpha\in \sproots$
there exists an element $I_{-\alpha} \in \sb_-$
such that
\begin{align}
 [f, I_{-\alpha}]=J_{-\alpha}.
\end{align}
The set
$\{P_i,\ 
I_{-\alpha};
i\in \sI,\ 
\alpha\in \sproots\}$
form a basis of $\sb_-$.
\subsection{Kostant's Theorems}
Define an element
$\schi_+\in \dual{\sn}_+$
by
\begin{align}
 \schi_+(x)=(f|x)\quad \text{ for $x\in \sn_+$}.
\label{eq;2006-10-14-03-27}
\end{align}
Then
$\schi_+$
is a character of $\sn_+$,
that is,
$\schi_+([\sn_+,\sn_+])=0$.
The character $\chi_+$
is non-degenerate
 in the sense of \cite{Kostant-Wh}:
%that is,
%$\schi_+$ is a character of $\sn_+$ 
%(i.e.$\schi_+([\sn_+,\sn_+])=0$)
%such that
 $\schi_+(J_{\alpha})\ne 0$
for all $\alpha\in \sPi$.
Set
\begin{align*}
 \dual{\schi}_+=-\schi_+.
\end{align*}
Then
 $\schi_+^*$
also defines a non-degenerate character of $\sn_+$.
Let
$\ker \schi_+^*=\Ker (\schi_+^*:
U(\sn_+)\ra \C)$
and set
\begin{align}
 \C_{\schi_+^*}=U(\sn_+)/\ker \schi_+^*.
\end{align}
This defines a one-dimensional representation of $\sn_+$.

Let $M$ be a $\sg$-module.
Define 
\begin{align}
 &\Wh (M)\teigi \{m\in M;  x m=\dual{\schi}_+(x)m\ 
\text{ for all $x\in \sn_+$}
\},\\
&\Wh^{\gen}(M)\teigi \{m\in M;  
(\ker\dual{\schi}_+)^r m=0\ \text{for }r\gg 0\}.
\end{align}
Then
$\Wh(M)\subset \Wh^{\gen}(M)$
and $\Wh^{\gen}(M)$ is a $\sg$-submodule of $M$.
An element of $\Wh(M)$ is called a {\em Whittaker vector}.

Let $H^{\bullet}(\sn_+, M)$ be the Lie algebra 
cohomology
of $\sn_+$ with coefficient in a $\sn_+$-module $M$.
By definition we have
\begin{align}
 \Wh(M)=H^0(\sn_+,M \* \C_{\schi_+}),
\end{align}
where 
%$\C_{{\schi}_+}$ is the one-dimensional $\sn_+$-module
%defined by the character ${\schi}_+$
%and
$\sn_+$ acts on $V\* \C_{\schi_+}$ by the tensor product action.

Define a $\sg$-module $\GG$
by
\begin{align}
 \GG\teigi U(\sg)\*_{U(\sn_+)}\C_{\dual{\schi}_+}.
\label{eq:2006-11-04-13-50}
\end{align}
Then $\GG=\Wh^{\gen}(\GG)$.
By the Frobenius reciprocity,  
we have an  isomorphism
\begin{align*}
  \Wh(M)\cong \Hom_{\g}(\GG,M)
\quad
\text{ for any $\sg$-module $M$.}
\end{align*}
In particular
$  \Wh(\GG)\cong\End_{\g}(\GG)$.
We consider 
 $\Wh(\GG)$  as a $\C$-algebra
through the identification
\begin{align}
\Wh(\GG)=\End_{\g}(\GG)^{\op}.
\end{align}

Let $\Center(\sg)$ be the center of the universal enveloping algebra
$U(\sg)$ of $\sg$.
The following theorem is well-known.
\begin{Th}[B. Kostant \cite{Kostant-Wh}]\label{Th:Kostant-WH}
The correspondence 
\begin{align*}
 \begin{array}{ccc}
   \Center(\sg)&\longrightarrow & \Wh(\GG)=H^0(\sn_+,\GG\* \C_{\schi_+}) 
\\
z&\longmapsto &z\*1
 \end{array}
\end{align*}
gives an algebra isomorphism.
\end{Th}

Let $\widetilde{\Cat}$
be the Serre full subcategory of the category of $\sg$-modules consisting
of objects $M$
such that $M=\Wh^{\gen}(M)$.
%Then $\GG$ belongs to $ \widetilde{\Cat}$.
\begin{Th}[{B. Kostant \cite{Kostant-Wh},
see also \cite[Appendix  by S.\ Skryabin]{Premet}}]\label{Th:Kostant-cat-eq}
$ $

\begin{enumerate}
 \item  Let $M$ be any object of $\wt{\Cat}$.
Then
  $H^i(\sn_+,V\* \C_{\schi_+})=\zero $
for all $i\ne 0$.
\item The functor  $M\mapsto \Wh(M)$
gives an equivalence of 
the category $\widetilde{\Cat}$ and the category of
	left $\Center(\sg)$-modules.
The inverse functor is given by the functor 
$E\mapsto \GG\*_{\Center(\sg)}E$.
\end{enumerate}
 \end{Th}
Let $\Cat$ be the full subcategory of $\widetilde{\Cat}$
consisting of objects $V$
satisfying
(1) $M$ is finitely generated over $U(\sg)$,
and (2) $M$ is locally finite over $\Center(\sg)$.
Let  
$\Fin\Center(\sg)$
be the category of finite-dimensional
$\Center(\sg)$-modules.
Then,
$\GG\*_{\Center(\sg)}E\in \Cat$
for $E\in \Fin \Center(\sg)$.
Therefore,
by Theorem \ref{Th:Kostant-cat-eq},
we have an equivalence of categories
\begin{align}
   \Cat\cong  \Fin \Center(\sg).
\end{align}

\subsection{Superspaces and Superalgebras}
A {\em superspace} 
 is a $\Z_2$-graded 
vector space $V=V_{\bar 0}\+ V_{\bar {1}}$.
We set $p(v)=\bar {j}$ if $v\in V_{\bar{j}}$.
The $p(v)$ is called the {\em parity} of $v$.
We use the convention that
if we write $p(v)$ then
the vector $v$ is assume to be 
 homogeneous.
A vector $v$ is called {\em even}
if $p(v)=\bar {0}$,
and
{\em odd} if $p(v)=\bar {1}$.
Also, in this paper 
we write
$V^{\even}$ for $V_{\bar 0}$
and $V^{\odd}$ for $V_{\bar 1}$.
A {\em supersubspace}
 $W$ of  $V$
is a $\Z_2$-graded subspace of $V$:
$W=W^{\even}\+ W^{\odd}$,
$W^{\even}=W\cap V^{\even}$,
$W^{\odd}=W\cap V^{\odd}$.

A {\em superalgebra}
 is a $\Z_2$-graded algebra $A=A^{\even}\+ A^{\odd}$.
For $a,b\in A$ we set
\begin{align*}
 [a,b]=ab-(-1)^{p(a)p(b)}ba.
\end{align*}
Any superalgebra can be viewed as a  
{\em Lie superalgebra}
with respect to this bracket.

The tensor product 
$A\* B$
of two superalgebras
$A$ and $B$
is considered as a superalgebra
by 
$p(a\* 1)=p(a)$,
$p(1\* b)=p(b)$,
\begin{align*}
(a_1\* b_1)\cdot (a_2\* b_2)=
(-1)^{p(b_1)p(a_2)}(a_1 a_2)\* (b_1 b_2).
\end{align*}If no confusion can arise
we often omit the tensor symbol for elements,
e.g. we write $a b$ for $a\* b$.

\subsection{The BRST Construction  of the Center
$\Center(\sg)$}
\label{subsection:fd-BSRT}
The form
$(\cdot|\cdot)$ of $\sg$
restricts to a
non-degenerate  symmetric bilinear form on 
$\sn_-\+ \sn_+$.
Let $\sCl$ be the associated Clifford algebra.
The
 $\sCl$
may be defined as the superalgebra
with % the following generators and relations:
\begin{align*}
 &\text{odd generators: $\psi_{\alpha}$ \quad with $\alpha\in \sroots$},\\
&\text{relations: $[\psi_{\alpha},\psi_{\beta}]=\delta_{\alpha+\beta,0}$
\quad
with $\alpha,\beta\in \sroots$}.
\end{align*}
Here $\psi_{\alpha}$ is regarded as the element of $\sCl$
corresponding to the root vector $J_{\alpha}\in \sg$.
The algebra $\sCl$
contains the Grassmann algebra
$\Lam(\sn_{\pm})$
of $\sn_{\pm}$ as its subalgebra:
\begin{align}
\Lam(\sn_{\pm})=\bra \psi_{\alpha};  \alpha\in \sroots_{\pm}\ket
\subset \sCl.
\end{align}
Also, there is an isomorphism of linear spaces
\begin{align*}
 \sCl\cong \Lam (\sn_+)\* \Lam(\sn_-).
\end{align*}
The space
$\Lam(\sn_+)$ and $\Lam(\sn_-)$ can be  considered as 
 $\sCl$-modules 
through the identification
\begin{align*}
 \Lam(\sn_+)=\sCl/(\sCl\cdot  \sn_-)
,\quad
\Lam(\sn_-)=\sCl/(\sCl\cdot \sn_+).
\end{align*}
These are irreducible $\sCl$-modules.

Define a superalgebra $\sCg$
by
\begin{align}
 \sCg\teigi U(\sg)\* \sCl,
\end{align}
where $U(\sg)$ is considered as a 
purely even superalgebra.
Let
$\sd_+\in \sCg$
be an odd element defined 
by
\begin{align}
\sd_+&=\sd^{\st}_++\schi_+
 , \label{eq:2006-08-31-14-30}
\\
&\text{where}\quad
\sd^{\st}_+\teigi \sum_{\alpha\in
\sproots}
J_{\alpha}{\psi}_{-\alpha}
-\frac{1}{2}\sum_{{\alpha,\beta,\gamma\in
\sproots}{ }}
c_{\alpha,\beta}^{\gamma}
{\psi}_{-\alpha}{\psi}_{-\beta}
\psi_{\gamma}.
\label{eq:2006-10-07-23-29}
\end{align}
Here
the character $\schi\in \dual{\sn}_+$
is   considered as an element of %identified with the element
$\sCl$, %via the isomorphism $\dual{\sn}_+\cong \sn_-\subset \sCl$,
that is,
\begin{align*}
 \schi_+
=\sum_{\alpha\in \sPi}\schi_+(J_{\alpha})\psi_{-\alpha}.
\end{align*}
The following assertion can be checked by a direct calculation. 
\begin{Pro}\label{Pro:squre=0-f-d}
We have
\begin{align*}
 [\sd_+^{\st},\sd_+^{\st}]=[\schi_+,\schi_+]=
[\sd_+^{\st},\schi_+]=0\quad \text{ in $\sCg$}.
\end{align*}
In particular,
$(\sd_+^{\st})^2=0$,
$\schi_+^2=0$
and
  $\sd_+^2=0$. 
\end{Pro}
By Proposition \ref{Pro:squre=0-f-d}
we have
\begin{align}\label{eq:fd-ad-q-2-0}
 \text{$(\ad\sd_+)^2=0$\quad
in 
$\End \sCg$,}
\end{align}
where $\ad \sd_+(c)=[\sd_+, c]$.
This follows from the relation
\begin{align}
 [\sd_+,[\sd_+,c]]=[[\sd_+,\sd_+],c]-[\sd_+,[\sd_+,c]]
\quad
\text{for $c\in \sCg$.}
\label{eq:2006-10-12-13-49}
\end{align}
Let
\begin{align}
 \text{$\deg \psi_{\alpha}=-1$,
$\deg {\psi}_{-\alpha}=1$ for $\alpha\in 
\sproots$
and $\deg u=0$ for $u\in U(\sg)$.
}
\label{eq:def-of-grading-cl}
\end{align}
Then this defines 
 a $\Z$-grading  of $\sCg$:
\begin{align}
 \sCg=
\bigoplus_{n\in \Z}\sC^n(\sg),\quad 
\sC^m(\sg)\cdot \sC^n(\sg)\subset 
\sC^{m+n}(\sg).
\end{align}
We  have
 $\deg \sd_+=1$.
Thus,
by \eqref{eq:fd-ad-q-2-0},
$(\sCg,\ad \sd_+)%=(U(\sg)\*\sCl,\ad \sd_+)
$
is a cochain complex.
The corresponding cohomology is denoted by
$H^{\bullet}(\sCg)$:
\begin{align}
H^{\bullet}(\sCg)%= H^{\bullet}(U(\sg)\*\sC,\ad \sd_+)
=\bigoplus_{i\in \Z} H^{i}(\sCg).
\end{align}
The space
$H^{\bullet}(\sC(\sg))$
inherits the $\Z$-graded superalgebra structure
from $\sCg$.
 \begin{Th}[B. Kostant and S. Sternberg \cite{KS}]\label{Th:zhu1}
$ $

\begin{enumerate}
\item The cohomology $H^i(\sCg)$ is zero 
for all $i\ne 0$.
 \item 
The correspondence
\begin{align*}
\begin{array}
{ccc}
\Center(\sg)&\rightarrow &H^0(\sCg)\\
z&\mapsto &z \*1
\end{array}
%\label{map:iso-th-center}
\end{align*}
gives an algebra isomorphism.
\end{enumerate}
 \end{Th}
\begin{proof}
We identify
$\Lam (\sn_-)$ with $\Lam (\dual{\sn_+})$
via $(\cdot|\cdot )$
and
set $\dual{\psi}_{\alpha}=\psi_{-\alpha}
\in \Lam(\sn_+^*)$
for $\alpha\in \sroots_+$.
For convenience   we put
$\sC=\sCg=U(\sg)\*\sCl$,
$\sC^n=\sC^n(\sg)$.
Then 
\begin{align}
\sC^n=\sum_{i-j
=n} U(\sg)\* \Lam^i (\dual{\sn_+})
\*
\Lam^j({\sn_+}).
\end{align}
 The following
assertion can be checked  
by a  direct calculation.
\begin{Claim}
We have
\begin{align}\label{eq:dec_bar_d}
\ad \sd_+=\bar{d}_++\bar{d}_-
\text{ }
\end{align}
on $\bar{C}$,
where
$\bar{d}_{\pm}\in \End \sCg$
is
defined by 
\begin{equation}\label{eq:bar-d-+}
\begin{aligned}
\bar{d}_+(u\*\omega_1\*
\omega_2)
=\sum_{\alpha
\in \sroots_+}
\left\{
(\ad J_{\alpha}
(u))
\* \dual{\psi}_{\alpha}\omega_1\*
\omega_2
+u\*\dual{\psi}_{\alpha}\omega_1\*
\ad J_{\alpha}(\omega_2)
\right\}\\
-\frac{1}{2}\sum_{\alpha,
\beta,\gamma\in \sroots_+}
c_{\alpha,\beta}^{\gamma} u\*
\dual{\psi}_{\alpha}\dual{\psi}_{\beta}
(\ad \psi_{\gamma}(\omega_1))
\* \omega_2,
\end{aligned}
\end{equation}
\begin{equation}\label{eq:bar-d--}
\begin{aligned}
(-1)^{i}\bar{d}_-(u\*\omega_1\*
\omega_2)
=\sum_{\alpha\in \sroots_+}
u(J_{\alpha}+\schi(J_{\alpha}))\*
\omega_1\*
\ad \dual{\psi}_{\alpha}
(\omega_2)
\\
-\frac{1}{2}\sum_{\alpha,
\beta,\gamma\in \sroots_+}
c_{\alpha,\beta}^{\gamma} u\*
\omega_1\*\psi_{\gamma}\ad \dual{\psi}_{\beta}
(\ad 
\dual{\psi}_{\alpha}
(\omega_2))
\end{aligned}
\end{equation}
for $u\in U(\sg)$,
$\omega_1\in \Lam^i(\sn_+^*)$,
$\omega_2\in \Lam(\sn_+)$.
\end{Claim}
Note that
\begin{equation}\label{eq:bar-d+-d-}
\begin{aligned}
\bar{d}_+(U(\sg)\* \Lam^i(\dual{\sn_+})\*
\Lam^{j}(\sn_+))
\subset
U(\sg)\* \Lam^{i+1}(\dual{\sn_+})\*
\Lam^{j}(\sn_+),
\\
\bar{d}_-(U(\sg)\* \Lam^i(\dual{\sn_+})\*
\Lam^{j}(\sn_+))
\subset
U(\sg)\* \Lam^{i}(\dual{\sn_+})\*
\Lam^{j-1}(\sn_+).
\end{aligned}
\end{equation}
Thus  we have
\begin{align}\label{ew:commutativity-dif-bar}
\bar{d}_+^2=\bar{d}_-^2=[\bar{d}_+,\bar{d}_-]=0.
\end{align}

Define
\begin{align}\label{eq:filtration-Kostant}
F^p \bar{C}\teigi \sum_{i\geq p}U(\sg)\* \Lam^i(\dual{\sn_+})\*
\Lam^{\bullet}(\sn_+)\subset \bar{C}.
\end{align}
Then
\begin{align}
&\bar{C}=F^0\bar{C}
\supset F^1\bar{C}\supset \dots
\supset F^{\dim \sn_+ +1}\bar{C}=\zero ,
\\
%&\bigcap_p F^p \bar{C}=\zero ,\\
&
\bar{d}_+ F^p \bar{C}
\subset F^{p+1} \bar{C},\quad
\bar{d}_- F^p \bar{C}\subset F^{p} \bar{C}.
\end{align}
Hence there is a  corresponding converging
spectral sequence
$E_r\Rightarrow H^{\bullet}(\sC)$.
By definition and 
\eqref{eq:bar-d--}, we have
\begin{align*}
E_1^{p,q}
=H^{p+q}(F^p \bar{C}/F^{p+1}\bar{C},\ad \sd)
=H^{q}(U(\sg)\* \Lam^p(\sn_+^*)\* \Lam^{\bullet}(\sn_+)
,\bar{d}_-).
\end{align*}

The formula \eqref{eq:bar-d--}
shows that
the cochain complex
$(\bar{C},\bar{d}_-)$
is identical to the Chevalley complex
for calculating the $\sn_+$-homology
$H_{\bullet}
\left(\sn_+,\left(U(\sg)\* \C_{\schi_+}\right)\*\Lam (\dual{\sn}_+)\right)
$
%\begin{align*}
%H_{\bullet}
%\left(\sn_+,\left(U(\sg)\* \C_{\schi_+}\right)\*\Lam (\dual{\sn}_+)\right)
%=
%H_{\bullet}\left(\sn_+,\left(U(\sg)\* \C_{\schi_+}\right)\right)
%\* \Lam (\dual{\sn}_+)
%\end{align*}
%with 
%coefficients in 
%$\left(U(\sg)\* \C_{\schi_+}\right)\*\Lam (\dual{\sn}_+)$
(equipped with the opposite homological gradation).
Here
$\left(U(\sg)\* \C_{\schi_+}\right)\*
\Lam (\dual{\sn}_+)$
is regarded as a right $U(\sn_+)$-module
on which
  $U(\sn_+)$
acts only on the first factor
$U(\sg)\* \C_{\schi}$.

Obviously
the right $\sn_+$-module
$U(\sg)$
is free over $\sn_+$,
and thus, so is $U(\sg)\* \C_{\schi_+}$.
Therefore
\begin{equation}
 \begin{aligned}
E_1^{p,q}&=\begin{cases}
\left(\left(U(\sg)\* \C_{\schi_+}\right)
/\left(U(\sg)\* \C_{\schi_+}\right)\sn_+
\right)
\*\Lam^p(\dual{\sn}_+)
&\text{for $q=0$,}\\
0&\text{for $q\ne 0$}
\end{cases}\\
&=\begin{cases}
\left(U(\sg)\*_{U(\sn_+)}\C_{\dual{\schi}_+}
\right)
\*\Lam ^p(\dual{\sn}_+)
&\text{for $q=0$,}\\
0&\text{for $q\ne 0$}.
\end{cases}
\end{aligned}\label{eq:E_1-Kostant}
\end{equation}

Next we calculate
the term $E_2$.
From
\eqref{eq:2006-11-04-13-50},
 \eqref{eq:bar-d-+}
and \eqref{eq:E_1-Kostant}
it follows that
\begin{align*}
E_2^{p,q}=\begin{cases}
	   H^{p}(\sn_+,\GG\* \C_{\schi_+})&\text{for $q=0$},\\
0&\text{for $q\ne 0$.}
	  \end{cases}
\end{align*}
But then
 by Theorem \ref{Th:Kostant-cat-eq}
(i)
we have
$E_2^{p,q}=0$ unless $(p,q)=(0,0)$
and
our spectral sequence
collapses at $E_2=E_{\infty}$.
This proves (i).
Because the 
term $E_{\infty}^{p,q}$
lies entirely in $E_2^{0,0}$
the corresponding filtration
of $H^{0}(\sCg)$ is trivial:
$H^0(\sCg)\cong \Wh(\GG)$.
Thus
Theorem \ref{Th:Kostant-WH}
proves (ii).
\end{proof}
Let 
\begin{align}
\sCl\ni a\mapsto a^t\in \sCl
\label{eq:2006-10-10-10-15}
\end{align}
be the anti-superalgebra automorphism 
defined by
$\psi_{\alpha}^t=\psi_{-\alpha}$
for  $\alpha\in \sroots$.
This induces the anti-superalgebra automorphism
\begin{align}\label{eq:acti-auto-of-sC(sg)}
 \sCg\ni c\mapsto c^t \in \sCg
\end{align}
defined by $(u\* a)^t=u^t\* a^t$
for $u\in U(\sg)$ and $a\in \sCl$.

Set
\begin{align}
 &\sd^{\st}_-\teigi 
(\sd^{\st}_+)^t=
\sum_{\alpha\in
\sroots_-}
J_{\alpha}{\psi}_{-\alpha}
-\frac{1}{2}\sum_{{\alpha,\beta,\gamma\in
\sroots_-}{ }}
c_{\alpha,\beta}^{\gamma}
{\psi}_{-\alpha}{\psi}_{-\beta}
\psi_{\gamma},
\label{eq:2006-10-07-23-30}
\\
& \schi_-:=
(\schi_+)^t=\sum_{\alpha\in \sroots_-}\schi_-(J_{\alpha})\psi_{-\alpha},
\\ &
\sd_-:=
(\sd_+)^t=\sd^{\st}_-+ \schi_-.
\end{align}
The element 
$\schi_-$ can be considered as a character of $\sn_-$
such that
\begin{align}
 \text{$\schi_-(J_{\alpha})=(e|J_{\alpha})$ for $\alpha\in \sroots_-$.}
\label{eq:2006-11-01-14-39}
\end{align}
It is clear that
$(\sd^{\st}_-)^2=(\schi_-)^2=
\sd_-^2
=0$,
$(\ad \sd_-)^2=0$
and
 $\ad \sd_-\cdot \sC^i(\sg)\subset \sC^{i-1}(\sg) $.
Thus
$(\sCg, \bd_-)$ is a chain complex.
Put
\begin{align*}
 H_{\bullet}(\sCg)\teigi H_{\bullet}(\sCg,\sd_-).
\end{align*}
The anti-automorphism \eqref{eq:acti-auto-of-sC(sg)}
induces an isomorphism $H^i(\sCg)\cong H_{-i}(\sCg)$
for all $i\in \Z$.
Therefore by Theorem \ref{Th:zhu1}
we have
$H_i(\sCg)=0$
for all $i\ne 0$ and the following algebra isomorphism:
\begin{align}
\begin{array}
{ccc}
\Center(\sg)&\isomap &H_0(\sCg)\\
z&\mapsto &z \*1
\end{array}
\end{align}
\subsection{Whittaker Functor}
\label{subsection:WhittakerFunctor}
Let
$\sBGG$
be the Bernstein-Gelfand-Gelfand category 
\cite{BGG}
of $\sg$,
that is,
a full subcategory
of the category of left $\sg$-modules
consisting
of objects $M$
such that the following hold:
\begin{itemize}
 \item  $M$ is finitely generated over $\sg$;
\item $\sh$ acts semisimply on $M$;
\item 
$\sn_+$
acts locally nilpotently on $M$.
\end{itemize}
Let $\sM(\slam)\in \sBGG$
be the Verma module of $\sg$ with highest weight $\bar{\lam}\in \dual{\sh}$,
 $\sL(\slam)\in \sBGG$ its unique simple
quotient.
Then
every simple object of 
$\sBGG$ is isomorphic to exactly one of 
the $\sL(\slam)$ with $\slam\in \dual{\sh}$.
Every object of $\sBGG$ has finite length.
It is known that
$\sM(\slam)=\sL(\slam)$
if and only if $\slam$ is {\em anti-dominant},
that is,
\begin{align}
 \text{$\bra \slam+\srho,\salpha\che\ket\not\in \{1,2,3,\dots\}$
for all $\salpha\in \sroots_+$.
}
\end{align}
Let $\BGG^{\tri}(\sg)$ be  the full subcategory of
$\sBGG$ consisting of objects
$M$ that admit a {\em Verma flag},
that is,
a finite filtration 
$M=M_0\supset M_1\supset \dots \supset M_r=\zero $
such that each successive subquotient $M_i/M_{i+1}$
is isomorphic to
some Verma module $\sM(\slam_i)$ with  $\slam_i\in \dual{\sh}$.
Let $\sP(\slam)$
be the projective cover of $\sL(\slam)$.
It is known that
$\sP(\slam)\in \BGG^{\tri}(\sg)$
for any $\slam\in \dual{\sh}$.

For a $\sg$-module $M$
set
\begin{align}\label{eq;def;sC}
 \sC(M)\teigi M\* \Lam(\sn_-).
\end{align}
The space
$ \sC(M)$
can be viewed naturally as a module
over the superalgebra
$\sCg=U(\sg)\* \sCl$.
Let $\sC_i(M)=M\* \Lam^i(\sn_-)$,
so that $\sC(M)=\bigoplus_{i\geq 0}\sC_i(M)$.
Then,
\begin{align*}
 \sd_- \cdot \sC_i(M)\subset \sC_{i-1}(M).
\end{align*}

One sees that the chain complex
$(\sC(M),\sd_-)$ is identically the Chevalley  complex for 
the Lie algebra homology
$H_{\bullet}(\sn_-,M\* \C_{\schi_-})$.
Put
\begin{align}\label{eq:def-of-sH}
 \sH_i(M)\teigi 
H_i(\sn_-,M\* \C_{\schi_-})=
H_i(\sC(M),\sd_-).
\end{align}
Note that the center
 $\Center(\sg)=H_0(U(\sg)\* \sCl,\ad \sd_-)$\
acts
on $\sH_i(M)$
naturally:
\begin{align*}
 \text{$[c]\cdot [m]=
[c\cdot  m]$\quad  for $c\in U(\sg)\* \sCl$, $v\in \sC_i(M)$
such that $[\sd_+,c]=\sd_+ m=0$.}
\end{align*}

By definition
we have
\begin{align}
 \sH_0(M)=M/(\ker\schi_-^* )M.
\end{align}
Here,
 $\schi_-^*=-\schi_-$
and
 $\ker\schi_-^*=\Ker
(\schi_-^*: U(\sn_-)\rightarrow \C)$.

\begin{Lem}\label{Lem:fini-dim-proj-fin}
 Let $M$ be any object of $\BGG^{\tri}(\sg)$.

\begin{enumerate}
 \item The homology $\sH_i(M)$
is zero  for all $i>0$.
\item
The space $M/(\ker\schi_-^*)^r M$
is finite-dimensional for all $r\in \Z_{\geq 1}$.
In particular
$\sH_0(M)$ is finite-dimensional.
\end{enumerate}
\end{Lem}
\begin{proof}
(i)
Since
an object 
$M$ of $\BGG^{\tri}(\sg)$
is free 
of 
finite rank over $U(\sn_-)$,
we have
\begin{align}\label{04-12-22-Jeri}
 \text{$M\* \C_{\schi_-}\cong M$ \quad as $\sn_-$-modules.}
\end{align}
Hence, $\sH_i(M)\cong H_i(\sn_-,M)$.
This gives
$\sH_i(M)=\zero $ for $i>0$
because $M$ is free over $\sn_-$.
(ii)
 By \eqref{04-12-22-Jeri}
we have
\begin{align}\label{04-12-22-Jeri-2}
 M/(\ker\dual{\schi}_-)^r M\cong M/ (\sn_-)^r M.
\end{align}
Hence the assertion follows. % from Lemma \ref{Lem:O-tri}.
\end{proof}
\begin{Lem}\label{Lem:f-d-H-0}
 Let $M$ be any object of  $\sBGG$.
Then the space $M/(\ker\schi_-^*)^r M$
is finite-dimensional for all $r\in \Z_{\geq 1}$.
In particular
$\sH_0(M)$ is finite-dimensional.
\end{Lem}
\begin{proof}
 The assertion follows from 
Lemma \ref{Lem:fini-dim-proj-fin}
because $\sBGG$ has enough projectives.
\end{proof}
Let  $M\che$
the full dual space of $M$:
\begin{align*}
 M\che\teigi \Hom_{\C}(M,\C).
\end{align*}
We regard $M\che$ as a $\sg$-module 
on which $\sg$ acts by
$(xf)(v)=f(x^t v)$ for $f\in M\che, \ v\in M$.
Then by definition
 we have
\begin{align}\label{eq:dual-Kostant}
 \Wh(M\che)=\Hom_{\C}(M/(\ker\schi_-^*)M,\C)
=\Hom_{\C}(\sH_0(M),\C). %\quad (M\in \sBGG).
\end{align}
Hence,
by Lemma \ref{Lem:f-d-H-0},
$\Wh(M\che)$ is finite-dimensional
for any object $M$ of $\sBGG$.
\begin{Rem}
By Theorem \ref{Th:Kostant-cat-eq} (ii)
we have
$\Wh^{\gen}(M\che)=\GG\*_{\Center(\sg)}\Wh(M^{\che})$
and
 $\Wh^{\gen}(M\che)$ belongs  to $\Cat$.
\end{Rem}
\begin{Th}[B.  Kostant \cite{Kostant-Wh}]\label{Th:Kostant-exact}
The cofunctor 
%defined by
\begin{align*}
 \begin{array}{ccc}
  \sBGG&\rightarrow  &\Fin \Center(\sg) \\
M&\mapsto &\Wh(M\che)
 \end{array}
\end{align*}
 is  exact.
\end{Th}
\begin{Th}\label{Th:excatness-H_0-bar}
The functor %defined by
\begin{align*}
 \begin{array}{ccc}
  \sBGG&\rightarrow  &\Fin \Center(\sg) \\
M&\mapsto &\sH_0(M)
 \end{array}
\end{align*}
 is  exact. 
\end{Th}
\begin{proof}
Let $M\in \BGG(\sg)$.
By  Lemma \ref{Lem:f-d-H-0}
and
 \eqref{eq:dual-Kostant},
$\sH_0(M)=\Hom_{\C}(\Wh(M\che), \C)$.
Hence 
by
 Theorem \ref{Th:Kostant-exact}
and Lemma \ref{Lem:f-d-H-0},
$\sH_0(?)$ is exact
because it  is 
 the composition of two
exact functors
$M\mapsto \Wh(M\che)$
and $M\mapsto \Hom_{\C}(M,\C)$.
 \end{proof}
\begin{Th}\label{Th:vanishing-finite-dim}
Let $M$ be any object of $\sBGG$.
Then
the homology 
$\sH_i(M) $ is zero
for all $i\ne 0$.
\end{Th}
\begin{proof}
 We prove the assertion
by induction on $i\geq 1$.
Because  the category  $\sBGG$
has enough projectives,
there exists a projective
object $P\in \sBGG$ 
and an exact sequence
\begin{align}
 \text{$
 0\rightarrow N\rightarrow P\rightarrow M\rightarrow 0
$
}
\end{align}in $\sBGG$.
Consider the corresponding long exact sequence
\begin{align}
%\begin{aligned}
  \dots \rightarrow 
&\sH_{i}(P)\rightarrow
\sH_{i}(M)
\rightarrow \sH_{i-1}(N)\rightarrow \cdots
\nno
\\
  \dots \rightarrow 
&\sH_1(P)
\rightarrow
\sH_1(M)
\rightarrow 
\sH_0(N)\rightarrow \sH_0(P)\rightarrow \sH_0(M)\rightarrow 0.
%\end{aligned}
\label{eq:long-exact-fd}
\end{align}
By Lemma \ref{Lem:fini-dim-proj-fin},
 we have \begin{align}\label{eq:in-proof:vanishing-fin-pro}
	  \sH_i(P)=\zero \text{  for $i>0$}
	 \end{align}
because $P$ is an object of $\BGG^{\Delta}(\sg)$.
Hence,
from Theorem \ref{Th:excatness-H_0-bar}
it follows that
$\sH_1(M)=\zero $.
Let $i\geq 2$.
Then from
\eqref{eq:long-exact-fd}
and \eqref{eq:in-proof:vanishing-fin-pro}
we see that  the vanishing of 
$\sH_{i-1}(N)$
for any object  $N$ of $\sBGG$
implies 
the vanishing of 
$\sH_{i}(M)$.
This completes the proof.
\end{proof}
Let 
\begin{align}
\gamma: \Center(\sg)\isomap S(\sh)^{\sW}
\label{eq:Harish-Chandra-cl}
\end{align}
be the Harish-Chandra
isomorphism with respect to the triangular decomposition
$\sg=\sn_-\+ \sh\+\sn_+$.
For $\slam\in \dual{\sh}$
let
\begin{align}
 \gamma_{\slam}:=(\text{evaluation at }\slam+\srho)
\circ \gamma:\Center(\sg)\rightarrow \C.
\label{eq:2006-11-03-00-39}
\end{align}
Then $\gamma_{\slam}$
is
 the infinitesimal  character of $\sM(\slam)$,
that is,  $z m=\gamma_{\slam}(z) m$ for $z\in \Center(\sg)$ and $m\in 
\sM(\slam)$.
We have
\begin{align}
\gamma_{\bar{\lam}}=\gamma_{\sw\circ \slam} 
\quad\text{for all $\sw\in \sW$.}
\end{align}Here,
\begin{align}
\sw\circ \slam=\sw(\slam+\srho)-\srho.
\end{align}
Let 
\begin{align}
 \C_{\gamma_{\slam}}:=\Center(\sg)/\ker \gamma_{\slam}.
\end{align}
%Let $\C_{\gamma_{\slam}}\in \Fin \Center(\sg)$
%be the one-dimensional representation of $\Center(\sg)$
%defined by the central character $\gamma_{\slam}$.
\begin{Th}[{\cite[Corollary 3.4.6, Theorem 3.4.7]{Ma}}]%
\label{Th:image-of-simple-f-d}
Let  $\slam\in \dual{\sh}$.
There are  the 
following isomorphisms of 
$\Center(\sg)$-modules:

 \begin{enumerate}
  \item $\sH_0(\sM(\lam))\cong \C_{\gamma_{\bar{\lam}}}$;
\item $\sH_0(\sL(\slam))=\begin{cases}
			  \C_{\gamma_{\bar{\lam}}}&
\text{if $\slam$ is anti-dominant},\\
0&\text{otherwise};
			 \end{cases}$
  \item $\sH_0(\sM(\lam)^*)\cong \C_{\gamma_{\bar{\lam}}}$.
 \end{enumerate}
\end{Th}
\begin{proof}
 (i) The assertion  follows from the fact that
$\sM(\slam)\* \C_{\chi_-}\cong \sM(\slam)$
as $\sn_-$-modules. (See the proof of Lemma \ref{Lem:fini-dim-proj-fin} (i).)
(ii) Suppose that
$\slam$ is anti-dominant.
% $\bra \slam+\srho,\alpha\che\ket\not\in \Z_{\geq 1}$
%for all $\alpha\in \sroots_+$.
Then $\sL(\slam)=\sM(\slam)$,
and hence the assertion follows from the first assertion.
Next suppose that
$\slam$ is not anti-dominant.
Then 
there exist  exact sequences
\begin{align*}
 0\rightarrow \sM(\bar{\mu})
\rightarrow
\sM(\slam)\rightarrow
\sM(\slam)/ \sM(\bar{\mu})
\rightarrow 0\\
\text{and}\quad
\sM(\slam)/ \sM(\bar{\mu})\rightarrow \sL(\slam)\rightarrow 0
\end{align*}
with some $\bar{\mu}\in \sW\circ \slam$.
By applying the exact functor $\sH_0(?)$
we get the exact sequences
\begin{align*}
 0\rightarrow \sH_0(\sM(\bar{\mu}))
\rightarrow
\sH_0(\sM(\slam))
\rightarrow
\sH_0(\sM(\slam)/ \sM(\bar{\mu}))
\rightarrow 0\\
\text{and}\quad \sH_0(\sM(\slam)/ \sM(\bar{\mu}))
\rightarrow \sH_0(\sL(\slam))
\rightarrow 0.
\end{align*}
But $\sH_0(\sM(\bar{\mu}))
\cong
\sH_0(\sM(\slam))$ by the first assertion.
Therefore,
$\sH_0(\sM(\slam)/ \sM(\bar{\mu}))=0$,
and hence,
$\sH_0(\sL(\slam))=0$.
(iii)
Let $K_0(\sBGG)$ and $K_0(\Fin \Center(\sg))$
denote the Grothendieck groups of $\sBGG$
and $\Fin \Center(\sg)$,
respectively.
Then,
by Theorem \ref{Th:excatness-H_0-bar},
$\sH_0(?)$
defines a well-defined map from 
$K_0(\sBGG)$
 to $K_0(\Fin \Center(\sg))$.
Because $[M(\lam)^*]=[M(\lam)]$
in  $K_0(\sBGG)$,
we have $[\sH_0(M(\lam)^*)]=[\sH_0(M(\lam))]=[\C_{\gamma_{\bar{\lam}}}]$
in $K_0(\Fin \Center(\sg))$.
Because $\C_{\gamma_{\slam}}$ is simple,
this completes the proof.
\end{proof}
By Theorem \ref{Th:image-of-simple-f-d},
we see in particular that
  any simple object of $\Fin \Center(\sg)$
is isomorphic to $\sH_0(\sL(\slam))$ for some anti-dominant 
weight $\slam$.
%\begin{Rem}
%By Theorem \ref{Th:image-of-simple-f-d}
%  any simple object of $\Fin \Center(\sg)$
%is isomorphic to $\sH_0(\sL(\slam))$ for some anti-dominant 
%weight $\slam$.
%\end{Rem}
\subsection{Identification with Soergel's Functor}
For $\slam\in \dual{\sh}$,
let
$\sW(\slam)\teigi \bra s_{\salpha};  \bra \lam+\srho,
\salpha\che\ket \in \Z\ket\subset \sW$.
It is known that
$\sW(\slam)$ is a Coxeter group 
and it
is called
 the {\em integral Weyl group}
of $\slam\in \dual{\sh}$.
Let $\sBGG^{[\slam]}$
 be the block of $\sBGG$ corresponding to
$\slam$,
that is,
$\sBGG^{[\slam]}$ is the Serre full subcategory of
$\sBGG$
whose objects have all their local composition factors isomorphic to
$\sL(w\circ \slam)$ with $w\in \sW(\slam)$.
Then
we have
\begin{align*}
 \sBGG=\bigoplus\limits_{\slam\in \dual{\sh}
\atop \slam\text{ is anti-dominant}}\sBGG^{[\slam]}.
\end{align*}

Let $\slam$ be anti-dominant.
For an object $M$ of $\sBGG^{[\lam]}$,
we regard $\Hom_{\sg}(\sP(\slam), M)$ as a $\Center(\sg)$-module 
by $(zf)(v)=z f(v)$  for $z\in \Center(\sg)$, $f\in
\Hom_{\sg}(\sP(\slam),M)$ and
$v\in M$.
Then, both $\Hom_{\sg}(\sP(\slam),?)$ 
and $\sH_0(?)$ 
define functors 
 from $\sBGG^{[\lam]}$
to  $\Fin \Center(\sg)$.
 They are both exact functors.
Indeed,
$\Hom_{\sg}(\sP(\slam),?)$ is  exact because $\sP(\slam)$ is
projective
and $\sH_0(?)$ 
is exact by 
Theorem \ref{Th:excatness-H_0-bar}.

The functor $\Hom_{\sg}(\sP(\slam),?):\sBGG^{[\slam]}\rightarrow
\Fin \Center(\sg)$
was studied by W. Soergel \cite{Soe}.
By a result of E. Backelin \cite{Bac}
it follows that
the functor $\sH_0(?)$
coincides the Soergel's functor $\Hom_{\sg}(\sP(\slam),?)$.
Indeed
one can describe the corresponding natural transformation
in the following manner:

Let
 $\tilde{v}_{\lam}
\in \sP(\slam)$ be any inverse image of the
highest
weight vector of  $\sL(\lam)$ by the canonical homomorphism
$\sP(\slam)\twoheadrightarrow \sL(\slam)$.
Define a natural transformation
\begin{align}
 \Phi: \Hom_{\sg}(\sP(\slam),?)\rightarrow \sH_0(?)
\end{align}
%of functors from $\sBGG^{[\lam]}$ to  $\Fin \Center(\sg)$
by 
\begin{align}
\begin{array}{cccc}
  \Phi_{M}:& \Hom_{\sg}(\sP(\slam), M)&\rightarrow &\sH_0(M)=M
/(\ker\schi_-^*)\\
&f &\mapsto &[f(\tilde{v}_{\lam})]
\end{array}
\end{align}
with $M\in \Obj\sBGG^{[\slam]}$.
(It is clear that $\Phi_M$ is a homomorphism of $\Center(\sg)$-module.)

% describes the
%corresponding natural transformation
%between two functors.
%the functor $\sH_0(?)$
%coincides the Soergel's functor $\Hom_{\sg}(\sP(\slam),?)$.
\begin{Th}[E. Backelin \cite{Bac}]
\label{TH:Soergel}
For each anti-dominant  $\slam\in \dual{\sh}$,
$\Phi$ defines a natural isomorphism 
$ \Hom_{\sg}(\sP(\slam), ?)\cong \sH_0(?)$
of functors
 from $\sBGG^{[\lam]}$
to  $\Fin \Center(\sg)$.
\end{Th}
\begin{proof}
We have to show that
$\Phi_M$ is isomorphism for each object $M$ of $\sBGG^{[\slam]}$.
We prove this by induction on the length $l(M)$
of the composition series of $M$.
%(It is known that each object of $\sBGG$ has a finite length.)
Let $l(M)=1$,
that is, $M=L(\bar{\mu})$ for some $\bar{\mu}\in \sW(\lam)\circ \slam$.
If $\bar{\mu}\ne \slam$,
then both $\Hom_{\sg}(\sP(\slam),\sL(\bar{\mu}))$
and $\sH_0(\sL(\bar{\mu}))$ are zero
(see  Theorem \ref{Th:image-of-simple-f-d}).
On the other hand,
if  $\bar{\mu}=\slam$,
then
the map $\Phi_{\sL(\slam)}$ is clearly an isomorphism
by Theorem \ref{Th:image-of-simple-f-d}.
Next let $l(M)>1$.
Then we have an exact sequence
$0\rightarrow M_1\rightarrow M\rightarrow M_2\rightarrow 0$
in $\sBGG^{[\slam]}$
such that $l(M_1), l(M_2)<l(M)$.
This induces a commutative diagram
\begin{align*}
\minCDarrowwidth1pc
 \begin{CD}
  0 @>>>\Hom_{\sg}(\sP(\slam),M_1) 
@>>> \Hom_{\sg}(\sP(\slam),M) 
@>>> \Hom_{\sg}(\sP(\slam),M_2) 
@ >>>  0\\
 @. @ VV \Phi_{M_1} V @ VV \Phi_{M} V @ VV \Phi_{M_2} V \\
 0@>>>  \sH_0(M_1) @ >>> \sH_0(M) @>>> \sH_0(M_2)@>>> 0.
  \end{CD}
\end{align*}
The upper row is obviously exact
and the lower row is exact by 
Theorem \ref{Th:excatness-H_0-bar}.
Hence we are done by the induction hypothesis and the five-lemma.
\end{proof} 
%By the well-known
%result of W. Soergel \cite{Soe}
%the natural map
%$\Center(\sg)\rightarrow \End_{\sg}(\sP(\slam))$
%is surjective.
%Hence we have proved the following result.

%\begin{Rem}
% Backelin's proof (\cite{Bac})
%of the  fact that $\Hom_{\sg}(\sP(\slam),?)\cong \sH_0(?)$
% is based on his characterization of the Soergel's functor
%$\Hom_{\sg}(\sP(\slam),?)$.
%\end{Rem}
\section{Filtration of Vertex (Super)Algebras
and BRST Cohomology}
\label{section:FIltration-VA}
\subsection{Some Notation}
\label{subsection:Some-Notation}
Let $V$ be any superspace.
We set
\begin{align}
L(V):=V \* \C[t,t\inv].
\label{eq:2006-11-04-23-16}
\end{align}
This is considered as a superspace
such that $L(V)^{\even}=L(V^{\even})$
and $L(V)^{\odd}=L(V^{\odd})$.
For  any subring $R$ of $\C[t,t\inv]$
we consider $V\*R$ as a supersubspace
of $L(V)$.

The 
{\em symmetric algebra}
$S(V)$ of a superspace $V$ is 
the superalgebra
 \begin{align}
  S(V)=S(V^{\even})\* \Lam(V^{\odd}),
\label{eq:super-symmetric^alg}
 \end{align}
where 
$S(V^{\even})$ is the
(usual)  symmetric algebra of $U$
and 
$\Lam(V^{\odd})$ is the  Grassmann algebra 
of $V^{\odd}$.

A filtration $F=\{F_p V; p\in \Z\}$ of 
a superspace $V$ 
is a 
$\Z_2$-graded filtration of the vector space $V$:
$F_p V=(F_p V)^{\even}\+ (F_p V)^{\odd}$,
where 
$(F_p V)^{\even}=F_p V\cap V^{\even}$,
$(F_p V)^{\odd}=F_p V\cap V^{\odd}$.

A filtration  $F$ of $V$
is called {\em exhaustive}
if $V=\bigcup_p F_p V$
and {\em separated}
if $\bigcap_p F_p V=\zero $.

For an increasing filtration $\{F_p V\}$ of 
a superspace
$V$,
we set
 $\gr^F V=\bigoplus_p \gr^F_p V$,
where $\gr^F_p V=F_{p}V/F_{p-1}V$.
For any filtration $F$ we write
$\sigma_p$ for the symbol map
$F_p V\ra \gr^F_p V=F_p V/F_{p-1}$.

\smallskip

Below up until the end of Section 
\ref{section:FIltration-VA}
we treat superspaces,
supersubspaces,
Lie superalgebras,
vertex superalgebras etc.,
 unless otherwise stated.
However for convention  we  shall
often drop the prefix ``super''.

\subsection{Vertex Algebras}\label{subsection:VA-sub}
(see \cite{KacBook2, FB}.)
Let $\V$ be a 
vertex algebra, 
\begin{align}
Y(v,z)=\sum_{n\in \Z}v_{(n)}z^{-n-1}\in (\End \V)[[z,z\inv]]
\label{eq:2006-09-28-12-19}
\end{align}
 the 
 field associated
with $v\in \V$,
$|0\ket$ the vacuum so that $Y(|0\ket,z)=\id_{\V}$,
$\TT  $
the translation operator:
\begin{align}
 \label{eq:vertex-axiom0}
[\TT  ,Y(v,z)]=Y(\TT   v,z)=\frac{d}{d z}Y(v,z)
\quad \forall v\in \V.
\end{align}
We have
\begin{align}
  &[u_{(m)},v_{(n)}]=\sum_{i\geq 0}
\begin{pmatrix}
 m\\ i
\end{pmatrix}(u_{(i)} v)_{(m+n-i)},\label{eq:vertex-axiom1}
\\
&(u_{(m)}v)_{(n)}\nonumber\\
&=
\sum_{i\geq 0}
\begin{pmatrix}
 m\\i
\end{pmatrix}
(-1)^i(u_{(m-i)}v_{(n+i)}-
(-1)^{p(u)p(v)}
(-1)^m v_{(m+n-i)}u_{(i)})
\label{eq:vertex-axiom2}
\end{align}
for $u,v \in \V$.
The sum in the right hand side of 
\eqref{eq:vertex-axiom1}
is finite because the axiom of vertex algebras
requires  that $u_{(n)}v=0$ for sufficiently large $n$.
 The relations \eqref{eq:vertex-axiom2}
is called the {\em Borchards identity}.

A subspace  $U$ 
of $\V$
is said to
  {\em strongly generate} $\V$
if 
$\V$ is spanned by the vectors 
\begin{align*}
a^{j_1}_{(-n_1)}\dots a^{j_r}_{(-n_r)}|0\ket
\end{align*}
with 
$a^{j_s}\in U$,
$r\geq 0$,
$j_s\in J$,
$n_s\geq 1$.

An even element $H$ of $\End \V$
is called a {\em Hamiltonian}
if
$H$ acts semisimply on $H$
satisfying
\begin{align}
[H, Y(v,z)]=Y(Hv,z)+zY(T v,z)
\quad \forall v\in \V.
\label{eq:Hamiltonial}
\end{align}
If $v\in \V$ is an
eigenvector of $H$
then its eigenvalue
$\Delta$ is called
the {\em conformal weight} of $v$
and denoted by
$\Delta_v$.
We use the convention
that
when
we write $\Delta_v$
the vector  $v$ is assumed to be an eigenvector of $H$.
By
\eqref{eq:Hamiltonial}
one has
\begin{align}
&\Delta_{|0\ket}=0,\\
& \Delta_{u_{(n)}v}=\Delta_u+\Delta_v-n-1.
\end{align}
Set
\begin{align}
 \V_{-\Delta}=\{v\in \V; H\cdot v=\Delta v\},
\end{align}
so that
$\V=\bigoplus_{\Delta\in \C}\V_{-\Delta}$.
The 
$\V$ is called {\em $\Z$-graded}
 (by $H$)
if $\V_{-\Delta}=\zero $
unless $\Delta\in \Z$;
$\V$ is called {\em $\Z_{\geq 0}$-graded}
if $\V_{-\Delta}=\zero $
unless $\Delta\in \Z_{\geq 0}$.

Let $\V$ be a $\Z$-graded vertex algebra,
$H$ the Hamiltonian.
We say that
 $\V$
 is 
{\em compatibly $\Z_{\geq 0}$-gradable}
if $\V$ admits a Hamiltonian $H'$
which 
gives a $\Z_{\geq 0}$-grading
of $\V$
and 
satisfies
$[H, H']=0$.

The vertex algebra $\V$ is called {\em commutative} if
$[u_{(m)},v_{(n)}]=0$ for all $u,v\in \V$,
$m,n\in \Z$.
It is well-known that
$\V$ is commutative if and only if
$u_{(m)}v=0$ for all $u,v\in \V$,
$m\geq 0$.

\subsection{M\"{o}bius conformal Vertex  Algebras
and Conformal Vertex  Algebras}
Let $\V$ be a $\Z$-graded vertex algebra,
$H$ the Hamiltonian.
Suppose there exists 
an even operator
$\Lone$ on $\V$
satisfying the following:
\begin{align}
 &
 [\Lone, H]=\Lone,\\
&[\Lone, Y(v,z)]=Y(\Lone v,z)+2z Y(H v,z)
+z^2 Y(\TT  v,z)\quad \forall v\in \V.
\end{align}
Then $\V$ is called a
{\em M\"{o}bious conformal} vertex  algebra  \cite[Section 4.9]{KacBook2}.
In this case the triplet
$\{\Lone, H, \TT  \}$
forms the Lie algebra
$\sl_2(\C)$ in $\End \V$:
\begin{align*}
[H,T^*]=-T^*,\quad  [T^*, T]=2H,\quad 
[H,T]=T.
\end{align*}
A vector which is annihilated by $T^*$ is called 
a {\em quasi-primary vector}.

A {\em conformal vector} 
of $\V$
is an even vector $\omega\in \V$
such that
the corresponding field $Y(\omega,z)=\sum_{n\in \Z}L(n)
z^{-n-2}$ 
has the following properties:
\begin{enumerate}
 \item 
$
[L(m),L(n)]=(m-n)L(m+n)+\frac{
m^3-m}{12}\delta_{m+n,0}c_{\V}\id_{\V}$
for $m, n \in\Z$,
where $c_{\V}\in \C$ (\em{the central charge});
\item $\TT  =L(-1)$;
\item the action of $L(0)$ 
is semisimple and  a Hamiltonian
(i.e.\ satisfies \eqref{eq:Hamiltonial})
\end{enumerate}
In this case 
 $\V$ has the 
M\"{o}bious conformal vertex algebra
structure
 with
$\Lone=L(1)$ 
and $H=L(0)$.
A  vertex algebra with a conformal vector is called
 {\em conformal}
(or {\em vertex operator algebra}).
%A vertex operator algebra $\V$ is
%considered to have the fixed Hamiltonian $L(0)$.
\subsection{Filtration of Vertex Algebras}
A vertex algebra $\V$ is called 
{\em filtered}
if there exists an increasing 
%separated,
%exhaustive
 filtration
$F=\{F_p \V\}$
of   $\V$ as a superspace
which is compatible with the  vertex algebra
structure in the following sense:
\begin{align}
&|0\ket\in F_0 \V\backslash F_{-1}\V,\label{eq:2006-06-20-11-49}
\\
%&F_{p} \V=\zero \text{ for $p\ll 0$;}
%\label{eq:2006-05-23-13-02}
%\\ 
&\TT   \cdot F_p \V\subset F_p \V;
%\label{eq:filt-partial}
\\
&v_{(n)}\cdot  F_q\V\subset F_{p+q}\V\quad \text{for  all $p,q\in \Z$,
$v\in F_p \V$ and $n\in \Z$ }.
\label{eq:2006-05-23-13-03}
\end{align}
Also,
unless otherwise stated,
we  require 
 a filtration of a vertex algebra $\V$
to be separated and exhaustive.

 A filtration of a $\Z$-graded vertex algebra $\V$
is a filtration $F$ of $\V$
which is  compatible with
the action of the Hamiltonian:
$H\cdot F_p \V\subset F_{p}\V$.

If $(\V, F)$ is a filtered vertex algebra
(resp.\ a filtered $\Z$-graded vertex algebra)
then
$\gr^F \V$ is naturally  a vertex algebra
(resp.\ a $\Z$-graded vertex algebra),
see e.g.\ \cite{Li}.
%Indeed one has
%\begin{align*} 
%\gr^F \V=\bigoplus_{\Delta \in \Z}(\gr^F \V)_{-\Delta},
%\end{align*}
%where $(\gr^F \V)_{\Delta}=\bigoplus_{p}F_p
%\V_{\Delta}/F_{p-1}\V_{\Delta}$.
%Also, from \eqref{eq:2006-05-23-13-03},
%$v_{(n)}$,
%with $v\in \gr \V$ and $n\in \Z$,
%defines a well-defined operator on $\gr^F \V$.
%In other words,
%the vertex operator
%$Y(?,z)$ is well-defined in $\gr^F V$.
%Also,
%the translation operator $\TT   $ is  well-defined on $\gr^F
%\V$.
%Further from \eqref{eq:2006-06-20-11-49},
%\eqref{eq:2006-05-23-13-02}
%and
%\eqref{eq:2006-05-23-13-03}
%we have
%\begin{align*}
%|0\ket \in F_0 \V\backslash F_{-1}\V,
%\end{align*}
%and it follows that
%the vacuum axion holds for $\gr^F \V$.
%Therefore,
%they define on
%$\gr^F \V$ a (graded) vertex algebra structure.
The vertex algebra $\gr^F \V$  is called the {\em graded vertex algebra
associated with the filtered vertex algebra $\V$}.

A filtered vertex (super)algebra 
$(\V, F)$ is called {\em quasi-commutative}
if $\gr^F \V$ is commutative.
If $(\V, F)$ is quasi-commutative
then $\gr^F \V$
is naturally a {\em vertex Poisson  algebra}
 (\cite[Proposition 4.2]{Li}).
However 
we do not use this fact
in this article.

\subsection{Standard Filtration}
\label{subsection:StandardFIltrations}
\begin{Th}[\cite{Li}]
\label{Th:Li}
  Let $\V$ be a 
vertex algebra,
which is
$\Z$-graded by $H$,
and 
compatibly $\Z_{\geq 0}$-gradable by $H'$.
Take a strongly generating
$H$, $H'$-invariant subspace $U$
of $\V$.
Let
$\{a^j; j\in J\}$
be a  basis of $U$,
$\Delta_i=\Delta_{a_i}$. 
%Set
%$G_{-1}\V=\zero $
%and
Let 
$G_p \V$, with $p\in \Z$,
be the subspace of $\V$ spanned
by all the vectors
\begin{align*}
a^{j_1}_{(-n_1)}a^{j_2}_{(-n_2)}\dots a^{j_r}_{(-n_r)}
|0\ket
\end{align*}
with $r\geq 0$, $j_s\in J$,
$n_s\geq 1$
%and $v\in V_0$
satisfying the relation 
\begin{align*}
 \Delta_{j_1}+\Delta_{j_2}+\dots+ \Delta_{j_r}
\leq p.
\end{align*}
Here by convention
$|0\ket \in G_0\V \backslash G_{-1}\V$.
Then
\begin{enumerate}
 \item $G=\{G_p \V\}$ 
gives a filtration of the 
$\Z$-graded vertex algebra $\V$;
\item $(\V, G)$ is quasi-commutative;
\item 
$G$ is the finest filtration of $\V$
such that
$V_{-\Delta}\subset G_{\Delta}\V$ for all $\Delta$.
In particular
the filtration $G$ is independent of the choice of
a
strongly generating $H$,
$H'$-invariant subspace $U$ of $\V$.
\end{enumerate}
\end{Th}
\begin{Rem}
  In \cite{Li}
it was assume that 
$\V$ is $\Z_{\geq 0}$-graded 
and $\V_0= \C |0\ket$.
However
it is easy to see that
Theorem \ref{Th:Li}
follows from the
the assertion in the case that $\V$ is $\Z_{\geq 0}$-graded.
Also,
%noting  that $G_0 \V$ is a supercommutative vertex subalgebra
in the case that  $\V$ is $\Z_{\geq 0}$-graded,
it is not difficult to remove the condition
that $\V_0=\C |0\ket$ .
%\item In the case of
%$\V=V_k(\sg)$ (see Section \ref{subsection:univ-affine}),
%Similarly,
%$G_p C_k(\sg)$ (see \eqref{eq:2006-05-17-13-15})
%is the subspace of
%$C_k(\sg)$ spanned by
%vectors
%of the form 
%$x_1x_2\dots x_r$ with $r\leq p$
%and $x_i\in \sg\* \C[t\inv]t\inv$,
%or $x_i=\psi_{\alpha}(-m)$,
%$\psi_{\beta}(-n)$
%with $\alpha\in \sroots_+$,
%$m\geq 1$,
%$\beta\in \sroots_-$,
%$n\geq 0$.
\end{Rem}
The filtration $G$
given in Theorem \ref{Th:Li} 
is called the {\em standard filtration}
of a $\Z$-graded,
compatibly $\Z_{\geq 0}$-gradable vertex algebra $\V$.

\subsection{PBW Basis of Vertex Algebras}
\label{subsection:PBW-VA-22-41}
Let $\V$ be a vertex algebra,
$F$  a filtration of $\V$
such that
$(\V,F)$ is quasi-commutative.
Let $U$ be a subspace
of $\V$ such that
the its image $\bar{U}$
in $\gr^F \V$
strongly generates $\gr^F \V$.
Then $U$
 strongly generates $\V$.
%Denote by 
%$ S(\bar{U}\* \C[t\inv]t\inv)$ the symmetric algebra
%of the superspace $\bar{U}\* \C[t\inv]t\inv$.
%Here
%we regard
%$(\bar{U}\* \C[t\inv]t\inv)^{\even}=
%\bar{U}^{\even}\* \C[t\inv]t\inv$,
%$(\bar{U}\* \C[t\inv]t\inv)^{\odd}=
%\bar{U}^{\odd}\* \C[t\inv]t\inv$.
We have the following surjective linear map
\begin{align}
\begin{array}{cccc}
 S(\bar{U}\* \C[t\inv]t\inv)
 &\longrightarrow & \gr^F \V &\\
(\bar a^1\* t^{-n_1})\dots (\bar a^r\* t^{-n_r})& 
\longmapsto & \bar a^1_{(-n_1)}\dots \bar a^r_{(-n_r)}|0\ket &
\end{array}
\label{eq:2006-09-15-PBW}
\end{align}
with $\bar a^i\in \bar{U}$, $a_i\geq 1$.

We say that 
a subspace
$U \subset \V$ {\em generates  a  PBW basis} of 
 $\V$
if there exists a
filtration 
$F$  of $\V$
such that $(\V, F)$ is quasi-commutative
and
the  map \eqref{eq:2006-09-15-PBW}
 is a linear isomorphism;
 In this case we say that {\em $\V$ admits a PBW basis.}

%Suppose $U$ generates a PBW basis 
%of a $\Z$-graded,
%compatibly $\Z_{\geq 0}$-gradable
%vertex algebra
%$\V$.
%Then the map
%\eqref{eq:2006-09-15-PBW}
%is an isomorphism
%for 
%$F=G$,
%the standard filtration of $\V$
%(of curse we take $\bar U$ as the image of $U$
%in $\gr^G \V$).
\subsection{Vacuum Subalgebras}
For a vertex algebra $\V$,
let
\begin{align}
 \Vac \V:=\{v\in \V; T\cdot v=0\}.
\end{align}
Then $|0\ket \in\Vac \V$.
The $\Vac \V$ is a vertex subalgebra of $\V$
and called the {\em vacuum subalgebra of $\V$}
(\cite[Remark 4.4b]{KacBook2}).

Suppose that 
$\V$ is filtered,
$F$ its filtration.
Then $F$ induces a filtration of $\Vac \V$:
$F_p \Vac \V=\Vac \V\cap F_p \V$.
There is a natural embedding of vertex algebras
\begin{align}
 \gr^F  \Vac \V\hookrightarrow  \Vac (\gr^F \V)
\label{eq:2006-10-10-13-38}
\end{align}
which is given by 
the correspondence
$\sigma_p(v)\mapsto \sigma_p (v)$
with $v\in F_p \Vac \V$.
\begin{Pro}\label{Pro:2006-10-10-15-11}
 Let $\V$ be a 
vertex algebra that admits a PBW basis.
Then $\Vac \V= \C |0\ket$.
\end{Pro}
\begin{proof}
By assumption 
there exist  a subspace $U$ 
of $\V$ and a filtration $F$ of $\V$
such that $\gr^F \V$  commutative and 
the map \eqref{eq:2006-09-15-PBW}
is a linear isomorphism.
By \eqref{eq:2006-10-10-13-38}
it is sufficient
to show that 
$\Vac (\Gr^F \V)
=\C |0\ket$.

By assumption
$\gr^F \V$
is isomorphic to
the space of 
 the following form:
\begin{align*}
R=\C[x^{(j)}_{(-n)}; j\in J_{\bar 0}, n\geq 1]\*
 \Lam(y^{(j)}_{(-n)};j\in J_{\bar 1},n\geq 1),
\end{align*}
where 
$ \Lam(y^{(j)}_{(-n)};j\in J_{\bar 1},n\geq 1)$
is the Grassmann algebra with generators
$y^{(j)}_{(-n)}$,
$j\in J_{\bar 1}$, $n\geq 1$.
Under this identification 
$T$ acts as the following even derivation:
\begin{align*}
 [T, x^{(j)}_{(-n)}
]=n x^{(j)}_{(-n-1)},\quad
 [T, y^{(j)}_{(-n)}
]=n y^{(j)}_{(-n-1)}
\end{align*}
Define the even derivations
$[H, ?]$,
$[T^*, ?]$
on $R$
by the following:
\begin{align*}
& [H, x^{(j)}_{(-n)}
]=n x^{(j)}_{(-n)},\quad
 [H, y^{(j)}_{(-n)}
]=n y^{(j)}_{(-n)}\\
& [T^*, x^{(j)}_{(-n)}
]=n x^{(j)}_{(-n+1)}
,\quad
 [T^*, y^{(j)}_{(-n)}
]=n y^{(j)}_{(-n+1)}.
\end{align*}
Here by convention
$x^{(j)}_{(0)}=y^{(j')}_{(0)}=0$
for all $j,j'$.
This gives the well-defined action of $\sl_2(\C)$ on $R$.
Therefore
by \cite[Proposition 4.9(a)]{KacBook2} 
we have
$\Vac (\gr^F \V)=\{v\in R; H\cdot v=0\}
=\C |0\ket$.
\end{proof}
\begin{Rem}
In general if $\V$ is a $\Z_{\geq 0}$-graded
M\"{o}ious conformal vertex algebra 
then
$\Vac \V\subset \V_0$,
see \cite[Proposition 4.9]{KacBook2}.
\end{Rem}
\subsection{Strict Filtration}
\label{subsecton-new:StrictFiltration}
Let $\V$ be a filtered vertex algebra,
$F$ its filtration.
Then
$\TT  F_p \V\subset F_p\V\cap \TT  \V$.
We call $F$ {\em{strict}}
if
\begin{align}\label{eq:condi-05-08}
  F_p\V\cap \TT   \V=\TT    F_p \V
\quad \text{ for all } p.
\end{align}

\begin{Lem}\label{Lem:06-05-16-16-06}
$ $

\begin{enumerate}
\item The filtration $F$ is strict 
if and only if 
the natural embedding 
$\gr^F \Vac \V\hookrightarrow \Vac (\gr^F \V)$
is an isomorphism of vertex algebras.
\item If $\Vac(\Gr^F \V)=\C |0\ket$,
then
$\Vac \V =\C |0\ket$
and $F$ is strict.
\end{enumerate}\end{Lem}
\begin{proof}
(i)
Let $\Vac (\gr^F _p \V)=\Vac (\gr^F \V)\cap \gr^F _p \V
=\{v\in F_p \V; Tv \in F_{p-1}\V\}/F_{p-1}\V$.
Then
\begin{align*}
  \Vac (\gr^F_p \V)/\gr^F_p \Vac \V=\{v\in F_p \V;
Tv\in F_{p-1}\V\}/
(F_{p-1}\V + F_p \Vac \V).
\end{align*}
One sees that
the correspondence $v\mapsto Tv$
induces an isomorphism
\begin{align}
 \Vac (\gr^F_p \V)/\gr^F_p \Vac (V)
\isomap
(F_{p-1}\V\cap T F_p \V)
/T F_{p-1}\V.
\end{align}
But as easily seen
$F$ is strict if and only if 
$F_{p-1}\V\cap T F_p \V=T F_{p-1}\V$ for all $p$.
Therefore the assertion follows.
(ii) 
The first assertion is obvious.
The second assertion
follows from (i).
%The first assertion follows from the fact that
% $|0\ket \in \Vac \V  $
%and we have the embedding
%$\ker_{\V} \TT  \hookrightarrow \Vac(\gr^F \V) \TT  $.
%The second assertion follows from (1).
\end{proof}
\begin{Pro}\label{Pro:2006-05-17-23-34}
 Let $\V$ be a 
vertex algebra,
$F$ its filtration.
Suppose that $\gr^F \V$ 
admits a PBW basis.
Then 
$\Vac \V  =\C |0\ket$
and
$F$ is strict.
\end{Pro}
\begin{proof}
By assumption 
and Proposition \ref{Pro:2006-10-10-15-11},
$\Vac(\gr^F \V)=\C |0\ket$.
Therefore Lemma \ref{Lem:06-05-16-16-06} (ii)
gives the assertion.
\end{proof}
\begin{Pro}\label{Pro:standard-is-strict}
Let  $\V$ be a
$\Z$-graded,
compatibly $\Z_{\geq 0}$-gradable
vertex algebra that admits a PBW basis.
Then
the standard filtration $G$ is strict.
\end{Pro}
\begin{proof}
%Let
%$G$ be  the standard filtration
%of $\V$.
Under the assumption of Proposition,
$\gr^G \V$ admits a PBW basis.
Thus the assertion follows from Proposition \ref{Pro:2006-05-17-23-34}.
\end{proof}
\subsection{Lie Algebras Attached to Vertex Algebras}
\label{subsection:LieAl-VA}
Let $\V$ be a  vertex algebra.
Define
\begin{align}
 \Lie{\V}\teigi L(\V)/\im \wpar
\end{align}
(see \eqref{eq:2006-11-04-23-16}),
where
\begin{align}
 \wpar
\teigi \TT   \* \id +\id \* \frac{d}{d t}.
\end{align}
Then $\Lie{\V}$ has the Lie algebra structure
\cite{Bor},
whose
the commutation relation is given by
\begin{align}
 [u_{\{m\}},v_{\{n\}}]=\sum_{r\geq 0}\begin{pmatrix} m\\r\end{pmatrix}
(u_{(r)}v)_{\{m+n-r\}},
\end{align}
where $v_{\{n\}}$ denotes the image of $v\* t^n$
with $v\in \V$ and $n\in \Z$
in $\Lie{\V}$.
By \eqref{eq:vertex-axiom1}
the correspondence $v_{\{n\}}\mapsto v_{(n)}$
defines a representation of
  $\Lie{\V}$ on $\V$.

Define the adjoint action $\ad T$ of
 $T$
on $\Lie{\V}$ by the following.
\begin{align}
\ad T\cdot v_{\{n\}}=
[T, v_{\{n\}}]:=(T v)_{\{n\}}=-nv_{\{n-1\}}.
\label{eq:2006-10-11-14-32}
\end{align}
If $\V$ is graded by a Hamiltonian $H$,
then there is a 
 adjoint  action $\ad H$ of
the Hamiltonian
$H$
on
$\Lie{\V}$ in view of  \eqref{eq:Hamiltonial}:
\begin{align}
\ad  H\cdot u_{\{n\}}
=[H, u_{\{n\}}]:=(H u)_{\{n\}}+ (Tu)_{\{n+1\}}.
\label{eq:2006-10-11-14-33}
\end{align}
This gives a $\Z$-grading of 
$\Lie{\V}$:
$\Lie{\V}=\bigoplus_{d\in \Z}\Lie{\V}_d$,
$[\Lie{\V}_d, \Lie{\V}_{d'}]\subset \Lie{\V}_{d+d'}$,
where 
\begin{align}
 \Lie{\V}_d=\{u \in \Lie{\V};\ad H \cdot u=-d u\}.
\end{align}
Set
\begin{align}\label{eq:def-of-degree-n-opeartor}
 v_n\teigi v_{\{n+\Delta_v-1\}}.
\end{align}
Then
$\ad H\cdot v_n=-n v_n$
and
$\Lie{\V}_d$ is spanned by 
$v_d$ with homogeneous vectors $v\in \V$.
One has
\begin{align}\label{eq:2006-05-18-11-12}
 [u_m, v_n]=\sum_{r\geq 0}
\begin{pmatrix}
 m+\Delta_u-1\\ r
\end{pmatrix}
(u_{(r)} v)_{m+n}.
\end{align}

Next suppose that $\V$ is M\"obius conformal.
Then we also have the 
following adjoint action
of $T^*$ on $\Lie{\V}$:
\begin{align}
 \ad T^* \cdot v_{\{n\}}=[T^*, v_{\{n\}}]:=(T^* v)_{\{n\}}+2 (H v)_{\{n+1\}}
+(T v)_{\{n+2\}}.
\label{eq:2006-10-11-14-34}
\end{align}
Thus
$\sl_2(\C)=\haru \{T^*, H, T\}$  
acts on $\Lie{\V}$
by even derivation.
% \begin{align}
%\ad X\cdot [u_m, v_n]=[\ad X\cdot u_m, \ad X \cdot v_n]
%\quad \text{for $X\in \{T^*, H, T\}$,
%}
% \end{align}
%$u,v\in \V$,
%$m,n\in \Z$.

The same proof of \cite[Proposition 4.1.1]{NT}
applies for the following assertion 
because only the M\"obius conformal
structure of $\V$ is used in the argument
(cf.\ \cite[Section 2.8]{FHL}).
\begin{Pro}[{\cite[Proposition 4.1.1]{NT},
cf.\ \cite[Section 5]{FHL}}]
\label{Pro:anti-auto}
Let $\V$ be a $\Z$-graded M\"obius conformal vertex algebra
on which $T^*$ acts locally nilpotently.
Set
\begin{align*}
 \theta(v_n)=
\begin{cases}
(-1)^{\Delta_v}(e^{T^*}v)_{-n}&\text{if $p(v)=\bar {0}$}\\
\sqrt{-1}(-1)^{\Delta_v}(e^{T^*}v)_{-n}&\text{if $p(v)=\bar 1$}
\end{cases}
\end{align*}
for 
 $v\in \V$,
$n\in \Z$.
Then $\theta$ defines a well-defined anti-Lie algebra
isomorphism of $\Lie{\V}$
such that $\theta^2(v)=(-1)^{p(v)}v$.
\end{Pro}
If $\V$ is purely even
then $\theta$ is an anti-Lie algebra involution
of $\Lie{\V}$.

\subsection{Filtration of $\Lie{\V}$}
Let 
$V$ be a filtered vertex algebra,
$F$ its filtration.
Define an increasing filtration on $L(\V)$
by 
\begin{align*}
F_p L(\V)
=L(F_p\V)
\end{align*}
(notation \eqref{eq:2006-11-04-23-16}).
This gives the quotient filtration 
$\{F_p \Lie{\V}\}$
of the Lie algebra
$\Lie{\V}$.
\begin{Pro}\label{Pro:fil-Lie}
Let $\V$ be a filtered vertex algebra,
$F=\{F_p \V\}$ its filtration.
Then 
there is a natural
surjective Lie algebra homomorphism
$\Lie{\gr^F \V}\ra \gr^F \Lie{\V}$
given by
the correspondence
\begin{align}
 (\sigma_p(v))_{\{n\}}\mapsto
\sigma_p(v_{\{n\}})\quad 
\text{for $p\in \Z$, $v\in \V$, $n\in \Z$}
\label{eq:2006-10-11-12-48}.
\end{align}
Moreover if 
 $F$ is strict then
this is an   isomorphism.
\end{Pro}
\begin{proof}
By definition
$\Lie{\gr^F \V}=
\bigoplus_{p\in \Z}\Lie{\gr^F_p \V}$,
where
\begin{align*}
\Lie{\gr^F_p \V}&=
\gr^F_p L(\V)/
\wpar(\gr^F_p L(\V))=F_p L(\V) 
/\left(F_{p-1}
L(\V)
+\wpar F_p L(\V )
\right).
\end{align*}
On the other hand
\begin{align*}
 \gr^F_p \Lie{\V}&=F_p\Lie{\V}/F_{p-1}\Lie{\V}
=F_p L(\V)
/
\left(F_{p-1}L(\V)
+
\wpar L(\V)
\cap F_p L(\V)
\right).
\end{align*}
Clearly
$
\wpar F_p L(\V)
\subset \wpar L(\V)
\cap F_p L(\V)
$.
Thus
\eqref{eq:2006-10-11-12-48}
is well-defined and surjective.
It is straightforward to see that it is  a
Lie algebra homomorphism.

Next assume that
$F$ is strict.
It is sufficient to show that
the opposite correspondence
\begin{align}
\sigma_p(v_{\{n\}})\mapsto  (\sigma_p(v))_{\{n\}}
\quad 
\text{for $p\in \Z$, $v\in \V$, $n\in \Z$}.
\label{eq:2006-10-11-12-50}
\end{align}
is well-defined.
For this
we need to show that
\begin{align}
\wpar L(\V)
\cap  
F_p L(\V)\subset 
\wpar F_p L( \V)
\label{eq:2006-05-16-11-28}.
\end{align}
Let $w\in \wpar L(\V)
\cap F_p L( \V)
$.
Then $w$ has the form
\begin{align*}
w=\sum_{i=1}^r\wpar(v_i\* f_i)=\sum_{i=1}^r
((\TT   v_i)\* f_i +v_i\* f_i')
\end{align*}
with $v_i\in \V$,
$f_i\in \C[t,t\inv]$.
We have to show that each $v_i$ belongs to $F_p \V$.
We may assume that
$\{f_i \}$ are linearly independent
homogeneous elements of $\C[t,t\inv]$,
and $\deg f_i\geq \deg f_{i+1}$ for all $i$.
Because
$f_1$ has the largest degree,
it follows that
$\TT   v_1$
 must belong to $F_p \V$.
Therefore  $v_1$ must belong to $F_p \V$
since $F$ is strict.
Repeating this argument,
we get $v_i\in F_p \V$ for all $i$.
\end{proof}

\subsection{Current Algebras  of Vertex Algebras and $\V$-Modules}
\label{subap:CurrentAl-of-VA}
Let $\V$ be a $\Z$-graded vertex algebra,
$H$ its Hamiltonian.
Following \cite{MNT},
we associate with $\V$
the {\em current algebra}
$\U(\V)$,
which is essentially the 
universal enveloping algebra of $\V$
in the sense of I. Frenkel and Y. Zhu \cite{FZ}:
Let $\UU(\V)$
be the 
 quotient of 
the universal enveloping algebra
$U(\Lie{\V})$
of $\Lie{\V}$
by the two-sided ideal
generated by $(|0\ket)_{0}-1$.
The action $\ad H$ naturally extends to $\UU(\V)$.
Let $\UU(\V)_d=\{u\in \UU(\V); \ad H\cdot u=-d u\}$.
This makes
$\UU(\V)$ a graded algebra:
\begin{align*}
 \UU(\V)=\bigoplus_{d\in \Z}\UU(\V)_d,\quad
\UU(\V)_d \cdot \UU(\V)_{d'}\subset \UU(\V)_{d+d'}.
\end{align*}
Denote by $\wt{\UU}(\V)=\bigoplus_{d\in \Z}\wt{\UU}(\V)_d$  
the 
{\em standard degreewise completion}
(see Section \ref{subsectuib:Compatible Degreewise Complete Algebra})
of $\UU(\V)$.
The $\wUU(\V)$ is equipped with 
the left linear topology
defined by the decreasing
sequence of the left ideals
\begin{align*}
\II_N(\wUU(\V))=\Ker\left(\wUU(\V)\twoheadrightarrow \UU(\V)/
(\bigoplus_{d\in \Z}\sum_{r>N}\UU(\V)_{d-r}\UU(\V)_{r})
\right)
\end{align*}with $N\geq 0$.
Let
 $\B(\V)=\bigoplus_{d\in \Z}\B(\V)_d$  be the 
ideal of 
$\wt{\UU}(\V)$
corresponding to the
Borchards identity,
i.e.\ 
the $\Z$-graded two-sided 
ideal
generated by
\begin{align}
& (u_{(m)}v)_{\{n\}}-\sum_{i\geq 0}
\begin{pmatrix}
 m\\i
\end{pmatrix}
(-1)^i(u_{\{m-i\}}v_{\{n+i\}}-(-1)^{p(u)p(v)}(-1)^m
 v_{\{m+n-i\}}u_{\{i\}})
\label{eq:2006-05-26-14-48}
\end{align}
with 
$u,v\in \V$, $m,n\in \Z$.
Let $\wt{\B}(\V)=\bigoplus_{d\in \Z}
\wt{\B}(\V)_d$ be the {\em degreewise closure}
(see Section \ref{subsectuib:Compatible Degreewise Complete Algebra})
of $\B(\V)$
in $\wUU(\V)$.
Now one can define
the current algebra $\U(\V)$
of $\V$ is
as 
 %the quotient algebra
\begin{align}
 \U(\V):=\wt{\UU}(\V)/\wt\B(\V).
\end{align}
By definition the $\U(\V)$ is graded by $\ad H$:
\begin{align*}
 \U(\V)=\bigoplus_{d\in \Z}\U(\V)_d,
\quad
\U(\V)_d=\{u\in \U(\V); \ad H\cdot u=-d u\}.
\end{align*}
The $\U(\V)$ is equipped with the left linear topology
induced by that of $\wUU(\V)$.
Its topology
is defined by the decreasing sequence of
left ideals $\II_N(\U(\V))$,
where $\II_N(\U(\V))$ is the image of 
$\II_N(\wUU(\V))$ in $\U(\V)$.
The space
$\II_N(\U(\V))$ coincides with  the degreewise closure 
of $\U(\V)\cdot \sum_{r>N}\U(\V)_r$
in $\U(\V)$.

The $\U(\V)$ is a {\em compatible degreewise complete 
algebra} \cite{MNT}:
Each $\U(\V)_d$ is complete with respect to the
relative topology:
\begin{align}
 \U(\V)_d=\lim\limits_{\leftarrow \atop N}\U_N(\V)_d,
\quad
\text{where }\U_N(\V)_d:=\U(\V)_d/\II_N(\U(\V))_d,
\end{align}
and 
 the multiplication map
$\U(\V)_d\times
\U(\V)_{d'}\ra
\U(\V)_{d+d'}$ is continuous,
see \cite{MNT} for details.
Here
$\II_N(\U(\V))_d=\{u\in \II_N(\U(\V));
\ad H\cdot u= -d u\}$.
We set $\U_N(\V)=\bigoplus_{d\in \Z}\U_N(\V)_d$.

The  image
of $v_n$, with $v\in \V$,
$n\in \Z$, by the natural
map $\Lie{\V}\rightarrow \U(\V)$
is denoted also by $v_n$,
and so is
 its image in $\U_N(\V)$.
\begin{Pro}[{\cite[Proposition 6.4.1]{MNT}}]
\label{Pro:2006-11-04-21-35}
The image of $\Lie{\V}$
in $\U(\V)$
is dense,
that is,
the natural map $\Lie{\V}
\ni v_n\mapsto v_n\in  \U_N(\V)$
is surjective for all $N\geq 0$.
\end{Pro}
Set
\begin{align}
 \U(\V)_{\geq 0}:=\bigoplus_{d\geq 0}\U(\V)_d,
\quad \U(\V)_{>0}:=\bigoplus_{d>0}\U(\V)_d.
\end{align}
These are subalgebras of $\U(\V)$.

\smallskip

Suppose that
$\V$  is M\"obius conformal.
Then,
by \cite[Section 5]{FLM},
the anti-Lie algebra
 isomorphism $\theta: \Lie{\V}\ra \Lie{\V}$ 
induces a anti-algebra 
isomorphism
$\theta: \U(\V)\ra \U(\V)$
such that $\theta(\U(\V)_d)=\U(\V)_{-d}$.
\subsection{$\V$-Modules and Zhu Algebras}
\label{subsection:MOdules-and-Zhu-algberas}
Let $\V$ be a $\Z$-graded vertex algebra,
$H$ its Hamiltonian.
A {\em $\V$-module} is by definition
a
$\U(\V)$-module.
A {\em graded $\V$-module}
 is a 
$\U(\V)$-module $M$
which carries
a $\C$-grading
$M=\bigoplus_{d\in \C}M_d$
such that
$\U(\V)_n\cdot M_d\subset M_{d+n}$
for all $n,d$.
Let $\V\GrMod$ be the category of graded $\V$-modules:
the objects of $\V\GrMod$ are graded $\V$-modules
and the morphisms of $\V\GrMod$ are
all the graded $\U(\V)$-module homomorphisms.
Here a $\U(\V)$-module homomorphism $\phi:M\ra N$ is called graded if
for each $d\in \C$
there exists $d'\in \C$
satisfying
$\phi(M_{d+n})\subset N_{d'+n}$ for all $n\in \Z$.

An {\em admissible $\V$-module} \cite{ABD}
is a graded $\V$-module such that
 $M_{d+n}=\zero $ for $n\gg 0$ with a fixed $d$.
If $M$ is a admissible $\V$-module then
the action
$\U(\V)\times M\ra M$
is continuous with respect to the 
topology on $\U(\V)$ and the discrete topology on $M$.
Let $\V\adMod$ be the full subcategory
of the category of graded $\V$-modules
consisting of admissible $\V$-modules.
If $\V$ is $\Z_{\geq 0}$-graded
then $\V$ belongs to $\V\adMod$
when considered as a $\V$-module.

%Then
%an object  $M $ of $ \V\GrMod$ is simple if and only if
%$M$ has no nonzero proper graded submodule.
%\begin{Rem}
% If $\V$ is a vertex operator algebra
%then a
%graded $\V$- module is called an admissible module \cite{ABD}.
%$\end{Rem}
%\subsection{The Zhu Algebra of Vertex Algebras}
%\label{subsection:Zhu's algebra}

Let $\V$ be a $\Z$-graded vertex algebra.
With it one associated 
the {\em Zhu Algebra}
 $\Zhu(\V)$
\cite{Zhu} of $\V$,
which is the unital  associative algebra
in the usual sense.
As remarked in \cite{FZ}
one may define $\Zhu(\V)$
as  
\begin{align}
 \Zhu(\V)&=\U_0(\V)_0
(=\{a\in \U_0(\V); \ad H\cdot a=0\}),
\label{eq:2006-11-04-23-29}
\end{align}
or equivalently,
\begin{align}
\Zhu(\V)=\U(\V)_0/
\overline{\sum_{r>0}\U(\V)_{-r}\U(\V)_r}
\end{align}
(see \cite[Theorem A.2.11]{NT}
for the proof of the equivalence with the usual definition).

\smallskip

For a graded $\V$-module  $M$,
let
 $M_{\tp}$ be the
sum of nonzero homogeneous 
subspace $M_d$
with $M_{d+n}=\zero $
for all $n> 0$.
Then 
$M_{\tp}$ is naturally a module over $\Zhu(\V)$.
If $M$ is a simple object of $\V\adMod$
then $M_{\tp}=M_d$ for some $d$.

\begin{Th}[{Zhu \cite{Zhu}}]
\label{Th:Zhu-08}
Let $\V$ be a $\Z_{\geq 0}$-graded vertex algebra.
The correspondence $M\mapsto M_{\tp}$
gives a  bijection
between  equivalence classes 
of
 simple objects of 
$\V\adMod$
and
 the equivalence classes of
 irreducible representations 
of $\Zhu(\V)$.
\end{Th}
\subsection{Filtration of $\U(\V)$}
Let $\V$ be a $\Z$-graded vertex algebra,
$F$ its filtration.
Then $F$ induces a filtration 
$\{F_p \UU(\V)\}$ of $\UU(\V)$:
\begin{align*}
 F_p \UU(\V)=
\sum\limits_{p_1+p_2+\dots+p_r\leq p}
F_{p_1}\Lie{\V}\cdot F_{p_2}\Lie{\V}\cdot \dots \cdot F_{p_r}\Lie{\V}.
\end{align*}
The following assertion follows from Proposition \ref{Pro:fil-Lie}.
\begin{Lem}\label{Lem:05-03-03-1}
There is
 a natural surjective
homomorphism 
 $ 
\UU(\gr^F \V)\rightarrow \gr^F \UU(\V)
$ of graded algebra,
given by the correspondence
\begin{align*}
 (\sigma_{p_1}(v^1))_{n_1}\cdot
\dots\cdot  (\sigma_{p_r}(v_r))_{n_r}
\mapsto 
\sigma_{p_1+\dots+p_r}(
v^1_{n_1}\cdot \dots \cdot v^r_{n_r})
\end{align*}
with $v^i\in \V$, $n_i\in \Z$.
This is an isomorphism
if $F$ is strict.
\end{Lem}
Let $F_p \wUU(\V)$ be the 
degreewise closure of the image of $F_p \UU(\V)$
in $\wUU(\V)$.
This
 gives
a filtration
of $\wUU(\V)$.
Let $\{F_p \U(\V)\}$ be the quotient filtration:
$F_p \U(\V)=F_p \wUU(\V)/F_p \wt{\B}(\V)$,
where 
 $F_p \wt{\B}(\V)=\wt{\B}(\V)\cap F_p\wUU(\V)$.
The union $\bigcup_p F_p \U(\V)$ is dense in $\U(\V)$.
Let $\{F_p \U_N(\V)\}$ with $N\geq 0$ be the
quotient filtration of  $\U_N(\V)$. 
This is an   exhaustive filtration.

\begin{Rem}
If $F_{-1}\V=0$
then
$\{F_p \U_N(\V)\}$ is obviously separated.
Let
$\V$ be a
$\Z$-graded,
compatibly $\Z_{\geq 0}$-gradable vertex algebra,
$G$ the standard filtration.
Then
from Proposition \ref{Pro:2006-11-04-21-35}
one sees that
the filtration
$\{G_p \U_N(\V)\}$ is  separated.
\end{Rem}

The 
following assertion follows from Proposition \ref{Pro:fil-Lie}.
\begin{Th}\label{Th;filtered-univesal-env}
Let 
$\V$ be
 a $\Z$-graded
vertex algebra,
$F$ its filtration. %such that $F_{-1}\V=0$.
Then there is a natural surjective
linear map
$\U_N(\gr^F \V)\ra \gr^F \U_N(\V)$
given by 
 the correspondence
\begin{align*}
 (\sigma_{p_1}(v^1))_{n_1}\cdot
\dots\cdot  (\sigma_{p_r}(v^r))_{n_r}
\mapsto 
\sigma_{p_1+\dots+p_r}(
v^1_{n_1}\cdot \dots \cdot v^r_{n_r})
\end{align*}
with $v^i\in \V$, $n_i\in \Z$
for each $N$.
Moreover if the filtration $F$ is strict
then this is an isomorphism.
 \end{Th}
Let $\wt{\gr}^F \U(\V)$ denote the degreewise completion of
$\gr^F \U(\V)$.
Then the 
surjection in Theorem \ref{Th;filtered-univesal-env}
extends to the
surjective homomorphism
$\U(\gr^F \V)\rightarrow \wt{\gr}^F \U(\V)$
of compatible degreewise complete algebras,
which is an isomorphism if $F$ is strict.

Let $F_p \Zhu(\V)=\Zhu(\V)\cap F_p\U_0(\V)
$.
This
gives the filtration 
$\Zhu(\V)$.
The following assertion is obvious from
Theorem 
\ref{Th;filtered-univesal-env}
and the definition \eqref{eq:2006-11-04-23-29}.
\begin{Th}\label{Th:2006-08-31013-10}
 Let $\V$ be a $\Z$-graded vertex algebra,
$F$ its filtration. %such that $F_{-1}\V=0$.
Then there is a natural 
surjective
algebra
homomorphism
$\Zhu(\gr^F \V)\ra \gr^F \Zhu(\V)$
given by 
the correspondence
\begin{align*}
 (\sigma_{p_1}(v^1))_{0}\cdot
\dots\cdot  (\sigma_{p_r}(v_r))_{0}
\mapsto 
\sigma_{p_1+\dots+p_r}(
v^1_{0}\cdot \dots \cdot v^r_{0})
\end{align*}
with $v^i\in \V$.
Moreover if the filtration $F$ is strict
then this is an isomorphism.
\end{Th}
\subsection{The PBW Theorem for Current Algebras and Zhu Algebras}
\label{subsection:PBW Thereom for Current Algebras and Zhu Algebras}
Let $\V$ be a $\Z$-graded vertex algebra,
$H$ its Hamiltonian,
$F$ its filtration such that
 $(\V, F)$ is 
quasi-commutative.
Take a $H$-invariant subspace $U$ of $\V$ 
such that its image $\bar{U}$ in $\gr^F U$
strongly generates $\gr^F \V$.

We define the adjoint action 
on  $L(\bar{U})$
(notation \ref{subsection:LieAl-VA})
by  $\ad H\cdot u\*t^n=(H u) \* t^n
-(n+1) u\* t^n$.
This action extends naturally to
the symmetric algebra $S(L(\bar{U}))$:
\begin{align*}
 S(L(\bar{U}))=\bigoplus_{d\in \Z}
S(L(\bar{U}))_d,
\quad
S(L(\bar{U}))_d=\{a\in S(L(\bar{U}));
\ad H\cdot a=-d a\}.
\end{align*}
Let $\S(\bar{U})$
denote
 the standard degreewise completion 
%(Section \ref{subsectuib:Compatible Degreewise Complete Algebra})
of 
$S(L(\bar{U}))$:
\begin{align}
& \S(\bar{U})=\lim_{\leftarrow \atop N}
\S_N(\bar{U}),\\
&\text{where }
\S_N(\bar{U})=
S(L(\bar{U}))/
\left(\bigoplus_{d\in \Z}\sum_{r>N}S(L(\bar{U}))_{d-r}S(L(\bar{U}))_r\right).
\end{align}

The space  $\S_N(\bar{U})$ can be identified with 
the subalgebra of 
$S(L(\bar{U}))$
spanned by the elements
of the form
\begin{align}
 \sigma_{p_1}(u^1)\* t^{n_1+\Delta_{u^1}-1}\cdot \dots
\cdot 
\sigma_{p_r}(u^r)\* t^{n_r+\Delta_{u^r}-1}
\label{eq:2006-10-14-00-48}
\end{align}
with $p_i\in \Z$,
$u^i\in U$,
$n_i\in \Z$
satisfying
\begin{align}
n_1+\dots +n_r\leq N.
\label{eq:2006-10-14-00-49}
\end{align}
 \begin{Th}
\label{TH:PBW-For-UV}
Let
$\V$,
$F$,
$U$,
$\bar{U}$ be as above.
Then
the map
\begin{align*}
 \begin{array}{ccc}
  \S_N(\bar{U})&\rightarrow  &\gr^F \U_N(\V) \\
\sigma_{p_1}(u^1)\* t^{n_1+\Delta_{u^1}-1}\cdot \dots
\cdot 
\sigma_{p_r}(u^r)\* t^{n_r+\Delta_{u^r}-1}
&\mapsto&
\sigma_{p_1+\dots +p_r}
(u^1_{n_1}\cdot \dots \cdot u^r_{n_r})
 \end{array}
\end{align*}
is surjective for each $N\geq 0$.
Further, 
if $\bar{U}$ generates a PBW basis of
$\gr^F \V$
then
  this is an isomorphism.
\end{Th}
\begin{proof}
First, the map
\begin{align*}
 \begin{array}{ccc}
  \S_N(\bar{U})&\rightarrow  & \U_N(\gr^F\V) \\
\sigma_{p_1}(u^1)\* t^{n_1+\Delta_{u^1}-1}\cdot \dots
\cdot 
\sigma_{p_r}(u^r)\* t^{n_r+\Delta_{u^r}-1}
&\mapsto&
\sigma_{p_1}
(u^1_{n_1})\cdot \dots \cdot \sigma_{p_r}(u^r_{n_r})
 \end{array}
\end{align*}
is surjective,
and 
bijective
if $\bar{U}$
 generates a PBW basis of $\gr^F \V$.
To show  this,
the proof of \cite[Lemma 4.3.2]{FB} applies.
(In fact the proof is easier because $\gr^F \V$ is supercommutative.)
Second,
if $\bar{U}$ generates a PBW basis of $\gr^F \V$
then $F$ is strict by Proposition \ref{Pro:2006-05-17-23-34}.
Therefore the assertion follows from 
Theorem \ref{Th;filtered-univesal-env}.
\end{proof}
The map in Theorem \ref{TH:PBW-For-UV}
extends to the 
surjective homomorphism
$\S(\bar{U})\ra \wt{\gr}^F \U(\V)$
of compatible degreewise complete algebra,
which is an isomorphism if $F$ is strict.

If $\V$ is compatibly $\Z_{\geq 0}$-gradable
then one can take the standard filtration 
as the filtration in 
Theorem \ref{TH:PBW-For-UV}.

The following immediately
follows from Theorem \ref{TH:PBW-For-UV}.
\begin{Th}
Let
$\V$,
$F$,
$U$,
$\bar{U}$ be as above.
Then 
the map
\begin{align*}
 \begin{array}{ccc}
 S(\bar{U}) &\rightarrow  &\gr^F \Zhu(\V) \\
\sigma_{p_1}(u^1)\cdot 
\dots \cdot
\sigma_{p_r}(u^r)&\mapsto
&\sigma_{p_1+\dots +p_r}
(u^1_0\cdot\dots \cdot u^r_0)
 \end{array}
\end{align*}
is a surjective algebra
homomorphism.
Further,
if $\bar{U}$ generates a PBW basis 
then
this is
 an isomorphism.
\end{Th}
\subsection{BRST Construction of Vertex Algebras}
\label{subsectin:BRST Construction of Vertex Algebras}
Let $\V=\bigoplus_{i\in \Z}\V^i$ be a 
a 
vertex algebra
with an additional $\Z$-gradation,
which is
shown by the upper index:
$
v_{(n)}\cdot \V^j\subset \V^{i+j}$
for $ v\in \V^i, n\in \Z$,
$
\V^i=\bigoplus_{\Delta\in \Z}\V^i_{-\Delta}$,
where $\V^i_{-\Delta}=\V^i\cap \V_{-\Delta}$.

Suppose there exists an odd operator
$\QO\in \End (\V)$ satisfying the following:
\begin{align}
&\QO^2=0,
\label{eq:2006-08-01-10-40}
\\
& Q\cdot\V^i\subset \V^{i+1},
\quad \QO |0\ket =0,
\quad [\QO,\TT  ]=0,
\label{eq:2006-09-01-15-49}
\\
&
[\QO,Y(v,z)]=Y(\QO v,z),\quad \forall v\in \V.
\label{eq:2006-08-01-10-41}
\end{align}
Then
$(\V, \QO)$ is a cochain complex.
Let  \begin{align*}
H^{\bullet}(\V):=H^{\bullet}(\V,\QO)
=\bigoplus_{i\in \Z}H^i(\V, \QO).
    \end{align*}
From \eqref{eq:2006-08-01-10-41}
it follows that 
the action of 
the field $Y([v],z)$
with $[v]\in H^{\bullet}(\V)$
is well-defined on $H^{\bullet}(\V)$.
Therefore
 $H^{\bullet}(\V)$ is a vertex algebra
and $H^0(\V)$ is a vertex subalgebra of $H^{\bullet}(\V)$.
%The vertex algebra $H^{\bullet}(\V)$ is 
%graded by cohomological degree.
We call $(\V, Q)$ a
{\em weak BRST complex} of  vertex algebras,
and $Q$ the corresponding {\em BRST operator}.
If
$Q$ is compatible with
a
Hamiltonian $H$ of $\V$,
that is, if
\begin{align}
 [\QO, H]=0,
\label{eq:2006-08-01-10-43}
\end{align}
then
$H^{\bullet}(\V)$ is graded by the Hamiltonian $H$.
In this case  $(\V, Q)$ is called
{\em  weak BRST complex of  $\Z$-graded vertex algebras}.
Further,
if $\V$ is M\"obius conformal  and
the  BRST operator $Q$ commutes also  with 
$T^*$ then
$H^{\bullet}(\V)$ is also 
M\"obius conformal.

A {\em BRST complex of vertex algebras}
is a weak BRST complex $(\V, Q)$ of 
%$\Z$-graded 
vertex algebras such that 
$Q$ coincides with the residue $A_{(0)}$ of 
the  filed $Y(A, z)$ associated with some 
odd element  $A\in \V^1$.
In this case
 the conditions 
\eqref{eq:2006-09-01-15-49} and \eqref{eq:2006-08-01-10-41}
are automatically satisfied.
\subsection{Filtration of 
Vertex Algebras
and Spectral Sequences}
\label{Subsection:Filtration and Spectral Sequences}
  Let $(\V, Q)$ be a weak BRST complex of
 vertex algebras.
Let  $F$ be
a filtration 
of $\V$
such that 
\begin{align}
F_{-1}\V=0\quad
\text{and}\quad
\text{$Q \cdot F_p \V\subset F_p \V$ for all $p$.}
\end{align}
Then the action of $Q$ on $\gr^F \V$ is well-defined
and $(\gr^F \V, Q)$ is a weak BRST complex of vertex algebras.
Also,  $F$ induces a filtration of $H^{\bullet}(V)$:
\begin{align}
F_p H^{\bullet}(\V)
=\im:(H^{\bullet}(F_p\V)\rightarrow
H^{\bullet}(\V)).
\end{align}
This is 
an increasing, 
separated, exhaustive filtration 
compatible with the vertex algebra structure
of $H^{\bullet}(\V)$.
Further,
 there is a natural vertex algebra homomorphism
\begin{align}
\gr^F H^{\bullet}(\V)\ra H^{\bullet}(\gr^F \V)
\end{align}
given by the
 correspondence
$
\sigma_p([v])\mapsto [\sigma_p(v)] $.

An increasing filtration $\{F_p \V\}$ 
can be transformed into a decreasing filtration by setting
$F^p\V=F_{-p}\V$.
Thus we have the following.
\begin{Pro}\label{Pro:2006-10-13-13-05}
 Let $\V$, $Q$, $F$ be as above.
Then  there exists a converging spectral sequence
$E_1^{p,q}\Rightarrow H^{\bullet}(\V)$
such that
    \begin{align*}
     E_1^{p,q}=H^{p+q}(\gr^F_{-p}\V),
\quad E_{\infty}^{p,q}=\gr^F_{-p}H^{p+q}(\V).
    \end{align*}
If $H^i(\gr^F \V)=0$ for all $i\ne 0$
then the spectral sequence collapses at $E_1=E_{\infty}$,
and consequently
$H^i(\V)=\zero $ for $i\ne 0$,
and the natural map
$\gr^F H^{0}(\V)\ra H^{0}(\gr^F
 \V)$
is  an   isomorphism 
of vertex algebras.
\end{Pro}

\subsection{BRST Cohomology of Attached Lie Algebras}
Let $(\V, Q)$
be a weak BRST complex of vertex algebras.
The action of $Q$ on $\V$
induces the  adjoint action $\ad Q$ on $\Lie{\V}$:
\begin{align*}
\ad Q \cdot v_{\{n\}}=(Q v)_{\{n\}} \quad \text{for }v\in \V, n\in \Z.
\end{align*}
We have $(\ad Q)^2=0$
and
$(\Lie{\V}, \ad Q)$ can be viewed as a cochain complex.
The cohomology
 $H^{\bullet}(\Lie{\V})=\bigoplus_{i\in \Z}H^i(\Lie{\V})$
is naturally a  Lie algebra.
If $\V$ is graded by $H$
and $Q$ is compatible with $H$
then 
 $H^{\bullet}(\Lie{\V})$ is $\Z$-graded by $\ad H$.

There is a natural
Lie algebra homomorphism
\begin{align}
 \Lie{H^{\bullet}(\V)}\ra H^{\bullet}(\Lie{\V})
\end{align}
given by
the correspondence
$[v]_{\{n\}}
\mapsto [v_{\{n\}}]$
with
$v\in \V$,
$n\in \Z$.
Here 
$[v]_{\{n\}}$
denotes the image 
of $[v]\* t^{n}$
in $\Lie{H^{\bullet}(\V)}$,
where
$[v]$ is
the cohomology class of a cocycle
$v\in \V$,
and  
$[v_{\{n\}}]$ denotes 
the cohomology class of a cocycle
$v_{\{n\}}\in \Lie{\V}$.
(We shall use similar convention throughout the paper).

\begin{Th}\label{Th:2006-05-25-00-33}
Suppose that
$\Vac \V  =\C |0\ket$,
$\Vac {H^{\bullet}(\V)} =\C |0\ket$,
$H^0(\V)\ne 0$
and
$H^i(\V)=0$ for all $i\ne 0$.
Then $H^i(\Lie{\V})=0$
 for $i\ne 0$
the natural map
$\Lie{H^0(\V)}\ra H^0(\Lie{\V})$
is an isomorphism of Lie algebras.
\end{Th}
The following assertion is easily seen.
\begin{Lem}\label{Lem:2005-05-25-00-33}
 Let $\V$ be a vertex algebra such that
$\Vac \V  =\C |0\ket$.
Then 
$\ker (\wpar: \LV\rightarrow \LV )=
\C (|0\ket \* 1)
$.
\end{Lem}
%\begin{proof}
% Let $F_p \LV$ be the subspace of
%$\LV$ spanned be the vectors $v\* f$
%with $v\in \V$,
%$f\in \C[t,t\inv]$ such that $\deg f\leq p$.
%Then this defines an increasing, separated, exhaustive
%filtration on $\LV$
%such that $\wpar F_p \LV\subset F_p \LV$.
%Under the identification
%$\gr^F \LV=\LV$,
%the operator $\wpar$ on $\gr^F \LV$
%is identified
%with $\TT     \* \id$. 
%Then from assumption  it follows that
%$\Ker (\wpar:\gr^F\LV\rightarrow \gr^F \LV)=\C (|0\ket \*1)$.
%This completes the proof.
%\end{proof}
\begin{proof}[Proof of Theorem \ref{Th:2006-05-25-00-33}]
First,
considering  $L(\V)$ as a cochain complex
 with the differential $Q\* \id$,
 we have
\begin{align}\label{eq-in-proof:Borcherds1}
 H^i(\LV)=\begin{cases}
			L(H^0(\V))
&\text{for $i=0$}\\
\zero &\text{for $i\ne 0$.}
			 \end{cases}
\end{align}

Second,
 consider the long exact sequence
associated with the
short exact sequence of cochain complexes
\begin{align*}
 0\rightarrow \C (|0\ket \* 1)
\rightarrow
\LV\rightarrow
\LV/\C (|0\ket \* 1)
\rightarrow 0.
\end{align*}
Then 
by \eqref{eq-in-proof:Borcherds1}
we obtain
%from
% \eqref{eq-in-proof:Borcherds1}
%we obtain 
%$
%H^i(\LV/\C (|0\ket \* 1))=\zero 
%$
%for $i\ne 0,-1$,
%and 
the exact sequences
\begin{align}
& 0\rightarrow H^{i}(\LV/\C (|0\ket \* 1))
\rightarrow 
H^{i+1}(\C (|0\ket \* 1))
\rightarrow 0
\quad \text{for $i\ne -1, 0$,}
\label{eq:2006-10-12-10-33}
\\
& 0\rightarrow H^{-1}(\LV/\C (|0\ket \* 1))
\rightarrow 
H^0(\C (|0\ket \* 1))
\rightarrow L (H^0(\V))
\nno\\
&\qquad \qquad\qquad 
\qquad \qquad\qquad 
\qquad \qquad\quad 
\rightarrow
H^0(\LV/\C (|0\ket \* 1))
\rightarrow 0.
\label{eq:2006-10-12-10-32}
\end{align}

We obviously  have 
that
$H^i(\C (0\ket \* 1))=0$ for $i\ne 0$
and
 \begin{align}
  H^0(\C (0\ket \* 1))
=\C (|0\ket \* 1).
\label{eq:2006-10-12-10-30}
 \end{align}
By \eqref{eq:2006-10-12-10-30}
it follows that  the middle map 
$H^0(\C (|0\ket \* 1))
\rightarrow L (H^0(\V))
$ in \eqref{eq:2006-10-12-10-32} is injective.
Therefore 
\eqref{eq:2006-10-12-10-33}
and \eqref{eq:2006-10-12-10-32} give
\begin{align}
 H^{i}(\LV/\C (|0\ket \* 1))=
\begin{cases}
L (H^0(\V))
/\C (|0\ket \* 1)&\text{for $i=0$}\\
 \zero &\text{for $i\ne 0$.}
\end{cases}
\label{eq:2006-05-25-00-58}
\end{align}

 Third,
by assumption and Lemma \ref{Lem:2005-05-25-00-33},
we have the following exact sequences:
\begin{align}
& 0
{\rightarrow}
\LV/
\C (|0\ket \* 1)
\overset{\wpar}{\rightarrow }
\LV\rightarrow
\Lie{\V}\rightarrow 0,
\label{eq:2006-05-25-00-38}
\\
&0 \rightarrow 
L (H^{0}(\V))
/\C (|0\ket \* 1)
\overset{\wpar}{\rightarrow }
L (H^{0}(V))
\rightarrow
\Lie{H^{0}(\V)}\rightarrow 0.
\label{eq:2006-05-25-01-02}
\end{align}
%Note that 
The \eqref{eq:2006-05-25-00-38}
is an exact sequence of cochain complexes.
Thus we have  the corresponding long exact sequence.
This together with
 \eqref{eq-in-proof:Borcherds1}, \eqref{eq:2006-05-25-00-58}
gives
 $H^{i}(\Lie{\V})=\zero $ for $i\ne 0,-1$
and  the exact sequence
\begin{align*}
 0\rightarrow H^{-1}(\Lie{\V})\rightarrow
L (H^0(\V))
/\C (|0\ket \* \1)
\rightarrow L (H^0(\V) )
\rightarrow H^0(\Lie{\V})
\rightarrow 0.
\end{align*}
However 
the middle map 
$L (H^0(\V))/\C (|0\ket \* \1)
\rightarrow L (H^0(\V))$
is the derivation $\wpar$.
Therefore by
 \eqref{eq:2006-05-25-01-02} this is exact.
Hence
$H^{-1}(\Lie{\V})$ is  zero and
that 
$\Lie{H^0(\V)}$
is isomorphic to $H^0(\Lie{\V})$.
This completes the proof.
\end{proof}
\subsection{BRST Cohomology of 
Current Algebras}
Let $(\V, Q)$
be a weak BRST complex of $\Z$-graded vertex algebras.
The
 action of $\ad Q$ extends to a 
degree-preserving derivation of $\U(\V)$
such that 
$(\ad Q)^2=$0
because $Q$ is odd (cf.\ \eqref{eq:2006-10-12-13-49}).
Thus each
$(\U_N(\V), \ad Q)$ 
 is a cochain complex,
which is 
$\Z$-graded by $\ad H$:
\begin{align*}
 &H^{\bullet}(\U_N(\V))=\bigoplus_{d\in \Z}
 H^{\bullet}(\U_N(\V))_d,\\
&
 H^{\bullet}(\U_N(\V))_d=
\{[v]\in  H^{\bullet}(\U_N(\V));
\ad H\cdot [v]=-d [v]\}.
\end{align*}
This  makes 
a
 projective system
$\{H^{\bullet}(\U_N(\V))_d\}
$
of linear spaces
for each $d$.
Set
\begin{align}
H^{\bullet}(\U(\V))=\bigoplus_{d\in \Z}H^{\bullet}
(\U(\V))_d,
\quad 
\text{where }
 H^{\bullet}(\U(\V))_d=\lim_{\leftarrow \atop N}
H^{\bullet}(\U_N(\V))_d.
\label{eq;2006-11-09-22-47}
\end{align}

There is a natural linear map
\begin{align}
   \U_N(H^{\bullet}(\V))\ra  H^{\bullet}(\U_N(\V)) 
\end{align}
for each $N$
given by the correspondence
\begin{align*}
\left[v^1\right]_{n_1}\dots \left[v^r\right]_{n_r} 
 \mapsto
\left[v^1_{n_1}\dots v^r_{n_r} \right].
\end{align*}

Let  $F$ a filtration 
of $\V$
which is
compatible with $H$
such that 
$F_{-1} \V=0$
and $Q\cdot F_p \V\subset F_p \V$ for all $p$.
Then the corresponding
filtration $\{F_p \U_N(\V)\}$
satisfies 
\begin{align}
 F_{-1}\U_N(\V)=0,\quad \ad Q\cdot F_p \U_N(\V)
\subset F_p \U_N(\V).
\end{align}
Thus we have the following.
\begin{Pro}\label{Pro:oo-2006-10-13-12-41}
 Let $\V$, $Q$, $F$ be as above.
 There exists a converging spectral sequence
$E_1^{p,q}\Rightarrow H^{\bullet}(\U_N(\V))$
such that
    \begin{align*}
     E_1^{p,q}=H^{p+q}(\gr^F_{-p}\U_N(\V)),
\quad E_{\infty}^{p,q}=\gr^F_{-p}H^{p+q}(\U_N(\V)).
    \end{align*}
If $H^i(\gr^F \U_N(\V))=0$ for all $i\ne 0$
then the spectral sequence collapses at $E_1=E_{\infty}$,
and consequently
$H^i(\U_N(\V))=\zero $ for $i\ne 0$,
and the natural map
\begin{align*}
\gr^F H^{0}(\U_N(\V))\ra H^{0}(\gr^F \U_N(\V))
\end{align*}
given by the
 correspondence
$
\sigma_p([a])\mapsto [\sigma_p(a)] $
is a linear isomorphism.
\end{Pro}
\begin{Pro}\label{Pro:hi--}
 Let $\V$, $Q$, $F$ be as above.
Assume that
\begin{itemize}
\item the filtration $F$ of $\V$ is strict;
 \item the filtration of $H^{\bullet}(\V)$
induced by $F$ is also strict;
\item the natural map
$\gr^F H^{\bullet}(\V)\ra H^{\bullet}(\gr^F \V)$
is an isomorphism of vertex algebras;
\item the natural map
$\gr^F H^{\bullet}(\U_N(\V))\ra H^{\bullet}(\gr^F \U_N(\V))$
is an linear isomorphism;
\item the natural map
$\U_N(H^{\bullet}(\gr^F \V))\ra
H^{\bullet}(\U_N(\gr^F \V))$
is a linear isomorphism.
\end{itemize}
Then
the natural map 
$\U_N(H^{\bullet}(\V))\ra H^{\bullet}(\U_N(\V))$
is also a linear isomorphism.
\end{Pro}
\begin{proof}
 The filtration
$\{F_p H^{\bullet}(\V)\}$ of $H^{\bullet}(\V)$
induces a filtration
$\{F_p \U_N(H^{\bullet}(\V))\}$
of $\U_N(H^{\bullet}(\V))$.
The natural map
$\U_N(H^{\bullet}(\V))\ra H^{\bullet}(\U_N(\V))$
preserves the filtration.
Therefore it is sufficient to show that
the induced homomorphism
$\gr^F \U_N(H^{\bullet}(\V))\ra \gr^F H^{\bullet}(\U_N(\V))$
is an isomorphism.

We have
\begin{align}
&\gr^F \U_N(H^{\bullet}(\V))
\cong \U_N(\gr^F H^{\bullet}(\V))\cong \U_N(H^{\bullet}(\gr^F \V)),\\
&\gr^F H^{\bullet}(\U_N(\V))
\cong H^{\bullet}(\gr^F \U_N(\V))
\cong H^{\bullet}(\U_N(\gr^F \V))
\end{align}
by assumption and Thereon \ref{Th;filtered-univesal-env}.
By chasing the isomorphisms one finds that
it is enough to show that
$\U_N(H^{\bullet}(\gr^F \V))\cong H^{\bullet}(\U_N(\gr^F \V))$.
But this is contained in the assumption.
\end{proof}

\subsection{Our Favorite Cases}
\begin{Th}
\label{Th:big-vanishing}
 Let $(\V, Q)$ be a weak BRST complex of
$\Z$-graded vertex algebras,
$H$ the Hamiltonian.
Assume that
there exists an $H$-invariant subspace
$U$
and a filtration $F$
of $\V$
(compatible with $H$)
satisfying the following:
\renewcommand{\theenumi}{\alph{enumi}}
\begin{enumerate}
\item $F_{-1} \V=0$;
\item $(\V, F)$ is quasi-commutative;
\item $Q\cdot F_p \V\subset F_p \V$ for all $p$;
\item  
The image $\bar{U}$ of $U$
in $\gr^F \V$
is preserved by the action of $Q$.
Further,
$H^i(\bar U)=0$ for $i\ne 0$;
\item 
$\bar{U}$ generates a PBW basis  of 
the commutative vertex algebra $\gr^F \V$,
 so that
 $U$ generates a PBW basis of
$\V$.
\end{enumerate}
Then we have the following:
\renewcommand{\theenumi}{\roman{enumi}}
\begin{enumerate}
\item $H^i(\V)=0$ and $H^i(\gr^F \V)=0$
for $i\ne 0$.
\item $F$ induces a strict filtration $\{F_p H^{0}(\V)\}$
of the $\Z$-graded vertex algebra $H^{0}(\V)$ 
such that $(H^{0}(\V), F)$ is quasi-commutative.
\item 
The natural map
$\gr^F H^0(\V)\ra H^0(\gr^F \V)$
is an isomorphism
of vertex algebras.
\item 
The natural map
$H^0(\bar U)\rightarrow H^0(\gr^F\V)$ is injective and 
its image 
 generates a PBW basis of $H^0(\gr^F \V)
=\gr^F H^0(\V)$.
%\item
%Any %$H$-invariant 
%subspace
%$W$ of $H^0(\V)$ whose image in $\gr^F H^0(\V)$
%coincides
%with $H^0(\bar U)$
%generates a PBW basis of $H^0(\V)$.
\item 
For each $N$,
$H^i(\U_N(\V))=0$ with $i\ne 0$.
\item
For each $N$,
the natural map
$\gr^F H^0(\U_N(\V)) \ra H^0(\U_N(\gr^F \V))$
is an  isomorphism.
\item
For each $N$, the natural map
$\U_N(H^0(\V))\ra H^0(\U_N(\V))$
is an isomorphism.
\end{enumerate}
\end{Th}
\begin{proof}
First,
we calculate
$H^{\bullet}(\gr^F \V)$.
By assumption,
there is 
the isomorphism
\begin{align}
 S(\bar{U}\* \C[t\inv]t\inv)\isomap \gr^F \V
\label{eq:2006-10-13-20-39}
\end{align}
 given by \eqref{eq:2006-09-15-PBW}.
Consider
$\bar{U}\* \C[t\inv ]t\inv$ as a cochain complex with the differential
$Q\* 1$.
Extend
 the action of $Q$
to an odd derivation 
of   $S(\bar U\* \C[t\inv]t\inv)$.
Then \eqref{eq:2006-10-13-20-39} 
is a cochain map.
Because
\begin{align}
 H^i(\bar U\* \C[t\inv]t\inv)=\begin{cases}
H^0(\bar U)\* \C[t\inv]t\inv&(i=0),
\\		     \zero &(i\ne 0)
		    \end{cases}
\label{eq:2006-10-13-10-23}
\end{align}
by assumption,
it follows that
\begin{align}
H^i(\gr^F \V)= H^{i}(S( \bar U\* \C[t\inv]t\inv))=
\begin{cases}
 S(H^0( \bar U) \* \C [t\inv]t\inv)
&(i=0),\\
\zero &(i\ne 0),
\end{cases}
\label{eq:2006-08-02-16-02}
\end{align}
see 
 \cite[Lemma 3.2]{KW2003}.
Applying Proposition \ref{Pro:2006-10-13-13-05},
the assertion  \eqref{eq:2006-08-02-16-02}
for $i\ne 0$ proves
(i) and (iii),
while
the assertion  \eqref{eq:2006-08-02-16-02} for $i=0$
proves (iv).
We have also proved (ii)
because the only non-trivial assertion is the
strictness of the filtration 
and this follows from (iv)
and Proposition \ref{Pro:2006-05-17-23-34}.
%(v) follows directly from (ii) and (iv).

The assertions (v), (vi) and (vii) remain to be proved.
For this
 we compute the cohomology
$H^{\bullet}(\gr^F \U_N(\V))$.
By Theorem \ref{TH:PBW-For-UV} 
and (ii),
we have
\begin{align}
 H^{\bullet}(\gr^F \U_N(\V))=
%H^{\bullet}(\U_N(\gr^F \V))=
H^{\bullet}(\S_N(\bar {U})).
\label{eq:2006-10-13-16-00}
\end{align}
Here,
$L(\bar{U})$ is considered as a complex
with the differential $Q\* \1$,
and its action is extended to $\S_N(\bar{U})$
as an odd derivation.
By definition,
we have the exact sequence of cochain maps
\begin{align}
 0\ra S(L(\bar{U}))\sum_{r>N}S(L(\bar{U})_r)
\ra S(L(\bar{U}))
\ra \S_N(\bar{U})
\ra 0.
\label{eq:2006-10-13-10-18}
\end{align}
But the differential is degree-preserving,
therefore by identifying
$\S_N(\bar{U})$ with the subalgebra
of $S(L(\bar{U}))$
as in 
Section \ref{subsection:PBW Thereom for Current Algebras and Zhu
 Algebras},
$\S_N(\bar{U})$ can be identified with 
a subcomplex of 
$S(L(\bar{U}))$.
Namely,
the exact sequence
\eqref{eq:2006-10-13-10-18} splits;
\begin{align}
H^{\bullet}(S(L(\bar{U})))=
H^{\bullet}(S(L(\bar{U}))\sum_{r>N}S(L(\bar{U})_r))
\+ H^{\bullet}(\S_N(\bar{U})).
\label{eq:2006-10-13-10-28}
\end{align}
On the other hand
we have
\begin{align}
H^{\bullet}(S(L(\bar{U})))=\begin{cases}
			    S(L(H^0(\bar{U})))&\text{for $i=0$,}
\\ \zero &\text{for $i\ne 0$.}
			   \end{cases}
\end{align}
in the same manner as \eqref{eq:2006-08-02-16-02}.
Together with \eqref{eq:2006-10-13-10-28},
this gives
\begin{align}
 H^i(\S_N(\bar{U}))=\begin{cases}
		     \S_N(H^0(\bar{U}))&\text{for $i=0$},\\
\zero &\text{for $i\ne 0$.}
		    \end{cases}
\label{eq:2006-11-05-01-30}
\end{align}
Thus we can  apply Proposition \ref{Pro:oo-2006-10-13-12-41}
(ii) to obtain
$H^i(\U_N(\V))=0$ for $i\ne 0$
and the isomorphism
\begin{align}
\gr^F H^0(\U_N(\V))\isomap 
H^0(\gr^F \U_N(\V))\isomap H^0(\U_N(\gr^F \V))
\isomap \S_N(H^0(\bar{U})).
\label{eq:2006-10-13-10-51}
\end{align}
Here 
the second isomorphism follows from 
(ii) and Theorem \ref{Th;filtered-univesal-env}
and the last isomorphism follows from
\eqref{eq:2006-11-05-01-30}
(iv) and Theorem \ref{TH:PBW-For-UV}.
We have proved (v) and (vi).
Finally,
because
$H^0(\bar{U})$ generates a PBW basis of
$H^0(\gr^F \V)$,
\begin{align*}
\U_N(H^0(\gr^F \V))\cong \S_N(H^0(\bar{U})),
\end{align*}
by Theorem \ref{TH:PBW-For-UV}.
This,
together with \eqref{eq:2006-10-13-10-51},
implies that
the natural map 
$\U_N(H^0(\gr^F \V))\rightarrow H^0(\U_N(\gr^F \V))$
is an isomorphism.
Therefore we conclude that 
all the conditions in 
Proposition \ref{Pro:hi--}
are satisfied.
Hence (vii) is proved.
This completes the proof.
 \end{proof} 
If the assumption of
Theorem \ref{Th:big-vanishing} 
is satisfied
then
the isomorphism in (viii)
extends to the isomorphism
\begin{align}
 \U(H^0(\V))\isomap H^0(\U(\V)).
\end{align}
In particular,
%in this case
$H^0(\U(\V))$ has the 
compatible degreewise complete algebra
structure.
\begin{Th}
\label{Th:2006-10-13-11-38}
 Let $(\V, Q)$ be a weak BRST complex of
$\Z_{\geq 0}$-graded vertex algebra,
$H$ the Hamiltonian.
Assume that
there exists an $H$-invariant subspace
$U$
satisfying the following:
\renewcommand{\theenumi}{\alph{enumi}}
\begin{enumerate}
  \item $U$ generates a PBW basis of $\V$;
\item $Q\cdot U\subset U$,
so that $U$ is a subcomplex of $V$.
Further,
$H^i(U)=0$ for $i\ne 0$.
\end{enumerate}
Let $G$ be the standard filtration of $\V$.
Then we have the following:
\renewcommand{\theenumi}{\roman{enumi}}
\begin{enumerate}
\item $H^i(\V)=0$ and $H^i(\gr^G \V)=0$
for $i\ne 0$.
\item The filtration $\{G_p H^{0}(\V)\}$
induced by $G$
is the standard filtration of $H^0(\V)$.
\item 
The natural map
\begin{align*}
\gr^G H^0(\V)\ra H^0(\gr^G \V)
\end{align*}is an isomorphism
of vertex algebras.
\item 
The natural map
$H^0(U)\rightarrow H^0(\V)$ is injective and 
its image 
 generates a PBW basis of $H^0( \V)$.
\item 
For each $N$,
$H^i(\U_N(\V))=0$ with $i\ne 0$.
\item
For each $N$,
the natural map
$\gr^G H^0(\U_N(\V)) \ra H^0(\U_N(\gr^G \V))$
is an linear isomorphism.
\item
For each $N$, 
the natural map
$\U_N(H^0(\V))\isomap H^0(\U_N(\V))$
is a linear isomorphism.
\end{enumerate}
\end{Th}
\begin{proof}
 If we take $G$
as the filtration $F$
in
Theorem \ref{Th:big-vanishing}
then $(U, G)$ satisfies the 
conditions (a)-(e) in 
Theorem \ref{Th:big-vanishing}.
Indeed,
because 
$Q\cdot U\subset  U$,
it follows that
$Q\cdot G_p \V\subset  G_p \V$,
from  the definition of the 
standard filtration
(see  Theorem \ref{Th:Li}).
So we only need to
show that the condition (d)
is satisfied,
that is,
\begin{align}
 H^{\bullet}(U)\cong H^{\bullet}(\bar{U}).
\label{eq:2006-10-12-16-04}
\end{align}
Because  $U$ generate a PBW basis of
$\V$,
$\sigma_{\Delta_u}(u)\ne 0$
for a nonzero  element of $u\in U$.
Because $\Delta_{Qu}=\Delta_{u}$,
 this implies that
$\sigma_{\Delta_u}(Qu)=0$ means
$Qu=0$.
Similarly,
if
$\sigma_{\Delta_u}(u)=Q\sigma_{\Delta_u}(u')=\sigma_{\Delta_u}
(Q u')$ for $u, u'\in U$,
then $u=Q u'$.
Thus \eqref{eq:2006-10-12-16-04} is proved.
Therefore all the assertions except for (ii)
follow directly from
Theorem \ref{Th:big-vanishing}.
But (ii) also,
follows from (iii), (iv)
and the definition of the standard filtration.
\end{proof}
\section{$\W$-Algebras}
\label{Section:definition-of-W-algebra}
Throughout this paper,
$k$ represents a complex number ($=$ the level of the affine Lie algebra
$\bg$
associated with $\sg$).
There is no restriction on $k$
unless otherwise stated.
\subsection{Affine Lie Algebras}
(see \cite{KacBook} for details)
We freely use the notation 
of Section 
\ref{section:quantum reduction for finite-dimensional Lie algebras}.
Let  $\bg$
be the (non-twisted) affine Lie algebra associated with
($\sg$,$(~,~))$.
This is the Lie algebra
given by
\begin{align*}
 \g= \sg\*\C[t,t^{-1}]\+\C K \+ \C \Dg
\end{align*}
with the commutation relations
\begin{align*}
 &[X(m), Y(n)]=[X, Y](m+n)+m\delta_{m+n,0}(X,Y)K,\\
&[\Dg, X(m)]=m X(m),\quad [K,\bg]=0
\end{align*}
for $X, Y\in \sg$,
$m, n\in \Z$.
Here
\begin{align}
 X(n)=X\* t^n
\text{ for $X\in \sg$,
$n\in \Z$.
}\end{align}
The subalgebra $\sg\* \C\subset \g$
is naturally identified with $\sg$.
The invariant symmetric bilinear form $(~|~)$
is 
extended from $\sg$ to $\g$
as follows:
\begin{align*}
 (X(m)|Y(n))=(X,Y)\delta_{m+n,0},
\quad (\Dg| K)=1,\\
(X(m)|\Dg)=(X(m)|K)=(\Dg|\Dg)=(K|K)=0.
\end{align*}

Fix
the  triangular decomposition
$\g=\g_-\+\h\+\g_+$
in the standard way.
That is,
\begin{align*}
&\h=\sh\+\C K\+\C \Dg,\\
&\text{$\g_-=\snn\* \C[t\inv]\+ \sh\* \C[t\inv]t\inv
\+ \snp\*\C[t\inv]t\inv$,}\\
&\text{$\g_+=\snn\* \C[t]t\+ \sh\* \C[t]t
\+ \snp\*\C[t]$.}
\end{align*}
Let
\begin{align}
\dual{\h}=\dual{\sh}\+\C \Lam_0\+\C \delta
\end{align}
be the  dual of $\h$.
Here,
$\Lam_0$ and
$\delta$ are dual elements of $K$ and $\Dg$,
respectively. For $\lam\in\dual{\h}$,
the number
$\bra \lam,K\ket$ is called the
{\em level of $\lam$}.
Let
$\dual{\h}_{\kappa}$
denote the
set of the
weights of level
$k$:
\begin{align}
\dual{\h}_{k}
\teigi\{\lam\in \dual{\h};\bra \lam,K\ket=k\}.
\end{align}
%$$\dual{\h}_{k}=\{\lam\in\dual{\h}
%;  \bra \lam+\rho,K\ket=k\}.$$
Let
$\bar{\lam}$ be the restriction of $\lam\in \dual{\h}$ to $\dual{\sh}$.
We refer to $\slam$ as the {\em classical part}
of $\lam$.

\smallskip

Let $\roots$ be the set of roots of $\g$,
$\roots_+$ the set of positive roots,
$\roots_-=-\roots_+$.
Then,
$\roots=\rroots\sqcup \iroots$,
where
$\rroots$
is
the set of real roots
and
$\iroots$ is the set of imaginary roots.
Let
$\Pi$ be the standard basis of $\rroots$,
$\rroots_{\pm}=\rroots\cap \roots_{\pm}$,
$\iroots_{\pm}=\iroots\cap \roots_{\pm}$.
Then
\begin{align*}
 \rroots_+=\{ \alpha+n\delta;  \alpha\in \sroots_+,\ n\geq 0\}\sqcup
\{ -\alpha+n\delta;  \alpha\in \sroots_+,\ n\geq 1\}.
\end{align*}
Let $Q$ be the root lattice,
$Q_+=\sum_{\alpha\in \proots}\Z_{\geq 0} \alpha\subset
Q$.
We define a partial 
ordering $\mu\leq  \lam$  on $\dual{\bh}$
by
$\lam-\mu\in Q_+$.

\smallskip 

Let
$W\subset GL(\dual{\h})$ be the Weyl group of $\g$
generated by the reflections $s_{\alpha}$
with $\alpha\in \rroots$,
where
$s_{\alpha}(\lam)=\lam-\bra \lam,\alpha\che\ket\alpha$
for $\lam\in \dual{\bh}$.
We have
$W=\sW\ltimes \sQ\che$.
Let $\eW$ be
the extended Weyl group of $\g$:
$\eW=\sW\ltimes \sP\che$.
We write $t_{\mu}$
for the element 
of $\eW$   corresponding to
$\mu\in \sP\che$:
\begin{align*}
t_{\mu}(\lam)=\lam+\bra\lam,K\ket\mu-
\left(\bra\lam,\mu\ket+\frac{1}{2}
|\mu|^2\bra\lam,K\ket
\right)\delta \quad \text{for $\lam\in \dual{\h}$.}
\end{align*}
Let \begin{align}
\eW_+:=\{w\in \eW;\prroots\cap
w\inv(\nrroots)=\emptyset  \}.
\label{eq:2006-11-06-21-48}
    \end{align}
Then
$\eW=\eW_+\ltimes W $.

The dot action of
$\eW$
on $\dual{\h}$ is defined by $w\circ \lam= w(\lam+\rho)-\rho$,
where
$\rho=\bar{\rho}+h\che
\Lam_0\in\dual{\h}$.

%
%
%For $\Lam\in \dual{\h}$,
%let
%\begin{align}
% \text{$R^{\Lam}=\{
%\alpha\in \rroots;
%\bra \Lam+\rho,\alpha\che\ket\in \Z\}$,
%}
%\end{align}$R^{\Lam}_+=R^{\Lam}\cap
%\prroots$.
%$\Pi^{\Lam}=\{\alpha\in R^{\Lam}_+;
%s_{\alpha}(R^{\Lam}_+\backslash \{\alpha\})
%\subset R^{\Lam}_+\}$.
%It is known that
%$R^{\Lam}$ is a subroot
%system of $\rroots$
%%with the basis $\Pi^{\Lam}$
%(\cite{MPBook,KT1}).
%Let 
%\begin{align}
% \text{$W^{\Lam}
%=\bra s_{\alpha};\alpha\in R^{\Lam}\ket\subset W$.
%}
%\end{align}
%The Coxeter group $W^{\Lam}$
%is called the {\em{integral Weyl group}} of $\Lam$.
%We have:
%\begin{align}
%\text{$
%R^{w\circ \Lam}=R^{\Lam}$
%for all $w\in W^{\Lam}$.}
%\end{align}
For $\lam\in \dual{\bh}$,
let
\begin{align}
& \roots(\lam)\teigi \{\alpha\in \rroots; 
\bra \lam+\rho,\alpha\che\ket\in \Z\},\\
&\roots_+(\lam)\teigi \roots(\lam)\cap \proots,\\
&\bW(\lam)\teigi \bra s_{\alpha};  \alpha\in \roots(\lam)\ket
\subset \bW.
\end{align}
One knows that $W(\lam)$ is a Coxeter subgroup of $W$
and it is called the {\em integral Weyl group} of $\lam\in \dual{\bh}$.

\subsection{Graded Duals}
Let $M$ be any
 semisimple
 $\h$-module.
We write
$M^{\lam}$ 
for the weight space of $M$ of weight $\lam$:
\begin{align}
 M^{\lam}\teigi \{m\in M;
hm=\lam(h)m \text{ for all }h\in \h\}.
\label{eq:2006-05-27-22-40}
\end{align}
Let $P(M)=\{\lam\in \dual{\h};M^{\lam}\not=
\zero \}$,
the set 
of weights of $M$.
If
$\dim M^{\lam}<\infty$ for all $\lam\in P(M)$,
then
we define the graded dual $M^*$
of $M$ by
\begin{align}\label{eq:graded-dual}
M^*=\bigoplus_{\lam}\Hom_{\C}(M^{\lam},\C)
\subset \Hom_{\C}(M,\C).
\end{align}
It is clear that $((M)^*)^*=M$.
\subsection{Universal Affine Vertex Algebras Associated with  
Lie Algebras}
\label{subsection:univ-affine}
Define
a $\bg$-module $\Vkg$ by
\begin{align}
 \text{$\Vkg\teigi
U(\g)\*_{U(\sg\*\C[t]\+\C K\+
\C \Dg)}\C %\in \BGG_{k}
$.
}
\end{align}Here,
$\C$
is 
considered as a
$\sg\*\C[t]\+\C K\+
\C \Dg$-module
on which
$\sg\*\C[t]\+
\C \Dg$
acts trivially
and $K$ acts as the multiplication by $k$.

It is
well-known that
the space $V_k(\sg)$
has the 
natural  vertex algebra
structure:
The vacuum vector is given by $|0\ket =1\* 1$;
The translation operator
is defined by the relations 
\begin{align}
 \text{$ \TT   |0\ket =0$,
\quad
$[\TT  , J(n)]=-n J(n-1)$ \quad for $J\in \sg$ and $n\in \Z$;
}
\end{align}
The filed corresponding to 
$J(-1)|0\ket$ with $J\in \sg$
is 
\begin{align}
 J(z):=\sum_{n\in \Z}J(n)z^{-n-1}
\label{eq:2006-11-05-09-31}
\end{align}
and
the fields corresponding to other vectors 
determined by \eqref{eq:2006-11-05-09-31}
in the sense of the Reconstruction Theorem
\cite[Theorem 4.5]{KacBook2}, \cite[Theorem 4.4.1]{FB}:
Explicitly,
they are given by
 \begin{align*}
  Y(J_{a_1}(-n_1-1)\dots J_{a_r}(-n_r-1)|0\ket,z)
 =\frac{1}{n_1!\dots n_r!}:\partial_z^{n_1}J_{a_1}(z)\dots
\partial_z^{n_r} J_{a_r}(z): 
 \end{align*}
for $a_i\in \sI\sqcup \sroots$,
$n_i\geq 0$.
Here $\partial_z=\frac{d}{dz}$
and $:a(z)b(z):$ is the normally ordered product
 \cite[Section 3.1]{KacBook2}.

The vertex algebra $\Vkg$ is M\"obius conformal:
The Hamiltonian
is $-\Dg$;
The operator $T^*$
is defined by
the relations
\begin{align}
 \text{$T^*|0\ket =0$,\quad 
$[T^*, J(n)]=-n J(n+1)$ \quad for  $J\in \sg$
and $n\in \Z$,
}
\end{align}
see \cite[Remark 5.7.d]{KacBook2}.

The M\"obius conformal  vertex algebra
$V_k(\sg)$ is
 called the {\em universal affine vertex algebra at level $k$}
associated with $\sg$.

It is clear that
$\Vkg$ is $\Z_{\geq 0}$-graded:
\begin{align}
 &V_k(\sg)=\bigoplus_{\Delta\in \Z_{\geq 0}}V_k(\sg)_{-\Delta}
\label{eq:grading-old-Vk(sg)} \\
& \Vkg_0=\C |0\ket,
\quad
 \dim V_k(\sg)_{-\Delta}<\infty \text{ for all $\Delta$}.
\end{align}

The standard filtration
of $\Vkg$ coincides with
the one induced by the
standard filtration 
of $U(\sg\* \C[t\inv ]t\inv)$ (in the usual sense).
The subspace
$\Vkg_{-1}$
generates a PBW basis of $\Vkg$ (see Section \ref{subsection:PBW-VA-22-41}).

For $k\ne -h\che$
$\Vkg$ is conformal
by the Sugawara construction:
The conformal vector 
is $\omega_{\sg}=1/2(k+h\che)\sum_{a\in \sI\sqcup
\sroots}J_a(-1)J^a(-1)|0\ket$,
where
$\{J^a; a\in \sI\sqcup \sroots\}$ is a basis
dual to $\{J_a\}$.
Its central charge is $k\dim \sg/(k+h\che)$.
We write 
\begin{align}
 Y(\omega_{\sg},z)=\Tg(z)=\sum_{n\in\Z}\Tg(n)z^{-n-2}.
\end{align}
\subsection{Affine Clifford Algebras and the  Corresponding  Vertex Algebras}
\label{subsection:Clifford}
Define the nilpotent subalgebra
$\Ln{\pm}$ of $\bg$
by
\begin{align}
\text{$\Ln{\pm}
=\sn_{\pm}\*\C[t,t\inv]$.
}\end{align}
We identify $\Ln{\pm}$ with $(\Ln{\mp})^*$ through
$(\cdot|\cdot)$.

Let
$\Cl$ be
  the {Clifford algebra}
associated with
$\Ln{+}\+\Ln{-}
$
and the restriction  of $(~|~)$ to $\Ln{-}\+ \Ln{+}$.
As $\sCl$  (see Section \ref{subsection:fd-BSRT}),
$\bCl$
may be defined as the superalgebra
with 
\begin{align*}
 &\text{odd generators: $\psi_{\alpha}(n)$ \quad with $\alpha\in
 \sroots$,
$n\in \Z$},\\
&\text{relations:
 $[\psi_{\alpha}(m),\psi_{\beta}(n)]=\delta_{\alpha+\beta,0}
\delta_{m+n,0}$
\quad
with $\alpha,\beta\in \sroots$,
$m,n\in \Z$}.
\end{align*}
Here $\psi_{\alpha}(n)$ is regarded as the element of $\bCl$
corresponding to the  vector $J_{\alpha}(n)\in \bg$.

The algebra 
$\Cl$   contains
the Grassmann algebra
$\Lam(\Ln{\pm})
$ of $\Ln{\pm}$ as its subalgebra:
$\Lam(\Ln{\pm})
=\bra \psi_{\alpha}(n);  \alpha\in \sroots_{\pm},n\in \Z\ket$.
One has
\begin{align}
 \Cl=\Lam(\Ln{+})\*\Lam(\Ln{-})
\label{eq:2006-07-31-14-06}
\end{align}
 as  linear spaces.

In view of \eqref{eq:2006-07-31-14-06},
the adjoint action of $\bh$
on $\Ln{\pm}$ induces an action of $\bh$ on $\Cl$.
In particular there is the adjoint action  of 
$\Dg$ on $\bCl$.

Let $\F$  be the irreducible
representation of $\Cl$ generated by a vector $\1$
such that
\begin{align}
 \text{$\psi_{\alpha}(n)\1=0$
\quad 
if
$\alpha+n\delta\in \prroots$.
}
\end{align}
Then
\begin{align}\label{eq:identi-FLn}
 \F=\Lam(\Ln{-}\cap \g_-)\*\Lam(\Ln{+}\cap \g_-)
\end{align}
as  linear spaces.
We regard $\F$ as an $\bh$-module
under this identification:
\begin{align}
\label{eq:qt-space-dec-of-F}
\begin{aligned}
  &\F=\bigoplus_{\lam\in -\bQ_+}\F^{\lam},
\\
& \1\in \F^0,
\quad \psi_{\alpha}(n)\F^{\lam}
\subset \F^{\lam+\alpha+n\delta}\quad
\text{ for $\alpha\in \sroots,\ n\in \Z$.}
\end{aligned}\end{align}
Then $\F^{\lam}$ is finite-dimensional for all $\lam$.

The space
$\F$ is naturally
a conformal vertex superalgebra:
The vacuum  vector  is $\1$;
The translation operator $T$
is defined 
by the relations
\begin{align*}
&T \1=0,\\
&[\TT  , \psi_{\alpha}(n)]=-n \psi_{\alpha}(n-1)
\quad 
\text{for $\alpha\in \sroots_+$, $n\in \Z$,}\\
&[\TT  , \psi_{\alpha}(n)]=-(n-1) \psi_{\alpha}(n-1)
\quad 
\text{for $\alpha\in \sroots_-$, $n\in \Z$;}
\end{align*}
The Hamiltonian is $-\Dg$;
The fields are determined by the Reconstruction Theorem 
and the following:
\begin{align}
& Y(\psi_{\alpha}(-1)\1,z)=
\psi_{\alpha}(z)\teigi \sum_{n\in \Z}\psi_{\alpha}(n)z^{-n-1}
\quad 
\text{for $\alpha\in \sproots$,
}\\
& Y(\psi_{\alpha}(0)\1,z)=
\psi_{\alpha}(z)\teigi \sum_{n\in \Z}\psi_{\alpha}(n)z^{-n}
\text{ for $\alpha\in \sroots_-$;
}
\end{align}
The conformal vector
is
chosen as
$
 \omega_{\F}=\sum_{\alpha\in \sproots}\psi_{-\alpha}(-1)
\psi_{\alpha}(-1)|0\ket
$.
We write
\begin{align}
 Y(\omega_{\F}, z)=L^f(z)=\sum_{n\in\Z}L^f(n)z^{-n-2}.
\end{align}

The vertex superalgebra $\F$ is $\Z_{\geq 0}$-graded:
\begin{align}
 \F=\bigoplus_{\Delta\n \Z_{\geq 0}}\F_{-\Delta},
\quad  \dim \F_{-\Delta}<\infty
\text{ for all $\Delta$}.
\label{eq:2006-11-06-02-41}
\end{align}
Denote by  $U_{\F}$ the subspace of $\F$ spanned by
$\psi_{\alpha}(-1)\1$, $\psi_{-\alpha}(0)\1$ with $\alpha\in\sproots$.
Then $U_{\F}$ generates a PBW basis of $\F$.

%The space  $G_p \F$,  for the standard filtration $G$,
%is spanned by the vectors
%\begin{align*}
% \psi_{\alpha_1}(-m_1)\dots \psi_{\alpha_r}(-m_r)
%\psi_{-\beta_1}(-n_1)\dots \psi_{-\beta_t}(-n_t)\1
%\end{align*}
%with 
%$r\leq p$, 
%$t\geq 0$ is any,
%$\alpha_i,\beta_i\in \sproots$,
%$m_i\geq 1$, $n_i\geq 0$.

\subsection{The BRST Complex of Quantized Drinfeld-Sokolov Reduction}
Recall that 
a tensor product of 
vertex superalgebras
is naturally a vertex superalgebra
(\cite[Lemma 1.3.6]{FB}).

Define a vertex algebra $\bCg$ by
\begin{align}
\bCg
\teigi \Vkg \* \F.
\label{eq:2006-05-17-13-15}
\end{align}
The vacuum vector
$|0\ket \* \1$
is also denoted by $|0\ket$.
The vertex algebra
$\bCg$
is naturally  M\"obius conformal
with the diagonal action of $\mathfrak{sl}_2(\C)$.

We consider $\bCg$
as an $\bh$-module by  the tensor product action.
Then  $-\Dg$ is the Hamiltonian.
The vertex superalgebra
$\bCg$ is clearly $\Z_{\geq 0}$-graded;
 The subspace  
\begin{align}
U=\Vkg_{-1} \+U_{\F} 
\label{eq:2006-10-14-21-09}
\end{align}
generates a PBW basis of $\bCg$.
Here %, as before,
 we have omitted  the tensor product
symbol:
$\Vkg_{-1}=\Vkg_{-1}\* \C \1$,
$U_{\F}=|0\ket \* U_{\F}$.
We have
\begin{align}
 \Vac \bCg=\C |0\ket
\label{2006-10-14-16-45}
\end{align}
by Proposition \ref{Pro:2006-05-17-23-34}.

Define
\begin{align}
 Q_+^{\st}(z)%=\sum_{n\in \Z}Q_{+}(n)z^{-n-1}
= Q_+^{\st}(z)+\chi_+(z)\in (\End \bCg)[[z,z\inv]],
\end{align}
by
\begin{align}
 Q_+^{\st}(z)
&=
\sum_{n\in \Z}Q_{+}^{\st}(n)z^{-n-1}
\nno \\ &
\teigi 
\sum_{\alpha\in \sroots_+}
J_{\alpha}(z)\psi_{-\alpha}(z)-\frac{1}{2}
\sum_{\alpha,\beta,\gamma\in \sroots_+}c_{\alpha,\beta}^{\gamma}
\psi_{-\alpha}(z)\psi_{-\beta}(z)\psi_{\gamma}(z), %\nonumber
 \\
 \chi_+(z)&=\sum_{n\in \Z}\chi_{+}(n)z^{-n}
\teigi 
\sum_{\alpha\in \sroots_+}\schi_+(J_{\alpha})
\psi_{-\alpha}(z),
\end{align}
where $\schi_+$ is the character of $\sn_+$
defined by \eqref{eq;2006-10-14-03-27},
and $c_{\alpha,\beta}^{\gamma}$ is the structure constant 
of $\sg$ as in
Section \ref{subsection:Setting}.
 The $\bd_+^{\st}(z)$ 
and $\chi_+(z)$ are fields corresponding to the vectors
$Q_+^{\st}(-1)|0\ket$
and
$\chi_+(0)|0\ket$,
respectively.

 By abuse of notation,
we set
\begin{align}
  Q_+^{\st}&\teigi (Q_{+}^{\st}(-1)|0\ket)_{(0)}=Q_+^{\st}(0)
\label{eq:def-of-Q-+-st}
\\
&=
\sum_{\alpha\in
\sproots,n\in \Z}
J_{\alpha}(-n){\psi}_{-\alpha}(n)
-\frac{1}{2}\sum_{\ud{\alpha,\beta,\gamma\in
\sproots}{ k+l+m=0}}
c_{ \alpha,\beta}^{\gamma}
{\psi}_{- \alpha}(k){\psi}_{ \beta}(l)
\psi_{\gamma}(m), \nonumber
\\
\chi_+&\teigi (\chi_{+}(0)|0\ket)_{(0)}
=\chi_+(1)=
\sum_{\alpha\in \sroots_+}\schi_+(J_{\alpha})\psi_{-\alpha}(1),
\label{eq:chi-++}\\
\bQ_+&\teigi 
(Q_+^{\st}(-1)+\chi_{+}(0)|0\ket)_{(0)}
=
\bQ_+^{\st}+\chi_+.
				      \end{align}
\begin{Lem}\label{Lem:bd-twice=0}
 We have $(\bd^{\st}_+)^2=\chi_+^2=[\bd^{\st}_+,\chi_+]=0$.
In particular  we have $\bd_+^2=0$.
\end{Lem}
\begin{proof}
 Direct calculation.
\end{proof}
Let $\F=\bigoplus_{i\in\Z}\F^i$ be an additional
$\Z$-gradation of the vertex algebra $\F$ defined 
by 
\begin{align*}
 \deg \1=0, 
\quad \deg \psi_{\alpha}(n)=\begin{cases}
						1&\text{for $\alpha\in
						\sroots_-$}\\ 
-1&\text{for $\alpha\in \sroots_+$}.
					       \end{cases}
\end{align*}
Set 
\begin{align*}
 C_k^i(\bg)=V_k(\bg)\* \F^i \quad \text{ with }i\in \Z.
\end{align*}
This gives a $\Z$-gradation
 of $\bCg$: $\bCg=\bigoplus_{i\in \Z}C_k^i(\bg)$.

By definition,
\begin{align*}
\bd_+\cdot C_k^i(\sg)\subset 
C_k^{i+1}(\sg).
\end{align*}Therefore,
by Lemma \ref{Lem:bd-twice=0},
$(\bCg, \bd_+)$
is a BRST complex of vertex algebras
in the sense of 
Section \ref{subsectin:BRST Construction of Vertex Algebras}.
This complex is called {\em the BRST complex of the
quantized Drinfeld-Sokolov (``$+$'') reduction} \cite{FF_W, FB}.
\begin{Rem} 
As in Section \ref{subsection:fd-BSRT},
$\chi_+$
may be 
 identified 
with the character of the Lie algebra
$\Ln{+}$
defined by
\begin{align*}
 \chi_+(J(n))=\delta_{n,-1}\schi_+(J) \text{\quad
 for $J\in \sn_+$ and $n\in \Z$.}
\end{align*}
\end{Rem}
\subsection{Change of Hamiltonian}
The BRST operator
$Q_+$  above
is not compatible with the Hamiltonian $-\Dg$
of $\bCg$.
Consider the 
weight space decomposition
$
\bCg=\bigoplus_{\lam\in \dual{\bh}}\bCg^{\lam}
$
with respect to the action of $\bh$.
Then
one 
 has
\begin{align}\label{eq:shift-of-weights_+}
 \bd^{\st}_+\cdot \bCg^{\lam}\subset \bCg^{\lam},\quad
\chi_+ \cdot \bCg^{\lam}\subset \sum_{\alpha\in \sPi}\bCg^{\lam-\alpha+\delta}.
\end{align}
Set
\begin{align}
 \Dg_{\new}\teigi \Dg+\srho\che\in \bh
\label{eq:2006-08-08-13-37}
\end{align}
where $\srho\che$ is as in Section \ref{subsection:Setting}.
\begin{Lem}\label{lem:d-commute-new}
The operator $-\DW$ defines a Hamiltonian
on $\bCg$.
The action of $\DW$ on $\bCg$
commutes with that of $\bd_+$.
\end{Lem}
\begin{proof}
The first assertion is easy to
check.
The second assertion follows from \eqref{eq:shift-of-weights_+}.
\end{proof}
We denote by
$\bCg_{\new}$
 the vertex algebra
$\bCg$ equipped with 
the new 
 Hamiltonian
$-\DW$.
We write
\begin{align}\label{eq:new-grad}
 \bCg_{\new}=\bigoplus_{\Delta\in \Z}\bCg_{-\Delta,\new},
\end{align}
where
\begin{align*}
\bCg_{-\Delta,\new}=\{c\in \bCg;  \Dg_{\new}\cdot c=-\Delta c\}.
\end{align*}
Note that  $\bCg_{\new}$ is no more $\Z_{\geq 0}$-graded.
It is compatibly $\Z_{\geq 0}$-gradable by $-\Dg$
(see Section \ref{subsection:VA-sub}).

The M\"obius conformal structure
is changed accordingly:
Let $T^*=T^* \*1 +1\* L^f(1)$. 
Set
\begin{align}
 T^*_{\new}:=T^*-2 \widehat{\srho\che}(1),
\end{align}
where $\widehat{\srho\che}(1)
$ is the operator
on $\bCg$ such that
\begin{align*}
&\widehat{\srho\che}(1) |0\ket=0,\\
&[\widehat{\srho\che}(1), J(n)]=[\srho\che, J](n+1)
+k\delta_{n,-1}(\srho\che|J)\id
\quad \text{for $J\in \sg$, $n\in \Z$},\\
&[\widehat{\srho\che}(1), \psi_{\alpha}(n)]=
\alpha(\srho\che)\psi_{\alpha}(n+1)
\quad\text{for $\alpha\in \sroots$,
$n\in \Z$.}
\end{align*}
The triplet 
$\{T^*_{\new}, -\DW, T\}$
 gives $\bCg$ a M\"obius conformal
structure
 (cf.\ Section \ref{subSection:The Virasoro Field at Non-Critical Level}).
We have
\begin{align}
 [\bd_+ ,T^*_{\new}]=0.
\end{align}
Also,
the action of $T^*_{\new}$ on $\bCg$ is locally nilpotent.

Put 
\begin{align}
 \r:=\C T^*_{\new}+ \C \DW+ \C T.
\end{align}
We consider $\r$
 as a Lie algebra %in $\End \bCg$
isomorphic to $\sl_2(\C)$.

Below
if no confusion can arise
 we  write just $\bCg$
for 
$\bCg_{\new}$.
We write $\bCg_{\old}$
for the  vertex algebra
$\bCg$
with the (old) Hamiltonian
$-\Dg$.
\subsection{Definition of $\W$-Algebras}
\label{subsection:def-of-W-algebras}
The following assertion
was proved by 
B. Feigin and E. Frenkel \cite{FF_W}  for generic $k$,
by J. de Boer and T. Tjin \cite{dBT} for all $k$
in the case that $\sg=\mathfrak{sl}_n$
and by E. Frenkel \cite{FB} for the general case.
\begin{Th}\label{Th:vanishing-vertex-2}
The cohomology
$\Hp{i}(\bCg)$ is zero for all $i \ne 0$.
\end{Th}

Define 
\begin{align}
 \Wg\teigi \Hp{0}(\bCg).
\label{eq:def-of-W-05-08}
\end{align}
This is a M\"obius conformal vertex algebra for all $k\in \C$,
because the action of $\r$ on $\bCg$ commutes with 
the BRST operator $\bd_+$.
The  vertex algebra  $\W_{k}(\sg)$
is called the
{\em $\W$-algebra associated with $\sg$
at level $k$}.

For the later purpose
we  recall the proof of 
Theorem \ref{Th:vanishing-vertex-2}
given by \cite{FB}
in next sections.
\subsection{The Tensor Product Decomposition of 
the Complex $\bCg$}\label{Subsection:Tensor Product Decomposition}
%We keep to write $J(n)$, with $J\in \sg$,
%$n\in \Z$, for the element
%$(J(-1)|0\ket)_{\{n\}}\in \Lie{\bCg}$
%or its image in $\End \bCg$, $\U(\bCg)$,
%$\U_N(\bCg)$.
%Similarly,
%we keep to write $\psi_{\alpha}(n)$,
%with $\alpha \in \sproots$,
%$n\in \Z$,
%for the element $(\psi_{\alpha}(-1)|0\ket)_{\{n\}}
%\in \Lie{\bCg}$,
%and $\psi_{-\alpha}(n)$,
%with $\alpha \in \sproots$,
%$n\in \Z$,
%for the element $(\psi_{-\alpha}(0)|0\ket)_{\{n+1\}}
%\in \Lie{\bCg}$.
Following \cite[Chapter 15]{FB},
we set
\begin{align*}
\widehat{J}_{a}(z)=
\sum_{n\in \Z}\widehat{J}_{a}(n)z^{-n-1}
\teigi J_a(z)-\sum_{\beta,\gamma\in \sproots}
c_{a,\beta}^{\gamma}:\psi_{-\beta}(z)\psi_{\gamma}(z)
:
\end{align*}
for $a\in \sI\sqcup \sroots$.
\begin{Lem}\label{Lem:relations-hat}
 We have the following relations for $m,n\in \Z$:
\begin{enumerate}
 \item $[\hatJ_a(m),\hatJ_b(n)]=\sum_{d}c_{a,b}^d\hatJ_d(m+n)
+(k+h\che)m(J_a,J_b)\delta_{m+n,0}$
if $a,b\in \sI\sqcup \sroots_+$ or 
$a,b\in \sI\sqcup \sroots_-$.
\item $[\hatJ_a(m),\psi_{\beta}(n)]=\sum_{\gamma}c_{a,\beta}^{\gamma}
\psi_{\gamma}(m+n)$
if $a\in \sI\sqcup \sroots_+$, $\beta\in \sroots_+$ or 
$a\in \sI\sqcup \sroots_-$, $\beta\in \sroots_-$.
\end{enumerate}
\end{Lem}
\begin{proof}
 Direct calculation.
\end{proof}
Let
 $\bCg'$
be the subspace of
$\bCg$ spanned by the elements
\begin{align*}
\widehat{J}_{\alpha_{i_1}}(-n_1)\dots
\widehat{J}_{\alpha_{i_p}}(-n_p)
{\psi}_{\alpha_{j_1}}
(-m_1)\dots {\psi}_{\alpha_{j_q}}(-m_q)
|0\ket
\end{align*}
with
$\alpha_{i_s},\alpha_{j_s}\in \sproots$,
$n_i, m_i\in \Z$.
Then,
by 
Lemma \ref{Lem:relations-hat},
$\bCg'$ is a vertex subalgebra
 of $\bCg$.
Similarly,
let $\bCgo$
be the vertex subalgebra of
$\bCg$ spanned by the elements
\begin{align*}
\widehat{J}_{a_1}(-n_1)\dots
\widehat{J}_{a_p}(-n_p)
{\psi}_{\alpha_{j_1}}
(-m_1)\dots {\psi}_{\alpha_{j_q}}(-m_q)
|0\ket
\end{align*}
with
$a_i\in \sI\sqcup \snroots$,
$\alpha_{j_i}\in \sroots_-$,
$n_i, m_i\in \Z$.
Then $\bCg''$ is a vertex subalgebra of $\bCg$
by 
Lemma \ref{Lem:relations-hat}.
\begin{Lem}\label{Lem:rel-2}
 We have the following relations for $n\in \Z$:
\begin{align*}
& 
[\bd_+^{\st},\psi_{\alpha}(n)]=\hatJ_{\alpha}(n)
\text{ and }[\bd_+^{\st},\hatJ_{\alpha}(n)]=0
 \text{ for $\alpha\in \sroots_+$};
\\
&[\chi_+, \psi_{\alpha}(n)]=\chi_+(J_{\alpha}(n))\text{ and }
[\chi_+, \hatJ_{\alpha}(n)]=0
\text{ for $\alpha\in \sroots_+$};
\\
&[\bd_+^{\st}, \psi_{-\alpha}(n)]=-\frac{1}{2}\sum_{\beta,\gamma\in
 \sroots_+\atop
k+l=n}
c_{\beta,\gamma}^{\alpha}\psi_{-\beta}(k)\psi_{-\gamma}(l)
\text{ and }[\chi_+, \psi_{-\alpha}(n)]=0
 \text{ for $\alpha\in \sroots_+$};
\\
&[\bd_+^{\st}, \hatJ_a(n)]=
\sum_{\alpha\in \sroots_+, b\in \sroots_-\sqcup
\sI\atop k+l=n}c_{\alpha, a}^b :\psi_{-\alpha}(l)\hatJ_b(k):n-k_{a}
\psi_{a}(n)
 \text{ for $a\in \sroots_-\sqcup \sI$;}
\\& [\chi_+,\hatJ_a(n )]=\sum_{\beta\in
 \sroots_+}([f,J_a],J_{\beta})\psi_{-\beta}(n+1)
\quad \text{for $a\in \sroots_-\sqcup \sI$,}
%\label{eq:fin-bra-chi_J}
\end{align*}
where $k_a=\begin{cases}
	    k-\sum_{\beta,\gamma\in \sroots_-}c_{a,\beta}^{\gamma}
c_{-a,-\beta}^{-\gamma}&\text{for $a\in \sroots_-$},\\
0&\text{for $a\in \sI$}.
	   \end{cases}$
\end{Lem}
\begin{proof}
 Direct calculation.
\end{proof}
By   Lemma \ref{Lem:rel-2} and the fact that
$ \bd_+ |0\ket =0$,
it follows that 
\begin{align}
&\bd_+^{\st}\cdot \bCg'\subset \bCg',
\quad \chi_+ \cdot \bCg' \subset \bCg'
\quad \text{so }\bd_+\cdot \bCg'\subset \bCg',\\
&\bd_+^{\st}\cdot \bCgo\subset \bCgo,
\quad \chi_+ \cdot \bCgo \subset \bCgo
\quad \text{so }\bd_+\cdot \bCgo\subset \bCgo.
\end{align}
Thus
both 
$\bCg'$ and $\bCGo$ are subcomplexes  of $\bCg$.
Hence
both
$(\bCg', \bd_+)$ and
$(\bCgo, \bd_+)$ are  weak BRST complexes
of vertex algebras
 (see Section \ref{subsectin:BRST Construction of Vertex Algebras}).
The cohomological
gradation takes only non-positive values on $\bCg'$
and only non-negative values
on $\bCgo$:
\begin{align}
 \bCg'=\bigoplus\limits_{i\leq 0}C_k^i(\sg)',
\quad
\bCgo=\bigoplus\limits_{i\geq 0}\bCgoo{i},
\label{eq;2006-09-22-31}
\end{align}
where  $C_k^i(\sg)'=C_k^i(\sg)\cap \bCg'$
and $\bCgoo{i}=C_k^i(\sg)\cap \bCgo$.

\begin{Pro}
\label{Pro:dec-tnsrot-kunneth}
$ $

\begin{enumerate}
 \item Let $U'$ be the subspace of
$\bCg'$ spanned by the vectors 
\begin{align*}
\hatJ_{\alpha}(-1)|0\ket
,\quad
\psi_{\alpha}(-1)|0\ket
 \quad
\text{with $\alpha\in \sroots_+$.}
\end{align*}Then $U'$ generates a PBW basis of $\bCg'$.
\item Let $U''$ be the subspace of $\bCgo$
spanned by the vectors
 \begin{align*}
\wJ_a(-1)|0\ket, \quad
\psi_{-\alpha}(0)|0\ket
\quad
\text{with $a\in \sI\sqcup \sroots_-$,
$\alpha\in \sroots_+$.
}
 \end{align*}
Then $U''$ generates a PBW basis of $\bCgo$.
\item 
The subspace $U'+U''=U'\+U''$
generates a PBW basis of $\bCg$,
and
the multiplication map
\begin{align*}
 \widehat{J}_{\alpha_{i_1}}(-n_1)\dots
{\psi}_{\alpha_{j_1}}
(-m_1)\dots
|0\ket
\*
\widehat{J}_{a_1}(-n_1)\dots
{\psi}_{\alpha_{k_1}}
(-m_1)\dots 
|0\ket
\\ \mapsto 
 \widehat{J}_{\alpha_{i_1}}(-n_1)\dots
{\psi}_{\alpha_{j_1}}
(-m_1)\dots
\widehat{J}_{a_1}(-n_1)\dots
{\psi}_{\alpha_{k_1}}
(-m_1)\dots 
|0\ket
\end{align*}
gives the isomorphism
of  cochain complexes
$\bCg'\* \bCg''\isomap  \bCg$ .
In particular
by 
the K\"{u}nneth theorem
one has
\begin{align*}
\Hp{n}(\bCg)=\bigoplus\limits_{p+q=n}
H^{p}(\bCg')
\* H^{q}(\bCgo)
\end{align*}
 for all $n\in \Z$.
\end{enumerate}
\end{Pro}
It is  clear that
$\bCg'$ and $\bCg''$
is preserved by the action of $\r$.
Therefore
both $\bCg'$ and $\bCg''$ are M\"obius conformal,
and thus so are
$H^{\bullet}(\bCg')$ and $H^{\bullet}(\bCgo)$.

\begin{Pro}[{\cite[15.2.6]{FB}}]\label{Pro:redu-tensor}
The cohomology 
$H^{i}(\bCg')$  is zero for all $i\ne 0$
and one has $H^{0}(\bCg')=\C$.
\end{Pro}
%\begin{proof}
%Let 
%$G_p \bCg'=G_p \bCg\cap \bCg'$.
%This is a standard filtration of $\bCg$
%with respect to the old grading.
%In particular  $\gr^G \bCg'$ is commutative
%and the image $\bar U'$ of $U'$
%in $\gr^G \bCg'$
%generates a PBW basis of $\gr^G \bCg'$.
%
%By \eqref{eq:2005-08-14-12-20}
%one has the corresponding converging
%spectral sequence
%such that
%\begin{align*}
%E_1^{p,q}=H^{p+q}(\gr^G_p \bCg').
%\end{align*}
%From Lemma \ref{Lem:rel-2}
%we have
%\begin{align*}
% \bd_+ \cdot \psi_{\alpha}(-1)|0\ket\equiv
%\hatJ_{\alpha}(-1)|0\ket\pmod{G_0 \bCg'},\quad
%\bd_+\cdot  \hatJ_{\alpha}(-1)|0\ket=0
%\end{align*}
%for all $\alpha\in \sproots$.
%This shows that
%\begin{align}
% \bd_+\cdot   \bar U'\subset  \bar U',\quad
%H^i(\bar U')={0}\quad \forall i.
%\label{eq:2006-08-30-15-30}
%\end{align}
%Hence 
%by  Theorem \ref{Th:2006-08-05-09-46}
%\begin{align}
%H^i(\gr^G \bCg')=\zero \quad \forall  i\ne 0,
%\quad H^0(\gr^G \bCg')=\C |0\ket.
%\label{eq:2006-08-28-14-08}
%\end{align}Therefore the spectral sequence collapses at $E_1=E_{\infty}$
%and  the assertion follows.
%\end{proof}
The following assertion follows directly from
 Propositions \ref{Pro:dec-tnsrot-kunneth}
and \ref{Pro:redu-tensor}.
\begin{Th}[{\cite[Lemma 15.2.7]{FB}}]\label{Th:realization_of_W_by_FB}
The natural embedding $\bCgo\hookrightarrow
\bCg$ induces the isomorphism
$\Hp{\bullet}(\bCgo)\isomap H^{\bullet}(\bCg)$
of  M\"obius conformal vertex algebras.
In particular
\begin{align*}
\Wg\cong H^0(\bCgo)%\cong \Wg %(=\Hp{0}(\bCg))
\end{align*}
as M\"obius conformal vertex algebras.
\end{Th}
\begin{Rem}
\label{Rem:Wg-is-a-sub-of-bCgo}
By \eqref{eq;2006-09-22-31}
we have
$
H^0(\bCgo)
=\{c\in \bCgoo{0};
\bd_+ c=0\}\subset \bCgo$.
Therefore,
by
Theorem \ref{Th:realization_of_W_by_FB}
$\Wg$ can be regarded as a vertex subalgebra of $\bCgo$,
and thus of $\bCg$.
\end{Rem}\subsection{The Filtration $F$ of $\bCgo$}
\label{subsection:Standard-Linearizing-filtration-of-W}
We keep the notation of previous sections.
The vertex algebra
$\bCgo$ is considered to be equipped with
the Hamiltonian $-\DW$.
We have %(with respect to the new grading)
\begin{align}
& \Delta_{\wJ_{-\alpha}(-1)|0\ket}=\height{\alpha}+1,\quad
\Delta_{\psi_{-\alpha}(0)|0\ket}=\height{\alpha}
\quad \text{for $\alpha\in \sproots$},\\
&\Delta_{\wJ_i(-1)|0\ket}=1\quad
\text{for $i\in \sI$.}
\end{align}
Thus
$\bCgo $ is $\Z_{\geq 0}$-graded:
\begin{align}
 \bCg''=\bigoplus_{\Delta\in \Z_{\geq 0}}\bCg''_{-\Delta},
\quad \bCg''_0=\C |0\ket,
\quad \dim \bCg''_{-\Delta}<\infty
\quad \forall \Delta.
\end{align}
It is clear that $\bCgo$ is an $\bh$-submodule of $\bCg$
and we have
\begin{align}
\bCgoo{n}=\bigoplus_{\bra \lam,\srho\che\ket \leq -n}
(\bCgoo{n})^{\lam}
\label{2006-10-14-10-51}
\end{align}for all $n$.
\begin{Pro}\label{Pro;The-Filtration-F}
Define an increasing filtration $F=\{F_p \bCgo\}$ of $\bCgo$
by
\begin{align}
&F_p \bCgo=\sum_{n}F_p \bCgoo{n},\quad
 F_p \bCgoo{n}\teigi 
\bigoplus_{\lam\atop \bra \lam,\srho\che\ket\geq -p-n}
(\bCgoo{n})^{\lam}.
\label{eq:2006-05-28-00-36}
\end{align}
Then we have the following:
\begin{enumerate}
 \item  $F$ is a
(separated, exhaustive,)
strict  filtration of the vertex algebra
$\bCg''$ 
such that
\begin{align*}
& \r\cdot F_p \bCgo\subset
F_p\bCgo\quad \forall p,\\
&F_{-1}\bCgo=0.
\end{align*} 
\item $\gr^F \bCgo\cong \bCgo$ as
M\"obius conformal  vertex algebras.
\item
For all $p$
we have $ \bd_+^{\st}\cdot   
F_p \bCgo\subset  F_{p-1} \bCgo
$ and $\chi_+ \cdot F_p \bCgo\subset  F_{p} \bCgo$,
and thus $\bd_+\cdot F_p \bCgo
\subset F_p \bCgo$.
\end{enumerate}
 \end{Pro}
\begin{proof}
For (i),
the strictness of $F$ follows from
(ii),
together with
Proposition \ref{Pro:dec-tnsrot-kunneth}
(ii) and Proposition \ref{Pro:2006-05-17-23-34}.
Other assertions are easily seen.
(ii)
 is easily seen.
 (iii) 
follows from
\eqref{eq:shift-of-weights_+}.
\end{proof}
The $(\bCgo,\chi_+)$ is also a weak BRST complex of vertex algebras.
By Proposition \ref{Pro;The-Filtration-F},
\begin{align}
 H^{\bullet}(\gr^F \bCg'')\cong H^{\bullet}(\bCgo,\chi_+)
\label{eq:2006-10-13-20-38}
\end{align}
as M\"obius conformal vertex algebras.
\subsection{The Vertex Algebra
$H^{\bullet}(\bCgo, \chi_+)$}
\label{subsection:VV-VV}
According to 
\cite[15.2.9]{FB},
the cohomology $H^{\bullet}(\bCg,\chi_+)$
is easy to calculate,
in view of Theorem \ref{Th:2006-10-13-11-38}:
Let
\begin{align*}
\widehat{P}_i(z)=\sum_{n\in \Z}\widehat{P}_i(n)z^{-n-1},
\end{align*}
with $i\in \sI$,
be the linear combination of
$\widehat{J}_a(z)$
corresponding to $P_i\in \sg^f$ (see  Section \ref{subsection:Setting}).
Similarly,
let  $\widehat{I}_{-\alpha}(z)=\sum_{n\in \Z}
\widehat{I}_{-\alpha}(n)z^{-n-1}$ with $\alpha\in \sproots$
be the 
linear combination of
$J_a(z)$ 
corresponding to $I_{-\alpha}$.
Then from the fifth formula of Lemma \ref{Lem:rel-2}
it follows that
\begin{align}
 &[\chi_+, \widehat{P}_i(n)]=0 \quad \forall i\in \sI,\ n\in \Z,
\label{eq:2006-10-13-14-09}\\
&[\chi_+, \widehat{I}_{-\alpha}(n)]=\psi_{-\alpha}(n+1)
\quad \forall \alpha\in \sproots,\ n\in \Z.
\label{eq:2006-10-13-14-10}
\end{align}

Let $L\sg^f=\sg^f \*\C[t,t\inv]$
(notation Section \ref{subsection:Setting}).
This is a commutative subalgebra of $\bg$.
Set $V(\sg^f)=U(L\sg^f)\cdot |0\ket \subset \Vkg$.
This is a commutative M\"obius conformal vertex subalgebra of $\Vkg$.
Here,
the M\"obius conformal structure
of $V(\sg^f)$ is considered to be given
by $\r$.
%We endow $V(\sg^f)$ with the grading
%defined by the Hamiltonian
%$-\DW$.
%Then $V(\sg^f)$ is $\Z_{\geq 0}$-graded
%and $V(\sg^f)_0=\C |0\ket$.

By \eqref{eq:2006-10-13-14-09}
there is an embedding of vertex algebras
of the following form:
\begin{align}
\begin{array}{ccc}
  V(\sg^f)%(\cong U(\sg^f\*\C[t\inv]t\inv))
&\rightarrow & H^0(\bCgo, \chi_+)\subset \bCgo\\
P_{i_1}
(-n_1)P_{i_2}(-n_2)\dots P_{i_r}(-n_r)|0\ket
&\mapsto &
\widehat{P}_{i_1}(-n_1)\widehat{P}_{i_2}(-n_2)
\dots \widehat{P}_{i_r}(-n_r)|0\ket.
\end{array}
\label{eq:map-E_1}
\end{align}
One sees that
\eqref{eq:map-E_1} is $\r$-equivalent.
\begin{Pro}\label{Pro:E_1}$ $
\begin{enumerate}
 \item \ban{\cite[Lemma 15.2.10]{FB}}
The cohomology
 $H^i(\bCgo,\chi_+)$ is zero 
for all $i\ne 0$
and the map \eqref{eq:map-E_1}
gives an isomorphism 
\begin{align*}
V(\sg^f)\isomap H^0(\bCGo,\chi_+)
%\label{eq:2006-10-13-14-26}
\end{align*}
of M\"obius conformal vertex algebras.
\item Let $N\geq 0$.
 The cohomology $H^i(\U_N(\bCgo),\ad \chi_+)$
is zero for all $i\ne 0$ 
and the isomorphism
$V(\sg^f)\isomap H^0(\bCgo,\chi_+)$ in \ban{i}
induces an isomorphism
\begin{align*}
\U_N(V(\sg^f))\isomap 
 H^0(\U_N(\bCgo),\ad \chi_+)
\end{align*} 
\ban{notation Section {\rm{\ref{subap:CurrentAl-of-VA}}}}.
\end{enumerate}\end{Pro}
\begin{proof}
The set
$\{\widehat{P}_i(-1)|0\ket,\
\widehat{I}_{-\alpha}(-1)|0\ket,\
\psi_{-\alpha}(0)|0\ket;
i\in \sI, \alpha\in \sproots\}$
forms a basis of $U''$ (defined in Proposition
 \ref{Pro:dec-tnsrot-kunneth}),
which generates a PBW basis of $\bCGo$.

By \eqref{eq:2006-10-13-14-09} and \eqref{eq:2006-10-13-14-10},
we have
\begin{align}
& \chi_+\cdot  \widehat{P}_i(-1)|0\ket=0\quad \forall i\in\sI,
\label{eq:2006-08-15-10-25}
\\
&\chi_+\cdot  \widehat{I}_{-\alpha}(-1)|0\ket=\psi_{-\alpha}(0)|0\ket,
\quad \chi_+\cdot \psi_{-\alpha}(0)|0\ket=0
\quad \forall \alpha\in \sproots.
\label{eq:2006-08-15-10-26}
\end{align}
Hence
\begin{align}
& \chi_+\cdot   U''\subset   U'',
\label{eq:2006-08-04-11-54}
\\
& H^i( U'',\chi_+)\cong \begin{cases}
		   \bigoplus_{i\in \sI}\C \widehat{P}_i(-1)|0\ket
&(i=0)\\
0&(i\ne 0).
		  \end{cases}
\label{eq:2006-08-04-11-55}
\end{align}
Therefore 
the  assertion follows 
from Proposition \ref{Pro:dec-tnsrot-kunneth} (ii)
and
Theorem  \ref{Th:2006-10-13-11-38}.
\end{proof}
\subsection{The Filtration $F$ of $\Wg$}
\label{subsection:Poorf-vanihing-FBB}
Let $\{F_p \Wg\}$ be the filtration of $\Wg
=H^0(\bCgo)$
induced by $F$:
\begin{align}
  F_p \Wg:=\im (H^0(F_p \bCgo)\rightarrow H^0(\bCgo)).
\label{eq:2006-10-14-11-50}
\end{align}
Then $F$ is compatible with the M\"obius conformal structure of $\Wg$:
$\r\cdot F_p \Wg\subset F_p \Wg$.

The following assertion 
proves Theorem \ref{Th:vanishing-vertex-2}:
\begin{Th}[{\cite[Theorem 15.1.9]{FB}}]\label{Th:vanishing-vertex-1}
The cohomology 
 $H^i(\bCgo)$  is zero for $i \ne 0$
and there is an isomorphism
$
\gr^F \Wg \cong V(\sg^f)
$
of
 M\"obius conformal vertex algebras.
\end{Th}
\begin{proof}
By 
Propositions \ref{Pro:2006-10-13-13-05},
the assertion follows from 
\eqref{eq:2006-10-13-20-38}
and 
Proposition 
\ref{Pro:E_1} (i).
\end{proof}
By Theorem \ref{Th:vanishing-vertex-1},
$(\Wg, F)$ is quasi-commutative.
\begin{Lem}\label{Lem:2006-11-05-14-45}
$ $

 \begin{enumerate}
  \item The filtration $F=\{F_p \Wg\}$
of $\Wg$ is strict.
\item We have $\Vac \Wg=\C |0\ket$.
 \end{enumerate}
\end{Lem}
\begin{proof}
Because 
\begin{align}
U_{V(\sg^f)}:=\bigoplus_{i\in \sI}P_i(-1)|0\ket
\label{eq:2006-11-05-14-49}
\end{align}
generates a PBW basis of $V(\sg^f)=\gr^F \Wg$,
we have that
%$ \Vac V(\sg^f)=\C |0\ket$,
$\Vac \Wg=\C |0\ket$
and $F$ is strict,
by  Proposition %\ref{Pro:2006-10-10-15-11} and 
\ref{Pro:2006-05-17-23-34}.
\end{proof}
\begin{Rem}
 The filtration $F$ differs from the standard filtration of $\Wg$.
One can characterize $F$
as the finest filtration of $\Wg$
such that $\Wg_{-\Delta}\subset F_{\Delta-1}\Wg$ for all
$\Delta>0$ (compare Theorem \ref{Th:Li} (iii), 
cf.\ Remark \ref{Rem:2006-10-31-10-58}).
The existence of such a filtration $F$
can be proved for any $\Z_{\geq 0}$-graded vertex algebra $\V$
with $\V_0=\C |0\ket$.
%(But $(\V, F)$ is not quasi-commutative for a
% general vertex algebra.
%The $F$ only linearizes $\V$,
%so that $\gr^F \V$ is generated by a Lie (super)algebra.)
 \end{Rem}
\subsection{Quasi-Primary Generators of $\Wg$}
Let $\{F_p \bCgo\}$ 
be the 
filtration
of $\bCgo$ 
defined in  Proposition \ref{Pro;The-Filtration-F},
$\{F_p \Wg\}$ the induced filtration \eqref{eq:2006-10-14-11-50}
of $\Wg$.

Recall the exponent 
 $d_i$
of $\sg$
 with $i\in \sI$,
see Section \ref{subsection:Setting}.
By definition $\widehat{P}_i(-1)|0\ket\in F_{d_i}\bCgo$.
\begin{Pro}\label{Pro:2006-10-14-16-23}
$ $

\begin{enumerate}
 \item  The action of $\r=\C T^*_{\new}+ \C \DW+ \C T$
on $\Wg$ is completely reducible.
 \item  For each $i\in \sI$,
there exists a 
quasi-primary vector $\WW_i\in F_{d_i}\Wg\subset \bCgo$
of conformal weight $d_i+1$
such that 
\begin{align*}
\WW_i=\widehat P_i(-1)|0\ket \pmod{F_{d_i-1}\bCgo}.
\end{align*}
%under the identification $\gr^F \Wg\cong V(\sg^f)$.
%Here $\sigma_p: F_p\Wg\rightarrow \gr^F_p\Wg$
%is the symbol map.
\item Let 
$\WW_1,\dots ,\WW_{l}$
be as in \rm{(}ii\rm{)}.
Then the subspace
$U_{\W}:=\bigoplus_{i\in \sI}\C \WW_i$
generates a PBW basis of $\Wg$.
\end{enumerate}
\end{Pro}
\begin{proof}
(i) Because 
the isomorphism in Theorem \ref{Th:vanishing-vertex-1} 
is $\DW$-equivalent,
it follows that
 there is no vector of conformal weight $1$ in $\Wg$.
This together Lemma \ref{Lem:2006-11-05-14-45} (ii)
shows that
$\Vac \Wg\cap  \im T^*_{\new}=0$.
Thanks to \cite[Proposition 4.9 (b)]{KacBook2},
this proves the assertion.
(ii)
Because 
 $F$ is compatible with the action of
$\r$,
\begin{align*}
 0\longra F_{d_i-1} \Wg
\longra F_{d_i} \Wg\overset{\sigma_{d_i}}{\longra} \gr_{d_i}^F \Wg
=V(\sg^f)
\longra 0
\end{align*}
is an exact sequence of $\r$-modules.
Since $\Wg$ is completely reducible over $\r$
 by (i),
there exist a $\r$-equivalent linear map
$s_{d_i}: \gr^F_{d_i} \Wg\ra F_{d_i} \Wg$
such that $\sigma_{d_i}\circ s_{d_i}=\id$.
The vector $\WW_i=s_{d_i}(P_i(-1)|0\ket)$
satisfies the desired property, because  $P_i(-1)|0\ket$
is quasi-primary.
(iii) is obvious from
the fact that
$U_{V(\sg^f)}$ (defined in  \eqref{eq:2006-11-05-14-49})
generates a PBW basis of $V(\sg^f)$.
\end{proof}
We fix quasi-primary generates
$\WW_1, \dots, \WW_{l}$ 
of $\Wg$ which appeared in Proposition \ref{Pro:2006-10-14-16-23}.
We write $\bar{\WW}_i$ for the image
of $\WW_i$ in $\gr^F \Wg$.
 \begin{Rem}\label{Rem:2006-10-31-10-58}
  From Proposition \ref{Pro:2006-10-14-16-23}
it follows that
\begin{align*}
\begin{aligned}
 & F_p \W_{k}(\sg)\\
&=
\haru\left \{(\WW_{i_1})_{(-n_1)}(\WW_{i_2})_{(-n_2)}\dots
( \WW_{i_r})_{(-n_r)}
|0\ket;  
 \begin{array}{l}
  n_i\geq 1,\ r\in \Z_{\geq 0},\\
d_{i_1}+d_{i_2}+\dots + d_{i_r}\leq p
 \end{array}\right\}.
\end{aligned}
%\label{eq:description-of-filraton-revisde}
\end{align*}
One can prove that
$F_p \Wg$
is spanned by vectors
\begin{align*}
w^1_{(-n_1)}\dots w^r_{(-n_r)}|0\ket
\end{align*}
with
$r\geq 0$, $n_i\geq 1$
and homogeneous vectors $w^1, \dots w^r\in \Wg$
satisfying 
\begin{align*}
 \Delta_{w^1}+\dots +\Delta_{w^r}\leq p-r.
\end{align*}
 \end{Rem}
\subsection{The Lie Algebra $\Lie{\Wg}$}
\begin{Th}$ $

\begin{enumerate}
\item The cohomology $H^i(\bCgo)$  is zero for all $i\ne 0$
and there is a natural  Lie algebra
 isomorphism
$\Lie{\Wg}\isomap H^0(\Lie{\bCgo})$
 \item The cohomology $H^i(\bCg)$  is zero for all $i\ne 0$
and there is a natural 
Lie algebra
 isomorphism
$\Lie{\Wg}\isomap H^0(\Lie{\bCg})$.
\end{enumerate}\end{Th}
\begin{proof}
 By applying Theorem \ref{Th:2006-05-25-00-33},
the assertion follows from 
\eqref{2006-10-14-16-45},
Lemma \ref{Lem:2006-11-05-14-45}
and Theorems \ref{Th:realization_of_W_by_FB}
and   \ref{Th:vanishing-vertex-2}.
\end{proof}
We write
$\WW_i(n)$ for $(\WW_i)_n\in \Lie{\Wg}=H^0(\Lie{\bCgo})=H^0(\Lie{\bCg})$.

\smallskip

Let $\{F_p \Lie{\Wg}\}$ be the filtration of $\Lie{\Wg}$
induced by  \eqref{eq:2006-10-14-11-50}.
Then by Lemma \ref{Lem:2006-11-05-14-45} (i)
and Proposition \ref{Pro:fil-Lie}
we have the isomorphism of Lie algebras
\begin{align}
 \Lie{\gr^F \Wg}\isomap \gr^F \Lie{\Wg}.
\end{align}
We write
$\bar \WW_i(n)$ for $(\bar \WW_i)_n\in \Lie{\gr^F\Wg}=\gr^F \Lie{\Wg}$.

\subsection{The Current Algebra of $\Wg$}
\begin{Th}\label{Th:vanishing-UBCGO}$ $

 \begin{enumerate}
 \item For each $N$,
the cohomology
$H^i(\U_N(\bCgo))$ is zero  for all $i\ne 0$
and the natural map 
$       \U_N(\Wg)\ra H^0(\U_N(\bCgo))$
is a linear isomorphism.
\item 
The cohomology
 $H^i(\Zhu(\bCgo))$ is zero for all  $i\ne 0$
and 
the natural map
$\Zhu(\Wg)\ra  H^0(\Zhu(\bCgo)
)$ is an algebra isomorphism.
 \end{enumerate}
\end{Th}
\begin{proof}
(i)
 By Proposition \ref{Pro:E_1} (ii)
and  \eqref{eq:2006-10-13-20-38},
one can
apply
Proposition \ref{Pro:oo-2006-10-13-12-41} (ii)
%to $\V=\bCgo$ and $F=\{F_p\ \bCgo\}$
to get:
\begin{align}
&H^i(\U_N(\bCgo))=0\quad \text{for $i\ne 0$},\\
&\gr^F H^0(\U_N(\bCgo))\isomap
H^0(\U_N(\gr^F \bCgo)).
\label{eq:2006-10-31-12-46}
\end{align}
Because $\{F_p \bCgo\}$ is strict,
we have
$\U_N(\gr^F \bCgo)\isomap \gr^F \U_N(\bCgo)$
by Theorem \ref{Th;filtered-univesal-env}.
Thus by \eqref{eq:2006-10-31-12-46}
\begin{align}
 \gr^F H^0(\U_N(\bCgo))\isomap H^0(\gr^F \U_N(\bCgo)).
\label{eq:2006-10-31-12-47}
\end{align}
Also, 
by Theorem \ref{Th:vanishing-vertex-1}
 and Proposition \ref{Pro:E_1} (ii),
\begin{align}
 \U_N(H^0(\gr^F \bCgo))
\cong \U_N (V(\sg^f))\cong
H^{0}(\U_N(\gr^F \bCgo)).
\label{eq:2006-10-31-12-50}
\end{align}
Therefore,
by Proposition \ref{Pro;The-Filtration-F} (i),
Lemma \ref{Lem:2006-11-05-14-45} (i),
Theorem \ref{Th:vanishing-vertex-1},
\eqref{eq:2006-10-31-12-47}
and \eqref{eq:2006-10-31-12-50},
 all the assumption  of 
Proposition \ref{Pro:hi--}
are satisfied.
Therefore one has
$\U_N(H^0(\bCgo))\cong H^{0}(\U_N(\bCgo))$.
(ii)
The assertion follows directly from (i).
Indeed,  
$\Zhu(\Wg)$ is a direct summand of the complex $\U_0(\Wg)$
by definition:
\begin{align*}
 \Zhu(\Wg)&=\{a\in \U_0(\Wg);\ad H\cdot a=0\}\\
&=\{a\in H^0(\U_0(\bCgo))
;\ad H\cdot a=0\}=H^0(\Zhu(\bCgo)).
\end{align*}
\end{proof}
By Theorem \ref{Th:vanishing-UBCGO} (i)
it follows that
\begin{align}
 \U(\Wg)\isomap H^0(\U(\bCgo))
\label{eq:2006-11-01-10-15}
\end{align}
(see \eqref{eq;2006-11-09-22-47}).
In particular $\U(\Wg)$ can be considered as s subalgebra
of $\U(\bCgo)$,
and thus of $\U(\bCg)$ by Theorem \ref{TH:PBW-For-UV}.
\begin{Th}\label{Th:2006-08-29-09-53}$ $

\begin{enumerate}
 \item For each $N$,
the cohomology
$H^i(\U_N(\bCg))$ is zero for all  $i\ne 0$
and 
the natural map       $
       \U_N(\Wg)\ra  H^0(\U_N(\bCg)))
$
is a linear isomorphism.
\item
The cohomology $H^i(\Zhu(\bCg))$ is zero for all $i\ne 0$
and the natural map
$\Zhu(\Wg)\ra H^0(\Zhu(\bCg)
)$
is an algebra isomorphism.
\end{enumerate}\end{Th}
\begin{proof}
Because
 (ii)  follows from (i) (see above),
it is sufficient to show (i).

Let $N\geq 0$.
The embedding
$\bCgo\hookrightarrow \bCg$
induces an embedding
\begin{align}
 \U_N(\bCgo)\hookrightarrow \U_N(\bCg),
\label{eq:2006-10-13-16-11}
\end{align}
by Theorem \ref{TH:PBW-For-UV}.
By Theorem \ref{Th:vanishing-UBCGO} (i),
it is sufficient to show that
the map \eqref{eq:2006-10-13-16-11}
induces an isomorphism
\begin{align}
H^{\bullet}(  \U_N(\bCgo))\isomap
H^{\bullet}(\U_N(\bCg)).
\label{eq:2006-10-13-16-29}
\end{align}
To prove this,
consider the standard filtration
$G=\{G_p \bCg\}$ of $\bCg$ 
with respect to the 
(old) grading defined by the Hamiltonian $-\Dg$.
Then 
\begin{align}
G_{-1}\bCg=0.
\end{align}
Also,
 $G_p \bCgo=G_p \bCg\cap \bCgo$ is the
standard filtration of
$\bCgo$ (with respect to the old grading).
%defined by the Hamiltonian $-\Dg$.
Because 
\begin{align}
 Q_+^{\st}(-1)|0\ket \in 
G_1 \bCg,
\quad \chi_+(0)|0\ket \in G_0 \bCg
\end{align}
and  $(\bCg, G)$ is quasi-commutative,
we have
\begin{align}
 Q_+^{\st}\cdot G_p \bCg\subset G_p \bCg,
\quad \chi_+\cdot G_p \bCg\subset G_{p-1}\bCg.
\label{eq:2005-08-14-12-20}
\end{align}
Therefore
there are converging  spectral sequences
\begin{align}
E_r^{p,q}\Rightarrow H^{\bullet}(\U_N(\bCg)),
\quad
 (E'')_r^{p,q}\Rightarrow H^{\bullet}(\U_N(\bCgo)),
\end{align}
see Proposition \ref{Pro:oo-2006-10-13-12-41}.
Because the restriction of
\eqref{eq:2006-10-13-16-11}
gives a cochain map
\begin{align}
G_p \U_N(\bCgo)\ra G_p \U_N(\bCg)
\end{align}
for each $p$,
\eqref{eq:2006-10-13-16-11} induces a map
$(E'')_r\rightarrow E_r$ of spectral sequences.
If
this is an isomorphism for $r=1$,
then this is an isomorphism for all $r\geq 1$,
inducing the desired isomorphism \eqref{eq:2006-10-13-16-29}.
So it is suffice to show the following assertion:
\begin{Pro}
The natural embedding
$\gr^G \bCgo\hookrightarrow  \gr^G \bCg$
induces an isomorphism 
 $H^{\bullet}(\U_N(\gr^G \bCgo))\isomap  H^{\bullet}(\U_N(\gr^G \bCg))$
for all $N$.
\end{Pro}
\begin{proof}
By Proposition \ref{Pro:dec-tnsrot-kunneth}
we  have
\begin{align*}
 \gr^G \bCg\cong \gr^G \bCg'\* \gr^G \bCgo
\end{align*}
as complexes and as (supercommutative) vertex algebras,
where $\gr^G \bCg'$ is the graded vertex algebra
associated with the standard filtration of $\bCg'$
(with respect to the old grading).
For convenience
we put $\V=\gr^G \bCg$,
$\V_1=\gr^G \bCg'$,
and $\V_2=\gr^G \bCgo$,
so that
$
 \V=\V_1\* \V_2
$.
Let $U'$ and $U''$
be as in 
Proposition 
\ref{Pro:dec-tnsrot-kunneth}.
Denote by  ${U}_1$ (resp.\ ${U}_2$) the image of $U'$
in $\V_1$
(resp.\ the image of $U''$ in $\V_2$).
Then $U_1$, $U_2$ and $U_1\+ U_2\subset \V$
generates PBW basis of
$\V_1$,
$\V_2$ and $\V$,
respectively.
By \eqref{eq:2005-08-14-12-20}
the differential
$\bd_+$ acts as $\bd_+^{\st}$ on $\V$,
$\V_1$ and $\V_2$.
%The subspaces
%$U_1$ are $U_2$  are preserved by the action 
%$\bd_+=\bd_+^{\st}$.

Let $F$ be  the filtration 
defined in Proposition \ref{Pro;The-Filtration-F}.
Let $F_p \V_2$ be the image of 
$F_p \bCgo$ in $\V_2$.
Then $\{F_2 \V_2\}$ gives a  filtration of the
supercommutative vertex algebra $\V_2$
such that 
\begin{align}
\bd_+ \cdot F_p \V_2\subset F_{p-1}\V_2
\label{eq:2006-10-13-22-23}
\end{align}
(see Proposition \ref{Pro;The-Filtration-F} (iii)).
Set
\begin{align*}
 F_p \V=\V_1\* F_p \V_2.
\end{align*}
Then this also gives
 a  filtration 
of the supercommutative vertex algebra $\V$
which is compatible with $\bd_+$.
By \eqref{eq:2006-10-13-22-23} one has
\begin{align}
 (\gr^F \V,\bd_+)=(\V_1\* \V_2, \bd_+^{\st}\* \id)
\label{eq:2006-08-30-15-31}
\end{align}
as weak BRST complexes of supercommutative vertex algebras.

Denote by  $\bar{U}_1$ (resp.\ $\bar{U}_2$) the image of $U_1$
in $\gr^F \V_1$
(resp.\ the image of $U_2$ in $\gr^F \V_2$).
Then 
$\bar{U}_2$ generates a PBW basis of $\gr^F \V_2$
and 
$\bar {U}:=\bar {U}_1+\bar{U}_2
=\bar{U}_1\+\bar {U}_2$
generates a PBW basis of the supercommutative
vertex algebra $\gr^F \V$.
By \eqref{eq:2006-10-13-22-23}
and the first formula in Lemma \ref{Lem:rel-2}
it follows that 
\begin{align*}
& \bd_+ \cdot \bar{U}_2=0,\quad
H^i(\bar{U}_2)=\begin{cases}
	      \bar{U}_2&(i=0)\\
\zero &(i\ne 0),
\end{cases}
\\&
 \bd_+ \cdot \bar {U}\subset \bar {U},\quad 
H^i(\bar{U})=\begin{cases}
	      \bar{U}_2&(i=0)\\
\zero &(i\ne 0).
	     \end{cases}
\end{align*}
Therefore,
 by Theorem \ref{Th:big-vanishing} (i), (iv),
it follows that
the $H^i(\V_2)=H^i(\V)=0$ for $i\ne 0$
and
the natural embedding
$\V_2\hookrightarrow \V$ induces an isomorphism
\begin{align}
H^0(\V_2)\isomap H^0(\V).
\label{eq:2006-10-14-01-10-1}
\end{align}
%because $H^0(U_2)$ generates PBW basis of both spaces.
Further, 
by
Theorem \ref{Th:big-vanishing} (v), (vii)
we have
\begin{align}
 &H^i(\U_N(\V_2))=H^i(\U_N(\V))=0\quad \text{ for $i\ne 0$},
\label{eq:2006-10-14-01-10-2}\\
&\U_N(H^0(\V_2))\cong H^0(\U_N(\V_2)),
\label{eq:2006-10-14-01-10-3}\\
&\U_N(H^0(\V))\cong H^0(\U_N(\V)).
\label{eq:2006-10-14-01-10-4}
\end{align}
Therefore,
by
\eqref{eq:2006-10-14-01-10-1},
\eqref{eq:2006-10-14-01-10-2},
\eqref{eq:2006-10-14-01-10-3},
and \eqref{eq:2006-10-14-01-10-4},
it follows that 
$\V_2\hookrightarrow \V$
induces an isomorphism
\begin{align}
 H^{\bullet}(\U_N(\V_2))\isomap H^{\bullet}(\U_N(\V)).
\end{align}
This completes the proof.
 \end{proof}
\end{proof}

\subsection{A Realization of $\U(\bCg)$
and
$\U(\Wg)$}
\label{subsection:theCurrentNew}
Let
\begin{align*}
  \bg'\teigi [\bg,\bg]=\sg\*\C[t,t\inv]\+ \C K\subset \g,
\end{align*}
and denote by  $U_k(\bg')$
 the quotient of $U(\bg')$ by the two-sided ideal generated by
$K-k\id$.
%The tensor product $U(\bg')\* \bCl$ 
%is regarded as a superalgebra
%naturally.  %(cf.\ Section \ref{subsection:fd-BSRT}).
There is a diagonal action of $\ad \DW$ on 
the superalgebra
$U(\bg')\* \bCl$.
Let
\begin{align*}
 (U_k(\bg')\*\bCl)_{d,\new}:=
\{u\in U_k(\bg)\* \bCl;[\DW, u]=d u\}.
\end{align*}
This gives 
$U_k(\bg')\*\bCl$
a graded algebra structure:
\begin{align}
 U_k(\bg')\* \bCl_{\new}
=\bigoplus_{d\in \Z}( U_k(\bg')\* \bCl)_{d,\new}.
\label{eq:2006-09-30-14-04}
\end{align}

There is a  natural 
homomorphism
$\nmap:U_k(\bg')\* \bCl\ra \U(\bCg)$ of graded algebras
induced by the correspondence
\begin{align}
& \bg\ni J(n)\mapsto J(n)(=(J(-1)|0\ket)_{\{n\}})\in \Lie{\bCg}\\
& \bCl\ni \psi_{\alpha}(n)\mapsto \psi_{\alpha}(n)
(=(\psi_{\alpha}(-1)|0\ket))_{\{n\}}\in \Lie{\bCg},\\
& \bCl\ni \psi_{-\alpha}(n)\mapsto \psi_{-\alpha}(n)
(=(\psi_{-\alpha}(0)|0\ket))_{\{n-1\}}\in \Lie{\bCg}
\end{align}
for $J\in \sg$,
$\alpha\in \sproots$,
$n\in \Z$.
Denote by
$\nmap_N$ the composite of $\nmap$
with the natural surjection
$\U(\bCg)\rightarrow \U_N(\bCg)$.

Set 
\begin{align}
\I_{N}:=U_k(\bg')\* \bCl_{\new}
\cdot
\sum_{d>N }(U_k(\bg')\* \bCl)_{d,\new}\subset U_k(\bg')\* \bCl_{\new}.
\end{align}
\begin{Pro}\label{Pro:2006-10-14-00-21}
For each $N\geq 0$,
there is an exact sequence
\begin{align*}
 0\longra \I_{N}\longra U_k(\bg')\* \bCl
\overset{\nmap_N}{\longra}\U_N(\bCg)\longra 0.
\end{align*}
\end{Pro}
\begin{proof}
 Let $\{G_p \U_N(\bCg)\}$ be the filtration of
$\U_N(\bCg)$ induced by the standard filtration $G$ of $\bCg$
(with respect to the old grading).
Let $U$ be the space as in \eqref{eq:2006-10-14-21-09},
which 
 generates a PBW basis of $\bCg$.
Then
by Theorem \ref{TH:PBW-For-UV}
we have 
\begin{align}
\gr^G \U_N(\bCg)\cong \S_N(U).
\label{eq:2006-10-14-00-41}
\end{align}
(Here $U$ is identified with its image in $\gr^G \bCg$.)

Let
$\{G_p (U_k(\bg')\*\bCl_{\new})
\}$ be the filtration of
$U_k(\bg')\*\bCl_{\new}$ defined by
\begin{align}
& G_{-1}(U_k(\bg')\*\bCl_{\new})
=0,\nno\\
&
 G_0(U_k(\bg')\*\bCl_{\new})=\Lam(\Ln{-})
=\bra \psi_{-\alpha}(n);
\alpha\in \sproots,\ n\in \Z\ket,\nno\\
&G_p (U_k(\bg')\*\bCl_{\new})\nno \\
 &=\sum_{n\in \Z, \,
J\in \sg}  J(n)\cdot G_{p-1}(U_k(\bg')\*\bCl_{\new})+
\sum_{n\in \Z,\, \alpha\in \sproots}
\psi_{\alpha}(n)
\cdot  G_{p-1}(U_k(\bg')\*\bCl_{\new})
\nno
\end{align}
for $p\geq 1$.
This defines a filtration 
of $U_k(\bg')\*\bCl_{\new}$ such that
\begin{align}
\gr^G (U_k(\bg')\*\bCl_{\new})
\cong S(\sg\* \C[t,t\inv])\* \Lam(\Ln{+})
\* \Lam(\Ln{-}).
\end{align}
This can be regarded as an isomorphism
\begin{align}
 \gr^G (U_k(\bg')\*\bCl_{\new})\cong S(L(U)).
\label{eq:2006-10-17-13-34}
\end{align}
By construction,
$\nmap$ preserves the filtration.
Therefore,
 $\nmap_N$ induces a map
\begin{align}
 \gr^G (U_k(\bg')\*\bCl_{\new})\ra \gr^G \U_N(\bCg).
\label{eq:2006-10-17-13-32}
\end{align}
But under the identification 
\eqref{eq:2006-10-14-00-41}
and \eqref{eq:2006-10-17-13-34},
the map
 \eqref{eq:2006-10-17-13-32} 
is identical to the natural surjection
$S(L(U))\ra \S_N(U)$.
Therefore \eqref{eq:2006-10-17-13-32} 
is surjective,
and its kernel
is
$S(L(U))\sum_{r>N}S(L(U))_r$,
which is exactly the image
of $\I_{N}$ in $\gr^G (U_k(\bg')\*\bCl_{\new}) $.
This completes the proof.
 \end{proof}
Let
$\wt{U_k(\bg')\* \bCl}_{\new}$
be the standard 
degreewise completion 
(Section \ref{subsectuib:Compatible Degreewise Complete Algebra})
of 
$ U_k(\bg')\* \bCl$ with respect to the grading
\eqref{eq:2006-09-30-14-04}.
%By definition
%\begin{align*}
% &\wt{U_k(\bg')\* \bCl}_{\new}
%=\bigoplus_{d\in \Z}
%(\wt{U_k(\bg')\* \bCl})_{d,\new},\\
%&
%(\wt{U_k(\bg')\* \bCl})_{d}
%=\lim_{\leftarrow \atop N}
%(U_k(\bg')\*\bCl)_{d, \new}/(\I_{N})_d,
%\end{align*}
%where $(\I_{N})_d=\{a\in \I_{N};
%[\DW, a]=-d a\}$. 

The following assertion immediately follows from
Proposition \ref{Pro:2006-10-14-00-21}
and \eqref{eq:2006-11-01-10-15}.
\begin{Pro}\label{Pro:2006-09-28-23-00}There
is an isomorphism
$\U(\bCg_{\new})\cong
\wt{U_k(\bg')\* \bCl}_{\new}
$ as 
compatible degreewise complete algebras.
Thus
there is an isomorphism
$\U(\Wg)\cong 
H^0(\wt{U_k(\bg')\* \bCl}_{\new}, \ad \bd_+)$.
%\item 
%The following map
%\begin{align*}
% \begin{array}{ccccccl}
%  U(\sg)\* \sCl&\ra & (U_k(\bg')\* \bCl)_{0,\new}
% &\overset{\nmap_0}{\ra}  &\Zhu(\bCg) & &\\
%J_{\pm \alpha}&\mapsto &J_{\pm \alpha}(\mp \height \alpha) & & & &
% (\alpha\in \sproots)\\
%\psi_{\pm \alpha}&\mapsto &\psi_{\pm \alpha}(\mp \height \alpha) & & & &
% (\alpha\in \sproots)
%\\
%J_i&\mapsto &J_i-k(\srho\che,J_i)1 & & & &
% (i\in \sI)
% \end{array}
%\end{align*}
%gives the algebra 
%isomorphism
%$U(\sg)\*\sCl\isomap \Zhu(\bCg)$.
%\end{enumerate}
\end{Pro}
%\begin{proof}
%(i) (ii)
%One can directly check 
%the fact that the map is an algebra homomorphism.
%The fact that it is bijective
% follows directly  from 
%Proposition \ref{Pro:2006-10-14-00-21}.
%\end{proof}
\subsection{New Grading vs Old Grading,
and a Realization of $\Zhu(\Wg)$}
\label{subsection:new-vs-old}
(cf.\ \cite[Section 2.2]{FKW})
Let
\begin{align*}
 (U_k(\bg')\*\bCl)_{d,\old}:=
\{u\in U_k(\bg)\* \bCg;[\Dg, u]=d u\}.
\end{align*}
This also gives 
$U_k(\bg')\*\bCl$
a graded algebra structure:
\begin{align}
 U_k(\bg')\* \bCl_{\old}
=\bigoplus_{d\in \Z}( U_k(\bg')\* \bCl)_{d,\old}.
\label{eq:2006-09-30-14-04-oo}
\end{align}
Let
$\wt{U_k(\bg')\* \bCl}_{\old}$
be the standard 
degreewise completion of 
$ U_k(\bg')\* \bCl$ with respect to the grading
\eqref{eq:2006-09-30-14-04-oo}.
\begin{Pro}\label{Pro:2006-10-17-13-01}
We have the following:
\begin{enumerate}
 \item  
$\U(\bCg_{\old})\cong
\wt{U_k(\bg')\* \bCl}_{\old}
$ as compatible degreewise complete algebras.
 \item  
The correspondence $\sg\ni J\mapsto J(0)\in \U_0(\bCg_{\old})$,
$\sCl\ni \psi_{\alpha}\mapsto \psi_{\alpha}(0)\in \U_0(\bCg_{\old})$
gives the algebra isomorphism
$U(\sg)\* \sCl\isomap \Zhu(\bCg_{\old})$.
\end{enumerate}
\end{Pro}
\begin{proof}
The proof of (1) is
exactly in the same manner as 
Proposition \ref{Pro:2006-09-28-23-00}.
(ii) follows from (i).
\end{proof}
For each $N\geq 0$,
$\U_N(\bCg_{\old})$ and $\U_N(\bCg_{\new})$
admit weight space decompositions
with respect to the adjoint action of $\bh$:
\begin{align}
& \U_N(\bCg_{\old})=\bigoplus_{\lam\in \dual{\bh}}
\U_N(\bCg_{\old})^{\lam}.\\
& \U_N(\bCg_{\new})=\bigoplus_{\lam\in \dual{\bh}}
\U_N(\bCg_{\new})^{\lam}
\end{align}

Let $t_{\slam}\in \eW$
with $\slam\in \sP\che$.
Then 
its action on $\bh$
extends to 
an automorphism 
of $U_k(\bg')\* \Cl$
by the following:
\begin{align*}
 &t_{\slam}(J_{\alpha}(n))=J_{\alpha}(n-\bra\alpha,\slam \ket)
\text{ for $\alpha\in \sroots$, $n\in \Z$};
\\
& t_{\slam}(h)=h-k (\slam,h)\text{ for $h\in \sh$};
\\ & t_{\slam}(J_i(n))=J_i(n)\text{ for $i\in \sI$, $n\ne 0$;}
\\ & t_{\slam}(\psi_{\alpha}(n))=\psi_{\alpha}(n-\bra
 \alpha,\slam\ket)
\text{ for $\slam\in \sroots$, $n\in \Z$}.
\end{align*}
One has:
\begin{align*}
&t_{\srho\che}((U_k(\bg')\* \Cl)^{\mu}_{d,\old})
\subset (U_k(\bg')\* \Cl)^{t_{\srho\che}(\mu)}_{d,\new},
\\
&t_{-\srho\che}((U_k(\bg')\* \Cl)^{\mu}_{d,\new})
\subset (U_k(\bg')\* \Cl)^{t_{-\srho\che}(\mu)}_{d,\old}
\end{align*}
This shows that
$t_{\srho\che}$
 and 
$t_{-\srho\che}$
extends to the mutually inverse 
isomorphisms
\begin{align}
& \wh t_{\srho}:\U_N(\bCg_{\old})\isomap \U_N(\bCg_{\new}),\\
& \wh t_{-\srho}:\U_N(\bCg_{\new})\isomap\U_N(\bCg_{\old})
\label{eq:2006-10-07-18-01-pp}
\end{align}
for all $N$,
inducing isomorphisms
\begin{align}
& \wh t_{\srho}:\U(\bCg_{\old})\isomap \U(\bCg_{\new}),
\label{eq:2006-10-07-18-01-oo}\\
& \wh t_{-\srho}:\U(\bCg_{\new})\isomap\U(\bCg_{\old})
\label{eq:2006-10-07-18-01}
\end{align}
of 
compatible degreewise complete algebras.
In particular,
we have
\begin{align}
& \wh t_{\srho}:\Zhu(\bCg_{\old})\isomap \Zhu(\bCg_{\new}),
\label{eq:2006-10-17-13-00}
\\
& \wh t_{-\srho}:\Zhu(\bCg_{\new})\isomap\Zhu(\bCg_{\old}).
\label{eq:2006-10-07-18-01-qq}
\end{align}
The following assertion follows immediately from
 Proposition \ref{Pro:2006-10-17-13-01},
\begin{Pro}
 The map $\wh t_{-\srho\che}$
induces the algebra isomorphism
$\Zhu(\bCg_{\new})\isomap U(\sg)\* \sCl$.
\end{Pro}

We have
\begin{align}
& \wh t_{-\srho\che}(Q_+^{\st})=Q_+^{\st},
\quad  \wh t_{-\srho\che}(\chi_+)=\chi_+',
\\
&\wh t_{-\srho\che}(Q_+)=Q_+':=Q_+^{\st}+
\chi_+',
\label{eq:2006-10-17-03-22}
 \end{align}
where 
\begin{align}
 \dchi_+'\teigi \sum_{\alpha\in \sPi}\psi_{-\alpha}(0).
\end{align}
Clearly,  we have
\begin{align}
& [\dchi_+', \bd^{\st}_+]=0,\quad
(\dchi_+')^2=0,\quad
	  (\bd_+')^2=0.
\end{align}
%which can be checked directly.
Hence each $(\U_N(\bCgo), \ad \bd_+')$ is a cochain complex
and one can define the corresponding cohomology:
\begin{align}
 H^{\bullet}(\U(\bCg_{\old}))
&=\bigoplus_{d\in \Z} H^{\bullet}(\U(\bCg_{\old}))_d,
\nno \\
& H^{\bullet}(\U(\bCg_{\old}))_d:=\lim_{\leftarrow \atop N}
 H^{\bullet}(\U_N(\bCg_{\old})).\nno
\end{align}
In view of 
Proposition \ref{Pro:2006-10-17-13-01} (i),
\begin{align*}
 H^{\bullet}(\U(\bCg_{\old}))
=H^{\bullet}(\oldWU, \ad \bd_+').
\end{align*}

By Proposition \ref{Pro:2006-10-17-13-01} (ii),
one finds that
the subcomplex  
$(\Zhu(\bCgo_{\old}), \ad \bd_+')
$
of
$\U_0(\bCgo_{\old})$
is identical 
to cochain complex 
$(\sC(\sg), \ad \sd_+)$
considered in Section \ref{subsection:fd-BSRT}.
Therefore
by  Theorem \ref{Th:zhu1} 
we have
\begin{align}
 H^0(\Zhu(\bCg)_{\old})\cong \Center(\sg).
\label{eq:2006-10-31-21-33}
\end{align}

By \eqref{eq:2006-10-07-18-01-pp} and \eqref{eq:2006-10-17-03-22},
it is clear that
\begin{align}
H^{\bullet}(\U_N(\bCg)_{\old},\ad \bd_+')=H^{\bullet}(\U_N(\bCg_{\new}),
\ad \bd_+)\quad \forall N.
\end{align}Therefore,
Theorem \ref{Th:2006-08-29-09-53}
and
\eqref{eq:2006-10-31-21-33}
give the following assertion.
\begin{Th}\label{Th:2006-10-19-30-31}
The $H^i(\U(\bCg_{\old}))$
is zero for all $i\ne 0$
and 
$\wh t_{-\srho\che}$ induces the following isomorphisms:
 \begin{enumerate}
  \item 
$\U(\Wg)\isomap H^0(\U(\bCg_{\old}))
=H^0(\oldWU, \ad \bd_+')$
as compatible degreewise topological algebras;
\item
\ban{{cf.\ \cite[Proposition 3.3 (1)]{FKW}}}
 $\Zhu(\Wg)\isomap H^0(\Zhu(\bCg_{\old}))=\Center(\sg)$
as algebras.
 \end{enumerate}
\end{Th}
Below we often identify $\Zhu(\Wg)$ with
$\Center(\sg)$ through the isomorphism $\wh t_{-\srho\che}$.
\smallskip

Similarly as above,
let 
 $\bCgo_{\old}$ be the vertex algebra $\bCgo$
equipped with the Hamiltonian $-\Dg$.
To avoid confusion we may write
$\bCgo_{\new}$
for the same vertex algebra with the Hamiltonian
$-\DW$.
By Theorem \ref{TH:PBW-For-UV}
and Proposition \ref{Pro:dec-tnsrot-kunneth},
$\U(\bCgo_{\old})$
and $\U(\bCgo_{\new})$ are subalgebras of $\U(\bCg_{\old})$
and $\U(\bCg_{\new})$,
respectively.

The  restriction of 
\eqref{eq:2006-10-07-18-01-oo}
and \eqref{eq:2006-10-07-18-01}
give
the mutually inverse  isomorphisms
\begin{align}
& \wh t_{\srho}:\U(\bCgo_{\old})\isomap \U(\bCgo_{\new}),\\
&\wh  t_{-\srho\che}: \U(\bCgo_{\new})\isomap
\U(\bCgo_{\old}).
\label{eq:2006-11-06-00-40}
\end{align}

The $\U(\bCgo_{\old})$ is a subcomplex
of $(\U(\bCg_{\old}),\ad \bd_+')$.
The following assertion follows from 
Theorem \ref{Th:vanishing-UBCGO}.
\begin{Th}\label{Th:2006-10-18-03-53}
The $H^i(\U(\bCgo_{\old}))$
is zero for all $i\ne 0$
and 
$\wh t_{-\srho\che}$ induces the following isomorphisms:
 \begin{enumerate}
  \item 
$\U(\Wg)\isomap H^0(\U(\bCgo_{\old}))$
as compatible degreewise topological algebras;
\item $\Zhu(\Wg)\isomap H^0(\Zhu(\bCgo_{\old}))$
as algebras.
 \end{enumerate}
\end{Th}
\begin{Rem}\label{Rem:2006-11-08-00-42}
 Let 
$\sa$ be the Lie superalgebra
such that
$\sa^{\even}$ is the Lie algebra $\sb_-$,
$\sa^{\odd}$ is the supercommutative Lie superalgebra $\sn_-$,
and $[x,y]=\ad (x)(y)$ for $x\in \sb_-=\sa^{\even}$,
$y\in \sn_-=\sa^{\odd}$.
Then
similarly as Proposition
\ref{Pro:2006-10-17-13-01} (ii)
one sees that
the correspondence
\begin{align*}
& \sa^{\even}\ni J_a\mapsto \wJ_a(0)\in \Zhu(\bCgo_{\old}),\\
&\sa^{\odd}\ni J_{\alpha}\mapsto \psi_{\alpha}(0)\in \Zhu(\bCgo_{\old})
\end{align*}
with $a\in \sI\sqcup \snroots$,
 $\alpha\in \snroots$,
give the isomorphism
\begin{align}
U(\sa)\isomap \Zhu(\bCgo_{\old})
\label{eq:2006-11-08-00-22}
\end{align}
By Theorem \ref{Th:2006-10-18-03-53}
$\Zhu(\Wg)=\Center(\sg)$ can be identified with the commutative
subalgebra of $U(\sa)$
consisting of elements
that commutes with the adjoint action
of $\bd_+'$
under the identification 
\eqref{eq:2006-11-08-00-22}.
\end{Rem}

\subsection{The Conformal Vector at Non-Critical Level}
\label{subSection:The Virasoro Field at Non-Critical Level}
Assume    that $k\ne -h\che$,
where $h\che$ is the dual Coxeter number of $\sg$.
Define
\begin{align*}
L(z)=\sum_{n\in \Z}
L(n)z^{-n-2}\teigi \Tg(z)+\Tf(z)+\frac{d}{d z}
\widehat{\srho\che}(z),
\end{align*}
where
\begin{align*}
&\widehat{\srho\che}(z)=
\sum_{n\in \Z}\widehat{\srho\che}(n)z^{-n-1}
=
{\srho\che}(z)+\sum_{\alpha\in\sproots}
\height \alpha:\psi_{\alpha}(z)\psi_{-\alpha}(z):.
\end{align*}

A short  calculation 
shows that
\begin{align}
 [\bd_+, L(z)]=0.\label{eq:comm_Virasoro-element}
\end{align}
Thus $L(-2)|0\ket$ 
defines an  element of $\Wg$.
Further, one sees that
\begin{align}
 T^*_{\new}=L(1),\quad -\DW=L(0),
\quad T=L(-1)
\label{eq:2006-10-15-00-11}
\end{align}
as elements of $\End \bCg$,
and that
\begin{align}
[L(m),L(n)]=(m-n)L(m+n)+\frac{
m^3-m}{12}\delta_{m+n,0}c(k)\id,
\label{eq:OPE-of-T}
\end{align}
where
\begin{align}
c(k)
&=l-12
\left((k+h\che)|\srho\che|^2
-2 \bra \srho,
\srho\che\ket+ \frac{1}{k+h\che}|\srho|^2
\right),
\label{eq:cc}
\end{align}
see \cite{FF_W, FKW, FB}
(Recall that $l$ is the rank of $\sg$). 
Therefore
$L(-2)|0\ket$ is a conformal 
vector of central charge $c(k)$.
Thus,
the vertex algebra
$\Wg$ 
is
conformal
for $k\ne -h\che$.

\begin{Rem}
Write $k$ as
 $k=-h\che +p/q$.
Then,
as explained in \cite{FKW},
the formula \eqref{eq:cc}
can be written as
\begin{align*}
 c(k)=l-12|q \srho-p\srho\che|^2/pq,
\end{align*}
which in the simply laced case becomes
\begin{align*}
 c(k)=l(1-h\che(h\che+1)(p-q)^2/pq).
\end{align*}
\end{Rem}

\begin{Pro}\label{Pro:image.of.L(0)}
Assume that $k\ne -h\che$.
Then
 under the isomorphism
$\wh t_{-\srho\che}:\Zhu(\W_{k}(\sg))\isomap \Center(\sg)$
in Theorem \rm{\ref{Th:2006-10-19-30-31}}
 \ban{ii},
the image of the Hamiltonian $L(0)$ in $\Zhu(\W_{k}(\sg))$
is mapped to the element
\begin{align*}
 \frac{1}{2(k+h\che)}\Omega
-\frac{1}{2}
\left(
{(k+h\che)}|\srho\che|^2-2\bra\srho,\srho\che\ket
%+\frac{|\srho|^2}{k}
\right)1\in \Center(\sg)
\end{align*} 
Here
$\Omega$ is the Casimir element of $U(\sg)$.
\end{Pro}
\begin{proof}
Direct calculation.
\end{proof}
%\begin{proof}
% \begin{align*}
%L(0)\equiv&\frac{1}{2\kappa}
%\left(\sum_{\alpha\in \sproots}J_{\alpha}(-\height
%\alpha)J_{-\alpha}(\height\alpha)+
%J_{-\alpha}(\height\alpha)J_{\alpha}(-\height\alpha)
%+\sum_i h_ih^i
%\right)\\
%&-\frac{k}{\kappa}h_{\srho\che}
%(0)+\frac{k^2-\kappa^2}{2\kappa}
%|\srho\che|^2+\bra\srho,\srho\che\ket.
%\end{align*}
%\begin{align*}
%\Ttot(0)-h_{\srho\che}(0)\equiv
%&\frac{1}{2\kappa}
%\left(\sum_{\alpha\in \sproots}J_{\alpha}(-\height
%\alpha)J_{-\alpha}(\height\alpha)+
%J_{-\alpha}(\height\alpha)J_{\alpha}(-\height\alpha)
%+\sum_i h_ih^i
%\right)\\&
%-h_{\srho\che}(0)+\frac{1}{\kappa}\sum_{\alpha\in
%\sproots}
%\height \alpha h_{\alpha}(0) -\frac{k h\che}{2\kappa}
%|\srho\che|^2
%\end{align*}
%\end{proof}
\section{Irreducible Highest Weight Representations  of 
 $\W$-Algebras}
\label{section:Simple module over}
\subsection{Verma Modules of $\W$-Algebras}
\label{subsection:VermaModulesofW-Algebras}
We identify $\Zhu(\Wg)$ with the center $\Center(\sg)$
of $U(\sg)$
through Theorem \ref{Th:2006-10-19-30-31} (ii).

Recall that
$\gamma_{\slam}:\Center(\sg)\ra \C$
with $\slam\in \dual{\sh}$
denotes the 
Harish-Chandra homomorphism 
(with respect to the triangular decomposition
$\sg=\sn_-\+ \sh\+\sn_+$), see \eqref{eq:2006-11-03-00-39}.
The one-dimensional 
$\Center(\sg)$-module $\C_{\gamma_{\slam}}$
can be 
naturally regarded as a $\U(\Wg)_{\geq 0}$-module
on which $\U(\Wg)_{>0}$
 acts trivially.

Define
\begin{align}
 \M(\gamma_{{\slam}})\teigi \U(\Wg)\*_{\U(\Wg)_{\geq 0}}
\C_{\gamma_{\slam}}.
\label{eq;def-of-Verma}
\end{align}
The $\Wg$-module $\M(\gamma_{\bar{\lam}})$
is called the {\em{Verma module of $\W_{k}(\sg)$
with highest weight $\gamma_{\bar{\lam}}$}}.
The canonical vector $1\* 1\in \M(\gamma_{\slam})$
is called the {\em highest weight vector} of $\M(\gamma_{\slam})$
and denoted by $|\gamma_{\slam}\ket$.

Define a subspace
$\U^-(\Wg)$ of $\U_0(\Wg)$ (notation Section \ref{subap:CurrentAl-of-VA})
by
\begin{align}
 \U^-(\Wg):=\haru \{ \WW_{i_1}(-n_1)\cdot \dots 
\cdot \WW_{i_r}(-n_r);
r\geq 0,\  i_s\in \sI,\
n_s\geq 1\}.
\label{eq:2006-11-08-14-56}
\end{align}
Then by Theorem \ref{TH:PBW-For-UV}
and Proposition \ref{Pro:2006-10-14-16-23} (iii)
there is a linear isomorphism
\begin{align}
 \begin{array}{ccc}
  \U^-(\Wg)&\isomap & \M(\gamma_{\slam}) \\
u&\mapsto & u|\gamma_{\slam}\ket.
 \end{array}
\label{eq:2006-10-18-16-30}
\end{align}

The Verma module $\M(\gamma_{\slam})$ 
is an admissible $\Wg$-module 
 (see Section \ref{subsection:MOdules-and-Zhu-algberas})
by the natural grading
such that
\begin{align}
 \M(\gamma_{\slam})_{\tp}%=\M(\gamma_{\slam})_0
=\C |\gamma_{\slam}\ket,
\end{align}
which is unique
up to a constant shift.
By \eqref{eq:2006-10-18-16-30}
each subspace $\M(\gamma_{\slam})_d$ is finite-dimensional.

Let $M$ be any graded $\Wg$-module.
Denote by
$S(M)$
the subspace of $M$
spanned by
the  homogeneous vectors $m$
such that $\WW_i(n)m=0$ for all $i\in \sI$,
$n>0$.
Then $S(M)$ is naturally
a $\Zhu(\Wg)$-module
containing
$M_{\tp}$ (notation  Section \ref{subsection:MOdules-and-Zhu-algberas}).
By definition we have the isomorphism
\begin{align}
 \Hom_{\Wg\GrMod}(\M(\gamma_{\slam}), M)
\cong \bigoplus_{d\in \C}\Hom_{\Center(\sg)}(\C_{\gamma_{\slam}}, S(M)_d),
\label{eq:2006-09-28-14-24}
\end{align}
where $S(M)_d=S(M)\cap M_d$.

Let $F=\{F_p \Wg\}$ be the
filtration  of $\Wg$ defined by \eqref{eq:2006-10-14-11-50},
$\{F_p \U(\Wg)\}$,
 the induced filtration
of
$\U(\Wg)$.
Set
\begin{align}
 F_p \M(\gamma_{\slam})\teigi F_p \U(\Wg)\cdot |\gamma_{\slam}\ket.
\end{align}
This defines an increasing,
separated, exhaustive  filtration of
$\M(\gamma_{\slam})$.
By definition,
$F_p \U(\Wg)\cdot F_q\M(\gamma_{\slam})\subset
F_{p+q}\M(\gamma_{\slam})$ for all $p,q$.
Thus
the graded space
$\gr^F \M(\gamma_{\slam})$
is naturally a module over the commutative vertex algebra
$\gr^F \Wg$.
\begin{Pro}\label{Pro:Verma-module-gr}
We have
\begin{align*}
 \gr^F \M(\gamma_{\slam})\cong \U(\gr^F \Wg)\*_{
\U(\gr^F \Wg)_{\geq 0}
}\C
\end{align*}as $\gr^F \Wg$-modules,
where $\C$ is the trivial 
$\U(\gr^F \Wg)_{\geq 0}$-module.
Namely
the
$\gr^F \Wg$-module
$ \gr^F \M(\gamma_{\slam})$
is isomorphic to 
the 
polynomial ring
\begin{align*}
\C[\bar \WW_i(-n_i);{i\in \sI, n_i\geq 1}]
\end{align*}
on which $\bar{\WW}_i(n)$, with $n<0$, acts by the multiplication,
and $\bar{\WW}_i(n)$, with $n\geq  0$, acts trivially.
\end{Pro}
\begin{proof}
 We have the surjective homomorphism
$\U(\gr^F \Wg)\*_{\U(\gr^F \Wg)_{\geq 0}}\C
\ra \gr^F \M(\gamma_{\slam})$
that sends $1\* 1$ to the image of $|\gamma_{\slam}\ket$
in $\gr^F \M(\gamma_{\slam})$.
From \eqref{eq:2006-10-18-16-30}
it follows that this is a bijection.
\end{proof}
\begin{Rem}\label{Rem:2006-11-01-17-14}
Let $\sM^{\dagger}(\slam)$
 be the 
lowest weight Verma module of $\sg$ with lowest weight
$\slam$.
Denote by
$\gamma'_{\slam}:
\Center(\sg)\ra \C$
be the evaluation at 
$\sM^{\dagger}(\slam)$.
Then  we have 
\begin{align}
\gamma'_{\slam}=\gamma_{w_0(\slam)},
\label{eq:2006-11-06-00-49}
\end{align}
and hence
$\M(\gamma'_{\slam})\cong \M(\gamma_{w_0(\slam)})$,
where $w_0$ is the longest element of $\sW$.
Indeed
let
\begin{align}
 \sP_{++}=\{\slam\in \sP;
\bra \lam,\salpha\che\ket
\in \Z_{\geq 0}\
\forall \salpha\in \sproots\}.
\label{eq:2006-11-01-12-47}
\end{align}
Then
for $\slam\in \sP_{++}$
the irreducible finite-dimensional
$\sg$-module $\sL(\slam)$ 
has the lowest weight $w_0(\slam)$.
Thus
\eqref{eq:2006-11-06-00-49}
holds
for all $\lam\in w_0(\sP_{++})$.
Hence it must hold 
for all $\lam\in \dual{\sh}$.
\end{Rem}
\subsection{A Realization of $\M(\gamma_{\slam})$}
\label{subsection-Old-coho-Varma}
Let $\C_{\lam}=\C 1_{\lam}$ with $\lam\in \dual{\bh}$  be the one-dimensional 
representation of $\U(\bCgo_{\old})_{\geq 0}$
such that
\begin{align*}
&\wJ_{-\alpha}(n)1_{\lam}=\psi_{-\alpha}(n)1_{\lam}=0
\quad \text{for $\alpha\in \sproots$,
$n\geq 0$}\\
&\wJ_i(n)1_{\lam}=0\quad \text{for $i\in \sI$, $n>0$},\\
&\wJ_i(0)1_{\lam}=\lam(J_i)1_{\lam}\quad \text{for $i\in \sI$.}
\end{align*}
We consider $\C_{\lam}$  as a  cochain complex
concentrated in degree $0$
on which  the differential $\bd_+'$  acts trivially.
Let $K(\lam)$ 
be a $\bCgo_{\old}$-module
defined by
\begin{align}
 K(\lam)=\U(\bCgo_{\old})\*_{\U(\bCgo_{\old})_{\geq 0}}\C_{\lam}.
\end{align}
 This  is  an admissible $\U(\bCgo)$-module with 
the natural grading such that
$K(\lam)_{\tp}=\C\* \C_{\lam}$.
We regard $K(\lam)$ as a cochain complex with
the differential $\bd_+'$,
whose action is defined by
the rule $\bd_+'(u\* 1)=[\bd_+', u]\* 1$ for $u\in \U(\bCgo)$.

By 
\eqref{eq:2006-11-01-10-15}
and \eqref{eq:2006-11-06-00-40},
$H^{\bullet}(K(\lam))$
can be viewed as a $\Wg$-module
by the following action:
\begin{align}
u\cdot [v]= [\wh t_{-\srho\che}(u)v]\quad
\text{for $u\in \U(\Wg)$
and $[v]\in H^{\bullet}(K(\lam))$.
}
\end{align}
Then
  $H^{\bullet}(K(\lam))$ is an admissible $\Wg$-module.
By construction,
\begin{align*}
H^{i}(K(\lam))_{\tp}=H^{i}(K(\lam)_{\tp})=\begin{cases}
					   \C &\text{$i=0$},\\
\zero &\text{$i\ne 0$}.
					  \end{cases}
\end{align*}
This
is a $\Zhu(\Wg)$-module.

%Recall that $w_0$ denotes the longest element
%of $\sW$ (see Section \ref{subsection:Setting}).
\begin{Pro}\label{Pro:2006-10-19-02-59}
 There is an isomorphism
$H^0(K(\lam))_{\tp}\cong \C_{\gamma_{\lam}}$
of  $\Zhu(\Wg)$-modules.
\end{Pro}
\begin{proof}
We shall use the identification 
\begin{align}
\Zhu(\Wg)&=
H^0(\Zhu(\bCgo_{\old}), \ad \bd_+')\\
&
=H^0(\Zhu(\bCg_{\old}),\ad \bd_+')
=H^0(U(\sg)\*\sCl, \ad \sd_+)=
\Center(\sg)
\label{eq:2006-11-01-10-29}
\end{align}
through Theorems \ref{Th:2006-10-19-30-31}
and
\ref{Th:2006-10-18-03-53}.

Observe
that  $\sM^{\dagger}(\slam +2\srho)$
has the central character
\begin{align}
 \gamma'_{\slam+2\srho}=\gamma_{w_0(\slam+2\srho)}
=\gamma_{w_0\circ \slam}=\gamma_{\slam}
\end{align}
(notation Remark \ref{Rem:2006-11-01-17-14}).

The superalgebra
$\Zhu(\bCg_{\old})=U(\sg)\* \sCl$
acts on $\sM^{\dagger}(\slam+2\rho)\* \Lam(\sn_+)$
as in Section \ref{subsection:WhittakerFunctor},
and by \eqref{eq:2006-11-01-10-29}
this induces the action of
$\Zhu(\Wg)$
 on the space 
\begin{align*}
H_{0}(\sn_+, \sM^{\dagger}(\slam+2\srho)\* \C_{\schi_+})=
H_0(\sM^{\dagger}(\slam+2\srho)\* \Lam(\sn_+), \sd_+).
\end{align*}
Let
$|\slam+2\srho\ket^{\dagger}
=v_{\slam+2\srho}^{\dagger}
\*1$,
where
$v_{\slam+2\srho}^{\dagger}$ is the lowest weight vector of 
$\sM^{\dagger}(\slam+2\srho)$.
Then $\C |\slam+2\srho\ket^{\dagger}$
is a subcomplex of
$(\sM^{\dagger}(\slam+2\srho)\* \Lam(\sn_+), \sd_+)$.
From (the proof of) Theorem \ref{Th:image-of-simple-f-d} (i),
it follows that
the natural embedding
$\C |\slam+2\srho\ket^{\dagger}
\hookrightarrow
\sM^{\dagger}(\slam+2\srho)\* \Lam(\sn_+)$
induces an isomorphism
\begin{align}
H_0(\C |\slam+2\srho\ket^{\dagger})
\cong 
H_0(\sn_+, \sM^{\dagger}(\slam+2\srho)\*
 \C_{\schi_+})=\C_{\gamma_{\slam}}
\label{eq:2006-11-06-01-11}
\end{align}as $\Center(\sg)$-modules.

%Therefore it is sufficient to show that
%\begin{align*}
%H^0(K(\lam))_{\tp}\cong 
%H_0(\C |\slam+2\srho\ket^{\dagger})
%\label{eq:2006-11-01-10-33}
%\end{align*}
%as $\Zhu(\Wg)$-modules.

Consider $\sM^{\dagger}(\slam+2\srho)\* \Lam(\sn_+)$
as a $\Zhu(\bCgo_{\old})$-module 
by the correspondence
\begin{align}
& \wJ_a(0)\mapsto J_a-\sum_{\beta,\gamma\in\sproots}
c_{a,\beta}^{\gamma}\psi_{-\beta}\psi_{\gamma}
\quad \text{for $a\in \snroots\sqcup \sI$},\\
& \psi_{\alpha}(0)\mapsto \psi_{\alpha}
\quad\text{for $\alpha\in \snroots$}
\end{align}
(see Remark \ref{Rem:2006-11-08-00-42}).
Then
a short calculation shows  that
 \begin{align}
  &\wJ_i(0)|\slam+2\srho\ket^\dagger=\lam(J_i)|\slam+2\srho\ket^{\dagger}
\quad  \text{for $i\in \sI$},\\
 & \wJ_{\alpha}(0)|\slam+2\srho\ket^\dagger=0
\quad  \text{for $\alpha\in \snroots$}.
 \end{align}
 Therefore
the correspondence
\begin{align}
 \sM^{\dagger}(\slam+2\srho)\* \Lam(\sn_+)
\ni |\slam+2\srho\ket^{\dagger}
\mapsto 1_{\lam}\in K(\lam)_{\tp}
\label{eq:2006-11-03-00-30}
\end{align}
gives 
a
$\Zhu(\bCgo_{\old})$-equivalent
cochain map.
Here 
$ \sM^{\dagger}(\slam+2\srho)\* \Lam(\sn_+)$
is considered as a cochain complex
by reversing the homological gradation.
Therefore,
it follows that
\eqref{eq:2006-11-03-00-30}
induces an isomorphism
\begin{align*}
H_0(\sM^{\dagger}(\slam+2\srho)\*
 \Lam(\sn_+),
\sd_+)
\isomap H^0(K(\lam))_{\tp}.
\end{align*}
By \eqref{eq:2006-11-06-01-11}
this completes the proof.
\end{proof}
\begin{Pro}\label{Pro:2006-10-19-03-43}
 We have $H^i(K(\lam))=0$ for $i\ne 0$
and $H^0(K(\lam))\cong \M(\gamma_{\slam})$
as $\Wg$-modules.
\end{Pro}
\begin{proof}
By \eqref{eq:2006-09-28-14-24}
and Proposition \ref{Pro:2006-10-19-02-59}
there is a unique $\Wg$-module homomorphism
\begin{align*} 
\phi:\M(\gamma_{\slam})\ra H^0(K(\lam))
\end{align*}
that sends $|\gamma_{\slam}\ket$
to $[1\* 1_{\lam}]\in  H^0(K(\lam))$.
%It is  sufficient to show that 
%$\phi$ is a bijection.

 Let $F$ be the filtration of 
$\bCgo$ defined 
in Section \ref{subsection:Standard-Linearizing-filtration-of-W},
$\{F_p \U(\bCgo_{\old})\}$
the  induced filtration.
Define
\begin{align}
 F_p K(\lam):=F_p \U(\bCgo_{\old})\cdot K(\lam)_{\tp}.
\end{align}
Then this defines an increasing,
separated, exhaustive
filtration of $K(\lam)$
such that
\begin{align}
& \bd_+^{\st}\cdot F_p K(\lam)\subset F_{p-1} K(\lam),
\quad
 \chi_+'\cdot F_p K(\lam)\subset F_{p} K(\lam),\quad
F_{-1}K(\lam)=0.
\label{2006-10-19-02-05}
\end{align}
Thus there is a converging spectral 
sequence $E_r^{p,q}\Rightarrow H^{\bullet}(K(\lam))$
such that
\begin{align}
 E_1^{p,q}=H^{p+q}(\gr^F_{-p} K(\lam)).
\end{align}
By construction,
we have
\begin{align}
 H^{\bullet}(\gr^F K(\lam))=H^{\bullet}(K(\lam),\chi_+').
\label{eq:2006-11-01-11-27}
\end{align}

From the
commutation relations
\begin{align}
& [\chi_+',\hatJ_a(-n )]=\sum_{\beta\in
 \sroots_+}([f,J_a],J_{\beta})\psi_{-\beta}(-n)
\quad \text{for $a\in \sroots_-\sqcup \sI$},
\label{eq:2006-10-18-04-10-a}
\\
&[\chi_+', \psi_{-\alpha}(-n)]=0
\quad \text{for $\alpha\in \sroots_+$,}
\label{eq:2006-10-18-04-10-b}
\end{align} 
one can show exactly in the same manner as 
 Proposition \ref{Pro:E_1}
that
\begin{align}
 H^i(K(\lam),\chi_+)=0\quad \text{ for $i\ne 0$}
\label{eq:2006-11-01-11-28}
\end{align}
and that there is a linear
 isomorphism
\begin{align}
 S(\sg^f\* \C[t\inv]t\inv)\isomap H^0(K(\lam),\chi_+')
\end{align}
given by 
 the correspondence
\begin{align*}
(P_{i_1}\* t^{-n_1})\dots (P_{i_r}\* t^{-n_r})
\mapsto
\widehat{P}_{i_1}(-n_1)\dots
\widehat{P}_{i_r}(-n_r)\*1.
\end{align*}

By \eqref{eq:2006-11-01-11-27} and 
\eqref{eq:2006-11-01-11-28}
it follows that
our spectral sequence collapses at $E_1=E_{\infty}$,
and consequently we have 
$H^{i}(K(\lam))=0$ for $i\ne 0$
and the isomorphism
\begin{align}
\gr^F H^0(K(\lam))\cong S(\sg^f\* \C[t\inv]t\inv).
\label{eq:2006-11-01-11-31}
\end{align}

Finally,
$\phi(F_p \M(\gamma_{\slam})\subset 
F_p H^0(K(\lam))$ by construction.
Thus we have a $\gr^F\Wg$-module homomorphism
\begin{align}
 \gr^F \M(\gamma_{\slam})
\ra \gr^F H^0(K(\lam)).
\label{eq;2006-10-19003020}
\end{align}
But from \eqref{eq:2006-11-01-11-31}
and Proposition \ref{Pro:Verma-module-gr}
it follows that \eqref{eq;2006-10-19003020} is a bijection.
Thus this completes the proof.
\end{proof}
\subsection{Irreducible Highest Weight Representations
of  $\W$-Algebras}\label{subsection:simple}
Let $\slam\in \dual{\sh}$.
Denote by
$\WN(\gamma_{\slam})$
 the 
the sum of all  
proper
graded
 $\W_{k}(\sg)
$-submodules of
$
\M(\gamma_{\slam})$.
Define 
the graded $\Wg$-module $\why(\gamma_{\slam})$
by
\begin{align}
 \why(\gamma_{\slam})
\teigi \M(\gamma_{\slam})/\WN(\gamma_{\slam}).
\end{align}
The image of $|\gamma_{\slam}\ket \in \M(\gamma_{\slam})$
in $\why(\gamma_{\slam})$
is also denoted by $|\gamma_{\slam}\ket$.
We have
$ \why(\gamma_{\slam})_{\tp}
=\C |\gamma_{\slam}\ket
$.
\begin{Th}\label{Th:classification-Simples-W}
The
$\Wg$-modules $\why(\gamma_{\slam})$ is a unique simple graded quotient of
$\M(\gamma_{\slam})$.
The set 
 \begin{align*}
\{\why(\gamma_{\slam});
\slam+\srho\in \sW\backslash \sh\}
 \end{align*}is the complete set of 
the isomorphism classes of simple objects of  $\Wg\adMod$.
\end{Th}
\begin{proof}
It is clear that
 $\WN(\gamma_{\slam})\not\ni |\gamma_{\slam}\ket$.
Thus it follows that
$\WN(\gamma_{\slam})$ is the unique maximal graded submodule 
of $\M(\gamma_{\slam})$.
The second assertion follows from
 Zhu's theorem (Theorem \ref{Th:Zhu-08}).
%applied to our case.
%which is easily seen from \eqref{eq:2006-09-28-14-24}.
\end{proof}
Note that we have
\begin{align}
 S(\why(\gamma_{\slam}))=\C |\gamma_{\slam}\ket
(\cong \C_{\gamma_{\slam}}),
\end{align}
because otherwise 
$\why(\gamma_{\slam})$
admits a non-trivial proper  graded submodule.
\subsection{The Simple Quotient of $\Wg$
is $\why(\gamma_{\vac})$}
\label{subsection:SimpleQuotinetVA}
It is clear that
$\Wg_{\tp}=\C |0\ket $ 
when $\Wg$ is considered as  a $\Wg$-modules.
Set
\begin{align}
 \vac:= \overline{t_{-\srho\che}\circ k \Lam_0}=
-(k+ h\che)\srho\che\in \dual{\sh}.
\end{align}
A short calculation
shows that
\begin{align*}
&\wh t_{\srho\che}(\wJ_i(0))|0\ket
=\bra \vac, J_i\ket |0\ket
\quad \text{ for $i\in \sI$,}\\
 & \wh t_{\srho\che}({\wJ_{\alpha}(0)})|0\ket =0,\quad
 \wh t_{\srho\che}({\psi_{\alpha}(0)})|0\ket =0
\quad \text{ for $\alpha\in \snroots$}
\end{align*}
(cf.\  \cite[Proposition 4.2]{A04}).
Hence 
by Proposition \ref{Pro:2006-10-19-02-59}
it follows that
\begin{align}
 \Wg_{\tp}\cong \C_{\gamma_{\vac}}
\end{align}
as $\Zhu(\Wg)$-modules.
Thus
there exist a unique surjection 
$\M(\gamma_{\vac})\twoheadrightarrow \W_{k}(\sg)$
of $\W_{k}(\sg)$-modules
that sends $|\gamma_{\vac}\ket$ to $|0\ket$.
Therefore
$\why(\gamma_{\vac})$
 is a %(unique in the case that $k\ne -h\che$)
simple quotient of $\Wg$.
By \cite[Remark 4.3.2]{FHL},
 $\why(\gamma_{\vac})$ is a simple
vertex algebra.
 \subsection{Duality Functor $D$}
\label{subsection:DualityFunctorD}
Let
$\theta:\Lie{\Wg}\rightarrow \Lie{\Wg}$
be the anti-Lie algebra
involution as in Proposition \ref{Pro:anti-auto}.
We have
\begin{align}
 \theta(\WW_i(n))=(-1)^{d_i+1}\WW_i(-n)
\end{align}
because $\WW_i$ is quasi-primary.

The $\theta$
 is 
induced by the anti-Lie superalgebra
isomorphism $\Lie{\bCg}\rightarrow \Lie{\bCg}$
which is also denoted by $\theta$.
We have:
\begin{align*}
 &\theta(J_{\alpha}(-\bra \alpha,\srho\che\ket ))
=
-(-1)^{\bra \alpha,\srho\che\ket}
J_{\alpha}(-\bra \alpha,\srho\che\ket)\text{ for $\alpha\in \sroots$,}
\\
&\theta(h(0))=-h(0)\text{ for $h\in \sh$},\\
&\theta(\psi_{\alpha}(-\bra \alpha,\srho\che\ket))
=-\sqrt{-1}
(-1)^{\bra \alpha,\srho\che\ket}
\psi_{\alpha}(-\bra \alpha,\srho\che\ket)\text{ for $\alpha\in
 \sroots$,}\\
&\theta(\psi_{-\alpha}(\bra \alpha,\srho\che\ket))
=\sqrt{-1}
(-1)^{\bra \alpha,\srho\che\ket}
\psi_{\alpha}(\bra \alpha,\srho\che\ket)\text{ for $\alpha\in \sroots$.}
\end{align*}
Therefore the following hold:
\begin{Lem}\label{Lem:2006-10-09-23-57}
The $\theta$ induces an anti-algebra 
isomorphism $\bar{\theta}: \Zhu(\bCg)\rightarrow \Zhu(\bCg)$
with $\bar{\theta}(a)^2=(-1)^{p(a)}a$ for $a\in \Zhu(\bCg)$,
which is given by the following formula
under the identification 
$\Zhu(\bCg)=U(\sg)\* \sCl$:
\begin{align*}
&\bar{\theta}(J_{\pm \alpha})=
-(-1)^{\height \alpha}J_{\pm \alpha}
\text{ for $\alpha\in \sproots$,}
\quad
\bar{\theta}(h)=-h~
\text{ for $h\in \sh$},\\
&\bar{\theta}(\psi_{\alpha})=-\sqrt{-1}
(-1)^{\height \alpha}\psi_{\alpha} \text{ for $\alpha\in\sproots$},\\
&\bar{\theta}(\psi_{-\alpha})=\sqrt{-1}
(-1)^{\height{\alpha}}
\psi_{-\alpha}
\text{ for $\alpha\in \sproots$.}
\end{align*}
\end{Lem}
Because
\begin{align}
 \theta(\bd_+)=-\sqrt{-1}\bd_+,
\end{align}
$\bar{\theta}$
induces 
anti-algebra
automorphism 
of $H^0(\Zhu(\bCg))=\Zhu(\Wg)=
\Center(\sg)$,
which 
 is also  denoted by $\bar{\theta}$.

\begin{Lem}\label{Lem:dual-center}
For $\bar{\lam}\in \dual{\sh}$,
we have
$\gamma_{\bar{\lam}}\circ
\bar{\theta}_{|\Center(\sg)}=\gamma_{-w_0(\bar{\lam})}$.
\end{Lem}
\begin{proof}
Let $\lam\in \sP_{++}$ (see \eqref{eq:2006-11-01-12-47}).
Then $\sL(\slam)$ has the central character $\gamma_{\slam}$.
Let $\sL(\slam)^{\bar {\theta}}$ denote
$\sL(\slam)$ viewed as a $\sg$-module
on which $J\in \sg$ acts by $J\cdot v=\bar{\theta}(J)v$.
Then 
the formulas in Lemma \ref{Lem:2006-10-09-23-57}
show that $\sL(\slam)$ has the central character $\gamma_{-w_0(\slam)}$.
This proves the assertion
for all $\slam\in \sP_{++}$.
Hence
it must hold for all $\slam\in \dual{\sh}.$
\end{proof}
Denote by  $\Wcat$  the Serre full subcategory of
$\Wg\adMod$
consisting of graded $\Wg$-modules
$M=\bigoplus_{d\in \C}M_d$
such that $\dim M_d<\infty$ for all $d$.
The Verma module
$\M(\gamma_{\slam})$ belongs to $\Wcat$.
Thus
by Theorem \ref{Th:classification-Simples-W}
any simple object of $\Wg\GrMod$ belongs to $\Wcat$.
 
 For an object $M$  of $\Wcat$,
define
\begin{align}
 D(M)\teigi \bigoplus_{d\in \C}\Hom_{\C}(M_{d},\C).
\end{align}
Then
the formula
\begin{align}
\bra u \cdot f,m\ket=\bra f, \theta(u)m\ket
\quad \text{for }u\in \U(\Wg),\
f\in D(M),
\ m\in M
\end{align}
gives a $\Wg$-module structure
on $D(M)$.
By definition 
$D(M)$ is an object of $\Wcat$.

One has
\begin{align}
D(D(M))=M
\label{eq:2006-11-06-01-53}
\end{align}
for $M\in \Wcat$.
The correspondence
$M\mapsto D(M)$ defines an
exact cofunctor from $\Wcat$
to itself (see Lemma \ref{Lem:elementary-exact}).

The following assertion is clear.
\begin{Lem}\label{Lem:2006-11-08-01-21}
 Let $M\in \Wcat$.
The correspondence
$N\mapsto D(N)$
gives a bijection 
between the set of graded submodules of $M$
and the set graded quotients of $D(M)$
\end{Lem}
Let $|\gamma_{\slam}\ket^*$
denote the homogeneous vector 
of $D(\M(\gamma_{\slam}))$
dual to $|\gamma_{\slam}\ket\in \M(\gamma_{\slam})$.
Then
$D(\M(\gamma_{\slam}))_{\tp}=\C |\gamma_{\slam}\ket^*$.
\begin{Th}\label{Th:Duality}
Let $\slam, \ \smu\in \dual{\sh}$. 

\begin{enumerate}
\item
There is an isomorphism 
$D(\why(\gamma_{\bar{\lam}}))\cong \why(\gamma_{-w_0(\bar{\lam})})$
of $\Wg$-modules.
\item
%Let $\slam\in \dual{\sh}$.
The $\Wg$-submodule of
$D(\M(\gamma_{\slam}))$
generated by 
$D(\M(\gamma_{\slam}))_{\tp}
=\C |\gamma_{\slam}\ket^*$
is isomorphic to $\why(\gamma_{-w_0(\slam)})$.
This is a
unique simple  graded submodule of
 $D(\M(\gamma_{\slam}))$.
 \item 
A (graded) homomorphic image of 
$\M(\gamma_{\slam})$
in  $D(\M(\gamma_{\smu})))$ is zero or 
isomorphic to $\why(\gamma_{\slam})$.
%Let $\slam$, $\smu\in \dual{\sh}$.
%For any nonzero element $\phi$ of
%  $ \Hom_{\Wcat}
%(\M(\gamma_{\slam}),
%D(\M(\gamma_{\smu})))$,
%we have 
%$\phi(\M(\gamma_{\slam}))
%\cong \why(\gamma_{\slam})$.
Hence
\begin{align*}
\dim \Hom_{\Wcat}
(\M(\gamma_{\slam}),
D(\M(\gamma_{\smu})))=
\begin{cases}
 1&\text{if $-w_0(\smu)\in \sW\circ \slam$,}\\
0&\text{otherwise.}
\end{cases}
\end{align*}
\end{enumerate}
\end{Th}
\begin{proof}
(i) %Because $D()$
By Lemma \ref{Lem:2006-11-08-01-21}
it is clear that
$D(\why(\gamma_{\slam}))$
is a simple object of $\Wcat$.
Thus it is sufficient to show that
$D(\why(\gamma_{\slam}))$ is a quotient of $\M(\gamma_{-w_0(\slam)})$.
Because its simple,
$S(D(\why(\gamma_{\slam})))
=D(\why(\gamma_{\slam}))_{\tp}
\cong 
 \C_{\gamma_{-w_0(\slam)}}$
by
Lemma \ref{Lem:dual-center}.
By \eqref{eq:2006-09-28-14-24},
this proves the assertion.
(ii) 
By Lemma \ref{Lem:2006-11-08-01-21}
and the first assertion of
Theorem \ref{Th:classification-Simples-W},
it follows that $D(\M(\gamma_{\slam}))$
has the unique simple graded submodule
$D(\why(\gamma_{\slam}))$,
which is
isomorphic to
$\why(\gamma_{-w_0(\slam)})$, by (i).
Further
we have
$|\gamma_{\slam}\ket^*\in D(\why(\gamma_{\slam}))$,
because it is the vector dual to $|\gamma_{\slam}\ket
\in \M(\gamma_{\slam})$.
Therefore the assertion follows.
(iii) 
Let
$\phi:\M(\gamma_{\slam})\ra
D(\M(\gamma_{\smu}))$ be a  homomorphism 
of
graded
$\Wg$-modules
and suppose that
$\phi(\M(\gamma_{\slam}))\ne \zero $.
Then
$\phi(\M(\gamma_{\slam}))$
must contain 
$|\gamma_{\smu}\ket^*$,
by (ii).
From this it follows that
\begin{align}
\phi(\M(\gamma_{\slam})_{\tp})=D(\M(\gamma_{\smu}))_{\tp}
\label{eq:2006-10-09-16-20}
\end{align}
because $\M(\gamma_{\slam})$ is generated by 
$\M(\gamma_{\slam})_{\tp}$
and $\phi$ is graded.
Therefore,
by \eqref{eq:2006-09-28-14-24},
$D(\M(\gamma_{\smu}))_{\tp}
\cong \C_{\gamma_{\slam}}$
as $\Center(\sg)$-modules.
But by (ii) we also have
$D(\M(\gamma_{\smu}))_{\tp}
=\C_{\gamma_{-w_0(\smu)}}
$,
thus $-w_0(\smu)\in \sW\circ \slam$.
Conversely,
if $-w_0(\smu)\in \sW\circ \slam$
then
\eqref{eq:2006-10-09-16-20} extends to the
 graded $\Wg$-module homomorphism 
$\M(\gamma_{\slam})\ra
D(\M(\gamma_{\smu}))$,
by \eqref{eq:2006-09-28-14-24}.
This completes the proof.
\end{proof}
\begin{Rem}\label{REm:2006-11-08-15-50}
For $k\ne -h\che$
 one has
$ \theta(L(0))=L(0)$.
Also,
it is easy to see that
 $\Delta_{\bar{\lam}}=\Delta_{-w_0(\bar{\lam})}$.
\end{Rem}
Because $w_0(\srho\che)=-\srho\che$,
we have the following assertion.
\begin{Pro}
 The simple vertex algebra
$\why(\gamma_{\vac})$ is self-dual,
that is,
there is an isomorphism
$D(\why(\gamma_{\vac}))\cong \why(\gamma_{\vac})$.
\end{Pro}
\subsection{Critical Level Cases vs Non-Critical Level Cases}
It is known \cite{FF_W} that
 the structure of $\Wg$ drastically
changes at the critical level $k=-h\che$.
Thus,
though our argument does not depend on the parameter $k$,
it makes sense to discuss  the difference
between the critical level case and non-critical level cases.

First,
let $k=-h\che$.
The following remarkable fact is well-known.
\begin{Th}[Feigin-Frenkel\cite{FF_W}]
The vertex algebra  $\W_{-h\che}(\sg)$ is commutative.
In fact it coincides with the center of $V_{-h\che}(\sg)$
(in the sense of vertex algebras).
\end{Th}
\begin{Rem}
 If $k\ne -h\che$
then
the center of $V_k(\sg)$ is trivial.
This follows from 
the fact that 
$\Vac {V_k(\sg)} =\C |0\ket$
%(because $V_k(\sg)$ admits a PBW basis, 
%see Proposition \ref{Pro:2006-05-17-23-34})
and \cite[Proposition 3.11.2]{LL}.
%(because 
%for $k\ne -h\che$
%$V_k(\sg)$ is a vertex operator algebra).
\end{Rem}
Since
$\W_{-h\che}(\sg)$ is commutative
it follows that 
$\WN(\gamma_{\slam})=\sum_{d<0}\M(\gamma_{\slam})_d$
and $\why(\gamma_{\slam})$ is one-dimensional.
Therefore  the following assertion immediately
follows.
\begin{Pro}\label{Pro:critiela-case-is-one-dim}
Let $k=-h\che$.
Then 
$\why(\gamma_{\slam})$
is the one-dimensional $\W_{-h\che}(\sg)$-module
on which $W_i(n)$,
with $i\in \sI$
and $n\ne 0$,
acts trivially
and $\Zhu(\Wg)=\Center(\sg)$
acts by the central character $\gamma_{\slam}$.
\end{Pro}

Next let $k\ne -h\che$.
In this case $\Wg$ has 
the conformal vector $L(-2)|0\ket$.
By \eqref{eq:cc} and Proposition \ref{Pro:image.of.L(0)},
we have:
\begin{align}
 L(0)|\gamma_{\bar{\lam}}\ket 
=\Delta_{\bar{\lam}}|\gamma_{\bar{\lam}}\ket,
\end{align}
%when $k \ne 0$,
where
\begin{align}
  \Delta_{\bar{\lam}}
=\frac{|\bar{\lam}+\srho|^2}{2(k +h\che)}
-\frac{\rank \sg}{24}+\frac{c(k)}{24}.
 \end{align}
Because
the grading of $\M(\gamma_{\slam})$
is determined by the  action of $-L(0)$
(up to some constant shift),
we have the following assertion.

\begin{Pro}\label{Pro:simple-def-follows}
Let  $k\ne -h\che$.
Then 
each $\why(\gamma_{\slam})$
is the unique irreducible quotient of 
the $\Wg$-module $\M(\gamma_{\slam})$.
\end{Pro}
\begin{Rem}
 In the case that  $k=-h\che$,
$\M(\gamma_{\slam})$ has (many)
non-graded irreducible quotients.
\end{Rem}
%\subsection{Category $\Wcat$
%and Normalized Characters}\label{subsection:categoryW}
%($k\ne -h\che$)
%Here and  until the end of Section \ref{subsection:DualityFunctorD}
% we assume that
%$k\ne -h\che$.
%the level $k$ is non-critical
For a $\Wg$-module $M$,
let
$M_{[h]}$ be the generalized eigenspace of
$L(0)$ of eigenvalue $-d\in \C$:
\begin{align*}
 V_{[d]}\teigi \{ v\in V;  
(L(0)+d)^r v=0\text{ for $r\gg 0$}\}.
\end{align*}
For an object $M$ of  $\Wcat$,
define the normalized formal character 
$\ch M$ of $M$
by
\begin{align}
 \ch M\teigi 
q^{-\frac{c(k)}{24}
}\sum_{d\in \C}q^{-d}\dim M_{[d]}.
\end{align}
One has
 $\ch M=q^{-\frac{c(k)}{24}}\tr_{M}q^{L(0)} $ if $L(0)$
acts on $M$ semisimply.

The following assertion is easily seen from Propositions
\ref{Pro:image.of.L(0)}
and \ref{Pro:Verma-module-gr}.
\begin{Pro}
\label{Pro:character-formula-of-Verma}
It holds that
 \begin{align*}
\ch \M(\gamma_{\bar{\lam}})=\frac{q^{
\frac{|\bar{\lam}+\srho|^2}{2(k+h\che)}
}}{
\eta(q)^{l}}
\end{align*}
for $\slam\in \dual{\sh}$,
where 
$
\eta(q)=q^{\frac{1}{24}}\prod_{i\geq 1}(1-q^i)$.
\end{Pro}

One of the main purpose
of this paper is to determine
each of $\ch \why(\gamma_{\slam})$.

\section{Functors $H_{+}^{\bullet}(?)$ and $H_-^{\bullet}(?)$}
\label{section:FunctorH}
%In this section 
%we recall the definition of
% the category $\BGG_k$
%and two reduction functors $H^{\bullet}_{+}(?)$
% and $H^{\bullet}_{-}(?)$.
%We also recall some  of the structure of
%the category $\BGG_k$%
%for the later purpose.
\subsection{Category $\BGG_k$}
\label{subsection:Category}
%of $\g$ at level $k\in \C$}%(See \cite{MPBook} for details.)
Let
$\BGG_{k}$ be the
Serre full subcategory of the category of left $\g$-modules
consisting of objects $M$
such that the following hold:
\begin{itemize}
 \item  $K$ acts as $k\id$ on $M$;
\item  $M=\bigoplus\limits_{\lam\in
\dual{\h}_{k}}M^{\lam}$ and
$\dim M^{\lam}<\infty$ for all $\lam$;
\item there exists
a finite subset
$\{\mu_1,\dots,\mu_n\}\subset
\dual{\h}_{k}$
such that
$P(M)\subset  \bigcup\limits_i \mu_i-Q_+$.
\end{itemize}
Let $M(\lam)\in \Obj \BGG_k$
%with
%$\lam\in \dual{\bh}_k$
be the Verma module with highest weight $\lam\in \dual{\bh}_k$,
$v_{\lam}\in M(\lam)$  its highest weight vector.
Let $L(\lam)$ be the unique simple quotient of $M(\lam)$.
Then every irreducible object of 
$\BGG_k$ is isomorphic to exactly one of 
the $L(\lam)$ with $\lam\in \dual{\bh}_k$.

The vertex algebra
$\Vkg$ is an object of $\BGG_k$
when considered as a $\bg$-module;
indeed
 $\Vkg$ is a quotient of $M(k\Lam_0)$.
%One knows the following assertion
% (\cite{FZ}, cf.\ Theorem \ref{Th:U-and-Zhu-of-bCg}):
%\begin{Lem}\label{Lem:BGG-cat-and-V}
% Every object of $\BGG_k$
%can be naturally regarded as a $V_k(\sg)$-module.
%\end{Lem}

The correspondence
$M\rightsquigarrow M^*$,
where $M^*$ is defined by \eqref{eq:graded-dual},
defines
the duality functor
in $\BGG_{k}$.
Here, $\g$ acts on $M^*$
by
\begin{align}
(Xf)(v)=f(X^tv)
\label{eq:2006-10-16-00-24}
\end{align}
where
$X\mapsto X^t$ is  the 
anti-automorphism of $\g$ define
by 
$K^t=K$, $\Dg^t=\Dg$ and $J(n)^t=J^t(-n)$ for $J\in \sg$, $n\in \Z$
(notation \eqref{eq:2006-10-10-10-14}).
\subsection{Category $\BGG_k^{\tri}$}
Let
$\BGG_k^{\tri}$
be the full subcategory of $\BGG_k$
consisting of objects $M$ that admit a Verma flag,
that is,
a finite filtration 
$M=M_0\supset M_1\supset \dots \supset M_r=\zero $
such that each successive subquotient $M_i/M_{i+1}$
is isomorphic to
some Verma module $M(\lam_i)$ with $\lam_i\in \dual{\bh}_k$.
The category $\BGG_k^{\tri}$
is stable under taking direct summands.
Dually,
let $\BGG_k^{\bigtriangledown}$  be 
the full subcategory of $\BGG_k$
consisting of objects $M$ such that
$\dual{M}\in \Obj\BGG^{\tri}_k$.
\subsection{Truncated Category $\BGG_{k}^{\leq \blam}$}
(see \cite[Section 2.10]{MPBook} for details.)
For $\lam\in \dual{\bh}_k$,
let
$\BGG_k^{\leq \lam}$
be the Serre full subcategory of $\BGG_k$
consisting of objects $M$ 
such that
$M=\bigoplus\limits_{\mu\leq \lam}M^{\mu}$.
Then
$\BGG_k^{\leq \lam}$
is  stable under taking (graded) duals.
It is known that every simple object $L(\mu)\in \Obj\BGG_k^{\leq \lam}$,
$\mu\leq \lam$,
admits an indecomposable projective
cover $P_{\leq \lam}(\mu)$
in $\BGG_{k}^{\leq \lam}$,
and hence,
every finitely generated object in $\BGG_k^{\leq \lam}$
is an image of a projective object
of the form
$\bigoplus_{i=1}^{r}
P_{\leq \lam}(\mu_i)$.
% for some $\mu_i\in \dual{\bh}$.
%such that $\mu_i\leq \lam$.
%Indeed,
%as in the Lie algebra case (see e.g.\ \cite{MPBOOK}),
%$P_{\leq \lam}(\mu)$
%can be defined as
% an indecomposable direct summand
%of
%\begin{align*}
% U(\bg)\*_{U(\bh\+ \bg_+)}\tau_{\leq \lam}\left(
%U(\bh\+ \bn_+)\*_{U(\bh)}\C_{\mu}\right)
%\end{align*}
%which has $L(\mu)$ as a quotient.
%Here,
%$\tau_{\leq \lam}(M)=M/\bigoplus\limits_\ud{\nu\in \dual{\bh}
%}{\nu\not\leq \lam}M^{\nu}$,
%and $\C_{\mu}$ is a one-dimensional $\bh$-module
%on which $h\in \bh$ acts as $\mu(h)\id$.
It is also known that
$P_{\leq \lam}(\mu)\in \Obj \BGG_k^{\tri}$.
%Moreover,A
%the BGG(Bernstein-Gelfand-Gelfand) reciprocity
%holds:
%\begin{align*}
% [P_{\leq \lam}(\mu):M(\mu')]
%=[M(\mu')
%:L(\mu)]
%\quad (\mu,\mu'\leq \lam).
%\end{align*}%
%Here,
%$[P_{\leq \lam}(\mu):M(\mu')]$
%is the multiplicity of $M(\mu')$
%in the Merma flag of $P_{\leq \lam}(\mu)$,
%and $[M(\mu'):L(\mu)]$
%is the multiplicity of $L(\mu)$
%in the local composition factor (\cite{KacBook}) of $M(\mu')$.
Dually,
$I_{\leq \lam}(\mu)=P_{\leq \lam}(\mu)^*$
is
the injective envelope of $L(\mu)$
in $\BGG_k^{\leq \lam}$.
In particular
 $M\in \Obj\BGG_k^{\leq \lam}$
is a submodule of an injective object of the form
$\bigoplus_{i=1}^r
I_{\leq \lam}(\mu_i)$
%for some $\mu_i\in \dual{\bh}$
if its dual $M^*$ is
finitely generated.

\subsection{Functor $H^{i}_{+}(?)$}
\label{subsection:H_+}
Let $M$ be an object of $\BGG_k$.
Define
\begin{align}
 C(\Ln{+},M)\teigi
M\* \F=\sum_{i\in \Z}
C^i(\Ln{+},M),
 \text{ where $C^i(\Ln{+},M)\teigi M\*
\F^i$}.
\label{eq:def-of-C-+(M)}
\end{align}
Then
$C(\Ln{+},M)$ 
is   naturally regarded as a module over
the vertex superalgebra
$\bCg=\Vkg\* \F(=C(\Ln{+},\Vkg))
$ (recall Section \ref{Section:definition-of-W-algebra}).
In particular,
$\bd_+$ acts 
on $C(\Ln{+},M)$
satisfying
\begin{align*}
\bd_+^2=0,\quad \bd_+\cdot C^i(\Ln{+},M)\subset 
C^{i+1}(\Ln{+},M).
\end{align*}Thus $(C(\Ln{+}, M), \bd_+)$ is a cochain complex.
Set
\begin{align}
H_{+}^{i}(M)
\teigi H^{i}
(C(\Ln{+},M),\bd_{+}) \quad \text{ for $i\in \Z$}.
\end{align}
The $H_{+}^{\bullet}(M)$
 is called the
{\em cohomology %$\Wcoho^{\bullet}(\n,M)$
of the BRST complex of the
quantized Drinfeld-Sokolov 
``$+$'' reduction
for $\Ln{+}$
%(^^ ^^ $+$" and ^^ ^^ $-$" reduction)
with coefficients in   $M$} (\cite{FF_W,FB,FKW}).
\begin{Rem}$ $

\begin{enumerate}
\item We have $H^{\bullet}(\bCg)=H^{\bullet}_+(\Vkg)$
and $\Wg=H^0_+(\Vkg)$.
 \item By definition
\begin{align*}
 H^{\bullet}_{+}(M)=H^{\semiinf+\bullet}(\Ln{+},M\* \C_{\chi_{+}}),
\end{align*}
where the right-hand-side is
%$H^{\bullet+\bullet}(\Ln{\pm},M)$ is 
the semi-infinite cohomology
 \cite{Feigin}
of the Lie algebra $\Ln{+}$ with coefficient in $M\* \C_{\chi_{+}}$.

\end{enumerate}\end{Rem}
The $\bCg$-module structure
of $C(\Ln{+},M)$ induces the $\Wg$-module structure 
on $H^i_+(M)$ with $i\in \Z$,
by Theorem \ref{Th:2006-08-29-09-53}.
Further,
there is a  natural action of the operator $\DW$
on $C(\Ln{+},M)$
commuting with 
$\bd_+$.
Thus $H^{\bullet}_+(M)$
is graded by $\DW$:
\begin{align}
 H^{\bullet}_+(M)=\bigoplus_{d\in \C}H^{\bullet}_+(M)_d,
\quad
H^{\bullet}_+(M)_d=\{v\in H^{\bullet}_+(M);
\DW v= d v\}.
\end{align}
Therefore
we have a functor 
form $\BGG_k$ to
$\Wg\GrMod$
defined
by 
\begin{align*}
M\mapsto H^i_+(M)
\end{align*}
for each $i\in \Z$.

\subsection{Functor $H^{i}_-(?)$}
\label{subsection:H_-}
Consider the automorphism of
$U_k(\bg)\* \bCl$ defined by 
\begin{align}
 \begin{aligned}
  &J(n)\mapsto -J^t(n)\quad \text{for $J\in \sg$, $n\in \Z$},\\
%&K\mapsto K,\\
&\psi_{\alpha}(n)\mapsto -(\psi_{\alpha}^t)(n)
\quad \text{for $\alpha\in \sroots$, $n\in \Z$}
%\\
%&\psi_{-\alpha}(n)\mapsto -\psi_{\alpha}(n)
%\quad \text{for $\alpha\in \sproots$, $n\in \Z$}
 \end{aligned}
\label{eq:2006-10-15-14-24}
\end{align}
(notation \eqref{eq:2006-10-10-10-14}
and  \eqref{eq:2006-10-10-10-15}).
Because it is degree-preserving,
\eqref{eq:2006-10-15-14-24}
induces an algebra involution of 
$\U(\bCg_{\old})=\oldWU$,
which we denote by $\sharp$.
Then 
\begin{align}
\wh t_{-\srho\che}^{\sharp}
:=\sharp\circ \wh t_{-\srho\che}
\end{align}
defines an isomorphism
$\U(\bCg_{\new})\isomap \U(\bCg_{\old})$
of compatible degreewise complete algebras.
We have:
\begin{align}
& \wh t_{-\srho\che}^{\sharp}(\bd_+)=\bd_-,
\quad \bd_-=\bd_-^{\st}+\chi_-,\\
& \wh t_{-\srho\che}^{\sharp}
(\bd^{\st}_+)=\bd_-^{\st},\quad
\wh t_{-\srho\che}^{\sharp}(\chi_+)=\chi_-,
\label{eq:y(Q)}
\end{align}
where 
$\bd_-^{\st}$
and $\chi_-$
are the following elements of $\U(\bCg_{\old})_0$,
respectively:
 \begin{align}
&\bd_{-}^{\st}=
\sum_{\alpha\in
\sroots_-,\ n\in \Z}
J_{\alpha}(-n){\psi}_{-\alpha}(n)
-\frac{1}{2}\sum_{\ud{\alpha,\beta,\gamma\in
\sroots_-}{ k+l+m=0}}
c_{ \alpha,\beta}^{\gamma}
{\psi}_{- \alpha}(k){\psi}_{- \beta}(l)
\psi_{ \gamma}(m),\\
& \chi_-=\sum_{\alpha\in \sroots_-}\schi_-(J_{\alpha})\psi_{-\alpha}(0),
\end{align}
where
$\schi_-$
is defined by \eqref{eq:2006-11-01-14-39}.
Clearly,
\begin{align}
(\bd_-^{\st})^2=\chi_-^2=
[\bd_-^{\st},\chi_-]=0,
\quad  \bd_-^2=0.
\end{align}

For $M\in \BGG_{k}$,
set
\begin{align}
C(\Ln{-},M)\teigi
M\* \F=\sum_{i\in \Z}
C^i(\Ln{-},M),
\text{ where $C^i(\Ln{-},M)\teigi M\*
\F^{-i}$}
\label{eq:def-of-C--(V)}
\end{align}
(compare \eqref{eq:def-of-C-+(M)}).
Then 
\begin{align*}
 \bd_- \cdot C^i(\Ln{-},M)\subset C^{i+1}(\Ln{-},M).
\end{align*}
Thus
$(C(\Ln{-},M), \bd_-)$ is a cochain complex.
The corresponding
cohomology is denoted by
$H^{\bullet}_-(M)$:
\begin{align}
 H^i_-(M)\teigi H^i(C(\Ln{-},M),\bd_-)
\quad \text{for $i\in \Z$}.
\end{align}
\begin{Rem}\label{Rem:-cohomology-as-semi-inf}
By definition
we have 
\begin{align*}
  H^{\bullet}_{-}(M)=H^{\semiinf+\bullet}(\Ln{-},M\* \C_{\chi_{-}}),
\end{align*}
where $\chi_-$ is considered as the character of $\Ln{-}$
defined by $\chi_-(J(n))=\schi_-(J)\delta_{n,0}$
with $J\in \sn_-$ and $n\in \Z$
(compare \eqref{eq:def-of-sH}).
\end{Rem}
By  \eqref{eq:2006-11-01-10-15}
and \eqref{eq:2006-11-06-00-40},
 each $H^i_-(M)$
has
a $\Wg$-module structure:
\begin{align*}
 u\cdot [v]=[\wh t^{\sharp}_{-\srho\che}(u)v]\quad \text{for $u\in
 \U(\Wg)$,
$v\in H^{i}_-(M)$.}
\end{align*}
Also, there is a
 natural action of $\Dg$ 
on $C(\Ln{-},M)$ commuting  with $\bd_-$.
Thus
$H^{\bullet}_-(M)$ is graded by $\Dg$:
\begin{align}
 H^{\bullet}_-(M)=\bigoplus_{d\in \C}H^{\bullet}(M)_d,
\quad H^{\bullet}(M)_d=\{ v\in H^{\bullet}(M);
\Dg v= d v\}.
\end{align}
%Therefore
The correspondence
\begin{align*}
M\mapsto H_-^i(M),
\end{align*}
with
$i\in \Z$,
defines %we have obtained 
a functor
from $\BGG_k$ to $\Wg\GrMod$,
which is refereed as 
the 
 ``$-$'' reduction functor.
The functor  $H^{i}_-(?)$ is essentially the same
functor studied by 
E. Frenkel, V. Kac and M. Wakimoto \cite{FKW}.

\section{Representation Theory of $\W$-Algebras Through the  Functor
$H_-^0(?)$}
\label{section:RepresentationTheoryofW-algebrasThRSRroughtheFunctor}
In this section we study
the representation theory of $\W$-algebras through the  
``$-$'' reduction functor
$H_-^0(?)$.
As in the previous sections,
the complex number $k$ is arbitrary
unless otherwise stated.
\subsection{Category $\dBGG_k$}
Let $M$ be an object of $\BGG_k$.
Regarded as a module over
$\Vkg$ with the Hamiltonian $\Dg$,
its grading is given by 
the operator $\Dg$:
\begin{align*}
 M=\bigoplus_{d\in \Z}M_d,
\quad  M_d\teigi \{v\in M;  \Dg v=d v\}
=\sum_{\lam\atop
\bra \lam, \Dg\ket=d}M^{\lam}.
\end{align*}
Each $M_d$ %with $d\in \C$
is naturally regarded as a module over $\sg$.
Let $\dBGG_k$
be the full subcategory of $\BGG_k$
consisting of objects $M$
such that
 each $M_{d}$,
with $d\in \C$,
belongs to the category $\BGG(\sg)$
(notation Section \ref{subsection:WhittakerFunctor}).
It is easily  seen that $\dBGG_k$ is a Serre 
full subcategory of $\BGG_k$.

\begin{Lem}\label{Lem:dot-BGG}
$ $

 \begin{enumerate}
\item Any Verma module $M(\lam)$,
with $\lam\in \dual{\bh}_k$,
belongs to $\dBGG_k$.
\item Any simple module $L(\lam)$,
with $\lam\in \dual{\bh}_k$,
belongs to $\dBGG_k$
  \item Any object of $\BGG_k^{\tri}$ belongs to $\dBGG_k$.
%In particular,
%$P_{\leq \lam}(\mu)$,
%$\mu\leq \lam$,
%belongs to $\dBGG_k$.
  \item Any object of $\BGG_k^{\bigtriangledown}$ belongs to $\dBGG_k$.
%In particular,
%$I_{\leq \lam}(\mu)$,
%$\mu\leq \lam$,
%belongs to $\dBGG_k$.
 \end{enumerate}
\end{Lem}
\begin{proof}
(i) 
Certainly
$\sn_+$ acts locally nilpotently and
$\sh$ acts semisimply  on $M(\lam)$.
We have to show that
each $M(\lam)_{d}$,
$d\in \C$,
is finitely generated over $\sg$.
But this follows from
 the PBW theorem.
(ii) and (iii) follow from (i).
(iv)
The category $\BGG(\sg)$
is closed under taking (graded) duals.
Hence $\dBGG_k$ is closed under taking (graded) duals.
Therefore (iv) follows from (ii).
\end{proof}

Let $\dBGG_k^{\leq \lam}$ be 
the Serre full subcategory of $\BGG_k$ consisting of objects 
that belong to both $\dBGG_k$ and $\BGG_k^{\leq \lam}$.
Then by Lemma \ref{Lem:dot-BGG},
$P_{\leq \lam}(\mu)$ and $I_{\leq \lam}(\mu)$,
with $\mu\leq \lam$, belong to $\dBGG_k^{\leq \lam}$.

From the following assertion 
it follows that
every object of $\BGG_k$
can be obtained as an injective limit of objects
of $\dBGG_k$:
\begin{Pro}\label{Pro:injective_limits}
 Let $M$ be an object of $\BGG_k^{\leq \lam}$
with
$\lam\in \dual{\bh}_k$.
Then there exists a sequence
$M_1\subset M_2\subset M_3\dots $
of objects of $\dBGG_k^{\leq \lam}$
such that $M=\bigcup_i M_i$.
\end{Pro}
\begin{proof}
Finitely generated objects of $\BGG_k$ belong to $\dBGG_k$
since each projective module $P_{\leq \lam}(\mu)$ does.
Let $\zero =M_0\subset 
M_1\subset M_2\subset M_3\dots $ be a highest weight filtration of
$M$,
so that
$M=\bigcup_i M_i$ 
and each successive subquotient  $M_i/M_{i-1}$ is a
 highest weight module.
In particular
each $M_i$ is finitely generated,
and hence
 belongs to $\dBGG_k$.
%since it is a quotient of
%a projective module of the form $P_)$
%by Lemma \ref{Lem:dot-BGG}.
%Proposition is proved.
\end{proof}

\begin{Lem}
 \label{Lem:finitely-gen}
Fix 
$d\in \C$.
Then for any 
object $M$ of $\dBGG_k^{\leq \lam}$,
with $\lam\in \dual{\bh}_k$,
there exists a 
finitely generated submodule $M'$ of $M$
such that $(M/M')_{d'}=\zero $ for all $d'\geq d$.
\end{Lem}
\begin{proof}
Let
$\P=\{v_1,v_2,\dots\}$
be a set of generators of $M$
such that 
(1) each 
$v_i$  is a weight vector of weight $\mu_i\in \dual{\bh}$,
%(2) $\mu_i \not\leq  \mu_j$ if $i<j$,
(2) if we set
$M_i=\sum_{r=1}^i U(\bg)v_r$ (and $M_0=\zero $),
then each $M_{i}/M_{i-1}$ is a highest weight module with highest weight
 $\mu_i$.
(so $M_1\subset M_2\subset \dots $ is a highest weight filtration
of $M$).
Then by definition $\sharp\{j\geq 1;  \mu_j=\mu\}\leq 
[M: L(\mu)]
$ for $\mu\in \dual{\bh}$,
where $[M: L(\mu)]$ is the multiplicity of $L(\mu)$
in the local 
composition factor of $M$.
Let
\begin{align}
 \P_{\geq d}
=\{ v_j\in \P;  \bra \mu_j,\Dg\ket\geq d\}\subset \P.
\end{align}
Then we see from the definition of $\dBGG_k$
that
$ \P_{\geq d}$
is a finite set.
%Now Lemma follows since
Thus $M'=\sum_{v\in \P_{\geq d}}U(\bg)v\subset M$
satisfies the desired properties.
\end{proof}

\subsection{Cohomology of Top Parts}
\label{Subsection:Cohomology of Top Degree}
Let $M$ be an object of $\BGG_k$
such that
$M=\bigoplus_{d\leq d_0}M_{d}$ for some $d_0\in \C$.
Then 
by \eqref{eq:2006-11-06-02-41}
we have $C(\Ln{-},M)=\bigoplus\limits_{d\leq d_0}C(\Ln{-},M)_{d}$,
where 
$C(\Ln{-},M)_{d}=\{c \in C(\Ln{-},M); \Dg c=d\}$.
Thus one has the following assertion.
\begin{Lem}
\label{Lem:09-wt-appearing}
 Let
$M$ be an object of $\BGG_k$.
Suppose that $M=\bigoplus\limits_{d\leq d_0}M_{d}$ for some $d_0\in \C$.
Then
$H^{\bullet}_-(M)=\bigoplus\limits_{d\leq d_0}H^{\bullet}_-(M)_{d}$.
In particular,
 $H^{\bullet}_-(M)=\bigoplus\limits_{d\leq \bra \lam,\Dg\ket}
H^{\bullet}_-(M)_{d}$
for any object $M$ of $ \BGG^{\leq \lam}_k$.
\end{Lem}

Recall that $\sH_i(M)=H_i(\sn_-,M\* \C_{\schi_-})=
H_i(M\* \Lam (\sn_-),\sd_-)$
(notation  Section \ref{subsection:WhittakerFunctor}).
\begin{Lem}\label{Lem:cohomology-of-zero-mode}
Let $M $ be an object of $\BGG^{\leq \lam}_k$
with
$\lam\in \dual{\bh}_k$.
Then,  we have
 \begin{align*}
  \text{$H^{i}_-(M)_{\bra \lam,\Dg\ket}=
\begin{cases}
 \sH_{-i}(M_{\bra \lam,\Dg\ket})&\text{for $i\leq 0$},\\
\zero &\text{for $i>0$}.
\end{cases}
$
}
 \end{align*}
\end{Lem}
\begin{proof}
Because $M$ is an object of 
$\BGG_k^{\leq \lam}$,
we have
\begin{align*}
 C(\Ln{-},M)_{\bra \lam,\Dg\ket}
=M_{\bra \lam,\Dg\ket}\* \F_{0}.
\end{align*}
Since
$\F_0=\haru\{ \psi_{-\alpha_1}(0)\dots \psi_{-\alpha_r}(0)\1
;  \alpha_i\in \sproots\}\subset \F$
and $\F_0$ is naturally identified with $\Lam(\sn_-)$.
Thus 
we have
\begin{align}
C(\Ln{-},M)_{\bra \lam,\Dg\ket}
=M_{\bra \lam, \Dg\ket}\* \Lam(\sn_-)
=\sC(M_{\bra \lam,\Dg\ket} )
\end{align}
(see \eqref{eq;def;sC}).
It is easily seen  that
the action of $\bd_-$
on $C(\Ln{-},M)_{\bra \lam,\Dg\ket}
=\sC(M_{\bra \lam,\Dg\ket} )
$ coincides with that of $\sd_-$.
This completes the proof.
\end{proof}
We identify
$\Zhu(\Wg)
$ with $\Center(\sg)$ through
Theorem \rm{\ref{Th:2006-10-19-30-31}}
\ban{ii}.
\begin{Pro}\label{Pro:coho-of-top-degree}
 For each $\lam\in \dual{\bh}_k$
we have the following
isomorphisms of $\Zhu(\Wg)$-modules
\begin{enumerate}
 \item $H^i_-(M(\lam))_{\tp}
=H^i_-(M(\lam))_{\bra \lam,\Dg \ket }
%=\C |\lam\ket
\cong \begin{cases}
					    \C_{\gamma_{-w_0(\bar{\lam})}}
&\text{if $i=0$},\\
\zero &\text{if $i\ne 0$;}
					   \end{cases}$
 \item $H^i_-(L(\lam))_{\bra \lam,\Dg\ket}
\cong \begin{cases}
					    \C_{\gamma_{-w_0(\bar{\lam})}}
&\text{if $i=0$ and $\slam\in \dual{\sh}$ is anti-dominant},\\
\zero &\text{otherwise;}
					   \end{cases}$
 \item $H^i_-(M(\lam)^*)_{\tp}
=H^i_-(M(\lam)^*)_{\bra \lam,\Dg\ket}
\cong \begin{cases}
					    \C_{\gamma_{-w_0(\bar{\lam})}}
&\text{if $i=0$},\\
\zero &\text{if $i\ne 0$.}
					   \end{cases}$
\end{enumerate}
\end{Pro}
\begin{proof}
First,
we have
\begin{align}
 \gamma_{\slam}\circ \sharp_{|\Center(\sg)}=\gamma_{-w_0(\slam)},
\end{align}
which can be shown in the same manner as
Lemma \ref{Lem:dual-center}.

Next
we have
$M(\lam)_{\tp}= M(\lam)_{\bra \lam,\Dg\ket}=\sM(\slam)$,
$L(\lam)_{\tp}= L(\lam)_{\bra \lam,\Dg\ket}=\sL(\slam)$
and $M(\lam)^*_{\tp}= M(\lam)_{\bra \lam,\Dg\ket}^*=\sM(\slam)^*$
as $\sg$-modules.
Therefore the assertion follows from
Theorems \ref{Th:vanishing-finite-dim}, \ref{Th:image-of-simple-f-d}
and Lemma \ref{Lem:cohomology-of-zero-mode}.
\end{proof}
\subsection{A Technical 
Lemma}
We need the following 
assertion for the later argument (because the category
$\BGG_k^{\leq \lam}$ is not Artinian in general).
\begin{Lem}
 \label{Lem:replace-by-finitely-gen}
Fix $\lam\in \dual{\bh}_k$
and $d\in \C$.
For any 
object $M$ of $\dBGG_k^{\leq \lam}$
we have the following:
\begin{enumerate}
 \item 
There exists a 
finitely generated submodule $M'$ of $M$
such that $H^{\bullet}_-(M)_{d}=H^{\bullet}_-(M')_{d}$.
\item
There exists a 
finitely generated submodule $N$ of $M^*$
such that $H^{\bullet}_-(M)_{d}=H^{\bullet}_-(N^*)_{d}$.
\end{enumerate}
\end{Lem}
\begin{proof}
(i)
Let $M'$ be a finitely generated 
submodule of $M$
as in Lemma \ref{Lem:finitely-gen},
so that 
\begin{align*}
 (M/M')_{d'}=\zero 
\text{ for all $d'\geq d$.}
\end{align*}
From  the exact sequence $0\rightarrow M'\rightarrow M\rightarrow
 M/M'
\rightarrow 0$,
we get the following long exact sequence:
\begin{align}\label{eq:10-28-exact-mornibg}
 \dots \rightarrow H^{i-1}_-(M/M')\rightarrow H^i_-(M')
\rightarrow H^i_-(M)\rightarrow H^i_-(M/M')\rightarrow \cdots.
\end{align}
Clearly,
 the restriction of \eqref{eq:10-28-exact-mornibg}
 to a $\Dg$-eigenspace
 remains exact.
We are done by 
%Lemma \ref{Lem:finitely-gen}
%and 
Lemma \ref{Lem:09-wt-appearing}.
(ii)
Similarly as above,
let $N$ be a finitely generated submodule of $M^*$
such that 
 $(M^*/N)_{d'}=\zero $ for all $d'\geq d$.
Then 
$(M^*/N)^*_{d'}=\zero $ for all $d'\geq d$.
Therefore 
$H_-^{\bullet}((M^*/N)^*)_{d'}=\zero $
for all $d'\geq d$
by Lemma \ref{Lem:09-wt-appearing}.
The assertion follows
by considering the
long exact sequence
corresponding to the
 exact sequence
$0\rightarrow (M^*/N)^*\rightarrow
M\rightarrow N^*\rightarrow 0$.
\end{proof}

\subsection{Estimate of  $\Dg$-Eigenvalues}
For
 $M\in \BGG_{k}$,
let
$\sH_{\bullet}(M)=\bigoplus_{d\in \C}\sH_{\bullet}(M)_d$.
where $\sH_{\bullet}(M)_d=\sH_{\bullet}(M_d)$.
%and $\sH(M)$ be the Lie algebra 
%homology of $\sn_-$
%is defined by \eqref{eq:def-of-sH}.
\begin{Pro}\label{Pro:reduction-to-Kostant-1}
$ $

\begin{enumerate}
 \item  Let $M$ be an object of $\dBGG_{k}$.
Then  $\sH_i(M)=0$ for $i\ne 0$
and $\sH_0(M)_d$ is 
finite-dimensional for all $d\in \C$.
\item 
For each object $M$ of $\BGG_k$
we have $\sH_i(M)=0$ with $i\ne 0$.
\end{enumerate}\end{Pro}
\begin{proof}
 (i)
follows from  Lemma \ref{Lem:f-d-H-0}
and Theorem \ref{Th:vanishing-finite-dim}.
(ii) follows from (i) and 
Proposition \ref{Pro:injective_limits}
since a homology functor commutes with injective limits.
\end{proof}
\begin{Th}[{cf.\ \cite[Proposition 7.6]{A04}}]\label{Th:estimate}
Let $\blam\in \dual{\bh}_{k}$.
\begin{enumerate}
 \item $($\cite{FKW}$)$  Let $M$ be an object or $\dBGG^{\leq \blam}_{k}$.
Then
each eigenspace $ H^i_-(M)_d$,
with $d\in \C$, is finite-dimensional.
\item  For each object $M$ of $\BGG_k^{\leq \lam}$
we have $ H^i_-(M)_d
=\zero $
unless
$
|i|\leq \bra \lam,\Dg\ket-d
$.\end{enumerate}
\end{Th}
\begin{proof}
Because
$
 \Ln{-}\cap \g_-=
(\sn_-\* \C[t\inv]t\inv)\+\sn_- 
$,
the space
$ \Lam^n (\Ln{-}\cap \g_-)$
decomposes as
\begin{align*}
 \Lam^n (\Ln{-}\cap \g_-)=
\bigoplus_{i+j=n}\Lam^i (\sn_-\* \C[t\inv]t\inv)\*
\Lam^j(\sn_-)\quad \forall n.
\end{align*}
Thus
by \eqref{eq:identi-FLn}
the space $\F=\bigoplus_{n\in \Z}\F^n$ decomposes as
\begin{align*}
& \F=\Lam^{\semiinf}(\Ln{-}/\sn_-)
\* \Lam(\sn_-),
\\&%\quad
 \F^n=
\bigoplus_{p+q=n
\atop p\in \Z,
q\in \Z_{\leq 0}}\Lam^{\semiinf+p}(\Ln{-}/\sn_-)\* \Lam^{-q}(\sn_-),
\end{align*}
where
\begin{align}
&\Lam^{\semiinf}(\Ln{-}/\sn_-)=\bigoplus_{p\in \Z}
\Lam^{\semiinf+p}(\Ln{-}/\sn_-),\\
& \Lam^{\semiinf+p}(\Ln{-}/\sn_-)\teigi 
\bigoplus_{i-j=p
\atop i,
j\in \Z_{\geq 0}}\Lam^i (\sn_+\* \C[t\inv]t\inv)
\* \Lam^j (\sn_-\* \C[t\inv]t\inv).
\label{eq:2006-10-04-12-35}
\end{align}
%see \eqref{eq:identification-of-F}.
Let
$M\in \BGG_k^{\leq \lam}$.
Put
$C^n=C^n(\Ln{-},M)$
and 
set
\begin{align*}
 F^p C^n\teigi \bigoplus_{i\geq p
\atop i-j=n}M\* \Lam^{\semiinf+i}(\Ln{-}/\sn_-)
\* \Lam^{j}(\sn_-).
\end{align*}
Then
\begin{align}
&\dots  \supset F^p C^n\supset F^{p+1} C^n\supset \dots
\supset F^{n+\dim \sn_- +1}C^n=0,\\
& C^n=\bigcup_{p}F^p C^n,\\
& \bd_{-} F^p C^n\subset F^p C^{n+1}.
\end{align}
Thus there is 
a corresponding converging spectral sequence
$E_{r}\Rightarrow H^{\bullet}_-(M)$.
This
is
the semi-infinite analogue of
Hochschild-Serre spectral sequence for the subalgebra
$\sn_-\subset \Ln{-}$.

%Because
%$F^p C^{p+q}/F^{p+1}C^{p+q}=M\* \Lam^{\semiinf +p}\* \Lam^{-q}(\sn_-)$,
By definition
we have
\begin{align}
 E_1^{p,q}=
\begin{cases}
 \sH_{-q}(M)\* \Lam^{\semiinf+ p}
(\Ln{-}/\sn_-)&\text{for $q\leq 0$},\\
\zero &\text{for $q>0$}.
\end{cases}\end{align}
Hence by Proposition \ref{Pro:reduction-to-Kostant-1}
it follows that
\begin{align}
 \begin{aligned}
 E_{1}^{p,q}=\begin{cases}
	      \sH_0(M)\* \Lam^{\semiinf+ p}(\Ln{-}/\sn_-)&
\text{for $q=0$}
,\\
\zero 
&\text{for $q\ne 0$}.
	     \end{cases}
\end{aligned}
\label{eq:05-09-03-4}
\end{align}
Therefore
we have
\begin{align}
 \dim H_-^i(M)_{\bra \lam,\Dg\ket-d}
&\leq \sum_{p+q=i}\dim (E_1^{p,q})_{\bra \lam,\Dg\ket-d}
=\dim (E_1^{i,0})_{\bra \lam,\Dg\ket-d}\nonumber \\
&=
\sum_{d'+d''=d
\atop d',d''\geq 0}\dim       \sH_0(M)_{\bra \lam,\Dg\ket-d'}
\cdot
\dim \Lam^{\semiinf+ i}(\Ln{-}/\sn_-)_{-d''},
\label{eq:05-09-03-2}
\end{align}
where
$\Lam^{\semiinf+\bullet}(\Ln{-}/\sn_-)_d
=
\{\omega\in
\Lam^{\semiinf+\bullet}(\Ln{-}/\sn_-);
\Dg \omega=d\omega \}$.
But
we have \begin{align}
&\dim  \Lam^{\semiinf+\bullet}(\Ln{-}/\sn_-)_{-d}<\infty \text{ for all $d$,}
\label{eq:05-09-03-3}
\\
&
 \Lam^{\semiinf+i}(\Ln{-}/\sn_-)_{-d}
=\zero 
\text{ if }d<|i|
\label{eq:05-09-03-5}
\end{align}
which  follows immediately
from the definition
\eqref{eq:2006-10-04-12-35}.

By \eqref{eq:05-09-03-2}
and \eqref{eq:05-09-03-5}
we have
$H^i(M)_{\bra \lam,\Dg\ket-d}=\zero $
for $d<|i|$.
This shows (ii).
Further,
if $M\in \dot{\BGG}_k^{\leq \lam}$
then 
by Proposition \ref{Pro:reduction-to-Kostant-1}
$\sH_0(M)_d$ is finite-dimensional for all $d$.
Therefore
 \eqref{eq:05-09-03-2} and 
\eqref{eq:05-09-03-3}
give the assertion (i).
\end{proof}
%\begin{Rem}
% A weaker version of  Theorem \ref{Th:estimate} (ii)
%was y proved in \cite[Proposition 7.6]{A04}.
%(It also has a similar (weaker)
%statement for the ``$+$'' reduction.
%)
%\end{Rem}
Recall the category $\Wcat$ defined 
in Section \ref{subsection:DualityFunctorD}.
\begin{Co}
The $\Wg$-module $H^0_-(M)$,
with $M\in \dBGG_k$,
belongs to $\Wcat$.
\end{Co}
\subsection{Images of Verma Modules
and their Duals
by the ``$-$'' Reduction}
The following assertion 
is essentially due to \cite{FB},
because its proof is the same as  
Theorem \ref{Th:vanishing-vertex-1}.
\begin{Th}[{\cite[Theorem 5.7]{A04}}]\label{Th:image-of-Verma--}
Let
   $\lam\in \dual{\bh}_k$.
The cohomology
$H^i_-(M(\lam))$ is zero for all $i\ne 0$
and 
there is an isomorphism 
$
H_-^0(M(\lam))\cong \M(\gamma_{-w_0(\bar{\lam})})
$.
\end{Th}
\begin{Rem}\label{Rem:2006-11-06-22-37}
$ $

\begin{enumerate}
 \item The fact that 
$H^0_-(M(\lam))=\M(\gamma_{-w_0(\slam)})$
is not explicitly
stated in \cite{A04},
but it easily follows from the proof of
{\cite[Theorem 5.7]{A04}}.
In fact 
by using exactly the same argument as
Section \ref{Subsection:Tensor Product Decomposition},
it
can be proved in the same manned as Proposition
       \ref{Pro:2006-10-19-03-43}.
\item 
We have
$H^0_-(M(w\circ \lam))\cong \M(\gamma_{-w_0(\slam)})$
for all $w\in \sW$.
\end{enumerate}
\end{Rem}
\begin{Th}\label{Th:vanishing-dual-Verma--}
Let  $\lam\in \dual{\bh}_k$.
The cohomology
 $H^i_-(M(\lam)^*)$ is zero  for all $i\ne 0$
and
there is an isomorphism 
$H^0_-(M(\lam)^*)\cong D(\M(\gamma_{\slam}))$.
\end{Th}
The proof of 
Theorem \ref{Th:vanishing-dual-Verma--}
is given in Section \ref{section:Proof-Vanishing}.

Theorem \ref{Th:image-of-Verma--} in particular implies the following 
assertion.
\begin{Th}[{\cite[Theorem 8.1]{A04}}]\label{Th:Vanishing-Projective}
Let $\lam,\mu\in \dual{\bh}_k$
such that $\mu\leq \lam$.
The cohomology
 $H^i_-(P_{\leq \lam}(\mu))$ is zero for all $i\ne 0$.
\end{Th}
Similarly,
the following assertion follows from
Theorem \ref{Th:vanishing-dual-Verma--}.
\begin{Th}\label{Th:Vanishing-Injective}
Let $\lam,\mu\in \dual{\bh}_k$
such that $\mu\leq \lam$.
The cohomology
 $H^i_-(I_{\leq \lam}(\mu))$ is zero for all
 $i\ne 0$.
\end{Th}

\subsection{Results for the ``$-$'' Reduction%
%Vanishing of Cohomology (for the ``$-$'' Reduction)
}
\label{subsection:REsults-for--}
The following
result 
generalizes \cite[Main theorem (2)]{A04}.
%is a generalization of 
%the result established in \cite{A04}.
\begin{Th}\label{Th:vanighing--}
For any object $M$  of $\BGG_k$,
the cohomology
 $H^i_-(M)$  is zero for all  $i\ne 0$.
\end{Th}
\begin{proof}
 We may assume that
$M\in \Obj \BGG_k^{\leq \lam}$ for some $\lam\in \dual{\bh}_k$.
Also, by Proposition \ref{Pro:injective_limits},
 we may further assume that
$M\in \Obj \dBGG_k^{\leq \lam}$, % for some $\lam\in \dual{\bh}_k$,
since the cohomology functor commutes with
injective limits.
Clearly,
 it is sufficient to show that $H_-^i(M)_{d}=\zero $ for $i\ne 0$
with each $d\in \C$.

First, we show that 
$H_-^i(M)_{d}=\zero $ for all $i> 0$.
By Lemma \ref{Lem:replace-by-finitely-gen} (1),
for a given $d$,
there exists a finitely generated submodule 
$M'$ of $M$ such that
\begin{align}
 H^{i}_-(M)_{d}=H^{i}_-(M')_{d}\quad \text{for all }i\in \Z.
\label{eq;in-the-proof-yale1}
\end{align}
Because $M'$ is finitely generated,
 there exists a 
projective object $P$ % of $\BGG_k^{\leq \lam}$
of the form $\bigoplus_{i=1}^r P_{\leq \lam}(\mu_i)$
and an exact sequence
\begin{align*}
  0\rightarrow N\rightarrow P\rightarrow M'\rightarrow 0
\end{align*}
in $\dBGG^{\leq \lam}_k$.
As the corresponding long exact sequence
 we obtain:
\begin{align}
 \dots \rightarrow H^i_-(P)\rightarrow H^i_-(M')\rightarrow H^{i+1}_-(N)
\rightarrow H^{i+1}_-(P)
\rightarrow \dots.
\end{align}
Hence, it follows that
$H^i_-(M')\cong H^{i+1}_-(N)$ for all $i>0$,
by Theorem  \ref{Th:Vanishing-Projective}.
Therefore,
we find
\begin{align}
 H^i_-(M)_{d}\cong  H^{i+1}_-(N)_{d} \quad \text{for all }i>0,
\end{align}
by \eqref{eq;in-the-proof-yale1}.
But then,
because $N\in \Obj\dBGG_k^{\leq \lam}$,
we can repeat this argument to find,
for each $r>0$,
some object $N_r$ of $\dBGG_k^{\leq \lam}$
such that \begin{align}
 H^i_-(M)_{d}\cong  H^{i+r}_-(N_r)_{d} \quad \text{for all }i>0.
\end{align}
By Theorem  \ref{Th:estimate} (ii),
this implies that
$H^i_-(M)_{d}=\zero $
for $i>0$.

Next, we show that
$H_-^i(M)_d=\zero $ for all $i<0$.
(The argument is parallel   to the above.)
By Lemma 
\ref{Lem:replace-by-finitely-gen} (ii)
there exists a finitely generated submodule $V$ of $M^*$
such that 
$
 H^{i}_-(M)_{d}=H^{i}_-(V^*)_{d}
$ for all $i\in \Z$.
Because $V$ is finitely generated,
 there exists an injective  object $I$ % of $\BGG_k^{\leq \lam}$
of the form $\bigoplus_{i=1}^r I_{\leq \lam}(\mu_i)$
and an exact sequence
\begin{align*}
  0\rightarrow V^*\rightarrow I\rightarrow L\rightarrow 0
\end{align*}
in $\dBGG^{\leq \lam}_k$.
Considering the corresponding long exact sequence,
Theorem  \ref{Th:Vanishing-Injective} gives
$H^i_-(V^*)\cong H^{i-1}_-(L)$ for all $i<0$.
This shows that
$ H^i_-(M)_{d}\cong  H^{i-1}_-(L)_{d}$
for all $i<0$.
By repeating this argument we find that
$H^i_-(M)_d=\zero $
for all $i<0$.
\end{proof}
%The fact that
%$H^{i}_-(L(\lam))=\zero $
%was conjectured in \cite{FKW}
%for an
%admissible weight $\lam$.
\begin{Co}
 \label{Co:exactness}
The correspondence $M\fmap H^0_-(M)$
defines an exact functor from $\BGG_k$
to  $\Wg\adMod$.
Its  restriction 
gives the exact functor
$H^0_-(?):\dBGG_k\rightarrow \Wcat$.
\end{Co}
%\subsection{Irreducibility (for  the ``$-$'' Reduction)}
\begin{Th}\label{Th:irr--}
Let $\lam\in \dual{\bh}_k$.
There is an isomorphism
\begin{align*}
 H^0_-(L(\lam))\cong \begin{cases}
		      \why(\gamma_{-w_0(\bar{\lam})})&\text{if
$\slam\in \dual{\sh}$ is anti-dominant,}\\
\zero&\text{otherwise.}
		     \end{cases}
\end{align*}
\end{Th}
\begin{proof}
Because 
$L(\lam)$ is a quotient of $M(\lam)$,
the exactness of the functor $H^0_-(?)$
(Corollary \ref{Co:exactness})
imply  that
there is an exact sequence of $\Wg$-modules
\begin{align}\label{eq:09-exact-proof-irr}
H^0_-(M(\lam))
\ra
H^0_-(L(\lam))
\ra 0.
\end{align}
By  Theorem \ref{Th:image-of-Verma--}
we have
\begin{align}
H^0_-(M(\lam))\cong  \M(\gamma_{-w_0(\slam)}).
\label{eq:2006-11-06-11-28}
\end{align}
Thus
$H^0_-(L(\lam))$ is generated by 
$H^0_-(M(\lam))_{\tp}
$ over $\Wg$.
By Proposition \ref{Pro:coho-of-top-degree}
(i)
$H^0_-(M(\lam))_{\tp}=H^0_-(M(\lam))_{\bra \lam ,\Dg\ket}$.
Because \eqref{eq:09-exact-proof-irr}
is $\Dg$-equivalent,
$H^0_-(L(\lam))\ne \zero$
if and only if
$H^0_-(L(\lam))_{\bra \lam,\Dg\ket}\ne \zero$.
Hence,
 by
Proposition \ref{Pro:coho-of-top-degree} (ii),
it follows  that
$H^0_-(L(\lam))$ is nonzero
if and only if 
its classical part
$\slam$ is anti-dominant.

Next
 suppose that
$\slam$ is anti-dominant,
so that $H^0_-(L(\lam))\ne \zero$.
Because
$L(\lam)$ is also a submodule of $M(\lam)^*$,
by Theorem \ref{Th:vanishing-dual-Verma--}
and Corollary \ref{Co:exactness},
 there is an exact sequence of $\Wg$-modules
\begin{align}
 0\ra H^0_-(L(\lam))
\ra D(\M(\gamma_{\slam})).
\label{eq:2006-11-06-11-10}
\end{align}
By Theorem \ref{Th:Duality} (iii),
it immediately follows from
\eqref{eq:09-exact-proof-irr},
\eqref{eq:2006-11-06-11-28} and
\eqref{eq:2006-11-06-11-10}
 that
$H^0_-(L(\lam))\cong \why(\gamma_{-w_0(\slam)})$.
 \end{proof}
Recall that
a weight $\lam\in \dual{\bh}$
is called {\em principal admissible} \cite{KW2}
if 
(1) $\lam$ is regular dominant,
that is,
$\bra \lam+\rho,\alpha\che\ket\not\in \Z_{\leq 0}$
for all $\alpha\in \rroots_+$,
and (2) $\Delta(\lam)\cong \Delta^{\re}$ as root systems.

Let $\Pr$ be the set of principal admissible weights 
of level $k$.
A principal admissible weight $\lam$ is called {\em non-degenerate} \cite{KW90}
if $\bra \lam+\rho,\alpha\che\ket \not\in \Z$
for all $\alpha\in \sroots$,
or equivalently,
if $\slam$ is anti-dominant.
Let $\Pr^{\nondeg}$ denote the 
set of non-degenerated principal admissible weighs
of level $k$.

\begin{Co}[{\cite[Conjecture 3.4${}_-$]{FKW}}]\label{Co:FKW-COng--}
 Let $\lam\in \Pr^{\nondeg}$.
\begin{enumerate}
 \item \ban{{\cite{A04}}} The cohomology $H^i_-(L(\lam))$
is zero for all $i\ne 0$.
\item There is an isomorphism
$H^0_-(L(\lam))\cong \why(\gamma_{-w_0(\slam)})$.
\end{enumerate}
\end{Co}
\begin{Rem}\label{Rem:2006-11-06-23-27}
$ $

\begin{enumerate}
 \item  It is clear from Theorems
\ref{Th:vanishing-+} and \ref{Th:Main-results-for-+}  
that $H^{\bullet}(L(\lam))=0$
for $\Pr\backslash \Pr^{\nondeg}$,
as stated in \cite{FKW} (without a complete proof).
\item
It is known \cite{KW90} that 
$\Pr^{\nondeg}\ne \emptyset$
if and only if
\begin{align}
k+h\che=p/q\in \Q,\quad
(p,q)=1,\quad
p\geq h\che,\quad  q\geq h, \quad
(q,r\che)=1,
\label{eq:2006-11-09-11-28}
\end{align}
where $h$ is the Coxeter number of $\sg$
and
$r\che$ is the lacing number of
       $\sg$.
%The set $\Pr^{\nondeg}$
%can be described as follows:
In this case
%Let $k$ be as \eqref{eq:2006-11-09-11-28}.
set
\begin{align*}
I_{p,q}=
(P_{++}^{p-h\che}\times P_{++}^{\vee, q-h})/\eW_+
\end{align*}where 
$w(\lam,\mu)=(w\lam,w\mu)$ for $w\in \eW_+$,
$P_{++}^{p-h\che}$
is the set of integral weights of
level $p-h\che$
and $P_{++}^{\vee, q-h}$ is the integral coweights of level
$q-h$:
\begin{align*}
 &P_{++}^{p-h\che}=
\{\lam\in \dual{\bh}_{p-h\che};
\bra \lam,\alpha\che\ket\in \Z_{\geq 0}
\text{ for all $\alpha\in \rroots_+$}\},\\
&P_{++}^{\vee q-h}=\{
\mu\in \dual{\bh}_{q-h};
(\alpha|\mu)\in \Z_{\geq 0}
\text{ for all $\alpha\in \rroots_+$}\}.
\end{align*}
Then 
there exists a  bijection
\begin{align*}
 \begin{array}{ccc}
\sW\times  I_{p,q}
&\isomap  &\Pr^{\nondeg} \\
(w, (\lam,\mu))&\mapsto & w\circ 
\Lam_{\lam, \mu},
 \end{array}
\end{align*}
where
\begin{align*}
\Lam_{\lam,\mu}= \slam-(k+h\che)(\smu+\srho\che)+(k+h\che)\Lam_0.
\end{align*}
Note that
$\vac=\Lam_{(p-h\che)\Lam_0, (q-h)\Lam_0}$.

By Corollary \ref{Co:FKW-COng--}
we have
$H^0_-(L(w\circ \Lam_{\lam,\mu}))\cong
H^0_-(L( \Lam_{\lam,\mu}))$
for all $w\in \sW$, $(\lam,\mu)\in I_{p,q}$.
Moreover
all the simple modules
\begin{align}
\{ H^0_-(L(\Lam_{\lam,\mu})); (\lam,\mu)\in I_{p,q}\}
\label{eq:2006-11-06-22-57}
\end{align}
are non-isomorphic.
This is all
 explained in \cite{FKW}.

It is natural to expect that
$\why(\gamma_{\vac})$ is rational and
$C_2$-cofinite, 
and the above set \eqref{eq:2006-11-06-22-57}
are exactly the isomorphism classes of
simple
$\why(\gamma_{\vac})$-modules.
 \item It is clear that 
our results also apply  to non-principal 
Kac-Wakimoto admissible representations,
see \cite[Remark 3.4.(c)]{FKW}.
\end{enumerate}
\end{Rem}

\subsection{Characters}
For an object $M$ of $\BGG_{k}$,
define the integer 
$[M: M(\mu)]$
by
\begin{align}
 \ch M=\sum_{\mu \in  \dual{\bh}_k}
[M: M(\mu)]\ch M(\mu),
\end{align}
where
$\ch M$ is the formal character of $M$:
$M=\sum_{\lam\in \dual{\bh}}
e^{\lam}\dim M^{\lam}$.
It is well-known 
\cite{KTI, Cas, KTII, KT1}
that
the coefficient
$[L(\lam): M(\mu)]$ is  expressed in terms of Kazhdan-Lusztig
polynomials provided that  $k\ne -h\che$,
see the formula \cite[Theoreom 1.1]{KT}.

The following assertion  follows immediately
 from 
Proposition \ref{Pro:character-formula-of-Verma},
Theorems \ref{Th:image-of-Verma--}, \ref{Th:irr--},
 Corollary \ref{Co:exactness}
and Remark \ref{REm:2006-11-08-15-50}.
\begin{Th}\label{Th:ch-formula}
 Let $k\ne -h\che$ and $\lam\in \dual{\bh}_k$. 
Suppose that 
the  classical part $\slam$ of $\lam$
is anti-dominant.
Then
we have
\begin{align*}
 \ch \why(\gamma_{-w_0(\slam)})&=\sum_{\mu}
[L(\lam): M(\mu)]
\ch \M(\gamma_{-w_0(\bar{\mu})}),
\end{align*}
or equivalently,
\begin{align*}
 \ch \why(\gamma_{\slam})&=\sum_{\mu}
[L(\lam): M(\mu)]
\ch \M(\gamma_{\bar{\mu}}),\\
&= \sum_{\mu}[L(\lam): M(\mu)]q^{|\smu+\srho|^2/2(k+h\che)}\eta(q)^{-l}.
\end{align*}
\end{Th}
Theorem \ref{Th:ch-formula}
gives
the character formula of all
 $\why(\gamma_{\slam})$.
%In particular
% we have proved
% the character formula 
%of the ``minimal series representations''
%of $\W$-algebras
%conjectured
%by E. Frenkel,
%V. Kac and M. Wakimoto \cite{FKW}.
\begin{Rem}$ $

\begin{enumerate}
 \item 
Using the vanishing of cohomology
one may also apply 
the Euler-Poincare character method \cite{FKW}
to obtain
Theorem \ref{Th:ch-formula}.
\item By setting the formal character
$\wt \ch H^0_-(M):=\sum_{d\in \C}q^d \dim H^0_-(M)_d$
for $M\in \dBGG_k$ 
we obtain  the formula
\begin{align*}
   \wt\ch H^0_-(L(\lam))
=\sum_{\mu}
[L(\lam): M(\mu)]q^{\bra \mu,\Dg\ket}\prod_{i\geq 1}(1-q^{-i})^{-l}
\end{align*}which is valid
even at $k=-h\che$.
If $k=-h\che$
then by Proposition \ref{Pro:critiela-case-is-one-dim}
we have\footnote{For an
 application of the formula \eqref{eq:2006-11-06-11-50},
see [T. Arakawa,
{\em A new proof of the Kac-Kazhdan conjecture},
International Mathematics Research Notices 
Volume 2006 (2006), Article ID 27091, 5 pages],
which appeared 
while the present article was being refereed.
}
\begin{align}
\sum_{\mu}
&[L(\lam): M(\mu)]q^{\bra \mu,\Dg\ket}\prod_{i\geq 1}(1-q^{-i})^{-l}
%\nno %\\&
=q^{\bra \lam,\Dg\ket} 
\label{eq:2006-11-06-11-50}
\end{align}
if if $\slam$ is
  anti-dominant,
and $\sum_{\mu}
[L(\lam): M(\mu)]q^{\bra \mu,\Dg\ket}\prod_{i\geq 1}(1-q^{-i})^{-l}
=0
$ if $\slam$ is not
  anti-dominant.
\end{enumerate}\end{Rem}
\section{Proof of Theorem \ref{Th:vanishing-dual-Verma--}}
\label{section:Proof-Vanishing}
In this section we prove 
Theorem \ref{Th:vanishing-dual-Verma--}.
%The argument is  similar to that of \cite{A05Duke}
%(though  the
% ``$+$'' reduction is treated in \cite{A05Duke}).
Fix $k\in \C$ and $\lam\in \dual{\bh}_k$.
%We keep the notation  throughout Section \ref{section:Proof-Vanishing}.
\subsection{Preliminaries}
The following elementary fact
has been already used in this article
and will be frequently used
in the sequel.
\begin{Lem}\label{Lem:elementary-exact}
The correspondence
$M\mapsto \Hom_{\C}(M,\C)$
defines an exact cofunctor from
the category of finite-dimensional
vector spaces to itself.
\end{Lem}
%By Lemma \ref{Lem:elementary-exact}%
%For instance the correspondence
%$M\mapsto M^*$ (notation \eqref{eq:graded-dual})
%gives 
%an exact cofunctor from the
%full subcategory of the  category
%of $\bh$-modules consisting of objects 
%that are direct sum of finite-dimensional%
%weight spaces to itself.
The following assertion will be also often used.
\begin{Lem}[the universal coefficient theorem]
\label{Lem:elementary-duality-of-cohomology}
 Let $C=\bigoplus_{i\in \Z}C^i$
be a  cochain complex of vector spaces,
$Q$ its differential.
Let $C^{\vee}=\bigoplus_{i\in \Z}(C^{\vee})^i$
with $(C^{\vee})^i=\Hom (C^{-i},\C)$.
Then $C^{\vee}$ can be regarded as a cochain complex with
the differential $Q^{\vee}$,
where $(Q^{\vee} f)(c)=f(Q c)$
with $f\in C^{\vee}$, $c\in C$.
If $C$ is finite-dimensional
then 
\begin{align}
H^i(C^{\vee})=\Hom (H^{-i}(C),\C)\label{eq:2006-10-06-00-37}
\end{align}
for all $i$.
\end{Lem}

For
 a semisimple $\C \Dg$-module $M$
we set 
\begin{align}
M_d=\{m\in M;
\Dg m=dm\}.
\end{align}If $\dim M_d<\infty$ for all $d$
we write
\begin{align}
 \D(M)\teigi \bigoplus_{d\in \C}\Hom_{\C}(M_{d},\C).
\end{align}

Recall the  cochain complex
 $(K(\lam), \bd_+')$ with   $\lam\in \dual{\bh}$
 (see Section \ref{subsection-Old-coho-Varma}).
In view of Lemma \ref{Lem:elementary-duality-of-cohomology},
we have
\begin{align}
 H^i(\D(K(\lam)))=\D(H^{-i}(K(\lam)))
=\begin{cases}
  \D(H^0(K(\lam)))&\text{for $i=0$},\\
0&\text{for $i\ne 0$}
 \end{cases}
\end{align}
by Proposition \ref{Pro:2006-10-19-03-43},
where 
$\D(K(\lam))$ is considered as a 
subcomplex of the dual complex $(K(\lam)\che, (\bd_+')\che)$.

The space 
$\D(K(\lam))$ 
is a $\bCgo$-module
%can be regarded a 
%$\bCgo_{\old}$-module 
through the 
anti-superalgebra
isomorphism
$ 
\U(\bCgo_{\old})\isomap
\U(\bCgo_{\old})$
induced by the anti-Lie superalgebra
isomorphism in Proposition \ref{Pro:anti-auto},
with respect to the 
(old)
M\"obius conformal structure
$\{T^*, -\Dg, T\}$.
We  write  $\theta_{\old}$
for this anti-superalgebra isomorphism.
To avoid confusion 
we write $\theta_{\new}$
for the anti-superalgebra
isomorphism $\U(\bCgo_{\new})\isomap \U(\bCgo_{\new})$
with respect to the 
M\"obius   conformal structure $\{T^*_{\new}, -\DW,
T\}$.

By \eqref{eq:2006-11-01-10-15}
the $\bCgo_{\old}$-modules structure
on $\D(K(\lam))$ gives rise to the
$\Wg$-module structure
on $H^0(\D(K(\lam)))=\D(H^0(K(\lam)))$
through the map \eqref{eq:2006-10-07-18-01}.
Because we have
\begin{align}
 \theta_{\old}\circ \wh t_{-\srho\che}=
\wh t_{-\srho\che}\circ \theta_{\new},
\end{align}
it follows that
\begin{align}
 \D(H^0(K(\lam)))\cong D(\M(\gamma_{\slam}))
\label{eq:2006-11-06-15-14}
\end{align}
as $\Wg$-modules.
\subsection{Identification with Duals}
For any  semisimple $\bh$-module $M$
such that
$M^{\lam}$ is finite-dimensional for all $\lam$,
we write $M^*$ for the linear space
$\bigoplus_{\lam}\Hom_{\C}(M^{\lam},\C)$,
as in 
\eqref{eq:graded-dual}.

Let
\begin{align}
 (\cdot|\cdot)_{\F}:
\F \times \F \rightarrow \C
\end{align}
be the 
symmetric bilinear form defined by
\begin{align*}
&(\1|\1)_{\F}=1,\\
&(\psi_{\alpha}(n)\omega_1|\omega_2)_{\F}
=-\sqrt{-1}( \omega_1|\psi_{-\alpha}(-n)\omega_2)_{\F}
\quad\text{for $\alpha\in \sproots,\
n\in \Z,\ \omega_1\in \F,\omega_2\in \F$,
}\\
&(\psi_{-\alpha}(n)\omega_1|\omega_2)_{\F}
=\sqrt{-1}( \omega_1|\psi_{\alpha}(-n)\omega_2)_{\F}
\quad\text{for $\alpha\in \sproots,\
n\in \Z,\ \omega_1\in \F,\omega_2\in \F$
}
\end{align*}
(recall Section \ref{subsection:Clifford}).
Then the form $(\cdot |\cdot )_{\F}$
is non-degenerate.
Indeed
its restriction to
$(\F^i)^{\mu}\times (\F^i)^{\mu}$
is non-degenerate for each $i\in \Z$ and $\mu\in \dual{\bh}$.

Let $M$ be an object of $\BGG_k$.
One may
identify $C(\Ln{-}, M^*)$ with 
$C(\Ln{+}, M)^*$
via the linear isomorphism
\begin{align}
 \begin{array}{ccc}
C(\Ln{-},M^*)=  M^*\* \F&\isomap &C(\Ln{+}, M)^*
=(M\*\F)^* \\
f\* \omega&\mapsto &(m\* \omega'\mapsto f(m)(\omega|\omega')_{\F}).
 \end{array}
\label{eq:2006-11-01-22-27}
\end{align}
Under the identification \eqref{eq:2006-11-01-22-27},
we have
\begin{align}
 (\sqrt{-1}\bd_- \phi)(c)=\phi(\bd_+' c)
\end{align}
for $\phi\in C(\Ln{+}, M)^*=C(\Ln{-},M^*)$
and $c\in C(\Ln{+},M)$,
where $\bd_+'$ is defined in \eqref{eq:2006-10-17-03-22}.
Thus,
in view of Lemma \ref{Lem:elementary-duality-of-cohomology},
 the cochain complex
$(C(\Ln{-},M^*), \bd_-)$
may be identified the subcomplex
$C(\Ln{+}, M)^*$
of 
$(C(\Ln{+,M})\che,(\bd_+')\che)$.
%We write this complex
%as $(C(\Ln{+}, M)^* ,\bd_-)$.

We would like to apply the 
formula \eqref{eq:2006-10-06-00-37}
%in
%Lemma \ref{Lem:elementary-duality-of-cohomology}
to the case that $M=M(\lam)$,
but unfortunately 
the complex
$(C(\Ln{+}, M(\lam)), \bd_+')$ is not 
a direct sum of finite-dimensional subcomplexes.
\begin{Rem}\label{Rem:acyclic1}
We have
$H^{\bullet}(C(\Ln{+},M),\bd_+')\ne H^{\bullet}_+(M)$ in general.
Indeed,
the cohomology  $H^i(C(\Ln{+},M),\bd_+')$
is zero
for all $i$ and $M\in \BGG_k$,
because $J_{\alpha}(0)$,
with
$\alpha\in \sPi$, acts locally nilpotently on
$M$ (\cite[Theorem 2.3]{FKW}).
\end{Rem}

Under the identification 
$C(\Ln{-}, M^*)=C(\Ln{+}, M)^*$,
the action 
of $\Wg$
of $H^{\bullet}_-(M^*)$ induced by the 
following $\bCgo_{\new}$-module structure of
$C(\Ln{+},M)^*$:
\begin{align}
 (u \cdot f)(c)= f( \theta_{\old}(\wh t_{-\srho\che}(u))c)
=f( \wh t_{-\srho\che}(\theta_{\new}(u))c)
\label{eq:2006-11-06-15-13}
\end{align}
for $u\in \U(\bCgo_{\new})$,
$f\in C(\Ln{+},M)^*$
and 
$c\in C(\Ln{+}, M)$.
\subsection{The First Reduction of the Computation
of $H^{\bullet}_-(M(\lam)^*)$}
\label{subsection:firlst-reduction-of-computation}
Below we  mainly  consider
$\bCgo_{\old}$,
so
if no confusion can arise
 we  write just $\bCgo$
for
$\bCgo_{\old}$.

The cochain complex $C(\Ln{+},M(\lam))= M(\lam)\* \F$
is endowed
with the grading defined by $\Dg$:
$C(\Ln{+}, M(\lam))=
\bigoplus_{d}C(\Ln{+}, M(\lam))_d$.
As a $\bCg$-module
\begin{align}
 C(\Ln{+},M(\lam))_{\tp}=
C(\Ln{+},M(\lam))_{\bra \lam,\Dg\ket}
\cong \sM(\slam)\* \Lam(\sn_-).\label{eq:2006-11-08-12-45}
\end{align}

Let $C_+(\lam)$ be the 
$\bCgo$-submodule of 
$C(\Ln{+},M(\lam))$ generated by the vector
$|\lam\ket$:
\begin{align}
\Clam=\U(\bCgo)\cdot |\lam\ket,
\end{align}
where $|\lam\ket =v_{\lam}\* \1$
and $v_{\lam}$ is the highest weight vector of $M(\lam)$.
By definition
$\Clam$ 
is linearly
spanned by the vectors
\begin{align*}
 \whJ_{a_1}(m_1)\dots \whJ_{a_r}(m_r)
\psi_{\alpha_1}(n_1)\dots \psi_{\alpha_s}(n_s)|\lam\ket
\end{align*}
with $a_i\in \sroots_-\sqcup \sI$,
$\alpha_i\in \sroots_-$,
$m_i,n_i\in \Z$.
Clearly,
\begin{align}
& C_+(\lam)=\bigoplus_{d\leq \bra \lam,\Dg\ket}C_+(\lam)_{d},\\
& C_+(\lam)_{\tp}=C(\Ln{+},M(\lam))_{\tp}=
\sM(\slam)\* \Lam(\sn_-).
\label{eq:2006-11-08-10-13}
\end{align}
\begin{Rem}
By considering the filtration induced by the standard filtration
of $\bCgo$,
one sees that
\begin{align}
 \Clam\cong 
S(\sn_-\* \C[t\inv]\+ 
\sh\* \C[t\inv]t\inv)\* \Lam(\sn_-\*\C[t\inv]).
\label{eq:2006-10-06-23-08}
\end{align}
%In fact,
%because it is generated
%by a Lie superalgebra,
%the standard filtration
%of $\bCGo$ coincides with the one induced
%by the standard filtration of Lie superalgebras.
\end{Rem}
\begin{Pro}
 We have $\bd_+^{\st}\cdot  C_+(\lam)\subset C_+(\lam)$
and $\dchi_+' \cdot C_+(\lam)\subset C_+(\lam)$.
Hence $
       \bd_+' \cdot C_+(\lam)\subset C_+(\lam),
      $
that is,
 $(C_+(\lam),\bd_+')$ is a subcomplex of
$(C(\Ln{+},M(\lam)),\bd'_+)$.
\end{Pro}
\begin{proof}
The assertion that
$\bd_+^{\st} \cdot C_+(\lam)\subset C_+(\lam)$
follows from the third and the fourth 
formulas in Lemma \ref{Lem:rel-2}
and the fact that
$\bd_+^{\st}\cdot |\lam\ket=0$.
The fact that $\dchi_+' \cdot C_+(\lam)\subset C_+(\lam)$
is obvious because this is an inner action of $\bCgo$.
\end{proof}
\begin{Rem}
Note that 
 $\chi_+' \cdot |\lam\ket\ne 0$.
\end{Rem}
%Clearly,
%\begin{align}\label{eq:2006-05-12-14-13}
%C_+(\lam)=\bigoplus_{\mu\in \dual{\bh}}C_+(\lam)^{\mu},
%\quad \dim C_+(\lam)^{\mu}<\infty,\ \forall \mu.
%\end{align}
The space
$C(\Ln{+},M(\lam))$
is a direct sum of finite-dimensional weight spaces
$C(\Ln{+},M(\lam))^{\mu}$ with $\mu\in \dual{\bh}$,
and so is $\Clam$.
Thus
%under the identification 
%\eqref{eq:iden-comp-of-dual},
the inclusion $C_+(\lam)\hookrightarrow C(\Ln{+},M(\lam))$
induces the surjection
\begin{align}
C(L\sn_-,M(\lam)^*)= C(\Ln{+},M(\lam))^*
\twoheadrightarrow C_+(\lam)^*
\label{eq:fininal-surj-comp}
\end{align}
by Lemma \ref{Lem:elementary-exact}.
The map
 \eqref{eq:fininal-surj-comp}
is a cochain map
if we 
consider $C_+(\lam)^*$
as a cochain complex with the differential
$\bd_-$,
where
$(\bd_-\phi)(c)=\phi(\bd_+' c)$
for $\phi \in C_+(\lam)^*$ and $c\in C_+(\lam)$.
The map \eqref{eq:fininal-surj-comp}
is also a $\U(\bCgo_{\new})$-module homomorphism,
whose action is given by \eqref{eq:2006-11-06-15-13}.
Thus
the map \eqref{eq:fininal-surj-comp}
induces a homomorphism of 
$\Wg$-modules
\begin{align}
 H^{\bullet}(M(\lam)^*)\ra H^{\bullet}(\Clam^*).
\label{eq:2006-11-09-12-21}
\end{align}
 \begin{Pro}[{\cite[Proposition 6.3]{A04}}]\label{Pro:2006-05-12-10}
  The map \eqref{eq:2006-11-09-12-21}
gives   the    isomorphism  
$
 H^{\bullet}_-(M(\lam)^*)\isomap  H^{\bullet}(C_+(\lam)^*)
$
of 
$\Wg$-modules.
\end{Pro}
\begin{Rem}
 One can use the argument of 
Section \ref{Subsection:Tensor Product Decomposition}
to prove Proposition \ref{Pro:2006-05-12-10}.
\end{Rem}
\subsection{The Assertion to Be Proved}
\begin{Pro}\label{Pro:2006-11-08-09-34}
  There exists an exhaustive,
decreasing filtration
\begin{align*}
\FF:\quad
\dots \supset \FF^p \Clam^*\supset \FF^{p+1}\Clam^*
\dots \supset \FF^{-1}\Clam^*
\supset \FF^0 \Clam^*=0
\end{align*}
of 
$\Clam^*$
 such that
\begin{align*}
& \U(\bCgo_{\new})\cdot \FF^p \Clam^* \subset \FF^p \Clam^*,\\
&\bd_-\cdot \FF^p \Clam^* \subset \FF^p \Clam,\quad
\Dg\cdot \FF^p \Clam^* \subset \FF^p \Clam^*
\end{align*}
for all $p$,
and the associated spectral sequence
$E_r\Rightarrow H^{\bullet}(\Clam^*)=H^{\bullet}_-(M(\lam)^*)$
satisfies the following:
there exists a $\Wg$-equivalent
isomorphism  of cochain maps
\begin{align}
 E_r^{p,q}\cong D(\M(\gamma_{\slam}))\* (E_r^{p,q})_{\bra \lam,\Dg\ket}
\quad
\text{for all $r\geq 1$.}
\label{eq:2006-11-08-09-49}
\end{align}
Here
$D(\M(\gamma_{\slam}))
\* (E_r^{p,q})_{\bra \lam,\Dg\ket}$
is considered as a cochain complex
with the differential $1\* d_r$
and as a $\Wg$-module
on which 
 $\Wg$ acts on the first factor 
$D(\M(\gamma_{\slam}))$.
\end{Pro}
Let us show that
Proposition \ref{Pro:2006-11-08-09-34}
implies Theorem \ref{Th:vanishing-dual-Verma--}.

Here and after,
for a decreasing filtration $F=\{F^p M\}$
of a superspace
$M$ we write $\gr_F M=\bigoplus_{p\in \Z}\gr_F^p M$,
$\gr_F^p M=F^p M/F^{p+1}M$.

Because  $\FF$ is compatible with the action of $\Dg$,
each $E_r$ is a direct sum of
$\Dg$-eigenspaces:
\begin{align}
 E_r=\bigoplus_{d\in \C}(E_r)_{d}.
\end{align}
Thus $(E_r)_d$
converges to $H^{\bullet}_-(M(\lam)^*)_d$:
\begin{align}\label{eq:convergenceof-sp-seqence-04}
( E_r)_{d}
\Rightarrow H^{\bullet}(C_+(\lam)^*)_{d}=H^{\bullet}_-(M(\lam)^*)_{d}
\quad \text{for each }d\in \C.
\end{align}
Also,
because 
$\FF$ is a filtration of $\bCgo$-modules,
it follows that each $E_r^{p,q}$ is a $\Wg$-module.
In particular
\begin{align}
E_{\infty}\cong \gr_{\FF} H^{\bullet}(\Clam^*)
\quad \text{as $\Wg$-modules.}
\label{eq:2006-11-08-09-45}
\end{align}
Here, $\FF^p H^{\bullet}(\Clam^*)
=\im (H^{\bullet}(\FF^p \Clam^*)\ra
H^{\bullet}(\Clam^*))$.

By \eqref{eq:convergenceof-sp-seqence-04}
and \eqref{eq:2006-11-08-09-45},
we have  the linear isomorphism
$\bigoplus\limits_{p+q=n}(E_{\infty}^{p,q})_d
\cong H^{n}_-(M(\lam)^*)_d$ for all $d$.
In particular by Proposition \ref{Pro:coho-of-top-degree}
(iii)
we have
\begin{align}
\bigoplus_{p} (E_{\infty}^{p,n-p})_{\bra \lam,
\Dg\ket}\cong H^{n}(M(\lam)^*)_{\tp}
=\begin{cases}
\C&\text{if $n=0$},
\\
  \zero&\text{otherwise.}
 \end{cases}
\label{eq:2006-10-08-04-10}
\end{align}
Because $\bigoplus_{p+q=0}(E^{p,n-p}_{\infty})_{\bra \lam,\Dg\ket}$
is one-dimensional,
it follows that
there exists $p_0\in \Z$
such that
\begin{align}
%(E^{p,q}_{\infty})_{\tp}
%=
(E^{p,q}_{\infty})_{\bra \lam,\Dg\ket}
=\begin{cases}
			   \C &\text{if $p=p_0$ and $q=-p_0$}\\
0&\text{otherwise.}
			  \end{cases}
\label{eq:2006-10-08-10-57}
\end{align}
But by \eqref{eq:2006-11-08-09-49}
we have
\begin{align}
 E_{\infty}^{p,q}
\cong D(\M(\gamma_{\slam}))\* (E_{\infty}^{p,q})_{\bra \lam,\Dg\ket},
\end{align}
and therefore 
there is the following isomorphism of $\Wg$-modules:
\begin{align}
 E_{\infty}^{p,q}
\cong 
\begin{cases}
 D(\M(\gamma_{\slam}))
&\text{if $(p,q)=(p_0,-p_0)$},\\
0&\text{otherwise.}
\end{cases}
\end{align}
This
 immediately
gives $H^i_-(M(\lam)^*)=0$ for $i\ne 0$.
Also,
because the term $E_{\infty}$ lies entirely in 
$E_{\infty}^{p_0,-p_0}$,
the corresponding filtration of
$H^0_-(M(\lam)^*)$ is trivial.
Therefore we conclude that 
\begin{align*}
H^{0}_-(M(\lam)^*) \cong
D(\M(\gamma_{\slam}))
\end{align*}
as $\Wg$-modules.

In next sections
we shall construct 
the filtration $\FF$
satisfying  the conditions of Proposition
\ref{Pro:2006-11-08-09-34}.
\subsection{The Filtration $\FF$}
%As a graded $\bCgo$-module,
%\begin{align*}
%C_+(\lam)_{\tp}=C_+(\lam)_{\bra \lam,\Dg \ket}.
%\end{align*}
The space
$C_+(\lam)_{\ttp}=\sM(\slam)\* \Lam(\sn_-)$
 is spanned by the vectors 
 \begin{align*}
 \whJ_{\alpha_1}(0)\dots \whJ_{\alpha_r}(0)
\psi_{\beta_1}(0)\dots \psi_{\beta_s}(0)|\lam\ket
\end{align*}
with $\alpha_i, \beta_i\in \sroots_-$.

Note there is a $\U(\bCgo)$-module isomorphism
\begin{align}
 \Clam\cong \U(\bCgo)\*_{\U(\bCgo)_{\geq 0}}\Clam_{\tp}.
\label{eq:2006-11-09-12-29}
\end{align}

We shall define a filtration $\{\FF_p C_+(\lam)_{\ttp}\}$
of $C_+(\lam)_{\ttp}$,
and then extend it to a filtration
$\{\FF_p C_+(\lam)\}$ of the whole space $C_+(\lam)$.
After that,
 we define the filtration
$\{\FF^p \Clam^*\}
$ 
of $\Clam^*$.

 Set
\begin{align}\label{eq:def-of-Gp-bar}
 \FF_p C_+(\lam)_{\ttp}
\teigi \sum_\ud{\mu \in \dual{\bh}}{ \bra \mu-\lam,\srho\che\ket\leq p}
C_+(\lam)_{\ttp}^{\mu}\subset C_+(\lam)_{\ttp}\quad
\text{for }p\leq 0.
\end{align}
Then,
$ \FF_p C_+(\lam)_{\ttp}$
is
a subspace of $C_+(\lam)_{\ttp}$
linearly spanned by the vectors
 \begin{align*}
 \whJ_{\alpha_1}(0)\dots \whJ_{\alpha_r}(0)
\psi_{\beta_1}(0)\dots \psi_{\beta_s}(0)|\lam\ket
\end{align*}
with $\sum_i \bra \alpha_i,\srho\che\ket+
\sum_j\bra \beta_j,\srho\che\ket\leq p(\leq 0)$.

The following assertion is a direct consequence of the
definition.
\begin{Lem}\label{Lem:2006-10-07-15-12}
Each $\FF_p C_+(\lam)_{\ttp}$ is preserved by the action
of $\widehat J_{a}(n)$,
$\psi_{\alpha}(n)$
with 
$a\in \sroots_-\sqcup \sI$,
$\alpha\in \sroots_-$,
$n\geq 0$.
Thus
$\FF_p\Clam_{\tp}$ is a $\U(\bCgo)_{\geq0}$-submodule of $\Clam_{\tp}$. 
\end{Lem}
\begin{Pro}\label{Pro:2006-05-11-1}
The following hold:
\begin{enumerate}
 \item $\{\FF_p C_+(\lam)_{\ttp}\}
$
is an increasing filtration
of $C_+(\lam)_{\ttp}$:
\begin{align*}
\dots \subset \FF_p C_+(\lam)_{\ttp}
\subset \FF_{p+1} C_+(\lam)_{\ttp}
\subset \dots \subset \FF_0 C_+(\lam)_{\ttp}
=C_+(\lam)_{\ttp}; %, \label{eq:filt-C--(lam)-1}.
\end{align*}
\item $\bigcap_p \FF_p C_+(\lam)_{\ttp}
=\zero$;
\item
For each $p$,
\begin{align*}
 \bd_+^{\st}\cdot  \FF_p C_+(\lam)_{\ttp}
\subset \FF_p C_+(\lam)_{\ttp},
\\\chi_+' \cdot \FF_p C_+(\lam)_{\ttp}\subset \FF_{p-1}
C_+(\lam)_{\ttp}.
\end{align*}
Hence
 $\FF_p C_+(\lam)_{\ttp}$ is a subcomplex of 
$(C_+(\lam)_{\ttp},\bd_-')$.
\item
 The space $C_+(\lam)_{\ttp}/\FF_p C_+(\lam)_{\ttp}$
is finite-dimensional for each $p\leq 0$.
\end{enumerate}
\end{Pro}
\begin{proof}
(i)-(ii) are easily seen.
(iii)
The first inclusion follows from the fact
$\bd_+^{\st}$ preserves weight spaces
with respect to the action
of $\bh$.
The second  assertion has been already seen in Lemma
\ref{Lem:2006-10-07-15-12}.
(iv) is obvious.
\end{proof}
Let
$\FF_p C_+(\lam)$
with
$p\leq 0$
be the $\bCgo$-submodule of
$\Clam$ generated by 
$\FF_p \Clam_{\tp}$:
\begin{align}
 \FF_p \Clam:=\U(\bCgo)\cdot \FF_p \Clam_{\tp}.
\end{align}
By 
Lemma \ref{Lem:2006-10-07-15-12},
 $\FF_p\Clam$
is 
spanned by the vectors
\begin{align*}
 \whJ_{a_1}(-m_1)\dots \whJ_{a_r}(-m_r)
\psi_{\beta_1}(-n_1)\dots \dots \psi_{\beta_s}(-n_s)v,
\end{align*}
with $a\in \sroots_-\sqcup \sI$,
$\beta_i\in \snroots$
$m_i,n_i>0$,
$v\in  \FF_p C_+(\lam)_{\ttp}$.
Further,
by 
\eqref{eq:2006-11-09-12-29},
 there is a $\bCgo$-module isomorphism
\begin{align}
 \FF_p\Clam\cong \U(\bCgo)\*_{\U(\bCgo)_{\geq 0}}
\FF_p \Clam_{\tp}.
\label{eq:2006-10-07-13-32}
\end{align}
The following assertion follows from
Proposition \ref{Pro:2006-05-11-1}
and \eqref{eq:2006-10-07-13-32}.
\begin{Pro}\label{Pro:2006-05-12-2}
 The following hold:
\begin{enumerate}
\item
$\dots \subset \FF_p C_+(\lam)
\subset \FF_{p+1} C_+(\lam)
\subset \dots \subset \FF_0 C_+(\lam)
=C_+(\lam)$;
\item $\bigcap_p \FF_p C_+(\lam)
=\zero$.
\item
Each $C_p C_+(\lam)$ is a subcomplex of 
$(C_+(\lam),\bd_+')$:
\begin{align*}
\bd_+^{\st} \cdot \FF_p C_+(\lam)
\subset \FF_p C_+(\lam),\quad
\chi_+' \cdot \FF_p C_+(\lam)\subset \FF_p C_+(\lam);
\end{align*}
\item The space  $(\Clam/\FF_p\Clam)_d$ 
is finite-dimensional for all $p$ and $d$.
\end{enumerate}
\end{Pro}
Define
\begin{align}
 \FF^p C_+(\lam)^*
\teigi \left(C_+(\lam)/\FF_p C_+(\lam)
\right)^*.
%\subset 
%C_+(\lam)^*\quad \text{for }p\leq 0.
\end{align}
This is a $\bCgo$-submodule of $\Clam^*$.
By Proposition \ref{Pro:2006-05-12-2} (iii),
we have
\begin{align}
 \bd_-\cdot \FF^p \Clam^*\subset \FF^p \Clam^*
\quad \text{for all $p$}.
\label{eq:2006-11-09-12-41}
\end{align}
Also, it is clear that
$\FF^p \Clam^*$
is preserved by the action of $\Dg$.

The following assertion
is a  consequence of
Proposition \ref{Pro:2006-05-12-2} (iv).
\begin{Lem}\label{Lem:2006-10-06-13-22}
The space
$\FF^p \Clam^*_d$ is finite-dimensional for all $p$,
$d$ and
we have $ \FF^p C_+(\lam)^*
=\D(C_+(\lam)/\FF_p C_+(\lam))
$ 
as subspaces of
$\Hom (\Clam, \C)$.
\end{Lem}
By definition
we have exact sequences
\begin{align}
& 0\rightarrow  \FF^p \Clam^*\ra
\Clam^*\ra (\FF_p\Clam)^*\rightarrow 0,\\
&0\ra \FF^p \Clam^*\ra \FF^{p-1}\Clam^*
\ra (\FF_p \Clam/\FF_{p-1}\Clam)^*\ra 0.
\label{eq:2006-10-06-14-22}
\end{align}
These are cochain maps by 
Proposition \ref{Pro:2006-05-12-2} (iii).
\begin{Pro}\label{Pro:2006-05-12-4}
The 
$\FF^p=\{\FF^p C_+(\lam)^*\}$
gives a 
decreasing,
separated,   exhaustive
filtration
of the $\bCgo$-module $C_+(\lam)^*$:
\begin{align*}
 & \dots \supset \FF^p C_+(\lam)^*
\supset \FF^{p+1} C_+(\lam)^*
\supset \dots \supset \FF^0
 C_+(\lam)^*
%\supset \FF^1 C_+(\lam)^*
 = \zero,\\
&C_+(\lam)^*=\bigcup_p \FF^p C_+(\lam)^*.
\end{align*}
\end{Pro}
\begin{proof}
We  need only to show that
 that $\FF^p$ is exhaustive.
For this,
 it is enough to show
that
\begin{align}\label{eq:2006-05-12-14-51}
\left(C_+(\lam)^*\right)^{\mu}
=\bigcup\limits_p (\FF^p
 C_+(\lam)^*)^{\mu}
\end{align}
for each $\mu$.
By Proposition \ref{Pro:2006-05-12-2} (ii),
$\bigcap_p \FF_p C_+(\lam)^{\mu}=\zero$.
Therefore,
there exists $p$
such that $\FF_p C_+(\lam)^{\mu}=\zero$,
because  $C_+(\lam)^{\mu}$ is finite-dimensional.
Thus $(C_+(\lam)^*)^{\mu}=\left(\FF^p C_+(\lam)^*\right)^{\mu}
$,
and this proves  \eqref{eq:2006-05-12-14-51}.
 \end{proof}
By Proposition \ref{Pro:2006-05-12-4}
and \eqref{eq:2006-11-09-12-41}
we obtain
the  converging spectral sequence
corresponding to the filtration $\{\FF^p C_+(\lam)^*\}$:
\begin{align}
 E_r\Rightarrow H^{\bullet}(C_+(\lam)^*)=H^{\bullet}_-(M(\lam)^*)
\label{eq:2006-11--8-10-05}
\end{align}
By definition,
\begin{align}\label{eq:E-1-first-identity}
 E_1^{p,q}=H^{p+q}(\gr^\FF_p C_+(\lam)^*
,\bd_-).
\end{align}

We shall show our spectral sequence
\eqref{eq:2006-11--8-10-05}
satisfies the property 
\eqref{eq:2006-11-08-09-49}
in Proposition \ref{Pro:2006-11-08-09-34}
in next sections.
\subsection{Computation of the Term $(E_1)_{\ttp}$}
First we compute
$(E_1^{p,q})_{\bra \lam,\Dg\ket}$.
By definition,
\begin{align}
 (E_1^{p,q})_{\bra \lam,\Dg\ket}=H^{p+q}(\gr^{\FF}_p \Clam^*_{\tp}).
\end{align}
We have:
%By definition
%$\{\FF^p \Clam^*_{\tp}\}$
%is the filtration dual to the filtration
%$\{\FF_p \Clam_{\tp}\}$
%defined by \eqref{eq:def-of-Gp-bar}:
\begin{align}
&\FF^p \Clam^*_{\tp}=
(\Clam/ \FF_p \Clam)^*_{\tp}=
(\Clam_{\tp}/\FF_p \Clam_{\tp})^*,
\label{eq:2006-10-16--8-47}
\\
&\gr_{\FF}^p \Clam_{\tp}^*= (\gr^{\FF}_{p+1}\Clam_{\tp})^*,
\end{align}
by \eqref{eq:2006-10-06-14-22}.
By Proposition \ref{Pro:2006-05-11-1} (iii),
\begin{align}
\bd_-^{\st}\cdot  \FF^p \Clam^*_{\tp}
\subset \FF^{p} \Clam^*_{\tp}\quad
\chi_-\cdot \FF^p \Clam^*_{\tp}
\subset \FF^{p+1} \Clam^*_{\tp}.
\end{align}
Hence
\begin{align}
 (E_1^{p,q})_{\bra \lam,\Dg\ket}&=
%H^{p+q}(\gr_\FF^p\Clam_{\tp}^* ,\bd_-)=
\bigoplus_{\mu\atop
\bra \mu,\srho\che\ket=-p-1}
H^{p+q}(\Clam_{\tp}^* ,
\bd_-^{\st})^{\mu},
\label{eq:2006-11-08-03-37}
\end{align}
by \eqref{eq:def-of-Gp-bar},
because  $\bd_-^{\st}$ commutes with the action of $\h$.
Therefore,
we can apply %\eqref{eq:2006-10-06-00-37}
 Lemma  \ref{Lem:elementary-duality-of-cohomology}
 to get:
\begin{align}
& (E_1^{p,q})_{\bra \lam,\Dg\ket}
=
\bigoplus_{\mu\atop\bra \mu,\srho\che\ket=-p-1}
H^{-p-q}(\Clam_{\tp}^{\mu}, \bd_+^{\st})^*
\label{eq:2006-10-07-13-43}
\end{align}

On the other hand
we have the following
assertion,
which is easily seen from \eqref{eq:2006-11-08-10-13}.
\begin{Lem}\label{Lem:2006-11-08-10-20}
% We have $\Clam^*_{\tp}
%=C(\Ln{-},M(\lam)^*)_{\tp}=
%\sM(\slam)^*\* \Lam(\sn_-)$.
The complex 
$(\Clam_{\tp}, \bd_+^{\st})$
is identical to the Chevalley complex
for calculating
the Lie algebra cohomology
$H^{\bullet}(\sn_+, \sM(\slam))$:
\begin{align*}
 H^i(\Clam_{\tp}, \bd_+^{\st})^{\mu}=
\begin{cases}
 H^{i}(\sn_+, \sM(\slam))^{\mu}&\text{for $i\geq 0$},\\
0&\text{for $i<0$}.
\end{cases}
\end{align*}
\end{Lem}
As is well-known.
$ H^{\bullet}(\sn_+, \sM(\slam))^{\mu}=0$
unless $\mu\in \sW\circ \slam$.
Thus
by Lemma \ref{Lem:2006-11-08-10-20}
it follows that
\begin{align}
  H^{\bullet}(\Clam_{\tp}, \bd_+^{\st})^{\mu}=0
\quad\text{unless  $\mu\in \sW\circ \slam$.}
\label{eq:2006-11-08-10-26}
\end{align}
\subsection{Computation of the Term $E_1$}
By Lemma \ref{Lem:2006-10-06-13-22},
each $\FF^p \Clam^*$ is a direct sum 
of finite-dimensional subcomplexes
$\FF^p \Clam^*_d$,
and thus so is the complex
$\gr_{\FF}^p \Clam^*$.
Hence
we have
\begin{align*}
\gr_\FF^p C_+(\lam)^*=
\left( \gr_\FF^{p+1}C_+(\lam)\right)^*
%\left( \FF_{p+1}C_+(\lam)/\FF_{p}C_+(\lam)\right)^*
=\D\left( \gr^\FF_{p+1}C_+(\lam)\right),
\end{align*}
By \eqref{eq:2006-10-06-14-22}.
Therefore,
applying
Lemma \ref{Lem:elementary-duality-of-cohomology},
we obtain the following assertion.
\begin{Pro}
\label{Pro:duality-gr}
There is an isomorphism 
 \begin{align*}
E_1^{p,q}
%\left(= H^{p+q}(\FF^p C_+(\lam)^*/\FF^{p+1}C_+(\lam)^*)\right)
\cong \D(H^{-p-q}(\gr^\FF_{p+1}\Clam))
 \end{align*}
of $\Wg$-modules
for each $p$ and $q$.
\end{Pro}
By Proposition \ref{Pro:duality-gr}
the computation 
of the term $E_1$ 
is reduced to the 
computation of
the cohomology
of the cochain complex 
$(\gr^{\FF}\Clam, \bd_+')$.

\smallskip

By  \eqref{eq:2006-10-07-13-32}
we have
\begin{align}
\gr^{\FF}_p \Clam\cong
\U(\bCgo)\*_{\U(\bCgo)_{\geq 0}}
(\gr^{\FF}_p \Clam_{\ttp})
\end{align}
as $\bCgo$-modules.
Under this identification
the action of $\bd_+'$
is described by the following formula:
\begin{align}
 \bd_+'(u\* v)=[\bd_+',u]\* v+(-1)^{p(u)}u\* \bd_+^{\st}v
\label{eq:2006-10-07-15-40}
\end{align}
with
$ u\in \U(\bCgo)$,
$v\in \gr^{\FF}\Clam_{\tp}$.

Observe
that
the $\U(\bCgo)_{\geq 0}$-module
$\gr^{\FF}_p \Clam_{\ttp}$
is a direct sum of one-dimensional
representations.
Indeed,
for a weight vector $v\in \gr^{\FF}_p \Clam_{\ttp}^{\mu}$,
we have
\begin{align}
 &\wJ_{\alpha}(n)v
=\psi_{\alpha}(n)v=0
\text{ for $\alpha\in \snroots$, $n\geq 0$;}
\label{eq:2006-10-13-00-40-2}\\
&\wJ_i(n)v
=0
\text{ for $i\in \sI$, $n>0$;}
\label{eq:2006-10-13-00-40-3}\\
&\wJ_i(0)v
 =\bra \lam,J_i\ket v
\text{ for $i\in \sI$. }
\label{eq:2006-10-13-00-40-4}
\end{align}
Hence
it follows that
\begin{align}
 \gr^\FF_p \Clam\cong   
\bigoplus_{\mu\atop
\bra \mu,\srho\che\ket=-p}
K(\mu)\+ C_+(\lam)^{\mu}_{\tp}
\label{eq:revised1}
\end{align}
(notation Section \ref{subsection-Old-coho-Varma})
as cochain complexes and
$\bCgo_{\old}$-modules,
where  
the right-hand-side
is considered
as a 
direct sum of tensor products of cochain complexes
$(K(\mu), \bd_+')$
and $(\Clam_{\tp}^{\mu},\bd_+^{\st})$,
and as a $\bCgo_{\old}$-module
on  which
 $\bCgo_{\old}$ acts on the first factor $K(\mu)$.

By \eqref{eq:revised1}
the following assertion follows 
from 
Proposition \ref{Pro:2006-10-19-03-43}
and \eqref{eq:2006-11-08-10-26}.
\begin{Pro}\label{Pro:2006-11-02-00-16}
 There is an isomorphism
\begin{align*}
H^{i}(\gr^\FF_p \Clam)\cong 
\M(\gamma_{\slam})\*
\bigoplus_{\mu\in \dual{\bh}
\atop 
\bra \mu,\srho\che\ket=-p}H^i(\Clam_{\tp}^{\mu},
\bd_+^{\st})
\end{align*}
of $\Wg$-modules
for all $i$.
Here $H^{i}(\gr^\FF_p \Clam)$
is considered as a $\Wg$-module
via the map $\wh t_{-\srho\che}$.
\end{Pro}
By 
\eqref{eq:2006-11-06-15-14},
\eqref{eq:2006-11-06-15-13},
\eqref{eq:2006-10-07-13-43},
\eqref{eq:2006-11-08-10-26},
Propositions \ref{Pro:duality-gr}
and \ref{Pro:2006-11-02-00-16}
 it follows that
we have
\begin{align}
 E_1^{p,q}\isomap 
D(\M(\gamma_{\slam }))\*  (E_1^{p,q})_{\bra \lam,\Dg\ket}.
\label{eq:2006-10-08-03-48}
\end{align}
as $\Wg$-modules 
 for any $p,q$.
\subsection{Computation of the Term $E_{r}$}
The $\Wg$ is a vertex subalgebra of $\bCgo_{\new}$
(see Remark \ref{Rem:Wg-is-a-sub-of-bCgo}).
Therefore
$\Wg$ acts on $\Clam$ itself,
rather than its cohomology,
and $\Clam_{\tp}$ is a $\U(\Wg)_{\geq 0}$-submodule
of $\Clam$.
By \eqref{eq:2006-11-01-10-15}
we have the
cochain map
\begin{align}
\begin{array}{ccc}
  \U^-(\Wg) \* (\Clam_{\tp}/\FF_p\Clam_{\tp})
&\ra &\Clam/\FF_p\Clam \\
u\* v&\mapsto &t_{-\srho\che}(u)\cdot v
\end{array}
\label{eq:2006-11-02-00-40}
\end{align}
(notation \eqref{eq:2006-11-08-14-56}),
where
$  \Um \* (\Clam_{\tp}/\FF_p\Clam_{\tp})$
is regarded as a cochain complex 
with the differential $1\* \bd_+'$.
By Proposition \ref{Pro:2006-05-11-1}
(iv),
\begin{align*}
 \D\left(\Um\* (\Clam_{\tp}/\FF_p \Clam_{\tp})\right)
\cong
\D(\Um)\* \FF^p C(\lam)_{\tp}^*.
\end{align*}
Here 
$\Um$ is considered as a $\Dg$-module
by $\Dg\cdot u=[\DW,u]$.
Therefore
\eqref{eq:2006-11-02-00-40}
 induces
a cochain map
\begin{align}
 \FF^p \Clam^*\ra 
\D(\Um)\*\FF^p\Clam_{\tp}^*,
\label{eq:2006-10-08-02-06}
\end{align}
where $\D(\Um)\*\FF^p\Clam_{\tp}^*$
is considered as a cochain  complex with the
differential $1\* \bd_-$.

It is clear the following diagram commutes:
\begin{align}
 \begin{CD}
  \Um \* (\Clam_{\tp}/\FF_{p}\Clam_{\tp})
@>>> \Clam/\FF_p\Clam\\
@VV  V @ V{} VV \\
  \Um\* (\Clam_{\tp}/\FF_{p+1}\Clam_{\tp})
@>>> \Clam/\FF_{p+1}\Clam.
 \end{CD}
\end{align}
Thus we have the
following commutative diagram 
of cochain complexes:
\begin{align}
 \begin{CD}
 \FF^p \Clam^*
@>>>  \D(\Um)
\* \FF^p \Clam_{\tp}^*
\\
@AA  A @ A AA \\
\FF^{p+1} \Clam^*
@>>>  \D(\Um)
\* \FF^{p+1} \Clam_{\tp}^*.
 \end{CD}
\end{align}
Therefore
the map \eqref{eq:2006-10-08-02-06}
induces the  following homomorphism of spectral sequences:
\begin{align}
 (E_r, d_r)
\rightarrow (\D(\Um)\* (E_r)_{\bra \lam,\Dg\ket},
1\* d_r)
\label{eq:2006-10-08-04-03}
\end{align}

We claim that
the map \eqref{eq:2006-10-08-04-03}
 is an isomorphism for $r\geq 1$.
To prove this 
it is sufficient to show the assertion  for $r=1$.
But for $r=1$
\eqref{eq:2006-10-08-04-03}
is identical to the isomorphism
\eqref{eq:2006-10-08-03-48},
under the isomorphism \eqref{eq:2006-10-18-16-30}.
Moreover
for $r\geq 1$
we have
\begin{align}
 E_r^{p,q}\cong D(\M(\gamma_{\slam}))\* (E_r^{p,q})_{\bra \lam,\Dg\ket}
\end{align}
as $\Wg$-modules,
where $(E_r^{p,q})_{\tp}$ is considered as the space 
of multiplicity.
Indeed this is true for $r=1$ by
\eqref{eq:2006-10-08-03-48},
and thus must be true for all $r\geq 1$.

We have shown that the filtration
$\{\FF^p \Clam^*\}$ satisfies the desired
property.
Proposition  \ref{Pro:2006-11-08-09-34} is proved.
Hence
 Theorem \ref{Th:vanishing-dual-Verma--} is proved.
\qed.

\section{Results for the Functor $H^0_+(?)$}
\label{section:results-for-+}
In this section we assume that $k\ne -h\che$.
\subsection{Results for the ``$+$'' Reduction}
For $\lam\in \dual{\bh}_k$
let $\BGG^{[\lam]}_k$
be the Serre full subcategory of $\BGG_k$
whose objects have all their local composition factors isomorphic to
$L(w\circ \lam)$ with $w\in \bW(\lam)$.
We have
\begin{align}
 \BGG_k=\bigoplus\limits_{[\lam]\in \dual{\bh}_k/\sim}\BGG^{[\lam]}_k.
\end{align}
where $\sim$
is  the equivalent relation defined by
$\lam\sim \mu \iff \mu\in W(\lam)\circ \lam$
(\cite{Kumer}).

Consider the following condition on $\lam\in \dual{\bh}_k$:
\begin{align}\label{eq:condition+}
 \roots_+(\lam)%\cap \proots
\cap t_{\srho\che}(\nrroots)
=\emptyset.
\end{align}
Since 
\begin{align}
 \prroots\cap t_{\srho\che}(\nrroots)= \{-\beta+n\delta; 
\beta\in \sroots_+,\ n=1,2,\dots,\height{\beta}\},
\end{align}
this condition is equivalent to
\begin{align}\label{cond:condition-+2}
 \bra \lam+\rho,\alpha\che\ket\not\in \Z
\text{ for  all $\alpha\in \{-\beta+n\delta; 
\beta\in \sroots_+,\ n=1,2,\dots,\height{\beta}\}
$}.
\end{align}\
Because $\roots(w\circ \lam)=\roots(\lam)$ for $w\in W(\lam)$,
the condition \eqref{eq:condition+}
can be regarded as a condition on the block $\BGG^{[\lam]}_k$.

The following assertion is straightforward from definition.
\begin{Lem}\label{Lem:roots-eq-t-rho}
 Suppose that $\lam\in \dual{\bh}_k$
satisfies the condition \eqref{eq:condition+}.
\begin{enumerate}
 \item The correspondence
$\alpha\mapsto t_{-\srho\che}(\alpha)$
gives a bijective map from $\roots_+(\lam)$ to 
$\roots_+(t_{-\rho\che}\circ \lam)$.
\item We have $\roots(t_{-\rho\che}\circ \lam
)\cap \sproots=\emptyset$.
That is,
the classical part 
\begin{align*}
\overline{t_{-\rho\che}\circ \lam}
=\slam-(k+h\che)\srho\che\in \dual{\sh}
\end{align*}
 is anti-dominant.
\end{enumerate}
\end{Lem}
\begin{Lem}\label{Lem:preserve-wt-+}
Assume that $\lam\in \dual{\bh}_k$
satisfies the condition \eqref{eq:condition+},
and let
$M(w\circ \lam)\hookrightarrow M(\lam)$,
with $w\in W(\lam)$,
be a non-trivial homomorphism
of $\bg$-modules.
Then
we have
$\bra w\circ \lam,\DW\ket<\bra \lam,\DW\ket$.
\end{Lem}
\begin{proof}
 See  \cite[Lemma 7.5]{A04}.
\end{proof}
The following result was 
established in \cite{A04}.
\begin{Th}\label{Th:vanishing-+}
Assume that $\lam\in \dual{\bh}_k$
satisfies the condition \eqref{eq:condition+}.
Let $M$ be an object of $\BGG_k^{[\lam]}$.
Then the cohomology 
$H^i_+(M)$ is zero  for all  $i\ne 0$.
%In particular
The correspondence
$M\mapsto H^0_+(M)$
defines an exact functor from $\BGG_k^{[\lam]}$
to $\Wg\adMod$.
\end{Th}

\begin{Th}\label{Th:Main-results-for-+}
Assume that $\lam\in \dual{\bh}_k$
satisfies the condition \eqref{eq:condition+}.
Then
there is an isomorphism
\begin{align*}
H_+^0(L(\lam))
\cong \why(\gamma_{\overline{t_{-\srho\che}\circ \lam}})
=\why(\gamma_{\slam-(k+h\che)\srho\che}).
\end{align*}
\end{Th}
\begin{Rem}$ $

\begin{enumerate}
 \item  If
$\lam\in \dual{\bh}_k$
satisfies the condition \eqref{eq:condition+},
then  
any element of $W(\lam)\circ \lam$ also satisfies 
the condition \eqref{eq:condition+}.
Thus 
$H_+^0(L(w\circ \lam))
\cong \why(\gamma_{\overline{t_{-\srho\che}w\circ \lam}})$
for any $w\in W(\lam)$
under the assumption of Theorem \ref{Th:Main-results-for-+}.
\item For $\sg=\mathfrak{sl}_2(\C)$
we do not need any   condition on the
block, see \cite{A05Duke}.
\end{enumerate}
\end{Rem}
\begin{Co}\label{Co:2006-11-08-21-11}$ $

 \begin{enumerate}
  \item \ban{\cite{FF_W}} Let  $k\not\in \Q$.
Then the cohomology $H^i_+(M)$ is zero for all $i\ne 0$
with any object $M$ of $\BGG_k$
and there is an isomorphism
$H^0_+(L(\lam))\cong \why(\gamma_{\slam-(k+h\che)\srho\che})$
for all $\lam\in \dual{\bh}_k$.
\item Suppose that
\begin{align*}
r\che j k\not\in \Z
\quad \text{for $j=1,2,\dots, h$},
\end{align*}
where $h$ is the Coxeter number of $\sg$.
Then the cohomology $H^i_+(L(k\Lam_0))$ is zero
for all $i\ne 0$ and 
$H^0_+(L(k\Lam_0))$ is  isomorphic to
$\why(\gamma_{\vac})$
as vertex algebras.
 \end{enumerate}
\end{Co}

Suppose that $\Pr^{\nondeg}\ne \emptyset$.
Then $k=p/q-h\che$
with some $p,q\in \N$ satisfying \eqref{eq:2006-11-09-11-28}.
Let
\begin{align}
  \dot{\Pr}\teigi \{ \slam-(k+h\che)\smu+k\Lam_0
;  \lam\in P^{p-h\che}_{++}
,\
\mu\in P_{++}^{\vee, q-h}
\}
\end{align}
(notation Section \ref{subsection:REsults-for--}).
Then
$\dot{\Pr}$ is a  subset of $\Pr$
containing $k \Lam_0$.

One checks that
any element of $\dot{\Pr}$ satisfies
the condition \eqref{eq:condition+}.
Thus we have the following
assertion, which proves 
\cite[Conjecture 3.4${}_+$(a)]{FKW} partially.
\begin{Co}\label{Co:CONG-FKW-PLUS}
 Suppose that 
there exists a non-degenerate principal admissible
weight of  $\bg$ of level $k$.
Let $\lam\in \dot{\Pr}$.
Then the cohomology $H^i_+(L(\lam))$ is zero for all $i\ne 0$
and there is an isomorphism 
$H^0_+(L(\lam))\cong \why(\gamma_{\slam-(k+h\che)\srho\che})$. 
\end{Co}

\begin{Rem}\label{Rem:2006-11-07-01-07}
 From Remark \ref{Rem:2006-11-06-23-27}
it follows that
there is 
the   surjective map
\begin{align*}
 \begin{array}{ccc}
  \dot{\Pr}&\twoheadrightarrow  &\Pr^{\nondeg} \\
\Lam&\mapsto &t_{-\srho\che}\circ \Lam.
 \end{array}
\end{align*}
Hence
$ H^0_-(L(\lam))\cong H^0_+(L(t_{\srho\che}\circ \lam))$
for all
$\lam\in \Pr^{\nondeg}$
and
all the isomorphism classes
in the set \eqref{eq:2006-11-06-22-57}
can be obtained as the image of the functor
$H^0_+(?)$.
\end{Rem}

\subsection{Proof of Theorem \ref{Th:Main-results-for-+}
}
For $M\in \BGG_k$,
let
\begin{align}
 H^{\bullet}_+(M)_d\teigi \{ [m]\in H^{\bullet}_+(M)
;  \DW [m]= d [m]\}.
\end{align} 
Then we have
$H^{\bullet}_+(M)=\bigoplus_{d\in \C}H^{\bullet}_+(M)_d$.

First,
we have the following assertion.
\begin{Pro}\label{Pro:image-of-Verma-+}
For each $\lam\in \dual{\bh}_k$
we have the following:
\begin{align*}
&\text{$H^i_+(M(\lam))=\zero$ for $i\ne 0$;}\\
&H^0_+(M(\lam))_{\tp}=H^0_+(M(\lam))_{\bra \lam,\DW\ket}=\C |\lam\ket,\\
& \ch H^0_+(M(\lam))=\ch \M(\gamma_{{ \slam-(k+h\che)\srho\che}}),
\end{align*}
where $|\lam\ket=v_{\lam}\* \1$.
\end{Pro}
\begin{proof}
 See \cite[Theorem 5.8]{A04} and its proof.
\end{proof}
\begin{Rem}
  It is not true that
$ H^0_+(M(\lam))\cong 
 \M(\gamma_{\slam-(k+h\che)\srho\che}
)$
in general (compare Theorem \ref{Th:image-of-Verma--}).
However 
this is probably  true for
$\lam$ that satisfies
the condition  \eqref{eq:condition+}.
\end{Rem}
\begin{Pro}\label{Pro:top-Zhu-+-05/02}
For any $\lam\in \dual{\bh}_k$
we have
$H^0_+(M(\lam))_{\tp}
\cong 
\C_{\gamma_{\slam-(k+h\che)\srho\che }}$
as $\Zhu(\Wg)$-modules.
\end{Pro}
\begin{proof}
The assertion follows
by observing
that
$\C |\lam\ket \cong K(t_{-\srho\che}\circ \lam)_{\tp}$
as $\Zhu(\Wg)$-modules,
see Proposition \ref{Pro:2006-10-19-02-59},
Section \ref{subsection:SimpleQuotinetVA}
and
  \cite[Proposition 4.2]{A04}.
\end{proof}

\begin{Lem}\label{Lem:+-hw-L}
Assume that $\lam\in \dual{\bh}_k$
satisfies the condition \eqref{eq:condition+}.
Then
$H^0_+(L(\lam))$ belongs to 
$\Wcat$ and
there is the following  isomorphism
of $\Zhu(\Wg)$-modules:
\begin{align*}
 H^0_+(L(\lam))_{\tp}=
H^0_+(L(\lam))_{\bra \lam,\DW\ket}
\cong \C_{\gamma_{\slam-(k+h\che)\srho\che}}.
\end{align*}
\end{Lem}
\begin{proof}
 By the exactness of the functor $H^0_+(?)$ (Theorem \ref{Th:vanishing-+}),
we have the exact sequence
of $\Wg$-modules
\begin{align*}%\label{eq:05/2/2-14-15}
 0\rightarrow H^0_+(N(\lam))\rightarrow
H^0_+(M(\lam))\rightarrow H^0_+(L(\lam))\rightarrow 0,
\end{align*}
where $N(\lam)$ is the maximal proper submodule of $M(\lam)$.
Hence $H^0_+(L(\lam))$
belongs to $\Wcat$
because $H^0_+(M(\lam))$ 
does, by Proposition \ref{Pro:image-of-Verma-+}.
To show the rest of the assertion,
we need only to show that
$H^0_+(N(\lam))_{\bra \lam,\DW\ket}=\zero$,
by Proposition \ref{Pro:top-Zhu-+-05/02}.
But this follows from
Lemma \ref{Lem:preserve-wt-+}.
\end{proof}
\begin{Pro}\label{Pro:ch-coincieds}
 Assume that $\lam\in \dual{\bh}_k$
satisfies the condition \eqref{eq:condition+}.
Then we have
$\ch H^0_+(L(\lam))=\ch \why(\gamma_{\slam-(k+h\che)\srho\che})
$.
\end{Pro}
\begin{proof}
 By the exactness of the functor $H^0_+(?)$
and Proposition \ref{Pro:image-of-Verma-+}
we have
\begin{align*}
 \ch H^0_+(L(\lam))&=\sum_{\mu}[L(\lam): M(\mu)]
\ch H^0_+(M(\mu))\\
&=\sum_{\mu}[L(\lam): M(\mu)]\ch
 \M(\gamma_{\smu-(k+h\che)\srho\che}).
\end{align*}
On the other hand
one knows from 
Theorem \ref{Th:ch-formula}
and
Lemma \ref{Lem:roots-eq-t-rho}
that
\begin{align*}
 \ch \why(\gamma_{\slam-(k+h\che)\srho\che})=
\sum_{\mu'}[L(t_{-\srho\che}\circ \lam):M(\mu')]\ch \M(\gamma_{\smu'}).
\end{align*}
Hence one needs only to show
that
\begin{align}
 [L(\lam):M(\mu)]=[L(t_{-\srho\che}\circ \lam): M(t_{-\srho\che}\circ
\mu)] \quad \text{for any $\mu$}
\end{align}
under  the condition  \eqref{eq:condition+} on $\lam$.
But this follow from Lemma \ref{Lem:roots-eq-t-rho}
and the character formula  \cite[Theorem 1.1]{KT}
of $L(\lam)$.
\end{proof}
\begin{proof}[Proof of Theorem \ref{Th:Main-results-for-+}]
By Lemma \ref{Lem:+-hw-L}
there is a
non-trivial $\Wg$-module homomorphism 
from  $\M(\gamma_{\slam-(k+h\che)\srho\che})$
to $H^0_+(L(\lam))$
that sends $|\gamma_{\slam-(k+h\che)\srho\che}
\ket$
to $|\lam\ket$.
But then 
$H^0_+(L(\lam))$ must be isomorphic to
  $\why(\gamma_{\slam-(k+h\che)\srho\che})$,
by Proposition \ref{Pro:ch-coincieds}.
\end{proof}
%We remark that, using (the proof of) \cite[Theorem 6.8]{A04},
%the following assertion can be shown in the same manner as 
%Theorem \ref{Th:vanishing-dual-Verma--}.
%\begin{Pro}\label{Pro:dual-Verma-+}
%Let $\lam\in \dual{\bh}_k$.
%Suppose that
%\begin{align*}
% \bra \lam+\rho,\alpha\che\ket\not\in \{1,2,\dots\}
%\text{ for all $\alpha\in \prroots\cap t_{\srho\che}(\nrroots)$.}
%\end{align*}
%Then
%there is an isomorphism of $\Wg$-modules
%\begin{align*}
%H^0_+(M(\lam)^*)\cong 
%D(\M(\gamma_{-w_0(\slam)-(k+h\che)\srho\che})).
%\end{align*}\end{Pro}
\appendix
\section{Compatible Degreewise Complete Algebras}
\label{section:The Current Algebra of Vertex (Super)algebras}
\subsection{Linear Topology}
Let 
$\II_0\supset \II_1\supset \dots$
be a decreasing sequence of linear subspaces
of a vector space $V$.
A {\em linear topology define by }$\II_n$
is a topology on $V$ such that
the set
$\{v+ \II_N;N\in \N\}$ forms a fundamental system
of open neighborhoods of $v\in V$.
The closure of
a subspace $U\subset  V$ is given by
$\bigcap_N (U+\II_N)$.
A completion 
$\wt V$
of $V$
with respect to the linear topology
defined by $\II_n$ is the projective limit
\begin{align*}
 \wt V=\lim_{\leftarrow \atop N}V/\II_N
\end{align*}
with the projective limit topology
induced from the discrete topology
on each $V/\II_N$.
 The topology on $\wt V$
coincides with the linear topology defined by $\wt \II_N$,
where
 $\wt \II_N$ is the 
kernel of the natural map
$\wt{V}\rightarrow V/\II_N$.
The following assertion is well-known (see e.g. \cite{MatsuBook})
\begin{Th}\label{Th:2006-05-23-14:44}
 Let $V$ a vector space
equipped with the linear topology 
defined by a decreasing sequence $\{\II_n\}$
of subspaces.
Let $U\subset V$ and give
$U$
and $V/U$ the induced topology.
Then we have the exact sequence
$0\rightarrow \wt U\rightarrow \wt V\rightarrow
\wt {U/V}\rightarrow 0$
and $\wt U$ coincides the
closure of
the image  of $U$ in 
$\wt V$.
\end{Th}
\subsection{Compatible Degreewise Complete Algebras}
\label{subsectuib:Compatible Degreewise Complete Algebra}
Following \cite{MNT},
by
a {\em compatible degreewise topological algebra}
we mean a 
$\Z$-graded  algebra
$A=\bigoplus_{d\in \Z}A_d$ 
satisfying the following:
\begin{enumerate}
 \item  each $A_d$ equipped with a
topology  
and the multiplication
map $A_d\times A_{d'}\rightarrow A_{d+d'}$
is continuous with respect to this topology;
\item 
the topology of
$A_d$
coincides with the linear topology 
define by 
$\II_N(A_d)$ with $N\geq \max (0,d)$,
where 
$\II_N(A_d)$ is the closure
of $\sum_{r>N}A_{d-r}A_r$ in $A_d$.
\end{enumerate}
Further,
if $A$ is {\em degreewise complete},
that is,
if each $A_d$ is complete,
then
$A$ is called a
  {\em compatible degreewise complete  algebra}.
The {\em degreewise completion}
of a compatible degreewise topological algebra
is 
$\wt{A}=\bigoplus_{d\in \Z}\wt{A}_d$,
where each $\wt{A}_d$ is the completion
of $A_d$:
\begin{align*}
 \wt{A}_d=\lim_{\leftarrow \atop N}
%Q_N(A_d),\quad
%Q_N(A_d):=
A_d/\II_N(A_d).
\end{align*}
Then $\wt{A}$ is a compatible degreewise complete algebra.
By definition $A$ is degreewise complete if and 
only if $\wt{A}=A$ as compatible degreewise topological algebras.

Let $A=\bigoplus_{d\in \Z}A_d$
be any $\Z$-graded algebra.
The {\em standard degreewise topology}
of $A$
is the
 the linear topology
defined by
the sequence $\bigoplus_{d\in \Z}\left( \sum_{r\geq N}A_{d-r}A_r\right)
$.
This 
makes $A$ a 
compatible  degreewise  topological algebra.
The corresponding
degreewise completion $\wt{A}$ of 
$A$ is called the {\em standard degreewise completion}.

Let $J=\bigoplus_{d\in \Z}J_d$
be a graded subspace of a compatible degreewise 
 topological algebra $A$.
The sum $\wt{J}$ of the closures
$\wt J_d$
of $J_d$ in $A_d$
is called the 
{\em degreewise closure}
of $J$. 
If $J$ is an ideal then $\wt{J}$ is also an ideal of $A$.
Further,
if $A$ is complete then 
$\wt{J}$ is the degreewise completion of $J$
with respect to the relative topology induced from $A$,
and the quotient algebra
$A/\wt{J}$ 
with the quotient topology is also degreewise complete
(Theorem \ref{Th:2006-05-23-14:44}).

\end{document}